\documentclass{memo-l}

\usepackage{amsmath, amsthm, amssymb}
\usepackage{mathtools}

\usepackage{graphicx, color}
\usepackage{listings}  

\newtheorem{thm}{Theorem}[chapter]
\newtheorem{lemma}[thm]{Lemma}
\newtheorem{cor}[thm]{Corollary}
\newtheorem{fact}[thm]{Fact}
\newtheorem{facts}[thm]{Facts}
\newtheorem{conj}[thm]{Conjecture}

\theoremstyle{definition}
\newtheorem{dfn}[thm]{Definition}
\newtheorem{dfns}[thm]{Definitions}
\newtheorem{example}[thm]{Example}
\newtheorem{examples}[thm]{Examples}

\theoremstyle{remark}
\newtheorem{remark}[thm]{Remark}
\newtheorem{remarks}[thm]{Remarks}

\numberwithin{section}{chapter}
\numberwithin{equation}{chapter}
\numberwithin{figure}{section}

\newcommand{\wt}[1]{\widetilde{#1}}
\newcommand{\wh}[1]{\widehat{#1}}
\newcommand{\co}{\colon\thinspace}

\newcommand{\bound}{\partial}
\newcommand{\pie}{\pi_1}
\newcommand{\AND}{\qquad \text{and} \qquad}

\newcommand{\C}{\mathbb{C}}
\newcommand{\IH}{\mathbb{H}}
\newcommand{\K}{\mathbb{K}}
\newcommand{\Q}{\mathbb{Q}}
\newcommand{\R}{\mathbb{R}}
\newcommand{\Z}{\mathbb{Z}}
\newcommand{\X}{\mathbb{X}}

\newcommand{\calA}{\mathcal{A}}
\newcommand{\calD}{\mathcal{D}}
\newcommand{\calE}{\mathcal{E}}
\newcommand{\calF}{\mathcal{F}}
\newcommand{\calG}{\mathcal{G}}
\newcommand{\calH}{\mathcal{H}}
\newcommand{\calN}{\mathcal{N}}
\newcommand{\calP}{\mathcal{P}}
\newcommand{\calQ}{\mathcal{Q}}
\newcommand{\calR}{\mathcal{R}}

\newcommand{\calU}{\mathcal{U}}
\newcommand{\calV}{\mathcal{V}}
\newcommand{\calZ}{\mathcal{Z}}

\newcommand{\dQ}{d_\calQ}

\newcommand{\frg}{\mathfrak{g}}

\newcommand{\frp}{\mathfrak{p}}
\newcommand{\frx}{\mathfrak{x}}
\newcommand{\fry}{\mathfrak{y}}
\newcommand{\frz}{\mathfrak{z}}
\newcommand{\dev}{\operatorname{dev}}
\newcommand{\expect}{\operatorname{expect}}

\DeclareMathOperator{\tr}{Tr}
\DeclareMathOperator{\PSL}{PSL}
\DeclareMathOperator{\SL}{SL}
\newcommand{\PSLtwo}[1]{\PSL_2\!\left(#1\right)}
\newcommand{\PGLtwo}[1]{\PGL_2\!\left(#1\right)}
\newcommand{\SLtwo}[1]{\SL_2\!\left(#1\right)}
\newcommand{\pslC}{\PSLtwo{\C}}

\newcommand{\pslQx}{\PSLtwo{\Q(x)}}
\newcommand{\pglQx}{\PGLtwo{\Q(x)}}
\DeclareMathOperator{\PGL}{PGL}
\DeclareMathOperator{\GL}{GL}

\newcommand{\redden}[1]{\color{red}#1\color{black}}

\usepackage{longtable}
\usepackage{booktabs}
\usepackage[colorlinks,linkcolor=blue,citecolor=blue,pdfstartview=FitH]{hyperref}

\makeindex

\begin{document}

\frontmatter

\title{Geometry for Kleinian Groups Generated by a Parabolic Pair}

\author{Eric Chesebro}
\address{Department of Mathematical Sciences, University of Montana} 
\email{Eric.Chesebro@mso.umt.edu} 

\begin{abstract}
This paper develops a unified and explicitly computable framework for studying the hyperbolic geometry of $2$-bridge links and Kleinian groups generated by two parabolic elements.  The framework is built around Sakuma--Weeks complexes together with a Farey recursive polynomial function $\calQ$ with values in $\Z[x]$.  Although $\calQ$ is defined by a short recursion on the Farey graph, it encodes a surprising amount of the geometry of these spaces and the associated $\pslC$-representations.

For a rational parameter $\alpha\in(0,1/2)$, the gluing equations carried by the corresponding Sakuma--Weeks complex are solved explicitly over $\Q(x)$.  The tetrahedral shape parameters are given by rational functions $\calZ$ constructed from $\calQ$.  When $\alpha$ determines a hyperbolic $2$-bridge link, folding the complex at the bottom gives the classical Sakuma--Weeks triangulation, and the additional equations reduce to
\[ \calQ(\alpha)\ =\ 0. \]
Thus the complete $2$-bridge structure appears as a distinguished specialization of this generic family of solutions to the gluing equations.  In this case, $\calQ(\alpha)$ has a unique geometric root which determines the complete hyperbolic structure.  Thus, given $\alpha$ for a hyperbolic $2$-bridge link, one can recursively compute $\calQ(\alpha)$, locate its geometric root, and evaluate $\calZ$ there to determine all tetrahedral shapes.  A second Farey recursive function $\calN$, together with $\calQ$, gives explicit developed domains and holonomy representations.  These formulas identify the holonomies of natural fundamental group elements arising from tangle diagrams, including meridians, crossing loops, and curves of prescribed slope on the $4$-punctured spheres in the decomposition.

The same theory applies when the complex is not folded at its bottom, including to incomplete structures for which the associated group need not be discrete.  In particular, it applies throughout the Riley slice and to parameters corresponding to Heckoid orbifolds, which, together with hyperbolic $2$-bridge link complements, arise in the classification of nonfree Kleinian groups generated by a parabolic pair.  A simple identity relates $\calQ(\alpha)$ to the Riley polynomial.  Associated to $\calQ$ is a discriminant polynomial function $D_{\calQ}$; the distinguished geometric roots of the polynomials $D_{\calQ}(\alpha)$ are precisely the cusp points of the Riley slice in the first quadrant.

Further applications include formulas for cusp fields, symmetries arising from reversed and palindromic continued fractions, and substantial reductions in the degrees of exact polynomial equations satisfied by geometric roots of certain palindromic $2$-bridge links. The same methods describe algebraic and geometric limits obtained by increasing the number of crossings in a twist region, leading to explicit descriptions of certain augmented-link complements.  The paper concludes with examples in which $\calQ$ and Farey recursion give geodesic polyhedral decompositions of the convex cores of Heckoid orbifolds.
\end{abstract}

\maketitle

\tableofcontents

\chapter{Acknowledgments}\label{chap: ack}

Thanks to Jessica Purcell for two fruitful trips to Melbourne and countless good ideas and inspiring meetings.  Thanks also to Alex Elzenaar for several illuminating conversations and help with references.  

Stephan Tillmann showed me how to think about the Sakuma--Weeks triangulations as described in Section \ref{sec: SW}.  This approach was essential for the constructions in this paper.  Gaven Martin taught me many things about the Riley slice.  In particular, he introduced me to the Heckoid orbifolds discussed in Section \ref{sec: Heck} and taught me the importance of \cite{AHOPSY}.

Thanks also to the referee.  They made many thoughtful comments and corrections to an earlier draft and the paper is greatly improved because of their work.

After the mathematical content of the paper was finished, I used ChatGPT to help with  grammar, spelling, and typesetting issues. I also used it to help me generate the symbols list and the index.  I take full responsibility for the content of this paper.

\mainmatter
\chapter{Introduction}\label{chap: intro}

This paper develops an explicit bridge between the geometry of the Sakuma--Weeks triangulations of hyperbolic $2$-bridge links and the parameter space of groups generated by two parabolic elements.  \index{2-bridge link}\index{Sakuma--Weeks complex!triangulation} The central algebraic object is a Farey recursive polynomial function\index{Farey recursive function}
\[ \calQ \co \widehat{\Q}_0 \longrightarrow \Z[x]\]
organized by the Stern--Brocot diagram.\index{Stern--Brocot!diagram}  Together with the related functions $\calN$ and $\calZ$, the function $\calQ$ gives explicit tetrahedral shapes, developed domains, and holonomy representations, linking the combinatorics of the Sakuma--Weeks complexes to the two-parabolic parameter space.  In particular, $\calQ$ is related explicitly to the classical Riley polynomials, the distinguished roots of the associated discriminant polynomials are identified with the cusp points of the Riley slice, and the geometric constructions are used to study limiting and Heckoid geometries.

\section{Background}

\begin{figure}
   \centering
   \includegraphics[width=3.4in]{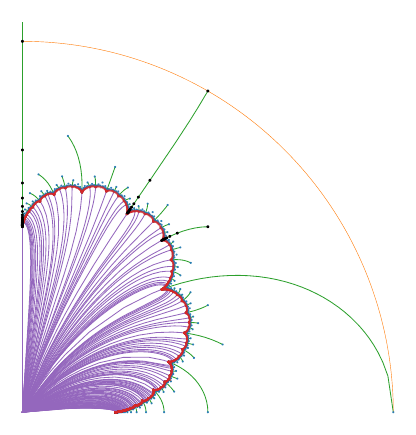} 
   \caption{The first quadrant of $\C-\{0\}$ provides a parameter space for the irreducible subgroups of $\pslC$ which are generated by a pair of parabolic elements.  The dots and curves in this figure help describe the anatomy of this parametrization.  For scale, the unit circle is shown in orange.}  \label{fig: pleats0}
\end{figure}

A Kleinian group is a discrete subgroup of $\pslC$; throughout this paper  we restrict attention to nonelementary Kleinian groups.  Kleinian groups act isometrically and properly discontinuously on hyperbolic 3-space, and the corresponding quotients are hyperbolic 3-orbifolds.  The study of Kleinian groups, initiated by Klein, Fricke, and Poincar\'e in the late nineteenth century, has continued to be a rich and central topic in mathematics.  A resurgence of interest followed the work of Ahlfors and Bers in the 1960s (see, for example, \cite{BersUniform}, \cite{AhlFinite},  \cite{Maskit}, \cite{MarRigid}, and \cite{JorIneq}).   Thurston's work in the late 1970s further deepened the subject by revealing its profound connections to the topology of $3$-manifolds.  The subsequent development is far too extensive to summarize here,  but some highlights include \cite{NR}, \cite{ATame}, \cite{CGTame}, \cite{MinskyLam}, and \cite{BCMLam}.

Among the simplest nontrivial Kleinian groups are those generated by two parabolic elements, and these form the focus of this paper.  Up to conjugacy in $\pslC$, any such group is generated by 
\begin{align*}  
\begin{bmatrix} 1&1\\0&1 \end{bmatrix} \qquad \text{and} \qquad  \begin{bmatrix} 1&0\\1/\frz&1 \end{bmatrix}\end{align*}
where $\frz$ is a nonzero complex number.  After perhaps inverting generators and conjugating in the full isometry group of $\IH^3$, the number $\frz$ can be assumed to lie in the first quadrant of $\C$.  (This parametrization is slightly nonstandard: the present choice simplifies the connection with the ideal triangulations constructed in this paper.)

For  $\frz \in \C^\ast = \C-\{0\}$, let $\Gamma_\frz$ be the subgroup of $\pslC$ generated by the matrices above.  In \cite{OM}, it is shown that, together with zero, the set
\[ \overline{\calR} \ =\ \left\{ \frz \in \C^\ast \, \big| \, \Gamma_\frz \text{ is discrete and free with rank 2} \right\}\]
is a closed topological disk whose boundary is a Jordan curve.   Symmetries of the generators imply that this disk is invariant under reflection in both coordinate axes.  The Jordan boundary curve is evident in Figure \ref{fig: pleats0}, which shows the first quadrant of $\C$.

The interior $\calR$ of $\overline{\calR}$ is the {\it Riley slice}\index{Riley slice} (Definition \ref{def: Riley slice}).  Classically, it is expressed as
\[ \calR \ = \ \left\{ \frz \in \overline{\calR} \, \big| \,  \calD(\Gamma_\frz)/\Gamma_\frz \text{ is a 4-punctured sphere}\right\}\]
where $\calD(\Gamma_\frz)$ is the domain of discontinuity for $\Gamma_\frz$ in the Riemann sphere.   A rich internal structure of $\overline{\calR}$ was revealed in \cite{KS}, where it was shown that it is foliated by {\it pleating rays}\index{pleating ray}; the rational rays form a dense subset parametrized by rational numbers and the endpoints of rational pleating rays are called {\it cusps}\index{Riley slice!cusp} and correspond to geometrically finite groups.  Several of these rational rays are shown in purple in Figure \ref{fig: pleats0}. \index{Riley slice!boundary}
  
According to the preface of \cite{ASWY}, the rational pleating rays should extend smoothly beyond the boundary of the Riley slice.  These extensions either terminate at a real point or they terminate at a point $\frz \in \IH^2$ for which $\Gamma_\frz$ is a Kleinian group which uniformizes a hyperbolic 2-bridge link complement.\index{2-bridge link!complement}   When this happens, the rational parameter of the pleating ray coincides with that of the 2-bridge link.  By \cite{Adams_2gen}, every torsion-free, finite-covolume Kleinian group generated by two parabolics arises in this way.  Many such extensions and points $\frz$ can be seen in Figure \ref{fig: pleats0}. 

Consider now a hyperbolic 2-bridge link $L$ with rational parameter $\alpha$.   As discussed below Definition \ref{def: SW}, one may assume that $\alpha \in (0,1/2)$ and that $\alpha$ is not the reciprocal of an integer.   Let $\frz \in \C^\ast - \overline{\calR}$ be the endpoint of the pleating ray for $\alpha$ so that $\IH^3/\Gamma_\frz$ is homeomorphic to the link complement $S^3-L$ and let $\hat{\frz}$ be the cusp point at the intersection of the pleating ray and the Jordan curve $\bound \calR$.  Along the ray, points between $\frz$ and $\hat{\frz}$ correspond to the hyperbolic cone manifolds $M_\theta$ obtained by assigning cone angle $\theta \in [0,2\pi]$ to the lower tunnel of $L_\alpha$ \cite{ASWY}.\index{core tunnel}   

At certain rational cone angles, these cone manifolds are hyperbolic orbifolds.  In the notation of Section~\ref{sec: Heck}, take $n\in\frac12\Z$ with $n\geq3/2$ and assign cone angle $2\pi/n$ to the lower tunnel.  Theorem~2.2 of \cite{Lee-Sakuma} shows that the associated Heckoid orbifold $S(\alpha;n)$ is geometrically finite and is uniformized by a group generated by two parabolics.\index{Heckoid orbifold}  If $n\in\Z$, then $S(\alpha;n)$ is an orbifold and is called an {\it even Heckoid orbifold}; if $n\notin\Z$, then $S(\alpha;n)$ is a $2$-fold quotient of the cone manifold and is called an {\it odd Heckoid orbifold} (Definition~\ref{def: Heckoid orbs}).\index{Heckoid orbifold!even}\index{Heckoid orbifold!odd}  The black dots along the green lines in Figure~\ref{fig: pleats0} correspond to Heckoid orbifolds on the extended pleating rays for $1/2$, $2/5$, and $3/8$.  The classification theorem \cite[Theorem~1.1]{AHOPSY} shows that every nonfree Kleinian group generated by a parabolic pair is, up to the standard symmetries, either a hyperbolic $2$-bridge link group or a Kleinian Heckoid group.\index{Heckoid group}

Although the location of these parameters $\frz$ relative to the Riley slice and its pleating rays gives a satisfying description of the corresponding Kleinian groups, it remains difficult to extract explicit
and computable geometric information about $\IH^3/\Gamma_\frz$ directly from the position of $\frz$.  The principal aim of this paper is to address this difficulty.  As indicated above, the approach here is to use Sakuma--Weeks complexes to associate explicit carried geometric structures to suitable points $\frz$ and, at complete structures, geometric triangulations, and to use these constructions to describe the geometry associated to $\Gamma_\frz$.  This will provide practical and detailed links between points $\frz \in \C^\ast$, the group $\Gamma_\frz$, and the geometry and topology of the associated hyperbolic spaces.  Often, the theory presented here will also provide geometric information about incomplete hyperbolic spaces associated to points $\frz$ for which $\Gamma_\frz$ is not discrete.

Two later applications are worth mentioning here.  Section~\ref{sec: limits} studies algebraic and geometric limits\index{algebraic limit}\index{geometric limit} obtained by increasing the number of crossings in a twist region\index{crossing loop!twist region} and gives explicit geometric descriptions of certain augmented link complements.\index{augmented link}  Section~\ref{sec: Heck} gives concrete convex core decompositions\index{convex core} for Heckoid orbifolds along extended pleating rays.  

\section{Outline and notable results}

The results below serve two complementary purposes: they give structural connections among triangulations, representations, and the Riley slice, and they provide effective procedures for computing the geometry of individual examples.  Broadly, the paper proceeds from recursive algebra to explicit hyperbolic geometry and then to applications in representation and deformation theory.

For readers interested primarily in using the constructions, the main computational route is short.  For moderate examples, compute $\calQ$ (and $\calN$ when developed domains are desired), obtain the
tetrahedral shape parameters from $\calZ$ using Theorem~\ref{thm: main1}, and locate the geometric root using
Theorem~\ref{thm: main} and Remark~\ref{rk: filter}.  SageMath implementations of the main steps are given in
Sections~\ref{subsec: code I}, \ref{subsec: code II}, and \ref{subsec: code III}, and early values of $\calQ$ and $\calN$ are listed in Appendix~\ref{chap: list}.  Much larger examples can be treated numerically without constructing $\calQ(\alpha)$ as a polynomial; see Figure~\ref{fig: Big}.  Once the geometric root is known, developed
domains and holonomy can be recovered from Theorem~\ref{thm: domain} and Corollary~\ref{cor: holonomy specialization}.

\subsection{The Sakuma--Weeks complexes and triangulations}\index{Sakuma--Weeks complex}\index{Sakuma--Weeks complex!triangulation} This subsection explains how the Sakuma--Weeks complexes arise naturally from the Stern--Brocot diagram. 

\begin{figure}
\setlength{\unitlength}{.1in}
\begin{picture}(50,51)
\put(0.1,0) {\includegraphics[width= 5in]{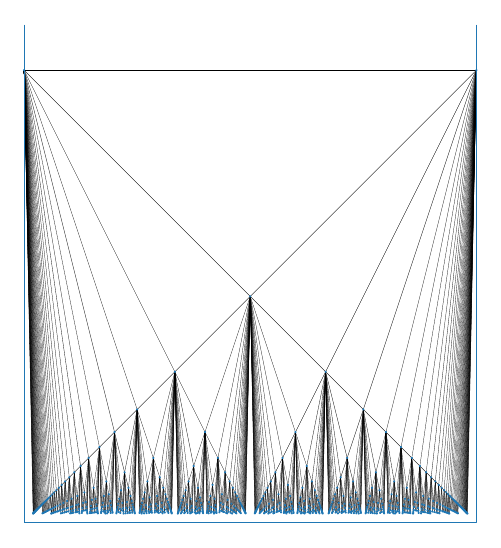}}
\put(3,48.5){$0/1$} 
\put(45.2,48.5){$1/1$}
\put(24.2,26.7){$1/2$} 
\put(17,19.4){$1/3$}
\put(31.2,19.4){$2/3$}
\end{picture}
 \caption{The Stern--Brocot diagram $\calG$.}   \label{fig: G}
\end{figure}

The Sakuma--Weeks complexes of Definition~\ref{def: SW} are obtained from the construction in \cite{SW}, which was inspired by J\o{}rgensen's ideal triangulations \cite{Jor} (also described in \cite{FH}) of hyperbolic 1-punctured torus bundles over $S^1$.  Their combinatorics are determined by a rational number $\alpha\in\Q\cap(0,1/2)$ with $\alpha\neq1/3$.  When $1/\alpha \notin \Z$, $\alpha$ is called {\it hyperbolic}\index{hyperbolic number} and the Sakuma--Weeks complex can be folded at its bottom to yield the well known Sakuma--Weeks triangulation of the 2-bridge link complement $M_\alpha$.  J\o{}rgensen, Sakuma, and Weeks use the Farey graph to encode these combinatorics.  In this paper, the {\it Stern--Brocot diagram}\index{Stern--Brocot!diagram} (see \cite{HO}, Figure~\ref{fig: G}, and Definition~\ref{def: SternBrocot}) is used instead, as many geometric features of this diagram are useful to the constructions introduced here.

The rational number $\alpha$ determines a {\it funnel} $F_\alpha$\index{funnel} (Definition~\ref{def: funnel} and Figure~\ref{fig: Funnels}) in the Stern--Brocot diagram.  Its interior edges $e_1,\ldots,e_n$\index{funnel!interior edges}, listed from top to bottom, index pairs of ideal tetrahedra $\Delta_1,\ldots,\Delta_n$.\index{Sakuma--Weeks complex!tetrahedral pair}  Gluing these pairs according to the labels inherited from $F_\alpha$ produces a triangulated thickened $4$-punctured sphere $K$.\index{Sakuma--Weeks complex!4-punctured spheres}  Folding the top boundary of $K$ gives the Sakuma--Weeks complex $M_\alpha^\circ$.  When $\alpha$ is hyperbolic, $M_\alpha^\circ$ is a manifold and folding the bottom boundary gives the classical Sakuma--Weeks ideal triangulation of $M_\alpha=S^3-L_\alpha$. If $1/\alpha\in\Z$, then $M_\alpha^\circ$ is instead a nonmanifold complex.  Figure~\ref{fig: SW link} relates this construction to the standard $4$-plat link diagram.\index{4-plat diagram}

\subsection{Solving gluing equations with Farey recursion}\index{gluing equation}  This subsection describes how a Farey recursive family of polynomials solves the $G$-invariant gluing equations carried by the Sakuma--Weeks complexes.

\begin{figure}[h] 
   \centering
   \includegraphics[width=5.0in]{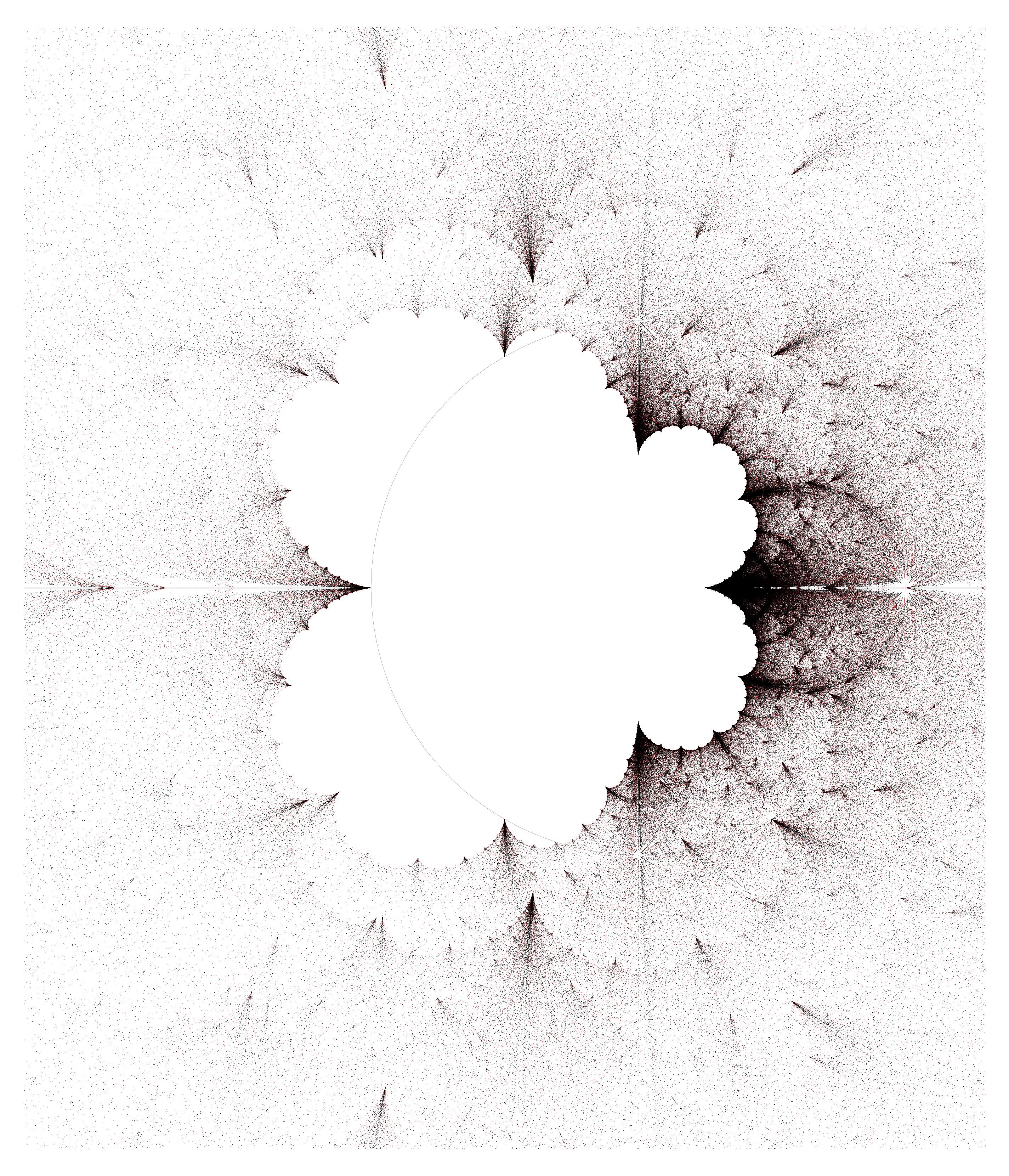}      
 \caption{A plot of over 400,000 points in $\C$ obtained as the roots of over 12,000 polynomials $\calQ(\alpha)$ and $D_\calQ(\alpha)$ with $\alpha \in [0,1/2]$.  For reference, the faint circle is the unit circle in $\C$.  If the roots of polynomials with $\alpha \in [0,1]$ are included, the picture becomes symmetric across the line $i \R$.  Theorem~\ref{thm: R is empty} shows that the white butterfly region contains no roots.}   \label{fig: allroots}
\end{figure}

A Klein $4$-group $G$ of simplicial involutions\index{involution} acts on $M_\alpha^\circ$ and preserves each tetrahedral pair $\Delta_j$ (Definition~\ref{def: Klein} and Fact~\ref{fact: isometry}).  Following \cite{BP}, the positive $G$-invariant solutions to the gluing equations carried by $M_\alpha^\circ$ form the space $\calH(M_\alpha^\circ)$ (Definition~\ref{def: hyperbolic structures}).  This definition also applies when $M_\alpha^\circ$ is not a manifold.  When $\alpha$ is hyperbolic, the lower hairpin fold\index{hairpin!fold} imposes additional equations and defines $\calH(M_\alpha)\subseteq\calH(M_\alpha^\circ)$.

Elements of $\calH(M_\alpha^\circ)$ can be identified by solving Thurston's gluing equations.  These equations are described in Section~\ref{sec: equations} and solved explicitly in Section~\ref{sec: solve} using a particular {\it Farey recursive function}\index{Farey recursive function} 
\[ \calQ \co \wh{\Q}_0 \longrightarrow \Z[x]\]
where $\wh{\Q}_0$ consists of $\{ \infty \}$ together with $\Q \cap [0,1]$.  As described in Definition \ref{def: FRF} and Example \ref{ex: Q}, the function $\calQ$ is determined succinctly from the rule $p/q \mapsto (-1)^p x^q$ along with the initial conditions
\begin{align*} 
\calQ\left(0\right)&=1 & \calQ\left(\infty\right)&=0 & \calQ\left(1\right)&=1.
\end{align*}
The function $\calQ$ can be thought of as a family of 2-term linearly recursive polynomial sequences whose recursion is governed by the geometry of the Stern--Brocot diagram.  Its values $\calQ$ are amenable to computer calculations.  Sample code for this is provided in Section \ref{subsec: code I}, and early values of $\calQ$ are listed in Appendix \ref{chap: list}.

Remark~\ref{rem: SW functions} identifies this formulation with the original Sakuma--Weeks one.  Although it is described differently, their polynomial $f$ is $\calQ$, and the tetrahedral parameters they obtain from $\psi$ agree with $\calZ$.  Remark~\ref{rem: Tsvietkova} compares $\calQ$ with the diagrammatic polynomials of Tsvietkova \cite{Tsv}.  Farey-recursive gluing equations also occur beyond the $2$-bridge setting.  Martinez \cite{Jose} uses them for hyperbolic Dehn fillings admitting triangulations by layered solid tori.

Also important are the {\it discriminant polynomials}\index{discriminant polynomial}
\[ D_\calQ(p/q) \ =\ \calQ(p/q)^2-4 \, (-1)^p x^q \]
from Definition~\ref{def: Disc}.  Figure~\ref{fig: allroots} plots roots of $\calQ(\alpha)$ and $D_\calQ(\alpha)$. Roots of polynomials $\calQ(\alpha)$ are shown as black points and roots of polynomials $D_\calQ(\alpha)$ are red.  The visible butterfly already suggests the connection with the Riley slice developed later in Section~\ref{sec: limits}.

\begin{figure}
\setlength{\unitlength}{.1in}
\begin{picture}(50,46)
\put(0,0) {\includegraphics[width= 5in]{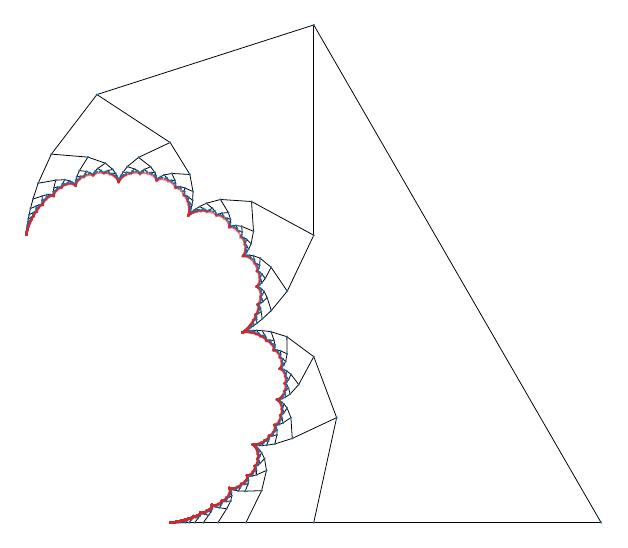}}
\put(26,42){$2/5$}
\put(5.2,36.7){$3/7$}
\put(1.4,31.8){$4/9$}
\put(46.5,0){$1/3$}
\put(24,0){$1/4$}
\put(25.8,24){$3/8$}
\put(27.8,9.3){$2/7$}
\end{picture}
\caption{The blue points are the geometric roots for the 2-bridge links.  The red points are the geometric roots for the discriminant polynomials $D_\calQ(\alpha)$.  To get a sense for the scale, the point for $1/3$ is one, the point for $1/4$ is one half, the point for $3/8$ is $\frac{1+i}{2}$, and the point for $2/5$ is $e^{i\pi/3}$. }  \label{fig: roots}
\end{figure}

Following Definition \ref{def: E and more}, let $\calE$ be the set of edges in the Stern--Brocot diagram whose vertices lie in the interval $[0,1]$.  Definition \ref{def: shape parameter function} provides a {\it shape parameter function}\index{shape parameter!function}
\[ \calZ \co \calE \longrightarrow \Q(x)\]
where the rational function $\calZ(e)$ depends on the slope of $e$ and the values of $\calQ$ on the endpoints of $e$.  As the name suggests, the edge parameters for the tetrahedra $\Delta_j$ under a structure in $\calH(M_\alpha^\circ)$ are given by the rational function $\calZ(e_j)$.  In fact, the gluing equations for $M_\alpha^\circ$ can be solved generically.  To make this precise, label the top/bottom edges of the tetrahedra in $\Delta_j$ with an indeterminate edge parameter $z_j$\index{shape parameter!edge parameter} and take $z_1=x$.

\begin{figure}
\setlength{\unitlength}{.1in}
\begin{picture}(50,41)
\put(0,0) {\includegraphics[width= 5in]{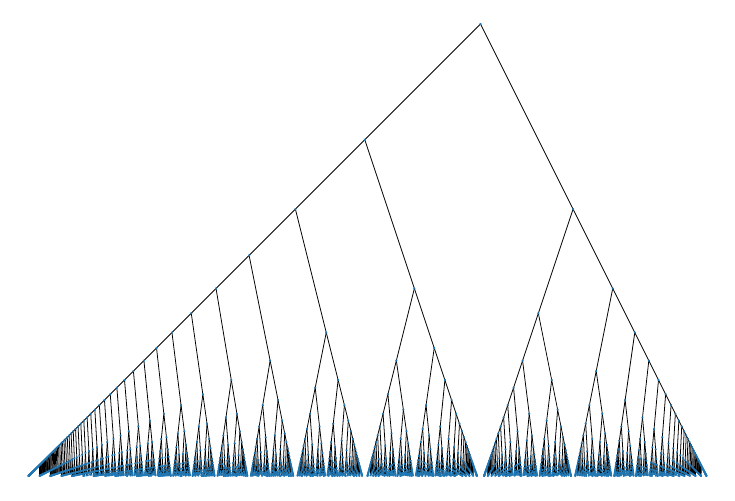}}
\put(1.2,0){$0$}
\put(31.5,33){$1/3$}
\put(22,25){$1/4$}
\put(31.5,0){$1/3$}
\put(47,0){$1/2$}
\put(40,20){$2/5$}
\put(42.7,15){$3/7$}
\put(29,15){$2/7$}
\end{picture}
\caption{The children of $[0;a_1, \ldots, a_k]$ in the Stern--Brocot tree are $[0;a_1, \ldots, a_{k-1}, 1+a_k]$ and $[0;a_1, \ldots, a_{k-1}, -1+a_k, 2]$.  Compare this tree to the tree embedded in Figure \ref{fig: roots}.  Here, brackets are used to denote continued fraction expansions.}  \label{fig: tree}
\end{figure}

\newcommand\GenericShapeSolution{Suppose $\alpha\in\Q\cap(0,1/2)$ and $\alpha\neq1/3$.  The $G$-invariant Thurston gluing equations for $M_\alpha^\circ$
determined by $F_\alpha$ have a unique solution in $\Q(x)^n$ satisfying $z_1=x$.  This solution is
\[ z_j\ =\ \calZ(e_j) \qquad \text{for} \qquad 1\leq j\leq n. \]}
\theoremstyle{plain}
\newtheorem*{GenericShapeSol}{Theorem \ref{thm: main1}}
\begin{GenericShapeSol}\GenericShapeSolution\end{GenericShapeSol}

\begin{figure}
   \centering
   \includegraphics[width=3in]{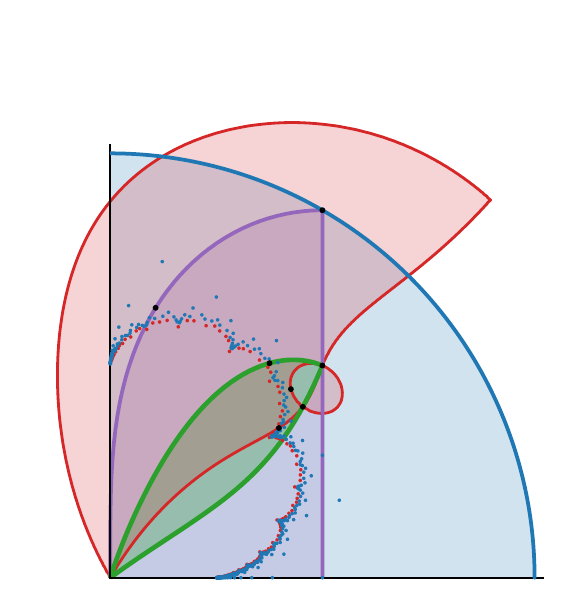} 
   \caption{This figure shows regions where $\calZ(e)$ takes values with positive imaginary part.  When $e$ is the edge $[1/2,1/3]$, the corresponding region is shaded blue, for $[1/3,2/5]$ it is purple, for $[1/3,3/8]$ it is green, and for $[4/11,3/8]$ it is red.  (See also Figure \ref{fig: RZ}.)}  \label{fig: fatregions}
\end{figure}

This shows that the edge parameter for the top pair of tetrahedra in the Sakuma--Weeks complex parametrizes $\calH(M_\alpha^\circ)$.  In particular, 
\[ \calH(M_\alpha^\circ)  \ =\ \left\{ \frz \in \IH^2 \, \big| \, \calZ(e_j)_\frz \in \IH^2 \text{ for every interior edge } e_j \text{ of } F_\alpha \right\}.\]
When $\alpha$ is hyperbolic, Theorem~\ref{thm: main} shows that the additional condition for a structure on $M_\alpha^\circ$ to descend through the lower hairpin fold is the single equation $\calQ(\alpha)_\frz=0$.  Equivalently,
 \[ \calH(M_\alpha)  \ =\ \left\{ \frz \in \calH\left(M_\alpha^\circ\right) \, \big| \, \calQ(\alpha)_\frz = 0  \right\}.\]
Figure~\ref{fig: fatregions} illustrates how the positivity conditions $\calZ(e_j)_\frz\in\IH^2$ cut out regions in the parameter plane and thereby determine the sets $\calH(M_\alpha^\circ)$.

Later in the paper, Theorem~\ref{thm: complete} shows that, when $\alpha$ is hyperbolic, every element of $\calH(M_\alpha)$ determines a complete hyperbolic structure on $M_\alpha$.  Together with the existence of the Sakuma--Weeks geometric structure and Mostow--Prasad rigidity, this implies that $\calH(M_\alpha)$ consists of exactly one point.  Its unique element is called the {\it geometric root} $\frz(\alpha)$ of $\calQ(\alpha)$ (Definition~\ref{def: geometric root}).\index{geometric root!2-bridge link}
 
Theorems~\ref{thm: main1} and~\ref{thm: main} give an effective procedure for locating the geometric root of a hyperbolic $\alpha$.  Sample code implementing this procedure is provided in Section~\ref{subsec: code II}.    For moderate examples, the sample code computes the roots of $\calQ(\alpha)$ and filters them by requiring that
$\calZ(e_j)_\frz\in\IH^2$ for every interior edge $e_j$.  Exactly one root survives this filtering process, namely the geometric root $\frz(\alpha)$.  For substantially larger examples, it is more efficient to evaluate the Farey-recursive data numerically at high precision and use Newton iteration to locate the geometric root without first constructing $\calQ(\alpha)$ as a polynomial.

Figure \ref{fig: roots} displays many geometric roots together with an embedded portion of the Stern--Brocot tree (Figure \ref{fig: tree}).  For more on the relationship between these two figures, see Conjecture 4 of \cite{Ober}.

\subsection{Developed domains and holonomies}\index{fundamental domain}\index{holonomy}  This subsection discusses the computation of explicit developed domains and holonomy representations, and explains when the developed domains are genuine embedded geometric fundamental domains.

\begin{figure}[p] 
   \centering
   \includegraphics[width=1.85in]{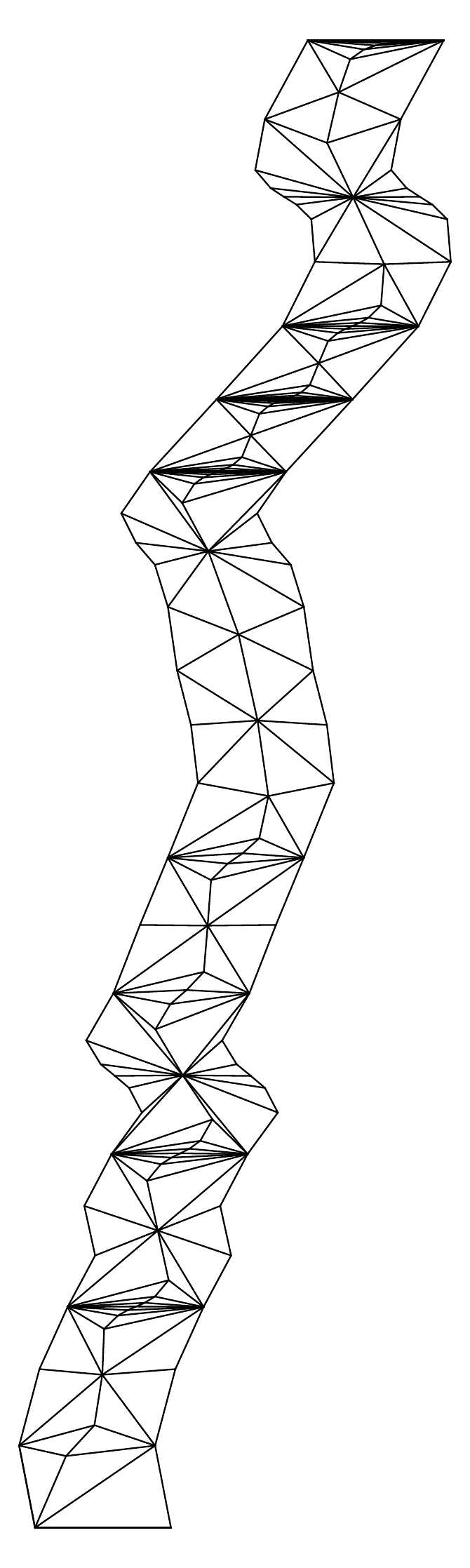} 
   \caption{\index{fundamental domain!specialized} This shows a geometric fundamental domain for the 2-bridge hyperbolic knot complement $M_\alpha$, where 
   \[ \alpha\ =\ \frac{7,243,082,905,463}{55,647,620,464,091}.\] 
The knot $L_\alpha$ has 104 crossings and the volume of $M_\alpha$ is approximately
\[105.110622606448.\]
Using Newton iteration together with the Farey recursive formulas, the geometric root, volume, and domain can be computed numerically without constructing $\calQ(\alpha)$ as a polynomial.}   
     \label{fig: Big}
\end{figure}

Thanks to the algorithmic definitions of $\calQ$, $\calN$, and $\calZ$, the constructions in Chapter~\ref{chap: FD} provide explicit developed domains for points of $\calH(M_\alpha^\circ)$ and geometric fundamental domains at complete $2$-bridge structures.  Figure~\ref{fig: Big}, for example, shows a geometric fundamental domain for the complete structure on the complement of the $2$-bridge link $L_\alpha$, where
\[ \alpha\ =\ \frac{7,243,082,905,463}{55,647,620,464,091}. \]
Although $\calQ(\alpha)$ has expected degree $27,823,810,232,045$, the numerical procedure described above readily produces the geometric root and the resulting $202$-tetrahedron fundamental domain.  In this case, $L_\alpha$ is a knot with $104$ crossings and
\[ \frz(\alpha)\ \sim \  0.29420206371521324+0.011160262055700585\,i. \]

More generally, if $F_\alpha$ has interior edges $e_1,\ldots,e_n$ and $\frz\in\calH(M_\alpha^\circ)$, then the domain $\Omega(\alpha)_\frz$ consists of a pair of ideal hyperbolic tetrahedra for each $e_j$, corresponding to the tetrahedral pair $\Delta_j$ of the Sakuma--Weeks complex (Definition~\ref{def: generic domain}).  These $2n$ tetrahedra are developed into $\IH^3$ with the combinatorics of the tetrahedral core of $M_\alpha^\circ$.

It is surprising that another Farey recursive function
\[ \calN \co \wh{\Q}_0 \longrightarrow \Z[x] \]
becomes relevant when computing domains.  Described in Example \ref{ex: N}, the definition for this function differs from that of $\calQ$ only in its initial conditions
\[ \calN\left(0\right)=\calN\left(\infty\right)=\calN\left(1\right)=1.\]
(See also Appendix \ref{chap: list}.) The ideal vertices of $\Omega(\alpha)_\frz$ are obtained by evaluating
\[ \calV = \frac{\calN}{\calQ} \]
at the vertices of the funnel $F_\alpha$ and specializing at $\frz$.

\newcommand\DomainThm{Suppose $\frz\in\calH(M_\alpha^\circ)$.  Then $\Omega(\alpha)_\frz$ is the image, under the developing map, of a connected fundamental domain for the action of $\pie(M_\alpha^\circ)$ on $\wt M_\alpha^\circ$.}
\theoremstyle{plain}
\newtheorem*{Domain}{Theorem \ref{thm: domain}}
\begin{Domain}\DomainThm\end{Domain}

For a general $\frz\in\calH(M_\alpha^\circ)$, Theorem~\ref{thm: domain} does not assert that the developing map\index{developing map} is injective on this fundamental domain, or that $\Omega(\alpha)_\frz$ is an embedded polyhedral region in $\IH^3$.  At the complete structure on a hyperbolic $M_\alpha$, the developing map is an isometry and this distinction disappears (Remark~\ref{rem: developed domain}).

Theorem~\ref{thm: domain} makes it possible to compute holonomy from $\Omega(\alpha)_\frz$.  This is done first over $\Q(x)$ and then specialized.  Definition~\ref{def: generic holonomy} identifies\index{holonomy!generic}\index{holonomy!generic representation} $\pie M_\alpha^\circ$ with the free group $\Gamma=\langle U_0,W_0\rangle\subset\pslQx$, where
\[ U_0=\begin{bmatrix}1&-1\\0&1\end{bmatrix}
   \qquad\text{and}\qquad
   W_0=\begin{bmatrix}1&0\\-1/x&1\end{bmatrix}. \]
For every $\frz\in\calH(M_\alpha^\circ)$, specialization\index{holonomy!specialization} gives the holonomy representation of the corresponding developing map (Corollary~\ref{cor: holonomy specialization}). At the geometric root of a hyperbolic $M_\alpha$, this is the complete holonomy representation.

Section~\ref{sec: holonomy} gives explicit elements $S_\omega,T_\omega\in\Gamma$, expressed in terms of $\calQ$ and $\calN$, which serve simultaneously as face pairings and as recognizable elements of the link group.  Theorem~\ref{thm: diagram} shows that, at a hub\index{funnel!hub}, $S_\omega$ is represented by a crossing loop\index{crossing loop}\index{crossing loop!holonomy}, that the relevant $T_\omega$ are distinguished Wirtinger meridians\index{holonomy!meridian}, and that $S_\omega$ represents the slope-$\omega$ curve on the embedded $4$-punctured spheres.\index{Sakuma--Weeks complex!4-punctured spheres}  Thus the same Farey-recursive formulas connect the developed domains directly to standard tangle diagrams.

\clearpage
\subsection{Recursion matrices, traces, and spirals}
Chapter~\ref{chap: logs} gives an algebraic connection between Farey recursion\index{Farey recursion} and the generic holonomy.  Take a number  $\omega \in \Q\cap(0,1/2]$ and let
\[ X_\omega=\begin{pmatrix}0&1\\-\dQ(\omega)&\calQ(\omega)\end{pmatrix} \]
be the recursion matrix for the sequences obtained by moving through the Farey neighbors of $\omega$.  Theorem~\ref{thm: spiral} constructs a conjugate $R_\omega=Y_\omega X_\omega Y_\omega^{-1}\in\PGL_2(\Q(x))$ and proves that, as M\"obius functions,
\[ R_\omega^2=S_\omega. \]
With the matrix representatives fixed in Section~\ref{sec: shapes and traces}, this gives the trace formula
\[ \tr S_\omega = 2-\frac{\calQ(\omega)^2}{\dQ(\omega)}. \]
Remark~\ref{rem: EMS} identifies the resulting trace polynomials with the Farey polynomials of \cite{EMS2} after the change of variable $x=-1/\mu$.

\newcommand\ShapeThm{Suppose $e \in \calE$ is an edge with endpoints in $\Q \cap [0,1/2]$.   Then
\[ \calZ(e) \ =\ -  \frac{2-\tr S_{e^L}}{2-\tr S_{e^R}}.\]}

Theorem~\ref{thm: traces and shapes} expresses the tetrahedral shape functions\index{shape parameter!function} directly in terms of these traces.
\theoremstyle{plain}
\newtheorem*{Shape}{Theorem \ref{thm: traces and shapes}}
\begin{Shape}\ShapeThm\end{Shape}

\begin{figure}[h] 
   \centering
   \includegraphics[width=4.2in]{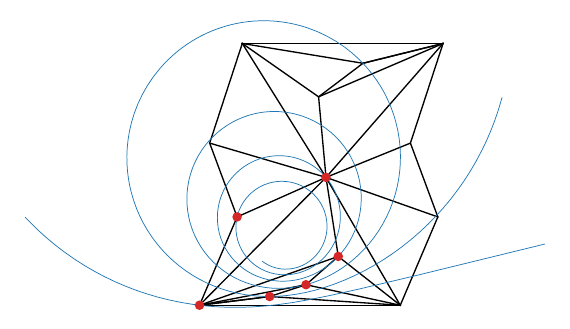} 
   \caption{\index{fundamental domain!specialized}  The blue curve is the logarithmic spiral for $(S_{3/13})_{\frz_1}$ through $\calV(1/4)_{\frz_1}$.  It lies nicely across the domain $\Omega(16/69)_{\frz_1}$, where $\frz_1$ is the geometric root for $\calQ(16/69)$.  Notice that the spiral passes through the ideal vertices of the domain (red) which are associated to the neighbors of $3/13$ in the funnel $F_{16/69}$.}        \label{fig: Spiral1669}
\end{figure}

When $(S_\omega)_\frz$ is loxodromic, its orbits on the Riemann sphere lie on logarithmic spirals\index{logarithmic spiral} about its fixed points.  The results of Section~\ref{sec: special} show how these spirals pass through the relevant ideal vertices of the developed domains.  At a complete $2$-bridge structure, when $\omega$ is a hub of $F_\alpha$, the axis of $(S_\omega)_\frz$ projects to the closed geodesic freely homotopic to the corresponding crossing circle (Corollary~\ref{cor: geodesic crossing circles}). Figure~\ref{fig: Spiral1669} shows a typical example.  A related M\"obius orbit phenomenon appears in Section~10 of \cite{GF}, where Gu\'eritaud writes certain ideal vertices in a fan in the form $\cot(a+sb)$.  

\subsection{Reversed fractions, cusp fields, and Riley polynomials} 
Sections~\ref{sec: reverse}--\ref{sec: riley} use the explicit shape and holonomy formulas in three different ways.  First, for a positive-termed expansion $\alpha=[0;a_1,\ldots,a_k]$ with $k\geq2$ and $a_1,a_k\geq2$, Section~\ref{sec: reverse} compares the Sakuma--Weeks complexes for $\alpha$ and the reversed continued fraction $\bar\alpha=[0;a_k,\ldots,a_1]$.\index{reverse}  Theorem~\ref{thm: reverse} gives a simplicial bijection $F_\alpha\to F_{\bar\alpha}$ reversing the order of the interior edges, and hence a simplicial homeomorphism $M_\alpha\to M_{\bar\alpha}$.  Theorem~\ref{thm: phi OP} describes the corresponding relation between their geometric shape parameters.  When $\alpha$ is palindromic,\index{palindrome} this becomes a self-symmetry of $M_\alpha$.

These symmetries lead to useful special values of the shape functions. In particular, for the palindromes $\alpha_e$ singled out in Lemma~\ref{lem: fixed edge}, $\calZ(e)_{\frz(\alpha_e)}=-1$. Clearing denominators therefore produces a polynomial $p_e$ satisfied by the geometric root which can have dramatically smaller degree than
$\calQ(\alpha_e)$.  Thus the symmetry can yield unexpectedly low degree exact algebraic equations for the geometric root.  For instance, $\calQ(46368/121393)$ has expected degree $60696$, whereas the
corresponding $p_e$ has degree $376$.

Section~\ref{sec: cusps} applies the same reversal and palindromic symmetries to cusp fields.\index{cusp field}  Theorem~\ref{thm: cuspfield} gives an explicit formula for the cusp field of $M_\alpha$ in terms of the values of $\calV$ at the bottom of the funnel.  Corollary~\ref{cor: field equality} recovers the relevant twist knot and Hoste--Shanahan cases and proves equality of cusp and trace fields\index{trace field} for several new infinite families of 2-bridge links.

Finally, Section~\ref{sec: riley} identifies each $\calQ(\alpha)$ with the classical Riley polynomial $\Lambda(\alpha)$.\index{Riley polynomial}  Specifically, with the notation of that section, Theorem~\ref{thm: QandL} proves
\[ \Lambda(\alpha)_x=\pm x^{\mu(\alpha)}\overline{\calQ}(\alpha). \]
Consequently, the nonzero roots of $\Lambda(\alpha)$ are precisely $-1/\frz$, where $\frz$ ranges over the roots of $\calQ(\alpha)$, with multiplicity preserved.  In particular, the geometric root $\frz(\alpha)$ gives the discrete faithful Riley root $-1/\frz(\alpha)$.  This explains the inversion picture in Figure~\ref{fig: inversion}.

\subsection{Limits, augmented links, and Heckoid orbifolds}
Section~\ref{sec: limits} studies what happens when the number of crossings in a twist region tends to infinity.  Theorems~\ref{thm: convergence I} and~\ref{thm: convergence II} show that the corresponding geometric roots approach roots of the discriminant polynomials $D_\calQ(\omega)$.\index{discriminant polynomial}  Theorem~\ref{thm: DQ uniqueness} singles out a unique such root $\hat{\frz}(\omega)$\index{geometric root!discriminant polynomial} for each $\omega\in\Q\cap[0,1/2]$.  These distinguished roots are exactly the cusp points of $\bound\calR$ in the first quadrant.\index{Riley slice!cusp}

The remainder of Section~\ref{sec: limits} works out the cusp slopes $1/3$ and $1/4$ as concrete examples of this limiting geometry.  Theorem~\ref{thm: core} identifies the closure of the quotient of the corresponding infinite developed domain with the convex core \index{convex core} of $\IH^3/\Gamma_{\hat{\frz}}$.  Its boundary consists of two totally geodesic $3$-punctured spheres.  Theorems~\ref{thm: augment example} and~\ref{thm: geometric limits} glue these pieces to obtain explicit uniformizing groups for the associated augmented link complements\index{augmented link} and identify them as geometric limits\index{geometric limit} of $2$-bridge link groups. These examples show that the Farey recursive machinery remains useful even when the limiting group is no longer generated by a parabolic pair. 

Section~\ref{sec: Heck} concludes the paper with Heckoid examples\index{Heckoid orbifold} for $\alpha\in\{1/2,2/5\}$.  These examples show how their convex cores\index{convex core} are assembled from geometric Sakuma--Weeks pieces and hyperbolic drums\index{drum} (Definition~\ref{def: drum}). In the case $\alpha=1/2$, the drums are regular with all dihedral angles $\pi/2$, so the corresponding volumes can be computed explicitly using Lobachevsky functions.  

The forthcoming paper \cite{CEP} develops the broader theory suggested by the limiting and Heckoid examples.

\subsection{Computer code}

Much of the work in this paper is straightforward to implement with a computer.    Sections \ref{subsec: code I}, \ref{subsec: code II}, and \ref{subsec: code III} include examples of computer code that runs in \lstinline$SageMath 10.1$ \cite{sagemath}.

\newcommand\ButterThm{
Suppose $\alpha \in \Q \cap [0,1/2]$.  
\begin{enumerate}
\item If $\frz$ is a root for $\calQ(\alpha)$, then $\frz \in \C^\ast - \overline{\calR}$.  
\item If $\hat{\frz}$ is a root for $D_\calQ(\alpha)$, then $\hat{\frz} \in \C^\ast - \calR$.
\end{enumerate}
}

\newcommand\RealThm{Let $e\in\calE$ be a nonhorizontal edge and suppose that $z_L$, $z_R$, $z_{\hat e}$, and
$z_\oplus$ are roots of $\calQ(e^L)$, $\calQ(e^R)$, $\calQ(\hat e)$, and $\calQ(e^-\oplus e^+)$, respectively.  Then
\[ \calZ(e)_{z_L}=0,\qquad \calZ(e)_{z_R}=\infty,\qquad \calZ(e)_{z_{\hat e}}=1,\AND \calZ(e)_{z_\oplus}=1. \]
If, in addition, the vertices of $e$ lie in $[0,1/2]$, $e\ne[0,1/2]$, and $1/\alpha_e \notin \Z$, then
\[ \calZ(e)_{\frz(\alpha_e)} \ =\ -1\]
where $\frz(\alpha_e)$ is the geometric root of $\calQ(\alpha_e)$.
}

\newcommand\CuspCor{
Suppose $m \in \Z_{\geq 1}$ and $\alpha$ is 
\[ \textstyle \frac{2}{2m+3}, \ \frac{m+1}{2m+3}, \ \frac{3}{3m+1}, \ \frac{m}{3m+1}, \ \frac{3}{3m+2}, \ \frac{m+1}{3m+2}, \ \frac{2m+1}{4m+4}, \ \frac{3m+1}{9m+6}, \ \text{or} \quad \frac{3m+2}{9m+3}  \]
then the cusp and trace fields of $M_\alpha$ coincide.
}

\chapter{Farey recursion}\label{chap: FR}
\index{Farey recursion}

\section{The Stern--Brocot diagram}  \label{sec: SBD}
In this section, we review the construction discussed in Section 1.4 of \cite{C}.  For more background, see also \cite{Hat}.  Of particular importance to this paper are the rational numbers in the unit interval.  Set
\[
\Q_0 = \Q \cap [0,1] \AND
\widehat{\Q}_0 = \Q_0 \cup \{\infty\}.
\]
Elements of $\Q_0 - \{0,1\}$ are always written in lowest terms with positive denominators.
The remaining three elements of $\widehat{\Q}_0$ are regarded as the fractions $0/1$, $1/0$, and $1/1$.

\begin{dfn} \label{def: Farey pair} \index{Farey pair|idxdef} \index{Farey sum|idxdef} \index{Farey triple|idxdef} A pair $\big\{ \frac{p}{q}, \frac{r}{s} \big\} \subset \wh{\Q}_0$ is called a {\it Farey pair} if $ps-rq=\pm 1$.  Its {\it Farey sum} is
\[ \frac{p}{q} \oplus \frac{r}{s} = \frac{p+r}{q+s}.\]
\end{dfn}
Farey sums have many remarkable properties.  For instance, when writing the Farey sum of a Farey pair according to the above formula, the right hand side gives a fraction in lowest terms.  Also, the Farey sum of a Farey pair makes a Farey pair with each of its summands, thus forming a {\it Farey triple}.  Each Farey pair belongs to exactly two Farey triples.  Hence, if $\{ \gamma, \alpha\}$ is a Farey pair, it is possible to recursively form new Farey sums 
\begin{align*}
\gamma \oplus^0 \alpha &= \gamma \\
\gamma \oplus^j \alpha &= \left(\gamma \oplus^{j-1} \alpha \right) \oplus \alpha.
\end{align*}
where $j$ is any nonnegative integer.

\begin{dfn} \label{def: SternBrocot} Farey sums inspire an attractive $2$-complex $\calG \subset \R^2$ whose $2$-simplices are Euclidean triangles.  See also Figure 1.3 of \cite{HO}.  As in \cite{C}, the complex $\calG$ is called the {\it Stern--Brocot diagram}.  The vertex set $\calG^{(0)}$ is defined to be the set $\Q_0 \subset \R^2$, where $\Q_0$ is identified with a subset of $\R^2$ via the injection
\[ \frac{p}{q} \mapsto \left( \frac{p}{q}, \frac{1}{q}\right).\]
The edges of $\calG$ are defined by connecting each Farey pair with a linear edge.  The result is shown in Figure \ref{fig: G}.  The 2-simplices of $\calG$ are the Euclidean triangles determined by Farey triples.  From the figure, it is evident that $\calG$ embeds nicely in $\R^2$; this is proven by Hatcher in \cite{Hat}.  The diagram $\calG$ is also often regarded as an embedded 2-complex in $\R \times i\R_{>0}$.
\end{dfn}
\index{Stern--Brocot!diagram|idxdef}

\begin{figure}
\setlength{\unitlength}{.1in}
\begin{picture}(12,38)
\put(1,0) {\includegraphics[width= 1.2in]{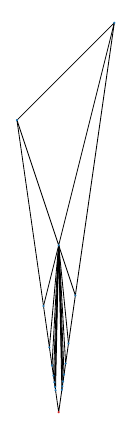}}
\put(1,28.5){$\frac{1}{4}$}
\put(2.8,11.5){$\frac{3}{11}$}
\put(8.8,12.3){$\frac{3}{10}$}
\put(5.8,20){$\frac{2}{7}$}
\put(12.2,37.5){$\frac{1}{3}$}
\end{picture}
\caption{The triangle $\Delta\left( \frac{2}{7} \right)$.}   \label{fig: Tri}
\end{figure}
  
\begin{dfn} \label{def: triangle} For $\alpha \in \Q_0$, the {\it triangle centered at} $\alpha$ is the subcomplex $\Delta(\alpha)$ of $\calG$ spanned by $\alpha$ together with the vertices which are Farey pairs with $\alpha$.  The triangle $\Delta(2/7)$ is shown in Figure \ref{fig: Tri}.  
\end{dfn}
\index{Stern--Brocot!triangle|idxdef}

Provided that $\alpha \notin \Z$, the triangle $\Delta(\alpha)$ is a Euclidean triangle with an ideal vertex at $\alpha$ on the $x$-axis and a pair of upper corner vertices $\kappa^L$ and $\kappa^R$ with $\kappa^L<\alpha<\kappa^R$.  If $q$ is the denominator of $\alpha$, the right side of $\Delta(\alpha)$ is the line segment of slope $q$ from the ideal vertex to the right corner $\kappa^R$.  The left side has slope $-q$ and connects the ideal vertex to $\kappa^L$.  

It is helpful to add an additional vertex for $\infty$ together with edges connecting this new vertex to its Farey partners $0$ and $1$.  The vertex $\infty$ becomes the corner for the right side of $\Delta(0)$ as well as the left side of $\Delta(1)$.  The two new edges are visualized as vertical rays emanating upwards from the vertices for $0$ and $1$.

\begin{dfns} \label{def: E and more} Let $\calE$ be the set of edges in $\calG$ whose vertices lie in $\Q_0$.  For $e \in \calE$, take 
\[\nu(e) \ =\ \begin{cases} -1 & \text{if the slope of } e \text{ is negative} \\
1 & \text{otherwise}. \end{cases}\]
The {\it opposite vertex} of $e$ is the number $\hat{e} \in \wh{\Q}_0$  with minimal denominator which makes a Farey triple with the endpoints of $e$.   Note that $e$ is an edge on a side of $\Delta(\hat{e})$.  If $e$ is not horizontal, its {\it high vertex} is the endpoint $e^-$ with smallest denominator and its other vertex $e^+$  is its {\it low vertex}.  Hence, $e^{\nu(e)}$ is always the leftmost vertex of $e$.  Write $e^L=e^{\nu(e)}$ and $e^R=e^{-\nu(e)}$.
\end{dfns}
\index{vertex!opposite|idxdef} \index{vertex!high|idxdef} \index{vertex!low|idxdef} \index{vertex!left|idxdef} \index{vertex!right|idxdef}

As explained in \cite{Hat}, there is a direct relationship between $\calG$ and the Farey graph.   To see this, allow the vertices of $\calG$ to slide straight down to the $x$-axis pulling their edges with them to form semicircles perpendicular to the $x$-axis.  The result is the Farey graph.

\section{Farey recursive functions} \label{sec: FareyRecursion}
The goal of this section is to review the notion of a Farey recursive function from Section 5.1 of \cite{C} and Section 3 of \cite{FRF}.

\begin{dfns} \label{def: FRF} Suppose $R$ is a commutative ring and that $d$ and $T$ are functions from $\wh{\Q}_0$ to $R$.  A function $\calF \co \wh{\Q}_0 \to R$ is a $(d, T)$-{\it Farey recursive function} (FRF) if, whenever $\{ \alpha, \gamma\}$ is a Farey pair, then
\begin{equation} \label{eq: defn} \calF(\gamma \oplus^2 \alpha) = -d(\alpha) \calF(\gamma) + T(\alpha) \calF(\alpha \oplus \gamma).\end{equation}

Said differently, if $\alpha \in \Q_0$, $\kappa$ is a corner for $\Delta(\alpha)$, and $\gamma_j = \kappa \oplus^j \alpha$ then $\{ \calF(\gamma_j) \}_0^\infty$ is linearly recursive with
\begin{equation} \label{eq: matrix}  \begin{pmatrix} 0 & 1 \\ -d(\alpha) & T(\alpha) \end{pmatrix} \, \begin{pmatrix} \calF(\gamma_{j-1}) \\ \calF(\gamma_j) \end{pmatrix} = \begin{pmatrix} \calF(\gamma_{j}) \\ \calF(\gamma_{j+1}) \end{pmatrix} .\end{equation}

 The $2\times 2$ matrix in \ref{eq: matrix} is called the {\it recursion matrix} at $\alpha$ for $\calF$ and the functions $d$ and $T$ are called the {\it determinant} and {\it trace} for $\calF$.   \label{def: matrix det trace} 
\index{Farey recursive function|idxdef}
\index{Farey recursive function!determinant|idxdef}
\index{Farey recursive function!trace|idxdef}
\end{dfns}

The expression \ref{eq: matrix} makes it clear that a linear combination of a pair of FRFs with the same determinant and trace is a FRF with the same determinant and trace.

It is often the case that the trace is taken to be $\calF$ itself.  The inductive definition continues to work when this is done because the denominator of $\alpha$ is smaller than that of $\gamma \oplus^2 \alpha$.  Excepting the next example, the ring $R$ will always be a subring of the field of rational functions $\Q(x)$.

\begin{example} \label{ex: Markoff}
In this example, the ring $R$ is the 3-variable polynomial ring $\Z[x,y,z]$.  Take $d \co \wh{\Q}_0 \to \Z[x,y,z]$ to be the constant function $\alpha \mapsto 1$.  Because the only value of $d$ is a unit in $\Z[x,y,z]$, Theorem 4.3 of \cite{FRF} shows that the assignments
\begin{align*} 
\calF\left(0\right)&=x & \calF\left(\infty\right)&=y & \calF\left(1\right)&=z
\end{align*}
extend uniquely to a Farey recursive function $\calF \co \wh{\Q}_0 \to \Z[x,y,z]$ with constant determinant $d=1$ and trace $\calF$.  In \cite{C}, this is referred to as  the {\em generic FRF}.  It is used there to describe the $\SLtwo{\C}$ and $\pslC$ character varieties for 2-bridge links.  

A specialization of the triple of variables $(x,y,z)$ to a point $\frp=(\frx, \fry, \frz) \in \C^3$ provides a specialization $\calF_\frp \co \wh{\Q}_0 \to \C$ which is an FRF with determinant $d=1$ and trace $\calF_\frp$.  If $\frp$ satisfies the Markoff equation
\[ \frx^2 + \fry^2 + \frz^2 \ =\ \frx \fry \frz\]
then the FRF $\calF_\frp$ is a {\em Markoff map} as defined in \cite{Bow}.  As shown by J\o rgensen in \cite{Jor} (explained also in 
\cite{Bow} and \cite{ASWY}), there is a natural bijection between the set of Markoff maps and the conjugacy classes of type-preserving representations of a once-punctured torus group into $\SLtwo{\C}$. 

This example shows that Farey recursive functions can be viewed as a generalization of Markoff maps.
\end{example}

The next three examples are central to this paper.

\begin{example} \label{ex: dQ}
The function
\[d_\calQ \co \wh{\Q}_0 \longrightarrow \Q(x)\] 
given by $d_\calQ(p/q)=(-1)^p x^q$ is a Farey recursive function with constant determinant $d=0$ and trace $T=\dQ$.  
\end{example}

Again, because the values of $d_\calQ$ are always units in $\Q(x)$, Theorem 4.3 of \cite{FRF} makes it easy to define FRFs with determinant $d_\calQ$.  In particular, if $T \co \wh{\Q}_0 \to \Q(x)$ is any function and  $a,b,c \in \Q(x)$ then there is a unique $(d_\calQ, T)$-Farey recursive function which takes the ordered triple $(0,\infty,1)$ to $(a,b,c)$.

\begin{example} \label{ex: Q}
Define the FRF
\[ \calQ \co \wh{\Q}_0 \longrightarrow \Q(x) \]
to be the unique FRF with determinant $d=d_\calQ$, trace $T=\calQ$, and
\begin{align*} 
\calQ\left(0\right)&=1 & \calQ\left(\infty\right)&=0 & \calQ\left(1\right)&=1.
\end{align*}
A primary goal of this paper is to show how $\calQ$ can be used to completely describe the geometry of all hyperbolic 2-bridge links as well as for many other Kleinian groups which are generated by a pair of parabolic elements in $\pslC$.

\index{Farey recursive function}
\end{example}

\begin{example} \label{ex: N}
Define the FRF
\[ \calN \co \wh{\Q}_0 \longrightarrow \Q(x) \]
to be the unique $(d_\calQ, \calQ)$-FRF with 
\[ \calN\left(0\right)=\calN\left(\infty\right)=\calN\left(1\right)=1.\]
\end{example}

Although the definitions for $\calQ$ and $\calN$ are remarkably simple, together they can be used to locate the ideal vertices of geometric fundamental domains for every hyperbolic 2-bridge link.  A list of polynomials $\calQ(\alpha)$ and $\calN(\alpha)$ where $\alpha$ has small denominators is given in Appendix \ref{chap: list}.

\begin{remark} \label{rem: constant term}
A quick inductive argument shows that, for every $\gamma \in \Q_0$,  the constant term in both $\calQ(\gamma)$ and $\calN(\gamma)$ is one.  
\end{remark}

In the next lemma, and the remainder of the paper, a subscript is used to indicate evaluation of a function.  The lemma seems interesting in its own right, but also becomes useful later in this paper. 

\begin{lemma} \label{lem: adjacent zeros}
For $e \in \calE$, $\calQ(e^-)$ and $\calQ(e^+)$ do not have a common zero. 
\end{lemma}

\proof
For a contradiction, assume that there is an edge $e$ and a number $\frz \in \C$ such that $\calQ(e^-)_\frz = \calQ(e^+)_\frz = 0$.

Choose $e$ among all such edges so that the denominator of $e^- \oplus e^+$ is as small as possible.  Since $\calQ=1$ on $\{ 0,1, 1/2 \}$, it must be true that the denominator of $e^-$ is at least three.  By Remark \ref{rem: constant term}, $\frz \neq 0$.  Hence, $\dQ(\beta)_\frz \neq 0$  for each $\beta \in \Q_0$.

Now, observe that $e^-$ must be a corner of $\Delta(\hat{e})$.  Otherwise, there is a number $\gamma \in \Q_0$ such that
\[ e^- = \gamma \oplus \hat{e} \qquad \text{and} \qquad e^+=\gamma \oplus^2 \hat{e}.\]
Then
\begin{align*}
0&= \calQ(e^+)_\frz \\
&= -\dQ(\hat{e})_\frz \calQ(\gamma)_\frz + \calQ(\hat{e})_\frz \calQ(e^-)_\frz \\
&=-\dQ(\hat{e})_\frz \calQ(\gamma)_\frz.
\end{align*}
Then $\calQ(\gamma)_\frz =0$ because $\dQ(\hat{e})_\frz \neq 0$.  But then, the edge with endpoints $\gamma$ and $e^-$ contradicts the minimality assumption for $e$.

Let $\kappa$ be the other corner of $\Delta(\hat{e})$.  Then
\begin{align*}
0&= \calQ(e^+)_\frz \\
&= -\dQ(e^-)_\frz \calQ(\kappa)_\frz + \calQ(e^-)_\frz \calQ(\hat{e})_\frz\\
&= -\dQ(e^-)_\frz \calQ(\kappa)_\frz
\end{align*}
and so, $\calQ(\kappa)_\frz =0$.  The edge with endpoints $\kappa$ and $e^-$ contradicts the minimality assumption for $e$ and this contradiction completes the proof.
\endproof

On the other hand, there are many pairs of numbers $\gamma_0, \gamma_1 \in \Q_0$ such that $\calQ(\gamma_0)$ and $\calQ(\gamma_1)$ have a common zero.  For instance, $\calQ(1/3)$ and $\calQ(1/6)$ are both zero at $x=1$.  Moreover, a similar version of Lemma \ref{lem: adjacent zeros} does not hold for $\calN$.  If $e$ is the edge with endpoints $1/3$ and $2/7$ then $x=1/2$ is a solution to $\calN(e^-)=0$ and $\calN(e^+)=0$.  (See Appendix \ref{chap: list}.)

\begin{remark} \label{rk: chebyshev}
The relationship between the polynomials $\calQ(1/n)$ and the Chebyshev polynomials of the second and fourth kind are discussed briefly in Example 5.6 of \cite{FRF}.  The Chebyshev polynomials of the second kind $\{ U_j \} \subset \Z[x]$ are determined by the second-order linear recursion relation $U_{j+1}=-U_{j-1}+2x\, U_j$ and initial conditions $U_0=1$ and $U_1=2x$.  The Chebyshev polynomials of the fourth kind are given by $W_j=U_j+U_{j-1}$.  In \cite{FRF}, it is mentioned that, for $n \in \Z_{\geq 1}$, 
\[
\calQ\left( \frac{1}{2n} \right)  =  x^{n-1} \,U_{n-1} \left( \frac{1-2x}{2x} \right) \quad \text{and} \quad \calQ\left(\frac{1}{2n+1}\right)  =  x^{n} \, W_n \left( \frac{1-2x}{2x} \right).
\]
A bit of trigonometry together with the identities 
\[ U_{j-1} (\cos \theta) \ = \ \frac{\sin(j \theta)}{\sin \theta} \qquad \text{and} \qquad W_j(\cos \theta) \ =\ \frac{\sin\left( \left( j+\frac12 \right) \theta \right) }{\sin \left(\frac{\theta}{2}\right)}
\]
shows that 
\[ \left\{ \frac14  \sec^2 \left( k\pi/n \right) \, \Big| \, k \in \Z \cap \left[1, (n-1)/2 \right] \right\} \]
is the set of roots of $\calQ(1/n)$.  In particular, every root of $\calQ(1/n)$ is real.
\end{remark}

\section{Computer Code} \label{subsec: code I}

The material from Sections \ref{sec: SBD} and \ref{sec: FareyRecursion} is straightforward to implement with a computer, as are many of the results in this paper.  Accordingly, some examples of computer code which work in \lstinline$SageMath 10.1$ \cite{sagemath} are included.  To keep things brief and transparent, the included code is not optimized for efficiency;  in practice it is a good idea to edit the given code to avoid redundant calculations as well as to manage precision and rounding errors. 

Given $\gamma \in \Q_0$, the block of code below can be used to find a triangle $\Delta(\alpha)$ with a corner $\kappa$ such that $\gamma = \kappa \oplus^j \alpha$ and $j\geq 2$.  The first two commands provide the maps $p/q \mapsto (p,q)$ and $(p,q) \mapsto p/q$.  The third command returns $(\alpha, \kappa, j)$ as desired.  This code utilizes the well-understood connection between $\calG$ and continued fractions \cite{Hat}.

\begin{lstlisting}[language=Python]
def vect(r):   # r is a rational number.
    return vector([r.numerator(), r.denominator()])

def quot(v):   # v is a vector.
    return QQ(v[0]/v[1])

def center_corner_index(v):   # v is a vector.
    r=quot(v)
    CF_center=list(r.continued_fraction())[:-1]
    center=vect(continued_fraction(CF_center).value())
    k2=vect(continued_fraction(CF_center+[1]).value())
    corner=k2-center
    ind=Integer((v[1]-corner[1])/center[1])
    return [center, corner, ind]
\end{lstlisting}

The above commands are used in the following code which computes $d_\calQ$, $\calQ$, and $\calN$.  Here, rational numbers $p/q$ are represented as vectors $v=(p,q)$.  It should be clear how to modify this code to compute other FRFs which have unit determinants.  

\begin{lstlisting}[language=Python]
x=var('x')

dQ = lambda v: (-1)^v[0]*x^v[1]

def Q(v):
    v = vector(v)
    if v == vector([0,1]): return 1
    elif v == vector([1,0]): return 0
    elif v == vector([1,1]): return 1
    else:
        (c,k,j) = center_corner_index(v)
        return (matrix([[0,1],[-dQ(c), Q(c)]])^j * 
        	vector([Q(k),Q(k+c)]))[0]    
    
def N(v):
    v = vector(v)
    if v == vector([0,1]): return 1
    elif v == vector([1,0]): return 1
    elif v == vector([1,1]): return 1
    else:
        (c,k,j) = center_corner_index(v)
        return (matrix([[0,1],[-dQ(c), Q(c)]])^j * 
        	vector([N(k),N(k+c)]))[0] 
\end{lstlisting}

Figure \ref{fig: allroots} was computed using a variant of this code.  It shows the complex roots (near zero) of the polynomials obtained by evaluating $\calQ$ at the points $p/q \in \Q_0$ with $q<250$.  This attractive image suggests a complicated structure and seems worthy of deeper study.
\chapter{Geometric triangulations}\label{chap: SakumaWeeks}

\begin{figure}[h] 
   \centering
   \includegraphics[width=2.5in]{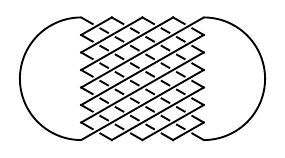} 
   \caption{The link $L_{4/7}$.  }  \label{fig: 47knot}
\end{figure}

\begin{dfn} \label{def: 2 bridge link} Fix $\alpha \in \Q \cup \left\{ \infty \right\}$.  The {\it 2-bridge link} $L_\alpha \subset S^3$ is obtained by connecting the corners of a square pillowcase with arcs of slope $\alpha$ as indicated in Figure \ref{fig: 47knot}.   Let $M_\alpha=S^3-L_\alpha$ denote the corresponding link complement. 
\end{dfn}
\index{2-bridge link|idxdef}
\index{2-bridge link!complement|idxdef}

These links became important examples after Schubert's thorough investigation in \cite{Sch}.  Although there is a 2-bridge link for every such $\alpha$, Schubert's work shows that any given 2-bridge link arises as $L_\alpha$ for $\alpha \in \left\{ \infty \right\} \cup \left( \Q \cap [0,1/2]\right)$.  

The links for $\alpha \in \{\infty, 0, 1/2\}$ correspond to the trivial 2-component link, the trivial knot, and the Hopf link, respectively. Since these are well-understood, this paper will primarily focus on the cases where $\alpha \in \Q \cap (0,1/2)$.  Menasco's work in \cite{Me} shows that, if $\alpha \in \Q \cap (0,1/2)$, then $M_\alpha$ admits a complete hyperbolic structure with finite volume if and only if $1/\alpha\notin\Z$.  

\begin{dfn}\label{def: hyperbolic alpha}
A rational number $\alpha\in\Q\cap(0,1/2)$ is called {\it hyperbolic} if $1/\alpha\notin\Z$.
Equivalently, $\alpha$ is hyperbolic precisely when $M_\alpha$ is a hyperbolic $2$-bridge link complement.
\end{dfn} 
\index{hyperbolic number|idxdef}

By Mostow--Prasad rigidity, if $\alpha$ is hyperbolic, then the complete finite volume hyperbolic structure on $M_\alpha$ is uniquely determined by $\alpha$.

For excellent background on 2-bridge links, see, for example, \cite{BZ} and \cite{Pur}.
 
\section{Sakuma--Weeks triangulations} \label{sec: SW}
Following \cite{GF} or \cite{SW}, the Stern--Brocot diagram $\calG$, from Definition \ref{def: SternBrocot}, can be used to construct the Sakuma--Weeks complexes associated to $2$-bridge links.  When $\alpha$ is hyperbolic, these are ideal triangulations of the link complements $M_\alpha$.  This first section mainly introduces definitions and notation.
\index{Sakuma--Weeks complex}
\index{Stern--Brocot!diagram}

\begin{dfns} \label{def: funnel} Assume $\alpha \in \Q \cap (0,1/2)$. Consider the 2-subcomplex of $\calG$ consisting of all simplices which meet the ray $\{ \alpha+it \, \mid \, t \in \R^+ \}$.  The {\it funnel} $F_\alpha$ is obtained from this complex by deleting all simplices which meet the vertex $\alpha$.  

The triangular simplices of $F_\alpha$ are denoted $s_j$, indexed as they are passed  falling along the ray $\{ \alpha+it \, \mid \, t \in \R^+ \}$ from above.  This sequence starts with $s_0$ and ends with $s_n$.  For $1\leq j\leq n$, the edge shared by $s_{j-1}$ and $s_j$ is denoted $e_j$.  The top edge of $s_0$ is denoted $e_0$ and the bottom edge of $s_n$ is denoted $e_{n+1}$. The edges $\{e_j\}_1^n$ are the {\it interior edges} of $F_\alpha$.
\end{dfns}
\index{funnel|idxdef} \index{funnel!interior edges|idxdef}

\begin{figure}
\setlength{\unitlength}{.1in}
\begin{picture}(50,25)
\put(0,0) {\includegraphics[width= 5in]{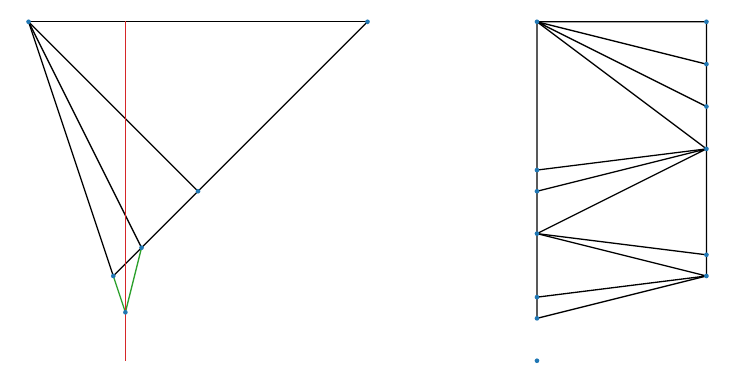}}
\put(0.4,24.2){$\frac{0}{1}$}
\put(25.6,24.2){$\frac{1}{1}$}
\put(9.2,20){\redden{$\frac{2}{7}+it$}}
\put(6,7){$\frac{1}{4}$}
\put(6.8,4){$\frac{2}{7}$}
\put(10.4,8.5){$\frac{1}{3}$}
\put(14.3,12.0){$\frac{1}{2}$}
\put(48.6,24.2){$\frac{1}{1}$}
\put(48.6,21.4){$\frac{1}{2}$}
\put(48.6,18.5){$\frac{1}{3}$}
\put(48.6,15.6){$\frac{1}{4}$}
\put(48.6,8.6){$\frac{4}{17}$}
\put(48.6,6.6){$\frac{7}{30}$}
\put(35.1,24.2){$\frac{0}{1}$}
\put(35.1,14.4){$\frac{1}{5}$}
\put(35.1,12.4){$\frac{2}{9}$}
\put(34.6,9.8){$\frac{3}{13}$}
\put(34.6,5.7){$\frac{10}{43}$}
\put(34.6,3.7){$\frac{17}{73}$}
\put(34.0,1.1){$\frac{24}{103}$}
\end{picture}
\caption{The funnels for $2/7$ and $24/103$.  Because the 2-simplices of the second funnel are extremely small, a labeled simplicial complex combinatorially isomorphic to the geometric funnel is shown instead.}  \label{fig: Funnels}
\end{figure}

The funnels for $2/7$ and $24/103$  are shown  in Figure \ref{fig: Funnels}.  If $\alpha=1/3$, the funnel is a triangle with no interior edges.  It is always true that $e_0=[0,1]$ and $e_1=[0,1/2]$.  

Now, assume that $\alpha \in \Q \cap (0,1/2)$ and $\alpha \neq 1/3$ so that $F_\alpha$ always has at least one interior edge.

\begin{dfns} \label{def: tetrahedral complex}
Let $\{ s_j\}_0^n$ and $\{ e_j\}_0^{n+1}$ be the triangles and edges of the funnel $F_\alpha$.  For each $s_j$,  associate a copy $\X_j$ of the punctured plane $\X=\C-\Z[i]$.  Let $\Theta$ be the group of Euclidean isometries of $\X$ generated by the order-2 rotations centered at the punctures of $\X$.  Each plane $\X_j$ has an ideal triangulation obtained by drawing all possible lines through the punctures of $\X_j$ whose slopes coincide with the rational numbers at the vertices of $s_j$.  When these triangulated planes are stacked top to bottom and identical lines are identified in adjacent planes, the result is an infinite complex $\wh K$. (The construction here is upside down from that given on p. 208 of \cite{SW}.)
\end{dfns}

By construction, the tetrahedra between $\X_{j-1}$ and $\X_j$ correspond to the edge $e_j$.  The edges in the triangulation of the top boundary $\X_0$ of $\wh{K}$ are labeled with the slopes $0$, $1$, and $\frac12$.  The edges in the bottom $\X_{n}$ of $\wh{K}$ are labeled with the numbers at the vertices of $s_n$.   The group $\Theta$ acts naturally on $\wh{K}$ as a covering space action preserving the labeled triangulation of each $\X_j$.

\begin{dfns} \label{def: quotient complex}
The quotient $K=\wh{K}/\Theta$ is an ideal simplicial complex whose edges are labeled by the rational slopes of their preimages in $\cup_j \X_j$.  The image $\Sigma_j$ of $\X_j$ in $K$ is homeomorphic to a 2-sphere with four punctures.  For $1\leq j\leq n$, the tetrahedra in $\wh{K}$ corresponding to the edge $e_j$ project to a pair of identically labeled ideal tetrahedra, denoted $\Delta_j$.  The vertex  $\hat{e}_{n+1}$ of $F_\alpha$ is called a {\it hairpin vertex}.  Similarly, at the top, $\frac12$ is a hairpin vertex.  The edges in $K$ which have hairpin labels are also called hairpins.  
\end{dfns}
\index{hairpin}
\index{hairpin!vertex|idxdef}

\begin{dfn} \label{def: SW}
Assume again that $\alpha \in \Q \cap (0,1/2)$ and $\alpha \neq 1/3$.  The {\it upper hairpin fold} identifies the ideal triangles in $\Sigma_0$ across the hairpins labeled $1/2$.  (See p. 2077 \cite{McShane1} for a more detailed description of this fold.) Let $M_\alpha^\circ$ denote the resulting ideal complex.

If $\alpha$ is hyperbolic (recall Definition~\ref{def: hyperbolic alpha}), then $M_\alpha^\circ$ is a manifold homeomorphic to the complement of a trivial tangle in a ball.  If $1/\alpha\in\Z$, then $M_\alpha^\circ$ contains at least one tetrahedral pair, but is not a manifold.  In this case, the same edges labeled $0$ occur in both $\Sigma_0$ and $\Sigma_n$, so points on these edges do not have $3$-ball neighborhoods.  Compare the construction preceding Theorem~II.2.4 and Remark~II.2.7 of \cite{SW}.

Now suppose that $\alpha$ is hyperbolic.  As proven in \cite{SW}, if the ideal triangles at the bottom of $M_\alpha^\circ$ are identified across the lower hairpins, the result is an ideal triangulation of the link complement $M_\alpha$.  Lemma II.2.5 of \cite{SW} shows that performing the upper and lower hairpin folds is topologically equivalent to attaching $2$-handles to $K$ along simple closed curves of slopes $\infty$ and $\alpha$, respectively.  The resulting ideal triangulation is referred to as the {\it Sakuma--Weeks triangulation} of $M_\alpha$.
\index{hairpin!fold|idxdef}
\index{Sakuma--Weeks complex!triangulation|idxdef}

An edge which is identified to another under a hairpin fold is called a {\it core tunnel}. \label{def: tunnel} The edges of $M_\alpha^\circ$ inherit rational labels from $K$.  When $\alpha$ is hyperbolic, the edges of the Sakuma--Weeks triangulation of $M_\alpha$ inherit the same labels.  An edge has a single rational label unless it is a core tunnel, in which case it is labeled by a pair of rational numbers.  In particular, the upper core tunnel of $M_\alpha^\circ$ is labeled by $(0,1)$ and, when $\alpha$ is hyperbolic, the lower core tunnel of $M_\alpha$ is labeled by $\big(e_{n+1}^-,e_{n+1}^+\big)$.
\index{core tunnel|idxdef}
\end{dfn}

\begin{remark} \label{rem: Sigma}
As before, assume that $\alpha\in\Q\cap(0,1/2)$ and $\alpha\neq1/3$. Given $j\in\{1,\ldots,n\}$, there is a deformation retraction from $K$ to $\Sigma_j$.  Hence, the inclusion induces a surjection
\[ \pie(\Sigma_j)\longrightarrow\pie(M_\alpha^\circ). \]
By the description of the upper hairpin fold in Lemma~II.2.5 of \cite{SW}, $M_\alpha^\circ$ is homotopy equivalent to a $4$-holed sphere with a disk attached along a simple loop of slope $\infty$.  Thus, $M_\alpha^\circ$ has the homotopy type of a genus-$2$ handlebody and $\pie(M_\alpha^\circ)$ is free of rank two.

When $\alpha$ is hyperbolic, the inclusion also induces a surjection
\[ \pie(\Sigma_j)\longrightarrow\pie(M_\alpha)\]
and $M_\alpha$ is homotopy equivalent to a $4$-holed sphere with disks attached along simple loops of slopes $\infty$ and $\alpha$.
\end{remark}
\index{Sakuma--Weeks complex!4-punctured spheres}

\section{4-plat link diagrams} \label{sec: diagrams}

Throughout this section, assume that $\alpha\in\Q\cap(0,1/2)$ and $\alpha\neq1/3$, so that $F_\alpha$ has at least one interior edge.

\begin{dfns} \label{def: hubs} Suppose that $e_0, \ldots, e_{n+1}$ are the edges of the funnel $F_\alpha$ that meet its defining ray.  If $j \in \{ 0, \ldots, n+1\}$, then $\omega = e_j^-$ is called a {\it hub} for $F_\alpha$.  A {\it left hub} is a hub on the left side of $F_\alpha$ and a {\it right hub} is one on the right.   If $\omega$ is a hub for $F_\alpha$ and $e_j^-=\omega$, then $e_j$ is called a {\it spoke} for $\omega$.  The vertex $0$ is always the first hub of $F_\alpha$ because, for every $\alpha$, $e_1$ is the edge between $0$ and $1/2$.  The horizontal edge between $0$ and $1$ is considered to be a spoke for $0$.  The number of spokes on a hub is called its {\it index}. 
\end{dfns}
\index{funnel!hub|idxdef}
\index{funnel!spoke|idxdef}

\begin{figure}
\setlength{\unitlength}{.1in}
\begin{picture}(48,20)
\put(0,1) {\includegraphics[width= 4.8in]{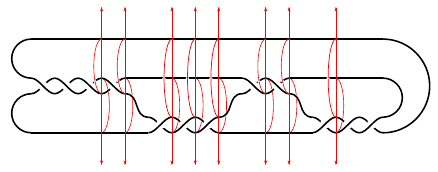}}
\put(10.5,0){$\Sigma_1$}
\put(9.6,14.3){$\scriptscriptstyle \frac12$}
\put(11.8,6.1){$\scriptscriptstyle \frac12$}
\put(11.2,12.8){$\scriptscriptstyle \frac13$}
\put(11.2,17){$\scriptscriptstyle 0$}
\put(10.3,7.4){$\scriptscriptstyle \frac13$}
\put(13.1,0){$\Sigma_2$}
\put(12.2,14.3){$\scriptscriptstyle \frac13$}
\put(14.4,6.1){$\scriptscriptstyle \frac13$}
\put(13.8,12.8){$\scriptscriptstyle \frac14$}
\put(13.8,17){$\scriptscriptstyle 0$}
\put(12.9,7.4){$\scriptscriptstyle \frac14$}
\put(18.2,0){$\Sigma_3$}
\put(17.6,14.3){$\scriptscriptstyle 0$}
\put(19.8,7.3){$\scriptscriptstyle 0$}
\put(18.9,12.8){$\scriptscriptstyle \frac14$}
\put(18.9,17){$\scriptscriptstyle \frac15$}
\put(20.8,0){$\Sigma_4$}
\put(20,14.3){$\scriptscriptstyle \frac15$}
\put(22.2,7.3){$\scriptscriptstyle \frac15$}
\put(21.5,12.8){$\scriptscriptstyle \frac14$}
\put(21.5,17){$\scriptscriptstyle \frac29$}
\put(23.4,0){$\Sigma_5$}
\put(22.6,14.3){$\scriptscriptstyle \frac29$}
\put(24.8,6.1){$\scriptscriptstyle \frac29$}
\put(24.1,12.8){$\scriptscriptstyle \frac14$}
\put(24.1,17){$\scriptscriptstyle \frac{3}{13}$}
\put(28.5,0){$\Sigma_6$}
\put(27.6,14.3){$\scriptscriptstyle \frac14$}
\put(29.8,6.1){$\scriptscriptstyle \frac14$}
\put(29.1,12.8){$\scriptscriptstyle \frac{4}{17}$}
\put(29.1,17){$\scriptscriptstyle \frac{3}{13}$}
\put(27.9,7.4){$\scriptscriptstyle \frac{4}{17}$}
\put(31.1,0){$\Sigma_7$}
\put(29.8,14.3){$\scriptscriptstyle \frac{4}{17}$}
\put(32.3,6.1){$\scriptscriptstyle \frac{4}{17}$}
\put(31.7,12.8){$\scriptscriptstyle \frac{7}{30}$}
\put(31.7,17){$\scriptscriptstyle \frac{3}{13}$}
\put(30.5,7.4){$\scriptscriptstyle \frac{7}{30}$}
\put(36.2,0){$\Sigma_8$}
\put(34.9,14.3){$\scriptscriptstyle \frac{3}{13}$}
\put(37.6,7.3){$\scriptscriptstyle \frac{3}{13}$}
\put(36.8,12.8){$\scriptscriptstyle \frac{7}{30}$}
\put(36.8,17){$\scriptscriptstyle \frac{10}{43}$}
\end{picture}
\caption{This figure is adapted from the first part of Figure II.3.3 of \cite{SW}.  Together with a standard 4-plat diagram $D_{24/103}$ for the knot $L_{24/103}$ it shows the edges in the Sakuma--Weeks triangulation.  Note that many edges are represented more than once in this diagram. This is done to make it easier to see the triangulations for the 4-punctured spheres $\Sigma_j$ in the knot complement.}  \label{fig: SW link}
\end{figure}

\begin{example} \label{ex: 24/103} If $\alpha = 24/103$ then $F_\alpha$ is pictured in Figure \ref{fig: Funnels} and $F_\alpha$ has four hubs, namely $0$, $1/4$, $3/13$, and $7/30$.  The spokes for $0$ are $e_0$, $e_1$, $e_2$, and $e_3$.  The next three edges are spokes for $1/4$ followed by the two spokes for $3/13$. The last two edges $e_9$ and $e_{10}$ are spokes for $7/30$.  This picture is closely related to the alternating 4-plat diagram for $L_\alpha$ shown in Figure \ref{fig: SW link}.  (This figure is adapted from the first part of Figure II.3.3 of \cite{SW}.)  Just as $F_\alpha$ has four hubs, this diagram has four twist regions.  With the exception of the rightmost twist region, the number of crossings in the twist region is equal to the index of the corresponding hub.  The rightmost twist region has one more crossing than the index of the last hub.  As in \cite{SW}, the figure is drawn so that it is easy to see the edges of the triangulations of the 4-punctured spheres $\Sigma_j$ in $M_\alpha$.
\end{example}

It is well-known that this picture holds in general.  

\begin{dfn} \label{def: diagram} If $F_\alpha$ has $k$ hubs then, as in Figure \ref{fig: SW link}, $L_\alpha$ has a corresponding alternating {\it 4-plat diagram} $D_\alpha$ with $k$ twist regions.
\end{dfn}
\index{4-plat diagram|idxdef}

 Let $a_j$ denote the number of crossings in its $j^\text{th}$ (from the left) twist region.  Let $\omega_1, \ldots, \omega_k$ be the sequence of hubs of $F_\alpha$ listed in order of increasing denominator.  Then, the index of $\omega_j$ is $a_j-1$ when $j$ is $1$ or $k$, and $a_j$ when $2 \leq j < k$.  This can be verified by combining Theorem 2.1 from \cite{Hat} with Proposition 12.13 from \cite{BZ}.  It follows that, in this situation, $a_1$ and $a_k$ are both at least two.  It is also well known (see, for example, \cite{Hat}) that $[0; \,a_1, \ldots, a_k]$ is a continued fraction expansion for $\alpha$.  The alternating diagram $D_\alpha$ is sometimes denoted as $D(a_1, \ldots, a_k)$.

\section{Gluing equations} \label{sec: equations} 
\index{gluing equation}

As in Section II.5 of \cite{SW}, this section systematically describes Thurston's gluing equations \cite{Th_notes} associated to the Sakuma--Weeks complexes.  Throughout the section, assume that $\alpha\in\Q\cap(0,1/2)$ and $\alpha\neq1/3$, so that $F_\alpha$ has at least one interior edge and $M_\alpha^\circ$ contains at least one tetrahedral pair. The gluing equations for $M_\alpha^\circ$ are defined under these assumptions.  Whenever the Sakuma--Weeks triangulation of $M_\alpha$ is considered,  always assume that $\alpha$ is hyperbolic as defined in Definition~\ref{def: hyperbolic alpha}.

\begin{dfns} \label{def: Klein}
Let $\wh{\Theta}$  be the group of Euclidean isometries of the punctured plane $\X=\C-\Z[i]$ generated by the order-2 rotations centered at the points of $\frac12 \, \Z[i]$.  The group $\Theta$ from Definition \ref{def: tetrahedral complex} is an index-4 normal subgroup of $\wh{\Theta}$ and the quotient $G=\wh{\Theta}/\Theta$  is a Klein 4-group.
\end{dfns}

The group $G$ acts simplicially on $K$, preserving its edge labeling, so this action descends to a simplicial action on $M_\alpha^\circ$.  When $\alpha$ is hyperbolic, it also descends to the Sakuma--Weeks triangulation of $M_\alpha$.

\begin{dfn} \label{def: orbifolds}
Define the quotient ideal complex
\[ O_\alpha^\circ\ =\ M_\alpha^\circ/G\]
and, when $\alpha$ is hyperbolic, also define
\[O_\alpha\ =\ M_\alpha/G.\]
In the hyperbolic case, $O_\alpha^\circ$ and $O_\alpha$ are orbifolds.  All of these complexes inherit well-defined edge labels from $M_\alpha$ and $M_\alpha^\circ$.
\end{dfn}

To be clear, this section computes the gluing equations on the quotient complex $O_\alpha^\circ$ and, when $\alpha$ is hyperbolic, on the orbifold $O_\alpha$.  Equivalently, these are the $G$-invariant gluing equations for $M_\alpha^\circ$ and, in the hyperbolic case, for $M_\alpha$.

To better understand the simplicial action of $G$ on these complexes, it is helpful to further discuss the anatomy of the tetrahedral pair $\Delta_j$.  

\begin{dfns} \label{def: edge types} Let $e_j$ be an interior edge of a funnel $F_\alpha$ as described in Definition \ref{def: funnel}.  Following Definition \ref{def: quotient complex}, take $\Delta_j$ to be the pair of tetrahedra in $K$ associated to $e_j$.  The edges of $\Delta_j$ which are labeled with the numbers $e^\pm_j$ at the endpoints of $e_j$ are called the {\it side edges} of $\Delta_j$.  (Notice that the side edges of $\Delta_j$ are shared by $\Delta_{j\pm 1}$.)  The {\it top edges} of $\Delta_j$ are those which are labeled by $\hat{e}_j$ and the {\it bottom edges} are labeled by $e^-_j \oplus e^+_j$.  
\end{dfns}
\index{Sakuma--Weeks complex!tetrahedral pair|idxdef}
\index{Sakuma--Weeks complex!tetrahedral pair!side edge|idxdef}
\index{Sakuma--Weeks complex!tetrahedral pair!top edge|idxdef}
\index{Sakuma--Weeks complex!tetrahedral pair!bottom edge|idxdef}

\begin{figure}
\setlength{\unitlength}{.1in}
\begin{picture}(25,14)
\put(0,0) {\includegraphics[width= 2.5in]{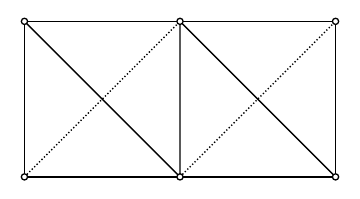}}
\put(6.5,0){$\scriptstyle \frac{1}{5}$}
\put(6.5,13.2){$\scriptstyle \frac{1}{5}$}
\put(17.5,0){$\scriptstyle \frac{1}{5}$}
\put(17.5,13.2){$\scriptstyle \frac{1}{5}$}
\put(0.5,6.8){$\scriptstyle \frac{1}{4}$}
\put(12.9,6.8){$\scriptstyle \frac{1}{4}$}
\put(23.7,6.8){$\scriptstyle \frac{1}{4}$}
\put(4.4,10.2){$\scriptstyle \frac{0}{1}$}
\put(15.3,10.2){$\scriptstyle \frac{0}{1}$}
\put(4.4,3.2){$\scriptstyle \frac{2}{9}$}
\put(15.3,3.2){$\scriptstyle \frac{2}{9}$}
\end{picture}
\caption{The pair of tetrahedra $\Delta_4$ for $\frac{24}{103}$.  The top edges are labelled $0$ and the bottom edges are labelled $\frac29$.}  \label{fig: Delta}
\end{figure}

Not only does this terminology reflect the combinatorics of the link complements $M_\alpha$, but also that of the funnels $F_\alpha$. Compare Figure \ref{fig: Funnels}.  Observe also that the opposite edge to a side edge of a tetrahedron in $\Delta_j$ is always a side edge in $\Delta_j$ with the same label, while the edge opposite a top edge is a bottom edge.  Figure \ref{fig: Delta} shows $\Delta_4$ when $\alpha= \frac{24}{103}$.

Consider a lift to $\wh{K}$ of a side edge of $\Delta_j$.  This is a line segment in $\X_j$ connecting a pair of punctures.  The rotation in $\wh{\Theta}$ which fixes the center of this segment represents a nontrivial element of $G$ which interchanges the tetrahedra in $\Delta_j$. Likewise, an element of $G$ which corresponds to a top or bottom edge does not interchange the tetrahedra in $\Delta_j$.  The next fact follows since every tetrahedron in $K$ has a side edge and a top edge. 

\begin{fact} \label{fact: isometry}
Given any $j \in \{ 1, \ldots, n\}$, $G$ contains an element which interchanges the two tetrahedra in $\Delta_j$.  $G$ contains another element which does not interchange them, but acts as a nontrivial involution on each.
\end{fact}

Every tetrahedron in $\Delta_j$ has six distinct edges in $K$; two pairs of side edges, a top edge, and a bottom edge.  If $g \in G$ corresponds to a rotation of the top or bottom edge, it cannot fix any side edge of this tetrahedron under its action on $K$.  In particular, if $E$ is an edge of $K$, then the only nontrivial element of $G$ which preserves $E$ setwise, acts on $E$ as a nontrivial involution.  The same is true for an edge of $M_\alpha^\circ$ which is not a core tunnel (Definition \ref{def: tunnel}) and, when $\alpha$ is hyperbolic, for an edge of $M_\alpha$ which is not a core tunnel. \index{core tunnel}

The situation is different for a core tunnel.  Under a hairpin fold, four edges of $K$ are identified to form the core tunnel, and there is
a nontrivial element of $G$ which fixes the resulting tunnel pointwise.  Since this element has order two and does not fix the
interior of any tetrahedron, it acts locally as a half turn about the tunnel.  Consequently, the edge cycle in $O_\alpha^\circ$ associated to a core tunnel is the quotient of the corresponding edge cycle by this order-$2$ rotational symmetry.  The same statement applies to $O_\alpha$ when $\alpha$ is hyperbolic.  

\begin{dfns} \label{def: involutions} The {\it upper involution} of $M_\alpha^\circ$ is the element of $G$ which fixes the points of the core tunnel labeled $(0,1)$ and acts by inversion on each of the hairpin edges labeled $1/2$.  When $\alpha$ is hyperbolic, the same involution acts on $M_\alpha$, which also has a {\it lower involution} in $G$ which fixes the points of the core tunnel labeled $\big(e_{n+1}^-, e_{n+1}^+\big)$ and acts by inversion on each of the hairpin edges labeled $\hat{e}_{n+1}$.   See Section 4 of \cite{ALSS} for illustrations of these involutions. 
\end{dfns}
\index{involution|idxdef}

For a pair $\{ v_1, v_2\}$ of distinct points in $\bound \IH^3$, let $[v_1, v_2]$ denote the geodesic whose ideal endpoints coincide with this pair.  Likewise, for a set of distinct points $\{ v_1, \ldots, v_k\} \subset \bound \IH^3$, let $[v_1, \ldots, v_k]$ denote the convex hull in $\IH^3$ of the set $\{[v_i, v_j] \, | \, 1 \leq i < j \leq k \}$.  

\begin{dfns} \label{def: parameters}
Suppose that $T=[v_0, \ldots, v_3]$ is an ideal tetrahedron and choose an edge $E=[v_i, v_j]$ of this tetrahedron.  Following \cite{Th_notes}, the edge $E$ of $T$ has a well-defined invariant $\frz \in \IH^2$ called the {\it edge parameter} and obtained by moving $T$ in $\IH^3$ by an isometry into a position $[0,1,\infty, \frz]$ where $[0,\infty]$ is the image of $E$.  
\index{shape parameter!edge parameter|idxdef}

If a tetrahedron $\delta$ in $\Delta_j$ is endowed with the structure of an ideal hyperbolic tetrahedron, this paper refers to the edge parameter $\frz$ of the top edge of $\delta$ as the {\it shape parameter} of $\delta$.  From Chapter 4 of \cite{Th_notes}, it is known that the edge parameter of the bottom edge of $\delta$ is also $\frz$.  It also follows from Thurston's work here that the edge parameters for the left and right edges of $\delta$ are $\zeta(\frz)$ and  $\zeta^2(\frz)$, where $\zeta$ is the M\"{o}bius transformation 
\[  \zeta(x) = \frac1{1-x}.\]
\end{dfns}
\index{shape parameter}

Suppose that $\alpha$ is hyperbolic.  By \cite{GF}, the Sakuma--Weeks triangulation of $M_\alpha$ is canonical in the sense of Epstein--Penner \cite{EP}.  In particular, it is geometric and each tetrahedron in the triangulation has a well-defined shape parameter in $\IH^2$.  By Mostow--Prasad rigidity, the simplicial action of $G$ is realized by isometries, so the two tetrahedra in each pair $\Delta_j$ are isometric and have the same shape parameter.
\index{geometric triangulation}

For any $\alpha$ satisfying the standing assumptions, a point ${\mathfrak Z}=(\frz_1,\ldots,\frz_n)\in(\IH^2)^n$ assigns, to each tetrahedral pair $\Delta_j$, a pair of positively oriented ideal hyperbolic tetrahedra by taking both tetrahedra to have shape parameter $\frz_j$.  Thus, the resulting assignment is invariant under the action of $G$.

The paragraph preceding Definition~\ref{def: involutions} explains why an ordinary quotient edge has edge parameter product $1$, while the quotient edge associated to a core tunnel has edge parameter product $-1$.  When $\alpha$ is hyperbolic, these conditions have the usual orbifold interpretation; the corresponding cone angles in the quotient are $2\pi$ and $\pi$, respectively.

\begin{dfns} \label{def: hyperbolic structures}
If $E$ is an interior edge of $M_\alpha^\circ$ which is not a core tunnel, then its image in $O_\alpha^\circ$ determines a {\it gluing equation} by setting the product of the edge parameters around this quotient edge equal to $1$.  The image in $O_\alpha^\circ$ of the upper core tunnel determines a gluing equation by setting the product of the edge parameters around this quotient edge equal to $-1$.  The collection of these equations is called the set of {\it $G$-invariant Thurston gluing equations} for $M_\alpha^\circ$.
\index{gluing equation|idxdef}

A solution ${\mathfrak Z}\in(\IH^2)^n$ to this system of equations determines a $G$-invariant {\it piecewise hyperbolic structure carried by} $M_\alpha^\circ$.  Let $\calH(M_\alpha^\circ)$ denote the set of all such solutions.

When $\alpha$ is hyperbolic, define $\calH(M_\alpha)$ analogously using the gluing equations on $O_\alpha$.  Every equation for $\calH(M_\alpha^\circ)$ is also an equation for $\calH(M_\alpha)$.  The lower hairpin fold introduces two additional equations: the ordinary edge equation associated to the lower
hairpin edge $\hat e_{n+1}$ and the tunnel equation associated to the lower core tunnel.  Hence,
\[ \calH(M_\alpha)\ \subseteq\ \calH(M_\alpha^\circ). \]
\end{dfns}
\index{piecewise hyperbolic structure|idxdef}

The equations defining $\calH(M_\alpha^\circ)$ are written on the quotient complex $O_\alpha^\circ$, but they imply the ordinary edge gluing equations on $M_\alpha^\circ$.  More precisely, for every interior edge of $M_\alpha^\circ$, the product of the edge parameters around the full edge cycle is $1$.  For an ordinary edge, this follows directly from the corresponding quotient edge equation.  For the upper core tunnel, the edge cycle in $O_\alpha^\circ$ is the quotient of the edge cycle in $M_\alpha^\circ$ by the order-$2$ half turn fixing the tunnel. Since the edge parameter product around the quotient edge is $-1$, the product around the full edge cycle in $M_\alpha^\circ$ is $1$.  This holds whether or not $\alpha$ is hyperbolic.

Completeness is not part of these definitions.  An element of $\calH(M_\alpha^\circ)$ records a compatible $G$-invariant assignment of positively oriented ideal hyperbolic tetrahedra to the tetrahedral pairs of $M_\alpha^\circ$.

For such a carried geometric structure, a {\it developing map} will mean a map
\[ \dev_{\mathfrak Z}\co
\wt M_\alpha^\circ\longrightarrow\IH^3 \] 
whose restriction to every lifted tetrahedron is an orientation-preserving isometry onto an ideal hyperbolic tetrahedron with the prescribed shape parameter.  A corresponding {\it holonomy representation} is a homomorphism
\[ \rho_{\mathfrak Z}\co \pie(M_\alpha^\circ)\longrightarrow\PSL_2(\C) \]
for which $\dev_{\mathfrak Z}$ is $\rho_{\mathfrak Z}$-equivariant.
\index{developing map|idxdef}
\index{holonomy|idxdef}

When $\alpha$ is hyperbolic, $M_\alpha^\circ$ is a manifold, so these edge gluing equations are the usual Thurston gluing equations on $M_\alpha^\circ$.  The construction in Section~E.6.i (especially Lemma~E.6.1) of \cite{BP} gives the usual developing map and holonomy representation.  When $1/\alpha\in\Z$, the ideal complex $M_\alpha^\circ$ is not a manifold. Nevertheless, in Chapter~\ref{chap: FD} a similar developing map and holonomy representation are constructed for every element of $\calH(M_\alpha^\circ)$.

When $\alpha$ is hyperbolic, an element of $\calH(M_\alpha)$ determines a (possibly incomplete) hyperbolic structure on $M_\alpha$.  Later, Theorem~\ref{thm: complete} will show that every element of $\calH(M_\alpha)$ is, in fact, complete.

An explicit description of the gluing equations for the Sakuma--Weeks complexes follows.  These equations can be read directly from the combinatorics of the funnel.

Recall that $\alpha \in \Q \cap (0,1/2)$ and $\alpha \neq 1/3$. The interior edges of $F_\alpha$ are $e_1,\ldots,e_n$ and $n\geq1$.  Suppose ${\mathfrak Z}=(\frz_1,\ldots,\frz_n)\in(\IH^2)^n$ and set $\frz_0=\frz_{n+1}=1$.   

An edge which is not a core tunnel has a single rational label.  This number also occurs as a vertex of the funnel $F_\alpha$.  The gluing equations can be labeled by the corresponding rational numbers.

Suppose that $\beta$ is a vertex of $F_\alpha$ which is not a vertex of $e_0$ or $e_{n+1}$ and that $\beta\neq\hat e_{n+1}$. Then the edge labeled $\beta$ is an ordinary interior edge of $M_\alpha^\circ$.  Let $e_k,\ldots,e_\ell$ be the consecutive interior edges of $F_\alpha$ which contain $\beta$.  The tetrahedral
pairs $\Delta_k,\ldots,\Delta_\ell$ are precisely those which contribute side edges labeled $\beta$ to the corresponding quotient edge cycle.  Each pair contributes the square of the appropriate side edge parameter, namely $\zeta(\frz_j)^2$ if $\beta$ lies on the left side of $F_\alpha$, and $\zeta^2(\frz_j)^2$ otherwise.  The tetrahedral pairs immediately above and below this string contribute $\frz_{k-1}$ and $\frz_{\ell+1}$, respectively.  Thus, the gluing
equation associated to $\beta$ is
\begin{align} \label{eq: mid glue}
1&= \begin{cases} \frz_{k-1} \frz_{\ell+1} \left( \prod_{k}^{\ell} \zeta \left( \frz_{j} \right) \right)^2 & \text{if } \beta \text{ is on the left side of } F_\alpha \\
\frz_{k-1} \frz_{\ell+1} \left( \prod_{k}^{\ell} \zeta^2 \left( \frz_{j} \right) \right)^2 & \text{otherwise.} \end{cases}
\end{align}

The same edge cycle bookkeeping applies to the core tunnels after the corresponding hairpin folds.  The difference is that, as explained above, the quotient edge associated to a core tunnel has edge parameter product $-1$.  For the upper core tunnel, let $\ell$ be the largest integer such that $e_\ell^-=0$.  The upper tunnel gluing equation is
\begin{align} \label{eq: top glue}
-1&= \frz_1 \frz_{\ell+1} \left( \prod_1^{\ell} \zeta(\frz_j) \right)^2
\end{align} 
\index{core tunnel}

Equations \ref{eq: mid glue} and \ref{eq: top glue}, as $\beta$ ranges over the labels described above, constitute the full set of $G$-invariant Thurston gluing equations for $M_\alpha^\circ$.

The three labels at the vertices of the bottom triangle $s_n$ do not determine gluing equations for $M_\alpha^\circ$.  The edges labeled $e_{n+1}^-$ and $e_{n+1}^+$ lie on the unfilled lower boundary, as does the hairpin edge labeled $\hat e_{n+1}$.  When $\alpha$ is hyperbolic, the lower hairpin fold folds $M_\alpha^\circ$ into the lower part of the Sakuma--Weeks triangulation of $M_\alpha$ and introduces two additional gluing equations.

The first additional equation is the ordinary edge equation associated to the lower hairpin edge $\hat e_{n+1}$.  This is an equation of the form \ref{eq: mid glue}, with $\beta=\hat e_{n+1}$.  The second is the equation associated to the lower core tunnel.  Let $\beta$ be the high vertex $e^-_{n+1}$ and take $k$ to be the smallest number such that $e_k$ has $\beta$ as a vertex.  Then the equation for the bottom tunnel is 
\begin{align} \label{eq: bottom glue}
-1&= \begin{cases} \frz_{k-1} \frz_{n} \left( \prod_{k}^{n} \zeta \left( \frz_{j} \right) \right)^2 & \text{if } \beta \text{ is on the left side of } F_\alpha \\
\frz_{k-1} \frz_{n} \left( \prod_{k}^{n} \zeta^2 \left( \frz_{j} \right) \right)^2 & \text{otherwise.} \end{cases}
\end{align}
Thus, when $\alpha$ is hyperbolic, the equations defining
$\calH(M_\alpha)$ consist of the equations defining
$\calH(M_\alpha^\circ)$ together with the ordinary edge equation for
$\hat e_{n+1}$ and the lower tunnel equation
\ref{eq: bottom glue}.

\section{Solving the gluing equations using Farey recursion} \label{sec: solve}
\index{Farey recursion}

As in the preceding sections, assume $\alpha \in \Q \cap (0,1/2)$ and $\alpha \neq 1/3$.  This section uses the Farey recursive functions $d_\calQ$ and $\calQ$ defined in Examples \ref{ex: dQ} and \ref{ex: Q} to solve the $G$-invariant gluing equations described in Section \ref{sec: equations}.   It will show that the gluing equations for $M_\alpha^\circ$ can be solved over the field $\Q(x)$ and that, in the hyperbolic case, imposing the additional equations arising from the lower hairpin fold is equivalent to the equation $\calQ(\alpha)=0$.
\index{Farey recursive function}

\begin{remark} \label{rem: SW functions}
The functions used here are closely related to the functions $f$ and $\psi$ introduced in Section~II.5 of \cite{SW}.  

Associate to a Sakuma--Weeks vector $\mathbf a = (a_1,\ldots,a_k)$ the rational number $r(\mathbf a) = [0;a_1,\ldots,a_k]$.  If $r(\mathbf a)=p/q$, written in lowest terms, the integer $\alpha(\mathbf a)$ of \cite{SW} is $q$ and the Sakuma--Weeks sign $\epsilon(\mathbf a)$ is $(-1)^p$.  Hence
\[ \epsilon(\mathbf a)\, x^{\alpha(\mathbf a)} \ =\ \dQ\big(r(\mathbf a)\big). \]
Under this identification, the recursion defining the Sakuma--Weeks function $f$ is precisely the Farey recursion defining $\calQ$, and the initial conditions agree as well.  Therefore
\[ f(\mathbf a)\ =\ \calQ\big(r(\mathbf a)\big). \]

There is a similar relationship between $\psi$ and the shape parameter function $\calZ$ from Definition~\ref{def: shape parameter function}.  Under the correspondence between the Sakuma--Weeks layering and the interior edges of $F_\alpha$, their tetrahedral parameters $w(i,1)$ agree, up to the indexing convention, with the parameters $\calZ(e_j)$ used here.  

Thus, the Farey recursive function $\calQ$ packages the collection of polynomials $f(\mathbf a)$ from \cite{SW} into a single function on the Stern--Brocot diagram, while $\calZ$ packages the tetrahedral parameters obtained from $\psi$.  
\end{remark}

\begin{remark}\label{rem: Tsvietkova}
Tsvietkova associates another polynomial to a hyperbolic $2$-bridge link in \cite{Tsv}.  Her recursion expresses the crossing and edge labels of a reduced alternating diagram in terms of the first crossing label $w_1$, and the final region relation gives a polynomial equation $P_\alpha(w_1)=0$.  She notes that the relationship between this polynomial and the polynomial arising from the Sakuma--Weeks construction is unclear in general.

From Remark~\ref{rem: SW functions}, to understand the relationship between the Sakuma--Weeks polynomials and the Tsvietkova polynomials, it suffices to compare her polynomials to $\calQ$.  The crossing labels in \cite{Tsv} carry orientation-dependent signs $\epsilon_j$. Write $\wt{P}_\alpha(x)=P_\alpha(\epsilon_1x)$ and omit factors introduced by clearing intermediate denominators in her recursion.  Then, a computation of the polynomials $\wt{P}_\alpha$, for every hyperbolic $\alpha$ with denominator at most 50, shows that $\wt{P}_\alpha$ always divides $\calQ(\alpha)$ in these examples.  In $315$ of the $338$ cases, the two polynomials agree up to sign.  In each of the remaining cases, the quotient $\calQ(\alpha)/\wt{P}_\alpha$ is, up to sign, one of 
\begin{align*}
   \calQ(1/3)&=1-x & \calQ(1/3)^3&=(1-x)^3 & \calQ(1/4)&=1-2x\\
   \calQ(1/5)&=1-3x+x^2& \calQ(2/5)&= 1-x+x^2.
\end{align*}
These computations suggest that, after the above normalization,  Tsvietkova's polynomial is, in general, a factor of $\calQ(\alpha)$.
\end{remark}

Recall from Definition \ref{def: E and more} that $\mathcal{E}$ is the set of edges in the Stern--Brocot diagram $\calG$ whose vertices lie in $\Q_0$ and that the function $\nu \co \calE \to \{ \pm 1\}$ is determined by the slope of the given edge.  

\begin{dfn} \label{def: shape parameter function}
Define the {\it shape parameter function}
\[ \calZ \co \calE \longrightarrow \Q(x) \]
by
\[ \calZ(e) \ =\ -\dQ(\hat{e})^{-\nu(e)} \cdot \left( \frac{\calQ(e^L)}{\calQ(e^R)}\right)^2. \]
\end{dfn}
\index{shape parameter!function|idxdef}

It follows immediately from the discussion in Definition \ref{def: E and more} that 
 \[ \calZ(e) = -\left( d_\calQ\left(\hat{e}\right) \, \frac{ \calQ(e^-)^2 }{\calQ(e^+)^2} \right)^{-\nu(e)}. \]
 If, as usual, $e_0, \ldots, e_{n+1}$ are the top, interior, and bottom edges of the funnel $F_\alpha$ (Definition \ref{def: funnel}), then it is always true that
\[ \calZ(e_0)=1 \qquad \text{and} \qquad \calZ(e_1)=x. \]

Let $z_1,\ldots, z_n$ be the coordinate functions on $\C^n$.  Since $\calH(M_\alpha^\circ)\subset(\IH^2)^n\subset\C^n$, these restrict to coordinate functions on $\calH(M_\alpha^\circ)$.  As such, $z_j$ returns the shape parameter for the tetrahedral pair $\Delta_j$.  When $\alpha$ is hyperbolic (Definition~\ref{def: hyperbolic alpha}), the same coordinate functions restrict to $\calH(M_\alpha)\subseteq\calH(M_\alpha^\circ)$.  From the equations in Section \ref{sec: equations}, these coordinate functions provide abstract gluing equations in variables $z_1, \ldots, z_n$.  

\begin{thm}\label{thm: main1}
\GenericShapeSolution
\end{thm}

Lemma \ref{lem: adjacent zeros} shows that any $\calZ(e)$ can be evaluated at any $\frz \in \C$ by evaluating its numerator and denominator separately.  Doing so provides a well-defined number $\calZ(e)_\frz \in \C \cup \{ \infty \} = \bound \IH^3$.

\begin{cor} \label{cor: main1}
Suppose $\alpha \in \Q \cap (0,1/2)$ and $\alpha \neq 1/3$.  A number $\frz \in \IH^2$ is the first shape parameter for an element of $\calH(M_\alpha^\circ)$ if and only if 
\[ \calZ(e_j)_\frz \ \in \ \IH^2 \]
for every interior edge $e_j$ of $F_\alpha$.   The corresponding element in $\calH(M_\alpha^\circ)$,  
\[ \Big( \calZ(e_1)_\frz, \ldots , \calZ(e_n)_\frz \Big)\]
is unique.
\end{cor}

\proof The recursive argument in the proof of Theorem~\ref{thm: main1} remains valid after setting $x=\frz$ since all of the shape parameters in question lie in $\IH^2$.  This forces the $j$th shape parameter to be $\calZ(e_j)_\frz$ for every $j$.   \endproof 

Using this corollary, the definition of $\calH(M_\alpha^\circ)$ may be replaced with
\[ \calH(M_\alpha^\circ)  \ =\ \left\{ \frz \in \IH^2 \, \big| \, \calZ(e_j)_\frz \in \IH^2 \text{ for every interior edge } e_j \text{ of } F_\alpha \right\}.\]
Now, if $x$ is interpreted as the coordinate on $\IH^2$, the rational functions $\calZ(e_j)$ provide the shape parameters for points of $\calH(M_\alpha^\circ)$.

\begin{thm} \label{thm: main}
Suppose $\alpha$ is hyperbolic.  A number $\frz \in \IH^2$ is the first shape parameter for an element of $\calH(M_\alpha)$ if and only if 
\[ \frz \in \calH(M_\alpha^\circ) \qquad \text{and} \qquad \calQ(\alpha)_\frz = 0.\]
\end{thm}

By Theorem \ref{thm: main}, the definition of $\calH(M_\alpha)$ may be replaced with 
\[ \calH(M_\alpha)  \ =\ \left\{ \frz \in \calH\left(M_\alpha^\circ\right) \, \big| \, \calQ(\alpha)_\frz = 0  \right\}.\]


Later, Theorem~\ref{thm: complete} will show that, when $\alpha$ is hyperbolic, every element of $\calH(M_\alpha)$ is complete.  Together with Mostow--Prasad rigidity and the existence of the Sakuma--Weeks geometric structure, this implies that $\calH(M_\alpha)$ consists of a single point.

The next corollary follows easily and will be used later in the paper.

\begin{cor} \label{cor: degree}
If $\alpha \in \Q \cap (0,1/2)$, then the degree of $\calQ(\alpha)$ is nonzero and zero is not a root of $\calQ(\alpha)$.
\end{cor}

\proof
Remark~\ref{rem: constant term} shows that zero is not a root of $Q(\alpha)$.

Now, suppose that $\alpha$ is hyperbolic.  By \cite{GF}, the Sakuma--Weeks triangulation of $M_\alpha$ is geometric, and hence determines an element of $\calH(M_\alpha)$.  Its first shape parameter $\frz$ cannot be real.  By Theorem~\ref{thm: main}, $\frz$ is a root of $\calQ(\alpha)$ and the degree of $\calQ(\alpha)$ must be at least two.

If $1/\alpha\in\Z$, the result follows from Remark~\ref{rk: chebyshev}.
\endproof

Theorems \ref{thm: main1} and \ref{thm: main} will be proved after a few lemmas.  Corollary \ref{cor: main1} follows immediately from Theorem \ref{thm: main1}.

\begin{lemma} \label{lem: det}
If $e$ is a nonhorizontal edge in $\calE$, then
\[ 
\begin{vmatrix} \calQ\left(e^-\right) & \calQ\left(e^+\right) \\ \calQ\left(e^+\right) & \calQ\left(e^+ \oplus \hat{e}\right)\end{vmatrix} \ =\ \dQ\left(e^-\right) \, \calQ\left(\hat{e}\right)^2.
\]
\end{lemma}

\proof
The proof is by induction on the denominator of $e^+$.  The smallest possible denominator is 2, and the two base possibilities are easily verified by hand.

Now, assume that the denominator $k$ of $e^+$ is larger than 2 and that the lemma holds for denominators smaller than $k$.  Let $\gamma$ be the opposite vertex for the edge with endpoints $e^-$ and $\hat{e}$.   

Suppose first that $e^-$ is not a corner of $\Delta(\hat{e})$.  Then $\gamma \oplus \hat{e}=e^-$ and, by the inductive assumption,
\[ \begin{vmatrix} \calQ\left(\gamma\right) & \calQ\left(e^-\right) \\
\calQ\left(e^-\right) & \calQ\left(e^+\right) \end{vmatrix} = d_\calQ\left( \gamma \right) \calQ(\hat{e})^2.\]
Hence,
\begin{align*}
\begin{vmatrix} \calQ\left(e^-\right) & \calQ\left(e^+\right) \\ \calQ\left(e^+\right) & \calQ\left(e^+ \oplus \hat{e}\right)\end{vmatrix} &=
\begin{vmatrix} 0 & 1 \\
-d_\calQ(\hat{e}) & \calQ(\hat{e}) \end{vmatrix} \ 
\begin{vmatrix} \calQ\left(\gamma\right) & \calQ\left(e^-\right) \\
\calQ\left(e^-\right) & \calQ\left(e^+\right) \end{vmatrix}  \\
&= d_\calQ(\hat{e}) \,  d_\calQ\left( \gamma \right) \calQ(\hat{e})^2.
\end{align*}
Since the determinant of $d_\calQ$ is zero, it is multiplicative on Farey sums.  Therefore, the last expression is equal to $d_\calQ \left( e^- \right) \calQ(\hat{e})^2$ as desired.

If, on the other hand, $e^-$ is a corner of $\Delta(\hat{e})$, then $\gamma \oplus^2 e^-=e^+$ and so
\[ \calQ\left( \gamma \right) = \frac{\calQ\left( e^- \right) \calQ(\hat{e}) - \calQ \left( e^+ \right)}{d_\calQ \left( e^- \right) }\]
by condition \ref{eq: defn}.  The inductive assumption and the multiplicativity of $d_\calQ$ give
\[ d_\calQ(\hat{e}) \, \calQ \left( e^- \right)^2 =
d_\calQ \left( e^- \right) d_\calQ \left( \gamma \right) \, \calQ \left( e^- \right)^2 = d_\calQ \left( e^- \right)
\begin{vmatrix} \calQ \left( \gamma \right) & \calQ (\hat{e}) \\ \calQ(\hat{e}) & \calQ \left( e^+ \right) \end{vmatrix}.\]

Then
\[ \calQ(\hat{e})\,  \calQ \left( e^- \right) \calQ \left( e^+ \right) - \calQ \left( e^+ \right)^2 - d_\calQ(e^-) \, \calQ \left( \hat{e} \right)^2 =  d_\calQ \left( \hat{e} \right) \calQ(e^-)^2\]
holds by substituting on the right for $\calQ \left( \gamma \right)$.  This, together with the expression
\[ \calQ \left( e^+ \oplus \hat{e} \right) = -d_\calQ(\hat{e}) \, \calQ \left( e^- \right) + \calQ(\hat{e}) \calQ \left( e^+ \right) \]
shows that
\[ \begin{vmatrix} \calQ\left(e^- \right) & \calQ\left(e^+\right) \\
\calQ\left(e^+\right) & \calQ\left(e^+ \oplus \hat{e}\right) \end{vmatrix} = d_\calQ \left( e^- \right) \, \calQ(\hat{e})^2 \]
as needed.  
\endproof

The next lemma utilizes the function $\nu \co \calE \to \{ \pm 1\}$ from Definition~\ref{def: E and more}.

\begin{lemma} \label{lem: zetas}
Suppose $e$ is a nonhorizontal edge in $\calE$.  Then
\[ d_\calQ(\hat{e}) \, d_\calQ(e^-\oplus e^+) = d_\calQ(e^+)^2.\]
If $\nu(e)=1$ then 
\[
\zeta \left( \calZ(e) \right) = \frac{d_\calQ(\hat{e}) \,\calQ(e^-)^2}{\calQ(\hat{e}) \,\calQ(e^-\oplus e^+)} \quad \text{and} \quad \zeta^2 \left( \calZ(e) \right) = \frac{\calQ(\hat{e}) \,\calQ(e^-\oplus e^+)}{\calQ(e^+)^2}.
\]
If $\nu(e)=-1$ then
\[
\zeta \left( \calZ(e) \right) = \frac{\calQ(e^+)^2}{\calQ(\hat{e}) \,\calQ(e^-\oplus e^+)} \quad \text{and} \quad \zeta^2 \left( \calZ(e) \right) = \frac{\calQ(\hat{e}) \,\calQ(e^-\oplus e^+)}{d_\calQ(\hat{e}) \,\calQ(e^-)^2}.
\]
\end{lemma}
\proof  By the multiplicativity of $d_\calQ$, 
\begin{align*}
d_\calQ(\hat{e}) \, d_\calQ(e^-\oplus e^+)&= d_\calQ(\hat{e}) \,d_\calQ(\hat{e}) \,d_\calQ(e^-)^2\\
&= \left( d_\calQ(\hat{e}) \,d_\calQ(e^-) \right)^2 \\
&= d_\calQ(e^+)^2.
\end{align*}

Suppose that $\nu(e)=1$.  Then \[ \calZ(e) = \frac{-\calQ(e^+)^2}{d_\calQ(\hat{e}) \, \calQ(e^-)^2}.\]
Now, for indeterminates $a$ and $b$, 
\[
\zeta\left( \frac{a}{b} \right) = \frac{b}{b-a} \quad \text{and} \quad \zeta^2\left( \frac{a}{b} \right) = \frac{a-b}{a}. 
\]
So, using Lemma \ref{lem: det},  
\[
\zeta \left( \calZ(e) \right) = \frac{d_\calQ(\hat{e}) \,\calQ(e^-)^2}{\calQ(\hat{e}) \,\calQ(e^-\oplus e^+)} \quad \text{and} \quad \zeta^2 \left( \calZ(e) \right) = \frac{\calQ(\hat{e}) \,\calQ(e^-\oplus e^+)}{\calQ(e^+)^2}.
\]
The proof is similar in the case that $\nu(e)=-1$.
\endproof

\begin{figure}
\setlength{\unitlength}{.1in}
\begin{picture}(40,27)
\put(.3,0) {\includegraphics[width= 4in]{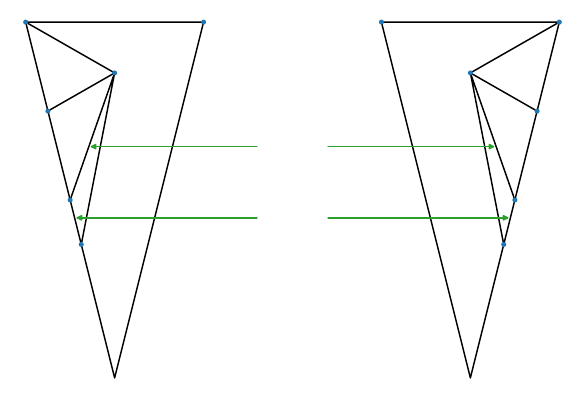}}
\put(0,25.8){$\gamma_k$}
\put(39.3,25.8){$\gamma_k$}
\put(15.1,25.8){$\hat{e}_k$}
\put(24.5,25.8){$\hat{e}_k$}
\put(19.0,17.2){$e_{n-1}$}\
\put(19.7,12.2){$e_{n}$}
\put(8.6,22){$\beta$}
\put(31.2,22){$\beta$}
\put(5.6,24.2){$e_k$}
\put(33.3,24.2){$e_k$}
\put(3.9,10.2){$\gamma_n$}
\put(35.5,10.2){$\gamma_n$}
\end{picture}
\caption{The labeling used in Lemma \ref{lem: product}.  On the left is the case $\gamma_k=\kappa^L$ and on the right is the case $\gamma_k=\kappa^R$.}  \label{fig: products}
\end{figure}

The next lemma provides formulas for certain products of shape parameter functions on edges of triangles $\Delta(\beta)$ in $\calG$.  Figure \ref{fig: products} shows the labeling used for the relevant edges and vertices.

\begin{lemma} \label{lem: product}
Take $\beta \in \Q \cap (0,1)$ and integers $k$ and $n$ with $1 \leq k < n$.   Take $\gamma_k$ to be a corner $\kappa$ of $\Delta(\beta)$ and define $\gamma_{k+j} = \gamma_k \oplus^j \beta$.  If $k\leq j<n$, let $e_j$ be the edge in $\calG$ between $\beta$ and $\gamma_j$.  Let $e_n$ be the edge between $\gamma_{n-1}$ and $\gamma_n$ and let $e_{k-1}$ be the top edge of $\Delta(\beta)$.  If $\gamma_k=\kappa^L$ then
\[ \prod_k^{n-1} \zeta^2 \left( \calZ(e_j) \right) \ = \ \frac{\calQ\left( \hat{e}_k \right) \calQ\left(\gamma_n \right)}{d_\calQ\left( \hat{e}_k \right) \calQ\left( \gamma_k \right) \calQ\left(\gamma_{n-1}\right)}\]
and
 \[\calZ\left( e_{k-1} \right) \left( \prod_k^{n-1} \zeta^2 \left( \calZ(e_j) \right) \right)^2 \ = \ \calZ(e_n)^{-1}.\]
 If $\gamma_k=\kappa^R$ then
 \[ \prod_k^{n-1} \zeta \left( \calZ(e_j) \right) \ = \ \frac{d_\calQ\left( \hat{e}_k \right) \calQ\left( \gamma_k \right) \calQ\left(\gamma_{n-1}\right)}{\calQ\left( \hat{e}_k \right) \calQ\left(\gamma_n \right)}\]
and
 \[\calZ\left( e_{k-1} \right) \left( \prod_k^{n-1} \zeta \left( \calZ(e_j) \right) \right)^2 \ = \ \calZ(e_n)^{-1}.\]
\end{lemma}

\proof
Suppose $\gamma_k=\kappa^L$.  Then, by Lemma \ref{lem: zetas},
\[ \zeta^2 \left( \calZ(e_k) \right) \ =\ \frac{\calQ \left( \hat{e}_k \right) \calQ(\gamma_{k+1})}{\dQ\left( \hat{e}_k \right) \calQ(\gamma_k)^2}\]
and, if $k<j<n$,
\[ \zeta^2 \left( \calZ(e_j) \right) \ =\ \frac{\calQ(\gamma_{j-1}) \calQ(\gamma_{j+1})}{\calQ(\gamma_j)^2}.
\]
Using these expressions, the left side of the first identity in the statement of the lemma collapses to the expression on the right side.  

Using the definition of $\calZ$ and the multiplicity of $\dQ$ from Lemma \ref{lem: zetas},
\[ \calZ(e_{k-1}) \ = \ - \frac{\dQ\left( \hat{e}_k \right) \calQ(\gamma_k)^2}{\dQ(\gamma_k) \calQ\left( \hat{e}_k \right)^2}.\]
Combined with the first formula stated in the lemma,
\begin{align*}
\calZ\left( e_{k-1} \right) \left( \prod_k^{n-1} \zeta^2 \left( \calZ(e_j) \right) \right)^2 & = - \frac{1}{\dQ\left( \hat{e}_k \right) \dQ(\gamma_k) } \cdot \frac{\calQ(\gamma_n)^2}{\calQ(\gamma_{n-1})^2}\\
&= - \left( \dQ(\beta) \, \frac{\calQ(\gamma_n)^2}{\calQ(\gamma_{n-1})^2} \right)^{\nu(e_n)} \\
&= \calZ(e_n)^{-1}
\end{align*}
since $\nu(e_n)=-1$.

The argument for the $\gamma_k=\kappa^R$ case is essentially the same.
\endproof

\proof[Proof of Theorem \ref{thm: main1}] The proof is by induction on the number $n$ of interior edges of $F_\alpha$.  It proves simultaneously that a solution with $z_1=x$ exists and that it is unique.

If $n=1$, then $\alpha\in\{1/4,2/5\}$ and $F_{1/4}=F_{2/5}$ has exactly one interior edge.  In this case there is one shape parameter and no ordinary interior edge gluing equation, so the result follows.

Now, assume $n>1$.  If $\alpha_0=e_n^- \oplus e_n^+$, then $F_{\alpha_0}$ is obtained from $F_\alpha$ by deleting its last simplex $s_n$.  The triangles $\{ s_j \}_0^{n-1}$ are common to both funnels.  By the inductive hypothesis, the gluing equations for $M_{\alpha_0}^\circ$ have the unique solution
\[ z_j\ =\ \calZ(e_j) \qquad \text{for} \qquad 1\leq j\leq n-1\]
subject to $z_1=x$.

Take $\beta$ to be the element of $\{ e_{n-1}^\pm \}$ which is not an endpoint of $e_n$.  The only $G$-invariant gluing equation for $M_\alpha^\circ$,
which is not a gluing equation for $M_{\alpha_0}^\circ$, is the equation labeled $\beta$.  Moreover, this is the only equation which involves the variable $z_n$.  Assign variables as in the statement of the theorem
\[\qquad z_j = \calZ(e_j) \]
for $j \in \{ 1, \ldots, n-1 \}$.  By the inductive assumption, it remains only to prove that the equation for $\beta$ is satisfied if and only if $z_n = \calZ(e_n)$.  There are three possibilities
\begin{enumerate}
\item $\beta$ is zero, 
\item $\beta$ is on the right side of $F_\alpha$, or
\item $\beta$ is not zero and is on the left side of $F_\alpha$.
\end{enumerate}

Suppose first that $\beta=0$ and note that, as in the proof of Lemma \ref{lem: product}, Lemma \ref{lem: zetas} provides
\[\prod_1^{n-1} \zeta\left( \calZ(e_j) \right) \ = \ \frac{\calQ(e_n^-)}{\calQ(e_n^+)}. \]
Consider equation \ref{eq: top glue} to see that the equation for $\beta$ is satisfied if and only if
\begin{align*}
z_n & = -\left( \calZ(e_1) \left( \prod_{1}^{n-1} \zeta\left( \calZ(e_j) \right) \right) ^2 \right)^{-1} \\
&= -\left( \calZ(e_1) \, \frac{\calQ(e_n^-)^2}{\calQ(e_n^+)^2} \right)^{-1}.
\end{align*}
Since $\calZ(e_1)=x=d_\calQ(0)$ and $\nu(e_n)=1$, this is equal to $\calZ(e_n)$.  Hence, the proof is finished in this case.

\begin{figure}
\setlength{\unitlength}{.1in}
\begin{picture}(26,16)
\put(5.6,0) {\includegraphics[width= 2in]{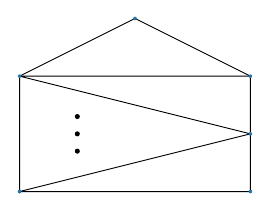}}
\put(15,0){$e_n$}
\put(14.7,10.5){$e_{k-1}$}
\put(0,.8){$e_n^-=\gamma_{n-1}$}
\put(5,10){$\gamma_k$}
\put(25.0,.8){$e_n^+$}
\put(25.0,5.2){$\beta$}
\put(24.8,10){$\hat{e}_k$}
\put(15,15){$\hat{e}_{k-1}$}
\end{picture}
\caption{Vertex labels for case (2) of Theorem \ref{thm: main1}.}  \label{fig: gamma labels}
\end{figure}

Now, assume that $\beta$ lies on the right side of $F_\alpha$.  For $k \leq j \leq n-1$, let $\gamma_j$ be the left vertex of $e_j$ as shown in Figure \ref{fig: gamma labels}.  Let $\gamma_n=e_n^+$.  Notice that $e_{n-1}$ is the last interior edge that has $\beta$ as an endpoint.  Take $k$ to be the index of the first such edge.  The gluing equation for $\beta$ is given by the second option listed in equation \ref{eq: mid glue}.  To prove this case, it is enough to prove 
\[ \calZ\left(e_n\right)^{-1} \ =\  \calZ\left(e_{k-1}\right) \left( \prod_k^{n-1} \zeta^2\left( \calZ\left( e_j \right)\right) \right)^2 .\]
This is true by Lemma \ref{lem: product}.  The same argument proves the final case.
\endproof

\proof[Proof of Theorem \ref{thm: main}] 
For $\frz \in \IH^2$ and $j \in \{ 1, \ldots, n\}$, write $\frz_j = \calZ(e_j)_\frz$.  By Theorem~\ref{thm: main1}, $\frz\in\calH(M_\alpha^\circ)$ if and only if $\frz_j \in \IH^2$ for each $j$.

The equations defining $\calH(M_\alpha)$ consist of the equations defining $\calH(M_\alpha^\circ)$ together with the ordinary edge equation for the lower hairpin edge $\hat e_{n+1}$ and the tunnel equation for the lower core tunnel.  Therefore, $\frz \in \calH(M_\alpha)$ if and only if $\frz \in \calH(M_\alpha^\circ)$ and both of these equations are satisfied by $(\frz_1, \ldots, \frz_n)$.

By Lemma~\ref{lem: product}, the gluing equation for $\hat e_{n+1}$ is equivalent to $\calZ\left(e_{n+1}\right) = 1$.  By the definition of $\calZ$, this is equivalent to
\[ 0 \ =\ \calQ\left(e_{n+1}^+\right)^2+\dQ\left( \hat{e}_{n+1} \right) \calQ\left(e_{n+1}^-\right)^2\]
which, by Lemma \ref{lem: det}, is equivalent to
\[ 0 \ = \ \calQ\left( \hat{e}_{n+1} \right) \calQ(\alpha).\]

The lower tunnel equation is redundant once all of the other gluing equations for $M_\alpha$ are satisfied.  Indeed, take the product of all of the
$G$-invariant gluing equations for $M_\alpha$.  For each tetrahedral pair $\Delta_j$, the product of its contributions is
\[ \Big( \frz_j\, \zeta(\frz_j)\, \zeta^2(\frz_j) \Big)^2\ =\ 1. \]
Thus, the product of the left sides of all the gluing equations is $1$.  The product of their right sides is also $1$ because the ordinary edge equations have right side $1$ and the upper and lower core-tunnel equations each have right side $-1$.  Therefore, any one of the gluing equations follows from all of the others.

Now, if $\frz \in \calH(M_\alpha)$ then, substituting $x=\frz$, it must be true that $\calQ\left( \hat{e}_{n+1} \right) \calQ(\alpha)$ is zero.  But $\frz_n = \calZ(e_n)$ is a nonzero complex number and this implies that $\calQ(\hat{e}_{n+1}) \neq 0$.  Hence $\calQ(\alpha)=0$.

On the other hand, if $\frz \in \calH(M_\alpha^\circ)$ and $\frz$ satisfies $\calQ(\alpha)$, then $\frz$ satisfies the equation for the lower hairpin edge $\hat e_{n+1}$ and the tunnel equation for the lower core tunnel.  Therefore, $\frz \in \calH(M_\alpha)$.
\endproof

Notice that the algebraic calculation in the proof does not require that $\alpha$ is hyperbolic.  

\section{Computer Code} \label{subsec: code II}

The code given in this section relies on the commands given in Section \ref{subsec: code I}.  As before, a rational number $\alpha$ is represented by a vector \lstinline$v$.  If $\{p/q,s/t\}$ is a Farey pair then there is an edge $e$ of $\calG$ with $e^-=p/q$ and $e^+=s/t$.  In the code below, $e$ is entered as a pair of vectors $[[p,q], [s,t]]$.  

The function \lstinline$e_sort(e)$ sorts the pair $e$ so that the vector representing $e^-$ is listed before $e^+$.  The functions \lstinline$s(e)$ and \lstinline$e_hat(e)$ return the slope $\nu(e)$ of $e$ and the vertex $\hat{e}$ for $e$.  The rational function $\calZ(e)$ is given by the command \lstinline$Z(e)$.  When \lstinline$v$ is the vector representing $\alpha$, \lstinline$funneledges(v)$ returns a list of the interior edges of the funnel $F_\alpha$ preceded with the top edge and ending with the bottom edge.

\begin{lstlisting}[language=Python]
def e_sort(e):   # e is a pair of vectors.
    return sorted(e, key = lambda v: v[1])

def s(e):   # e is a pair of vectors.
    (m,p) = e_sort(e)
    if quot(m) < quot(p): return -1
    else: return 1

def e_hat(e):   # e is a pair of vectors.
    (m,p) = e_sort(e)
    return vector(p) - vector(m)

def Z(e):   # e is a pair of vectors.
    (m,p) = e_sort(e)
    return -( dQ(e_hat(e)) * Q(m)^2/Q(p)^2 )^(-s(e))
    
def funnel_edges(v):   # v is a vector.
    v = vector(v)
    edges = [list(map(vector, [[0,1],[1,1]]))]
    e = edges[-1]
    Fsum = sum(e)
    while Fsum != v:
        if quot(v) < quot(Fsum):
            edges.append([e[0],Fsum])
        else:
            edges.append([Fsum,e[1]])
        e = edges[-1]
        Fsum = sum(e)
    return edges
\end{lstlisting}

The next functions provide standard information regarding the 2-bridge link determined by $\alpha \in \Q_0$.  In particular, they compute the number of tetrahedra in the Sakuma--Weeks complex for $F_\alpha$, the crossing number for $L_\alpha$, and the sequence of twist numbers in the diagram $D_\alpha$.

\begin{lstlisting}[language=Python]
def tet_number(v):   # v is a vector.
    return 2 * (len(funnel_edges(v)) - 2)

def crossings(v):   # v is a vector.
    return 3 + tet_number(v)/2

def braid(v):   # v is a vector.
    return list(quot(v).continued_fraction())[1:]
\end{lstlisting}

It is now possible to compute various geometric properties for $M_\alpha$.  The code below gives a simple implementation which works well for moderate examples.  In particular, it can be used to find decimal approximations for the shape parameters of the hyperbolic ideal tetrahedra in the Sakuma--Weeks triangulation and to compute the hyperbolic volume of $M_\alpha$.

The function \lstinline$roots(p)$ computes the complex roots of a polynomial \lstinline$p$.  The boolean function
\lstinline$shape_test(z)$ determines whether or not a complex number lies in $\IH^2$.  When $\alpha$ is hyperbolic, the function \lstinline$geometric_root(v)$ filters the roots of $\calQ(\alpha)$ by requiring all of the shape parameters $\calZ(e_j)$ to lie in $\IH^2$. By Theorem~\ref{thm: complete}, exactly one root satisfies these
conditions.  (If $1/\alpha\in\Z$, Remark~\ref{rk: chebyshev} explains that there are no such roots.)  To list the shape parameters for the pairs of tetrahedra $\Delta_j$, apply the function \lstinline$param_list(v)$.  Finally, if $\alpha$ is hyperbolic, the function \lstinline$volume(v)$ computes the volume of $M_\alpha$ by summing the volumes of its ideal tetrahedra. \index{geometric root!2-bridge link}

\begin{lstlisting}[language=Python]
import numpy

def roots(p):   # p is a polynomial.
    coef_list = p.coefficients(sparse=False)
    coef_list.reverse()
    return list(numpy.roots(coef_list))

def shape_test(z):   # z is a complex number.
    return imag(z) > 1.1e-25 

def geometric_root(v):   # v is a vector.
    rootlist = roots(Q(v))
    fedges = funnel_edges(v)[1:-1]
    for fe in fedges:
        if len(rootlist) > 1:
            rootlist = list(filter(lambda rt: 
            	shape_test(Z(fe).substitute(x=rt)), 
            	rootlist))
            if len(rootlist)==1: return rootlist[0]
        else: return rootlist[0]
        
def param_list(v):   # v is a vector.
    groot = geometric_root(v)
    fedges = funnel_edges(v)[1:-1]
    return [Z(fe).substitute(x=groot) for fe in fedges]
    
def tet_volume(z):  # z is a complex number.
    return dilog(z).imag() + ln(abs(z)) * arg(1 - z)

def volume(v):   # v is a vector.
    return 2 * sum(list(map(tet_volume, param_list(v))))
\end{lstlisting}

Similar code was used to compute the plot shown in Figure \ref{fig: roots} of all geometric roots for $\alpha \in \Q_0$ with denominators at most 128.  Notice the similarity with the Stern--Brocot tree shown in Figure \ref{fig: tree}.

For much larger examples, it is inefficient to construct the polynomial $\calQ(\alpha)$ and compute all of its roots.  Instead, the Farey recursive data can be evaluated numerically at high precision and Newton iteration can be used to follow the geometric root through the successive stages of the funnel.  At each stage, the preceding geometric root provides a good initial approximation for the next iteration, while the conditions
\[ \Im\bigl(\calZ(e_j)_\frz\bigr)\ >\ 0 \]
provide a check that the resulting value continues to determine the geometric structure.  Once the geometric root has been found, the same recursive evaluations of $\calZ$ and $\calV=\calN/\calQ$ give the tetrahedral shapes, volume, and developed domain without constructing $\calQ(\alpha)$ as a polynomial.  This is the numerical strategy used for the large example in Figure~\ref{fig: Big}.

\chapter{Developed domains and holonomy}\label{chap: FD} \index{fundamental domain} \index{holonomy}

This chapter uses the Farey recursive functions $\calQ$ and $\calN$ from Examples~\ref{ex: Q} and~\ref{ex: N} to give explicit geometric realizations of the tetrahedral pairs of $M_\alpha^\circ$ determined by points of $\calH(M_\alpha^\circ)$, defined in Definition~\ref{def: hyperbolic structures}.  As in the preceding chapter, assume that $\alpha\in\Q\cap(0,1/2)$ and $\alpha\neq1/3$.  This way, the funnel $F_\alpha$ always has at least one interior edge.

For every $\frz\in\calH(M_\alpha^\circ)$, the construction in this chapter gives an explicit collection $\Omega(\alpha)_\frz$ of ideal hyperbolic tetrahedra which geometrically realizes the tetrahedral core of $M_\alpha^\circ$.  It will be shown that $\Omega(\alpha)_\frz$ is the image, under a developing map, of a connected fundamental domain for the action of $\pie(M_\alpha^\circ)$ on $\wt M_\alpha^\circ$, and that this developing map has an associated holonomy representation
\[ \pie(M_\alpha^\circ)\longrightarrow\PSL_2(\C). \]
Now, recall from Definition~\ref{def: hyperbolic alpha} that $\alpha$ is called hyperbolic if $1/\alpha\notin\Z$.  When $\alpha$ is hyperbolic, $M_\alpha^\circ$ is a manifold and these are the usual developing map and holonomy representation of the hyperbolic structure carried by $M_\alpha^\circ$.  When $1/\alpha\in\Z$, they are understood in the piecewise-hyperbolic sense defined in Chapter~\ref{chap: SakumaWeeks}.
\index{developing map}

\begin{dfn} \label{def: calV}
Define
\[ \calV \co \wh{\Q}_0 \longrightarrow \Q(x) \cup \left\{ \infty \right\} \qquad \text{by} \qquad \calV=\calN/\calQ.\]
 \end{dfn}
 
This function will be used to locate the ideal vertices for these fundamental domains.   Recall that $\calQ(\gamma)$ is the zero polynomial only when $\gamma=\infty$, so $\calV \left( \Q_0 \right) \subset \Q(x)$.

Recall from Definition~\ref{def: hyperbolic structures} that $\calH(M_\alpha^\circ)$ consists of the positive $G$-invariant solutions to the gluing equations carried by the tetrahedral pairs of $M_\alpha^\circ$.  By Corollary~\ref{cor: main1}, it is viewed as a subset of $\IH^2$, with $\frz\in\calH(M_\alpha^\circ)$ giving shape parameter $\frz$ to the tetrahedra in $\Delta_1$ and shape parameter $\calZ(e_j)_\frz$ to those in $\Delta_j$.   For each $\gamma\in\Q_0$, $\calV(\gamma)$ is a rational function of the coordinate $x$.  These functions will be specialized at points $\frz\in\calH(M_\alpha^\circ)$; Lemma~\ref{lem: more zeros} below shows that these specializations are always well-defined as elements of $\C\cup\{\infty\}=\bound\IH^3$.

\begin{figure}
\setlength{\unitlength}{.1in}
\begin{picture}(41,14.5)
\put(4.6,2) {\includegraphics[width= 3.5in]{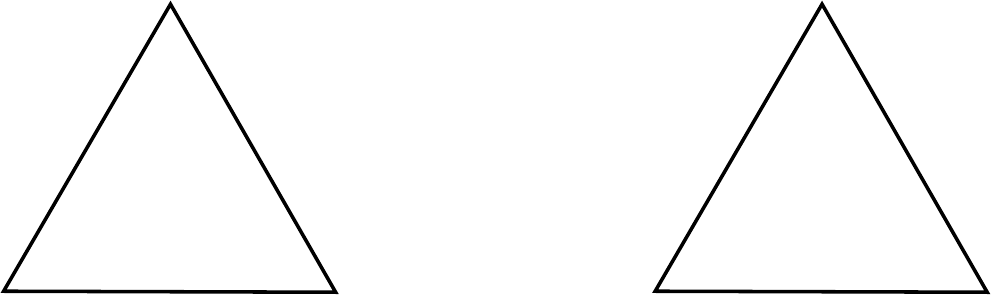}}
\put(0,0){$-1+\calV \left(e^L\right)$}
\put(8.8,13.5){$\varepsilon+\calV \left(\hat{e}\right)$}
\put(32,13.5){$\calV \left(e^R\right)$}
\put(38.2,0){$\calV \left(e^R\right)$}\
\put(14.5,0){$\calV \left(e^R\right)$}
\put(24,0){$\calV \left(e^L \oplus e^R\right)$}
\put(9.3,5.7){$\delta(e)$}
\put(32.2,5.7){$\delta'(e)$}
\put(5.5,6.5){$e^R$}
\put(14.5,6.5){$e^L$}
\put(28.4,6.5){$e^R$}
\put(37.6,6.5){$\hat{e}$}
\put(8.4,2.5){$e^L\oplus e^R$}
\put(32.8,2.5){$e^R$}
\end{picture}
\caption{This diagram shows how pairs in the generic tetrahedra $\delta(e)$ and $\delta'(e)$ are labeled by numbers in $\Q_0$.  This is also described in Definition~\ref{def: generic tetrahedra}.  The choices made here are made so that the specialized tetrahedra are labeled in agreement with the Sakuma--Weeks complexes.}  \label{fig: generic}
\end{figure}

\begin{dfns} \label{def: generic tetrahedra}
Let $e \in \calE$ be a nonhorizontal edge.  Then $\{ \hat{e}, e^L, e^R, e^L\oplus e^R\}$ is a subset of $\Q_0$.  Notice that,  if $\nu(e)=1$, then $\hat{e} < e^L, e^R$ and, if $\nu(e)=-1$, then $e^L, e^R < \hat{e}$. Set
\[\varepsilon \ =\ \begin{cases} -1 & \text{if } \nu(e)=1 \\
0 & \text{otherwise}. \end{cases}\] 
Define the {\it generic tetrahedra}  \index{generic tetrahedron|idxdef}
\[ \delta(e) \ =\ \Big\{ \varepsilon + \calV\left( \hat{e} \right), \ -1+\calV\left( e^L \right), \  \calV\left( e^R \right)\Big\} \ \subset \ \Q(x)
\]
and
\[ 
\delta'(e)= \Big\{ \calV \left( e^L \oplus e^R \right), \ \calV \left(e^L\right), \  \calV \left( e^R \right) \Big\}\ \subset \ \Q(x).\]
Think of these triples as the vertices of triangles.  Label their
edges as indicated in Figure~\ref{fig: generic}.  Equivalently, use the following rules. 
\begin{itemize}
\item The pair $\left\{ \varepsilon + \calV(\hat{e}), \, -1+\calV(e^L) \right\}$ in $\delta(e)$ is labeled $e^R$.
\item The pair $\left\{ \varepsilon + \calV(\hat{e}), \, \calV(e^R) \right\}$ in $\delta(e)$ is labeled $e^L$.
\item The pair $\left\{ -1+\calV(e^L), \, \calV(e^R) \right\}$ in $\delta(e)$ is labeled $e^L \oplus e^R$.
\item The pair $\left\{ \calV(e^L\oplus e^R), \, \calV(e^L) \right\}$ in $\delta'(e)$ is labeled $e^R$.
\item The pair $\left\{ \calV(e^L \oplus e^R), \, \calV(e^R)\right\}$ in $\delta'(e)$ is labeled $e^L$.
\item The pair $\left\{ \calV(e^L), \, \calV(e^R)\right\}$ in $\delta'(e)$ is labeled $\hat{e}$.
\end{itemize}

For $\frz\in\C$, define $\delta(e)_\frz$ and $\delta'(e)_\frz$ to be the convex hull in $\IH^3$ of $\infty$ together with the points obtained by evaluating the corresponding three vertices at $\frz$.  Thus, $\delta(e)$ and $\delta'(e)$ may be regarded as one parameter families of possibly degenerate ideal tetrahedra.  The triangle opposite $\infty$ is called the
{\it finite face}; it carries the edge labeling just defined. \index{finite face|idxdef}

Suppose now that $e=e_j$ is an interior edge of a funnel. The finite face labels agree with the edge labels on the corresponding faces of the tetrahedra in $\Delta_j$ from Definition~\ref{def: edge types}. The labels $e^L$ and $e^R$ are side labels, $\hat e$ is the top label, and $e^L\oplus e^R$ is the bottom label.  When $\delta(e)_\frz$ and $\delta'(e)_\frz$ are ideal tetrahedra, extend the labeling from the finite face to the three edges incident to $\infty$ so that the full edge labeling agrees with that of $\Delta_j$.  

Finally, these conventions extend to the horizontal edge $e=[0,1]$.  In this case,
\[ \delta(e)\ =\ \{0,1\} \AND \delta'(e)\ =\ \{1,-x+1\}.\]
Using that $\hat e=\infty$, $\calV(\infty)=\infty$, and  $\calV(0)=\calV(1)=1$,  take
\[ \delta(e)_\frz\ =\ [\infty,0,1] \AND \delta'(e)_\frz\ =\ [1,-\frz+1,\infty]\]
which are ideal triangles provided $\frz \neq 0$.  The edges of these triangles are labeled using the inclusions
\[ \delta(e)_\frz\ \subset\ \delta(e_1)_\frz \AND \delta'(e)_\frz\ \subset\ \delta'(e_1)_\frz\]
where $e_1=[0,1/2]$.
\end{dfns}
 
 \begin{lemma} \label{lem: delta parameters}
Assume $\alpha \in \Q \cap (0,1/2)$ and $\alpha \neq 1/3$.  Let $e_0$ and $e_{n+1}$ be the top and bottom edges of $F_\alpha$ and $\{ e_1, \ldots, e_n\}$ the interior edges of $F_\alpha$.  If $\frz \in \calH(M_\alpha^\circ)$ then $\delta(e_j)_\frz$ and $\delta'(e_j)_\frz$ are ideal hyperbolic tetrahedra.  The top and bottom edges of these tetrahedra each have edge parameter $\calZ(e_j)_\frz$. 
\end{lemma}

The proof  will utilize the following lemmas. 

\begin{lemma} \label{lem: det2}
If $e$ is a nonhorizontal edge in $\calE$, then
\[ 
\begin{vmatrix} \calN\left(e^-\right) & \calQ\left(e^-\right) \\ \calN\left(e^-\oplus^i \hat{e}\right) & \calQ\left(e^- \oplus^i \hat{e}\right)\end{vmatrix} \ =\ -\nu(e) \, \dQ\left(e^-\right) \, \calQ\left(\hat{e}\right)^i \]
provided $i \in \{ 1,2 \}$.
\end{lemma}

\proof  The following argument for the $i=1$ case is very similar to the inductive argument of Lemma \ref{lem: det}.   Here again, the proof is by induction on the denominator of $e^+$.  As before, the base case is easily verified by hand.  

Assume that the denominator $k$ of $e^+$ is larger than 2 and that the lemma holds for denominators smaller than $k$.  Let $\gamma$ be the opposite vertex for the edge with endpoints $e^-$ and $\hat{e}$. 

Suppose first that $e^-$ is not a corner of $\Delta(\hat{e})$.  Then $\gamma \oplus \hat{e}=e^-$ and, using the inductive assumption and that $\nu$ takes the same value on both $e$ and the edge between $\gamma$ and $e^-$,
\begin{align*}
\begin{vmatrix} \calN(e^-) & \calQ(e^-) \\ \calN(e^+) & \calQ(e^+) \end{vmatrix} &= 
\begin{vmatrix} 0 & 1 \\ -d_\calQ(\hat{e}) & \calQ(\hat{e}) \end{vmatrix}\, \begin{vmatrix} \calN(\gamma) & \calQ(\gamma) \\ \calN(e^-) & \calQ(e^-) \end{vmatrix} \\
&= -\nu(e) \, d_\calQ(\hat{e}) \, d_\calQ(\gamma) \, \calQ(\hat{e}) \\
&= -\nu(e) \, d_\calQ(e^-) \, \calQ(\hat{e}).
\end{align*}
Since $e^- \oplus \hat{e}=e^+$, this is the $i=1$ result.

Now, assume that $e^-$ is a corner of $\Delta(\hat{e})$.  Then the edge $f$ between $e^-$ and $\gamma$ is the top edge of $\Delta(\hat{e})$.  Moreover, $\gamma \oplus^2 e^- = e^+$, so
\begin{align*}
\calQ(e^+) &= -\dQ(e^-) \, \calQ(\gamma) + \calQ(e^-) \, \calQ(\hat{e}) \\
\calN(e^+) &= -\dQ(e^-) \, \calN(\gamma) + \calQ(e^-) \, \calN(\hat{e}).
\end{align*}
These formulas imply that
\begin{align} \label{eq: linear}
 \begin{vmatrix} \calN(e^-) & \calQ(e^-) \\ \calN(e^+) & \calQ(e^+) \end{vmatrix} & = \dQ(e^-) \, \begin{vmatrix} \calN(\gamma) & \calQ(\gamma) \\ \calN(e^-) & \calQ(e^-) \end{vmatrix} + \calQ(e^-) \, \begin{vmatrix} \calN(e^-) & \calQ(e^-) \\ \calN(\hat{e}) & \calQ(\hat{e}) \end{vmatrix}.
 \end{align}
Using the inductive assumption again,
\begin{align} \label{eq: f edge}
\begin{vmatrix} \calN(f^-) & \calQ(f^-) \\ \calN(f^+) & \calQ(f^+) \end{vmatrix} &= -\nu(f) \, \dQ(f^-) \, \calQ\left(\hat{f}\right)
\end{align}
and
\begin{align} \label{eq: hat}
 \begin{vmatrix} \calN(e^-) & \calQ(e^-) \\ \calN(\hat{e}) & \calQ(\hat{e}) \end{vmatrix} \ = \  -\nu(e) \, \dQ(e^-) \, \calQ(\gamma).
 \end{align}
Here, we have also used the fact that the value of $\nu$ agrees on $e$ and on the edge between $e^-$ and $\widehat e$.  Now, if
\[ \beta \ =\ \begin{cases} e^- & \text{if } \nu(e)=\nu(f) \\
\gamma & \text{otherwise} \end{cases},\]
there are four cases to check to see that identity \ref{eq: f edge} implies
\begin{align} \label{eq: beta}
 \begin{vmatrix} \calN(\gamma) & \calQ(\gamma) \\ \calN(e^-) & \calQ(e^-) \end{vmatrix} \ = \ \nu(e) \, \dQ(\beta) \, \calQ\left( \hat{f} \right).
 \end{align}
Combining identities \ref{eq: linear}, \ref{eq: hat}, and \ref{eq: beta} gives
\begin{align*}
\begin{vmatrix} \calN(e^-) & \calQ(e^-) \\ \calN(e^+) & \calQ(e^+) \end{vmatrix} &= \dQ(e^-) \, \nu(e) \, \dQ(\beta) \, \calQ\left( \hat{f} \right) - \calQ(e^-) \, \nu(e) \, \dQ(e^-) \, \calQ(\gamma) \\
&= -\nu(e) \, \dQ(e^-) \, \left( - \dQ(\beta) \, \calQ\left( \hat{f} \right) + \calQ(e^-) \, \calQ(\gamma) \right) \\
&= -\nu(e) \, \dQ(e^-) \, \calQ(\hat{e}).
\end{align*}
This proves the lemma in the case $i=1$.

Combined with the $i=1$ case, the equalities 
\begin{align*}
\calQ(e^- \oplus^2 \hat{e}) &= -\dQ(\hat{e}) \, \calQ(e^-) + \calQ(\hat{e}) \, \calQ(e^+) \\
\calN(e^- \oplus^2 \hat{e}) &= -\dQ(\hat{e}) \, \calN(e^-) + \calQ(\hat{e}) \, \calN(e^+)
\end{align*}
imply that
\begin{align*} \begin{vmatrix} \calN(e^-) & \calQ(e^-) \\ \calN(e^- \oplus^2 \hat{e}) & \calQ(e^- \oplus^2 \hat{e}) \end{vmatrix} &=
\calQ(\hat{e}) \, \begin{vmatrix} \calN(e^-) & \calQ(e^-) \\ 
  \calN(e^+) &   \calQ(e^+) \end{vmatrix} \\
  &= - \nu(e)\, \dQ(e^-) \, \calQ(\hat{e})^2.
\end{align*}
\endproof

\begin{lemma}  \label{lem: Vshape}
Suppose $e \in \calE$ and let $\beta$ and $\gamma_0$ be its left and right endpoints.  Then
\[ \frac{\calV(\gamma_0)-\calV(\gamma_0 \oplus \beta)}{\calV(\beta)-\calV(\gamma_0 \oplus \beta)} = \calZ(e).\]
\end{lemma}

\proof
Write $\gamma_1=\gamma_0 \oplus \beta$ and let $e_1$ and $e_2$ be the edges from $\gamma_0$ to $\gamma_1$ and from $\beta$ to $\gamma_1$.  Then, using Lemma \ref{lem: det2},
\[ \calV(\gamma_0)-\calV(\gamma_1) = \frac{\calN(\gamma_0)\calQ(\gamma_1)-\calN(\gamma_1) \calQ(\gamma_0)}{\calQ(\gamma_0)\calQ(\gamma_1)} = \frac{-\nu(e_1) d_\calQ(\gamma_0) \calQ(\beta)}{\calQ(\gamma_0)\calQ(\gamma_1)}\]
and
\[\calV(\beta)-\calV(\gamma_1) = \frac{-\nu(e_2) d_\calQ(\beta) \calQ(\gamma_0)}{\calQ(\beta)\calQ(\gamma_1)}.\]
Using that $\nu(e_1)$ and $\nu(e_2)$ have opposite signs and, after considering the two cases for the sign of $\nu(e)$, this gives
\[ \frac{\calV(\gamma_0)-\calV(\gamma_1)}{\calV(\beta)-\calV(\gamma_1)} = -\frac{d_\calQ(\gamma_0) \calQ(\beta)^2}{d_\calQ(\beta)\calQ(\gamma_0)^2} = \calZ(e).\]  
\endproof

Similar arguments prove the next lemma.  

\begin{lemma}  \label{lem: Vshape2}
Suppose $\{ \beta, \gamma_0\} \subset \Q_0$ is a Farey pair with $\beta < \gamma_0$.  Let $\gamma_1=\gamma_0 \oplus \beta$ and $e$ be the edge between $\gamma_0$ and $\gamma_1$.  Then
\[ \frac{\calV(\gamma_0)-\calV(\beta)+1}{\calV(\gamma_1)-\calV(\beta)} = \calZ(e).\]
\end{lemma}

\proof[Proof of Lemma \ref{lem: delta parameters}]
Suppose that $\{ v_1, v_2, v_3 \}$ are distinct points in $\C$ ordered in agreement with the right hand orientation of the triangle they span.  The edge parameter for the edge $[v_1, \infty]$ of the ideal tetrahedron $[ v_1, v_2, v_3, \infty ]$ is given by the cross-ratio formula $(v_3-v_1)/(v_2-v_1)$.  Using this, together with Lemmas \ref{lem: Vshape} and \ref{lem: Vshape2}, the result follows. 
\endproof

The next lemma is a useful extension of Lemma \ref{lem: adjacent zeros} which follows easily from Lemma \ref{lem: det2}.  

\begin{lemma} \label{lem: more zeros}
If $e \in \calE$, then none of $\calQ(e^-)$, $\calQ(e^+)$, and $\calQ(e^+\oplus \hat{e})$ share a zero.  Also, for $\gamma \in \Q_0$, $\calQ(\gamma)$ and $\calN(\gamma)$ do not have a common zero.
\end{lemma}

\proof
Write $\gamma_0=e^-$ and $\gamma_j = \gamma_0 \oplus^j \hat{e}$.  Suppose $n \in \{ 1, 2 \}$.  If $\calQ(\gamma_0)$ shares a zero $\frz$ with either $\calQ(\gamma_n)$ or $\calN(\gamma_0)$ then, by Remark \ref{rem: constant term}, $\dQ(\gamma_0)_\frz \neq 0$.  Hence, Lemma \ref{lem: det2} implies that $\calQ(\hat{e})_\frz=0$.  But then, the edge between $\hat{e}$ and $\gamma_0$ provides a contradiction to Lemma \ref{lem: adjacent zeros}.
\endproof

The last part of Lemma \ref{lem: more zeros} shows that, for any $\gamma \in \wh{\Q}_0$ and any $\frz \in \C$,
\[ \calV(\gamma)_\frz \ =\ \frac{\calN(\gamma)_\frz}{\calQ(\gamma)_\frz}. \]

For additional perspective with regard to Lemma \ref{lem: more zeros}, it is worth mentioning again the discussion following Lemma \ref{lem: adjacent zeros}.  It points out that if $e$ is the edge between $1/3$ and $1/4$, then $\calQ(e^-)$ and $\calQ(e^+ \oplus^2 \hat{e})$ have a common zero.  

\section{Domains} \label{sec: domains}

Assume $\alpha \in \Q \cap (0,1/2)$ and $\alpha \neq 1/3$.  As usual, let $e_0$ and $e_{n+1}$ be the top and bottom edges of the funnel $F_\alpha$ and let $\{e_1,\ldots,e_n\}$ be the sequence of interior edges of $F_\alpha$ listed from top to bottom.  Write $\delta_j=\delta(e_j)$ and $\delta'_j=\delta'(e_j)$ as in
Definition~\ref{def: generic tetrahedra}.

Fix $\frz\in\calH(M_\alpha^\circ)$.  This assigns to each tetrahedron in $\Delta_1$ an ideal hyperbolic tetrahedron with shape parameter $\frz$.  Choose one such tetrahedron $\delta$.  By definition, its top edge is labeled $1$, its side edges are labeled $0$ and $1/2$, and its bottom edge is labeled $1/3$. By Lemma~\ref{lem: delta parameters}, $(\delta_1)_\frz$ is an ideal hyperbolic tetrahedron in $\IH^3$ with the same edge labeling and whose top edge has parameter $\frz$.  Thus, there is a label preserving isometry from $\delta$ to $(\delta_1)_\frz$.

\begin{dfn} \label{def: generic domain}
Using the discussion above, define the {\it generic domain} for $\alpha$ to be the $2n$-tuple of triples,
\[\Omega(\alpha) =\left( \delta_1, \delta'_1, \ldots, \delta_n, \delta'_n \right) \ \in \ \left( \Q(x)^3 \right)^{2n}.\]
Also, given $\frz \in \calH(M_\alpha^\circ)$, 
define $\Omega(\alpha)_\frz$ to be the $2n$-tuple of ideal tetrahedra obtained by evaluating the triples $\delta_j$ and $\delta'_j$ at $\frz$.  
\end{dfn}
\index{fundamental domain!generic|idxdef}
\index{fundamental domain!specialized|idxdef}

\begin{lemma} \label{lem: delta pairs}
Suppose $\tau_0$ and $\tau_1$ are terms in $\Omega(\alpha)$.   The intersection $\tau_0 \cap \tau_1$ consists of a pair of elements if and only if one of the following holds.
\begin{enumerate}
\item $\{ \tau_0, \tau_1 \} = \{ \delta_j, \delta_{j-1} \}$ where $j \in \{ 2, \ldots, n \}$.  In each $\tau_i$, the shared pair is labeled by the number $e_{j-1}^- \oplus e_{j-1}^+$. 
\item $\{ \tau_0, \tau_1 \} = \{ \delta_j', \delta_{j-1}' \}$ where $j \in \{ 2, \ldots, n \}$.  In each $\tau_i$, the shared pair is labeled by the number $\hat{e}_j$. 
\item $\{ \tau_0, \tau_1 \} = \{ \delta_j, \delta_{j-1}' \}$ where $j \in \{ 2, \ldots, n \}$ and $e_{j-1} \cap e_j$ lies on the left side of $F_\alpha$. In each $\tau_i$, the shared pair is labeled by $e_j^-$, the number at the intersection of $e_{j-1}$ and $e_j$. 
\item $\{ \tau_0, \tau_1 \} = \{ \delta_1, \delta_1' \}$. The shared pair is $\{ 1, -x+1\}$ and is labeled $0$ in $\delta_1$ and labeled $1$ in $\delta_1'$.
\end{enumerate}
Likewise, the intersection $(1+\tau_0) \cap \tau_1$ consists of a pair of elements if and only if $\{ \tau_0, \tau_1 \} = \{ \delta_j, \delta_{j-1}' \}$, where $e_j$ is a spoke for a hub on the right of $F_\alpha$.  Here, the intersection is labeled by the hub.  
\end{lemma}

\proof
The intersections listed in the statement follow directly from the definitions.  It remains to show that there are no others.

Choose a hyperbolic $\beta\in\Q\cap(0,1/2)$ such that $F_\alpha\subset F_\beta$.  By \cite{GF}, the Sakuma--Weeks triangulation of $M_\beta$ is geometric.  Let $\frz\in\calH(M_\beta)$ be its first shape parameter and normalize the developing map so that a tetrahedron in $\Delta_1$ develops to
$(\delta_1)_\frz$, preserving labels.

By Lemma~\ref{lem: delta parameters}, the specialized generic tetrahedra have the same shapes and edge labels as the corresponding tetrahedra in the Sakuma--Weeks triangulation.  Using the intersections already verified above, they may be attached successively to $(\delta_1)_\frz$ across the prescribed faces. Thus, the specializations of the terms of $\Omega(\alpha)$ are the corresponding developed lifts of tetrahedra in the Sakuma--Weeks triangulation of $M_\beta$.  The unit translates in the statement are treated similarly using the simplicial symmetry induced by $z\mapsto z+1$ on the triangulated punctured plane.
\index{Sakuma--Weeks complex!triangulation}

If there were any additional intersection consisting of a pair of elements, either between $\tau_0$ and $\tau_1$ or between $1+\tau_0$ and $\tau_1$, then after specializing at $\frz$ the corresponding developed tetrahedra would have an additional common nonfinite face.  Since the developing map for the complete structure on $M_\beta$ is an isometry, this would give an additional face adjacency in the lifted Sakuma--Weeks triangulation, contrary to its combinatorics.
\endproof

\begin{figure}
\setlength{\unitlength}{.1in}
\begin{picture}(21,25)
\put(0,0) {\includegraphics[width= 2.1in]{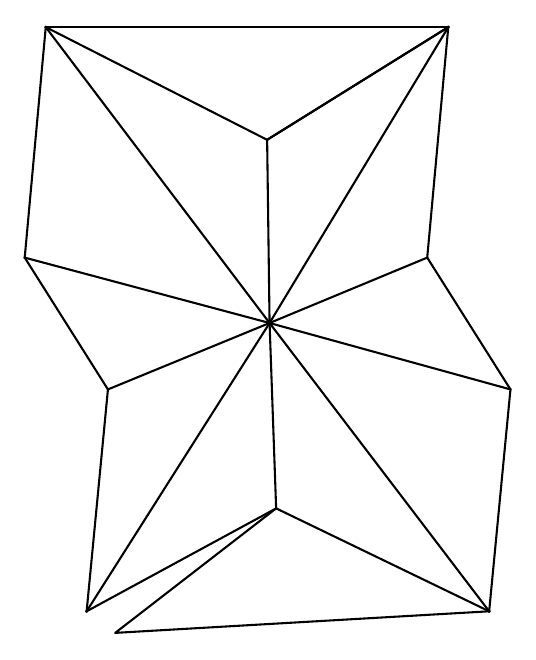}}
\put(10.1,22.3){$\delta_1$}
\put(7.8,19){$\delta_2$}
\put(3.5,17.5){$\delta_3$}
\put(4.6,12.6){$\delta_4$}
\put(5.1,8){$\delta_5$}
\put(8,6.8){$\delta_6$}
\put(12,19){$\delta'_1$}
\put(14.8,17.5){$\delta'_2$}
\put(14.9,12.8){$\delta'_3$}
\put(16.5,8){$\delta'_4$}
\put(12.2,6.8){$\delta'_5$}
\put(10.2,2.9){$\delta'_6$}
\end{picture}
\caption{\index{fundamental domain!specialized} The domain $\Omega(10/33)_\frz$ where $\frz=.45+.28i$.}  \label{fig: split domain}
\end{figure}

\begin{figure}
\setlength{\unitlength}{.1in}
\begin{picture}(31,26)
\put(3.5,.3) {\includegraphics[width= 2.5in]{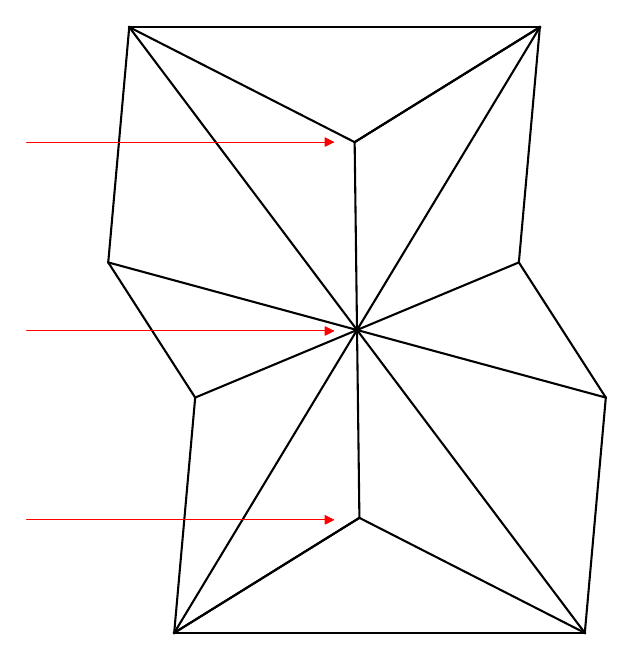}}
\put(.7,20.5){$\scriptstyle \calV(1/2)$}
\put(.7,13){$\scriptstyle \calV(1/3)$}
\put(.1,5.6){$\scriptstyle \calV(4/13)$}
\put(6.5,26){$\scriptstyle -1+\calV(0)$}
\put(24,26){$\scriptstyle \calV(1)$}
\put(25,16){$\scriptstyle \calV(1/4)$}
\put(1.7,16){$\scriptstyle -1+\calV(1/4)$}
\put(28,10.5){$\scriptstyle \calV(2/7)$}
\put(5,10.5){$\scriptstyle -1+\calV(2/7)$}
\put(25.2,0){$\scriptstyle \calV(3/10)$}
\put(5,0){$\scriptstyle -1+\calV(3/10)=\calV(7/23)$}
\end{picture}
\caption{\index{fundamental domain!specialized} The domain $\Omega(10/33)_\frz$ where $\frz$ is the unique element of $\calH\left(M_{10/33}\right)$.  Here, $\frz \sim 0.451069913 + 0.280155498 \, i$.
   }  \label{fig: domain}
\end{figure}

Recall from Definitions~\ref{def: generic tetrahedra} and \ref{def: generic domain} that, when passing from $\Omega(\alpha)$ to $\Omega(\alpha)_\frz$, each triple $\delta_j$ or $\delta_j'$ in $\Omega(\alpha)$ becomes an ideal hyperbolic tetrahedron whose edge labeling agrees with that of the tetrahedra in $\Delta_j$.  (See Figures~\ref{fig: split domain} and~\ref{fig: domain} for $\Omega(10/33)_\frz$ with two choices of $\frz$.)  In this setting, Lemma~\ref{lem: delta pairs} shows that certain nonfinite faces of the tetrahedra in $\Omega(\alpha)_\frz$ either coincide or differ by a unit translation.  With their edge labelings, these face identifications agree with the corresponding identifications among the tetrahedral pairs in the tetrahedral core of $M_\alpha^\circ$.

By the discussion following Definition~\ref{def: hyperbolic structures}, the shape parameters determined by $\frz$ satisfy the edge gluing
equations on $M_\alpha^\circ$.  Hence, after choosing a lift of the tetrahedron $\delta\subset\Delta_1$ fixed above, its label preserving isometry to  $(\delta_1)_\frz$ extends across faces to a well-defined developing map
\[ \dev_\frz\co\wt M_\alpha^\circ\longrightarrow\IH^3 \]
with an associated holonomy representation
\[ \rho_\frz\co\pie(M_\alpha^\circ)\longrightarrow\PSL_2(\C). \]
When $\alpha$ is hyperbolic, this is the usual developing map and holonomy representation.
\index{developing map}
\index{holonomy}

Lemma~\ref{lem: delta pairs} now gives the following theorem.

\begin{thm}\label{thm: domain}
\DomainThm
\end{thm}
\index{fundamental domain}

\begin{remarks}\label{rem: developed domain}
\[ \] 
\vspace{-.3in}
\begin{enumerate}
\item Theorem~\ref{thm: domain} asserts that $\Omega(\alpha)_\frz$ is the developed image of a connected fundamental domain in $\widetilde M_\alpha^\circ$.  For a general point $\frz\in\calH(M_\alpha^\circ)$, this does not imply that the restriction of the developing map to this fundamental domain is globally injective, or that the tetrahedra of $\Omega(\alpha)_\frz$
form an embedded polyhedral region in $\IH^3$.  When $\frz$ determines the complete structure on $M_\alpha$, the developing map extends to an isometry $\widetilde M_\alpha\to\IH^3$, and no such distinction arises.
\item The degenerate tetrahedra $(\delta_0)_\frz$ and $(\delta'_0)_\frz$
arise naturally in $\Omega(\alpha)_\frz$.  By
Definition~\ref{def: generic tetrahedra}, $(\delta_0)_\frz$ is the
geodesic triangle $[0,1,\infty]$, which is a face of
$(\delta_1)_\frz$, while $(\delta'_0)_\frz$ is the geodesic triangle
$[1,-\frz+1,\infty]$ in
$(\delta_1)_\frz\cap(\delta'_1)_\frz$.
\end{enumerate}
\end{remarks}

\section{Completeness} \label{sec: completeness}

Throughout this section, assume that $\alpha$ is hyperbolic.

Recall from Definition~\ref{def: hyperbolic structures} that $\calH(M_\alpha)$ consists of the positive $G$-invariant solutions to the gluing equations carried by the Sakuma--Weeks triangulation of $M_\alpha$.  Following Theorem~\ref{thm: main}, $\calH(M_\alpha)$ is identified with a subset of $\IH^2$ by interpreting $\frz\in\calH(M_\alpha)$ as the shape parameter of the tetrahedra in $\Delta_1$.  In particular,
\[ \calH(M_\alpha) \ =\ \left\{ \frz\in\calH(M_\alpha^\circ) \,\Big|\, \calQ(\alpha)_\frz=0 \right\}\]
where
\[ \calH(M_\alpha^\circ) \ =\  \left\{ \frz\in\IH^2 \,\Big|\, \calZ(e_j)_\frz\in\IH^2 \text{ for every interior edge $e_j$ of $F_\alpha$} \right\}. \]

\begin{thm} \label{thm: complete}
Suppose $\alpha$ is hyperbolic.  Every element of $\calH(M_\alpha)$ determines a complete hyperbolic structure on $M_\alpha$.  Consequently, $\calH(M_\alpha)$ consists of exactly one element.
\end{thm}

The proof for this theorem is given at the end of the current section.

\begin{dfn} \label{def: geometric root}
If $\alpha$ is hyperbolic, then the unique element of $\calH(M_\alpha)$ is called the {\it geometric root} $\frz(\alpha)$ of $\calQ(\alpha)$.  If $\alpha$ is clear by context, it may be denoted as $\frz_1$.
\index{geometric root!2-bridge link|idxdef}
\end{dfn}

\begin{remark} \label{rk: filter} In practice, approximations for geometric roots can be found by filtering the roots of $\calQ(\alpha)$ by the requirement that $\calZ(e_j)_\frz$ has positive imaginary part for each $j \in \{ 1, \ldots, n \}$.   Theorem \ref{thm: complete} implies that, if $\alpha$ is hyperbolic, this will always yield a single possibility. 
\end{remark}

After completeness is established, the final statement of  Theorem \ref{thm: complete} follows from Mostow rigidity.  The completeness argument uses the next lemma.  See Figure \ref{fig: domain}, to see how (at least in the case of $\alpha=10/33$) this lemma follows from the gluing equation for $4/13$.  The proof given below uses properties of Farey recursive functions.

\begin{lemma} \label{lem: difference}
If $\alpha$ is hyperbolic and $\frz \in \calH(M_\alpha)$, then
\[ -1 + \calV\left(e_{n+1}^L \right)_{\frz} \ =\ \calV\left(e_{n+1}^R \right)_{\frz}.\]
\end{lemma}

\proof
Write $\beta=\hat{e}_{n+1}$, $\gamma_0=e_{n+1}^-$ and $\gamma_1=e_{n+1}^+$.  By Lemma \ref{lem: det2}, 
\[ \calV(\gamma_1) -\calV(\gamma_0) = \nu(e_{n+1}) \frac{d_\calQ (\gamma_0) \calQ(\beta)}{\calQ(\gamma_0) \calQ(\gamma_1)}.\]
Hence, using Farey recursion,
\begin{align*}
\calV(\gamma_1) -\calV(\gamma_0) &= -\nu(e_{n+1}) \left( \frac{-d_\calQ(\gamma_0) \calQ(\beta) + \calQ(\gamma_0) \calQ(\gamma_1)}{\calQ(\gamma_0) \calQ(\gamma_1)} \right) + \nu(e_{n+1}) \\
&= -\nu(e_{n+1}) \frac{\calQ(\alpha)}{\calQ(\gamma_0) \calQ(\gamma_1)} + \nu(e_{n+1}).
\end{align*}
Now, $\calQ(\alpha)$ is satisfied by $\frz$ and Lemma \ref{lem: adjacent zeros} guarantees that neither $\calQ(\gamma_0)$ nor $\calQ(\gamma_1)$ is zero at $\frz$.  So, together with the definitions of $e^L$ and $e^R$, this proves the lemma.  
\endproof

\begin{cor} \label{cor: delta n+1}
Suppose $\alpha \in \Q \cap (0,1/2)$ and $\frz \in \calH(M_\alpha)$.  Then $(\delta_{n+1})_{\frz}$ is a geodesic face of the tetrahedron $(\delta_n)_{\frz}$, and $(\delta'_{n+1})_{\frz}$ is a geodesic face of the tetrahedron $(\delta'_n)_{\frz}$.  Moreover, if $\nu(e_{n+1}) = 1$, then
\[ (\delta_{n+1})_{\frz} \ =\ \left[ -1+\calV(e_n^L)_{\frz}, \, \calV(e_n^R)_{\frz}, \, \infty \right] \]
and 
\[ (\delta'_{n+1})_{\frz} \ =\ \left[ \calV(e_n^R)_{\frz}, \, \calV(e_n^- \oplus e_n^+)_{\frz}, \, \infty \right]. \]
Otherwise, $\nu(e_{n+1})=-1$ and 
\[ (\delta_{n+1})_{\frz} \ =\ \left[ \calV(e_n^R)_{\frz}, \, \calV(e_n^- \oplus e_n^+)_{\frz}, \, \infty \right] \]
and 
\[ (\delta'_{n+1})_{\frz} \ =\ \left[ \calV(e_n^- \oplus e_n^+)_{\frz}, \, \calV(e_n^L)_{\frz}, \, \infty \right]. \]

\end{cor}

\proof
By definition of $\calH(M_\alpha)$ and Lemma \ref{lem: delta parameters}, $(\delta_n)_{\frz}$ and $(\delta'_n)_{\frz}$ are ideal hyperbolic tetrahedra.  Also, 
\[\delta_{n+1} \ = \ \left\{ \epsilon + \calV\left( \hat{e}_{n+1} \right), \, -1+\calV\left(e^L_{n+1}\right), \, \calV\left(e^R_{n+1}\right) \right\} \]
and
\[ \delta'_{n+1} \ = \ \left\{ \calV\left(e^-_{n+1}\right), \, \calV\left(e^+_{n+1}\right), \, \calV(\alpha) \right\}\]
where $\epsilon$ depends on the slope of $e_{n+1}$.  Since $\frz \in \calH(M_\alpha)$, $\calQ(\alpha)_{\frz}=0$ and so 
\[ (\delta'_{n+1})_{\frz} \ = \ \left[ \calV\left(e^-_{n+1}\right)_{\frz}, \, \calV\left(e^+_{n+1}\right)_{\frz}, \, \infty \right].\]
By Lemma \ref{lem: difference}, 
\[(\delta_{n+1})_{\frz} \ = \ \left[ \epsilon + \calV\left( \hat{e}_{n+1} \right)_{\frz}, \, \calV\left(e^R_{n+1}\right)_{\frz}, \, \infty \right]. \]
The corollary now follows from the fact that, if $\nu(e_{n+1})=1$, then
\begin{align*}
\hat{e}_{n+1} &= e_n^L & e_{n+1}^L &= e_n^- \oplus e_n^+ & e_{n+1}^R &= e_n^R
\end{align*}
and, if $\nu(e_{n+1})=-1$, then
\begin{align*}
\hat{e}_{n+1} &= e_n^R & e_{n+1}^L &= e_n^L & e_{n+1}^R &= e_n^- \oplus e_n^+.
\end{align*}
\endproof

\begin{figure}
\setlength{\unitlength}{.1in}
\begin{picture}(49,13)
\put(0,0) {\includegraphics[width= 4.9in]{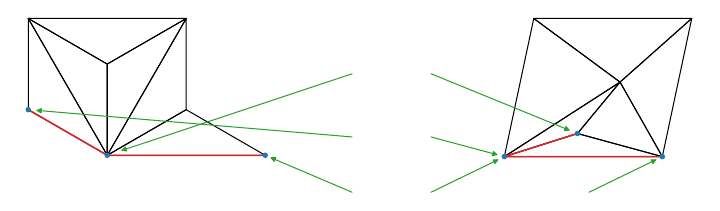}}
\put(24.2,.8){$\scriptscriptstyle \calV\left( e_n^- \oplus e_n^+\right)$}
\put(24.2,4.8){$\scriptscriptstyle -1+\calV\left( e_n^L\right)$}
\put(25.3,9.4){$\scriptscriptstyle \calV\left( e_n^R\right)$}
\put(37,.5){$\scriptscriptstyle \calV\left( e_n^L\right)$}
\end{picture}
\caption{\index{fundamental domain!specialized} This figure illustrates Corollary \ref{cor: delta n+1} and the comments following its proof.  Taking $\frz \in \calH(M_\alpha)$ for $\alpha \in \{ 3/10, 4/11 \}$, the left figure shows $\Omega(3/10)_{\frz}$ and the right shows $\Omega(4/11)_{\frz}$.  The geodesic triangles $(\delta_{n+1})_{\frz}$ and $(\delta'_{n+1})_{\frz}$ lie over the lines highlighted in magenta.  In both cases, the triangle $(\delta'_{n+1})_\frz$ lies over the horizontal magenta line segment.  The geodesic triangles $(\delta_0)_\frz$ and $(\delta'_0)_\frz$ are highlighted in green. 
   }  \label{fig: BottomFold}
\end{figure}

It is worth noticing that 
\[ \calV(e_n^- \oplus e_n^+)_{\frz} \ = \ 
\begin{cases}
\calV(e_n^R)_\frz+1 & \text{if } \nu(e_{n+1})=1 \\
\calV(e_n^L)_\frz -1& \text{otherwise}.
\end{cases}\]
So, if $\nu(e_{n+1})=1$, then $1+(\delta_{n+1})_{\frz}$ is a face of $(\delta'_n)_{\frz}$ and, if $\nu(e_{n+1})=-1$, then $(\delta_{n+1})_{\frz}$ is the intersection of $(\delta_n)_{\frz}$ and $(\delta'_n)_{\frz}$.  For $\alpha \in \{ 3/10, 4/11 \}$, these observations, together with Corollary \ref{cor: delta n+1}, are illustrated in Figure \ref{fig: BottomFold}.  

It is also worthwhile, at this point, to define specific elements of $\pglQx$ which specialize to lifts of the upper and lower involutions in $G$ (Definition \ref{def: involutions}).

\begin{dfn} \label{def: concrete involutions}
In $\pglQx$, define 
\[ \sigma_u \ = \ \begin{bmatrix} -1&1\\0&1 \end{bmatrix} \qquad \text{and} \qquad \sigma_\ell \ = \ \begin{bmatrix} -1 & \calV\left( e_{n+1}^- \right)+\calV\left( e_{n+1}^+ \right) \\ 0 & 1 \end{bmatrix}. \]
\end{dfn}

Both of these M\"obius transformations have order two.  The involution $\sigma_u$ is the unique order-2 isometry which fixes $\infty$ and interchanges the pair in $\delta_1$ labeled $1/2$ (see Definition \ref{def: generic tetrahedra}).  If $\frz \in \calH(M_\alpha^\circ)$, then $\sigma_u$ acts by isometry on the associated geometric Sakuma--Weeks triangulation as the upper involution.  The involution $\sigma_\ell$ is the unique order-2 isometry which fixes $\infty$ and interchanges the pair in $\delta_n'$ labeled $\hat{e}_{n+1}$.  If $\frz \in \calH(M_\alpha)$, it acts on the associated geometric Sakuma--Weeks triangulation as the lower involution.  These involutions act on the Riemann sphere $\bound \IH^3$ as
\[ z \mapsto 1-z \qquad \text{and} \qquad z \mapsto \calV\left( e_{n+1}^- \right)_\frz+\calV\left( e_{n+1}^+ \right)_\frz-z.\]

\proof[Proof of Theorem \ref{thm: complete}] By \cite{GF}, the Sakuma--Weeks triangulation of $M_\alpha$ is geometric, so $\calH(M_\alpha)$ is nonempty.  If $\frz \in \calH(M_\alpha)$ then the tetrahedra in $\Omega(\alpha)_\frz$ all share the ideal vertex $\infty$ so the intersection $\Upsilon(\alpha)_\frz$  of $\Omega(\alpha)_\frz$ with a horosphere $H$ centered at $\infty$ is a collection of Euclidean triangles.  The vertices for a given triangle are obtained from a corresponding triple $D \in \Omega(\alpha)$ by evaluating the terms of $D$ at $\frz$.

Let $\omega_1, \ldots, \omega_k$ be the sequence of distinct hubs for $F_\alpha$, listed in order with increasing denominators.  
By definition, hubs with even indices are on the right and hubs with odd indices are on the left.  It is clear that $0=\omega_1$ is a left hub.  Since $\alpha$ is hyperbolic, there must be at least two hubs and at least one right hub, namely $\omega_2$.

In $
\pslC$, define
\[ U_0 \ =\ \begin{bmatrix} 1&-1 \\ 0 & 1 \end{bmatrix}\]
with its action on $\Q(x) \cup \{ \infty \}$ by M\"obius transformation.  Then, for every right hub $\omega_j$ and every spoke $e_i$ for $\omega_j$, $U_0$ fixes $\infty$ and, by Lemma \ref{lem:  delta pairs}, takes the pair in $\delta_i$ labeled $\omega_j$ to the similarly labeled pair in $\delta'_{i-1}$.  By definition of the Sakuma--Weeks triangulation and because $F_\alpha$ has a right hub, this implies that $U_0$ belongs to the holonomy group for $\frz$ for every $\frz \in \calH(M_\alpha^\circ)$.  This is evident in Figures \ref{fig: split domain}, \ref{fig: domain}, and \ref{fig: BottomFold}.

Now, assume $\frz \in \calH(M_\alpha)$ and take involutions $\sigma_u$ and $\sigma_\ell$ as given in Definition \ref{def: concrete involutions}.  Then the upper and lower involutions in $G$, acting on $M_\alpha$, lift to $(\sigma_u)_\frz$ and $(\sigma_\ell)_\frz$.  The upper involution acts as an order-2 rotation about the upper core tunnel, the edge of $M_\alpha$ labeled by the pair $\{0,1 \}$.  It acts by inversion on each of the edges labeled $1/2$.  Similarly, the lower involution acts as an order-2 rotation about the lower core tunnel, the edge of $M_\alpha$ labeled by the pair $\{e_{n+1}^-, e_{n+1}^+ \}$ and acts by inversion on each of the edges labeled $\hat{e}_{n+1}$. 

The involutions $(\sigma_u)_\frz$ and $(\sigma_\ell)_\frz$ restrict to involutions on the horosphere $H$.  As such, $
(\sigma_u)_\frz$ inverts the horizontal edge $(\delta_0)_\frz \cap H$ between $0$ and $1$ and $(\sigma_\ell)_\frz$ inverts the edge $(\delta'_{n+1})_\frz \cap H$ between $\calV\left( e_{n+1}^- \right)_\frz$ and $\calV\left( e_{n+1}^+ \right)_\frz$ which, by Lemma \ref{lem: difference}, is also horizontal.

It follows that the holonomy group for the cusp of $O_\alpha$ is generated by $U_0$, $(\sigma_u)_\frz$, and $(\sigma_\ell)_\frz$.  The group generated by $U_0$ and the parabolic isometry 
\[ z \mapsto z+2 \, \calV\left( e_{n+1}^- \right)_\frz \]
is an index-2 subgroup of the holonomy group for the cusp of $O_\alpha$.  This shows that the torus cusps of $M_\alpha$ are complete and establishes the theorem.
\endproof

\section{Holonomy and crossing circles} \label{sec: holonomy}

As usual, assume $\alpha \in \Q \cap (0,1/2)$ and $\alpha \neq 1/3$.

This section describes a generic holonomy representation 
\[ \pie M_\alpha^\circ \longrightarrow \pslQx\]
which specializes to the holonomy for $\frz \in \calH(M_\alpha^\circ)$ by setting $x=\frz$.  This will help to make concrete connections between the algebra, topology, and geometry for elements of $\calH(M_\alpha^\circ)$.  This representation can be considered as a monomorphism of the free group of rank two to $\pslQx$ and is independent of $\alpha$.

\subsection{Definitions, notation, and lemmas} \label{sec: defs and lems}

\begin{figure}
\setlength{\unitlength}{.1in}
\begin{picture}(48,20)
\put(0,12.3){$k_0$}
\put(0,6){$k_1$}
\put(1.7,1) {\includegraphics[width= 4.5in]{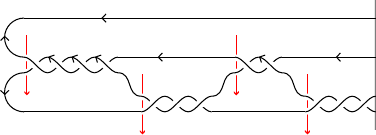}}
\put(5.7,6){$\redden{S_{0}}$}
\put(6.5,11.4){$T_{1}$}
\put(9.1,11.4){$T_{1/2}$}
\put(12.9,11.4){$T_{1/3}$}
\put(13.7,16){$U_0$}
\put(19.5,1.3){$\redden{S_{1/4}}$}
\put(20.6,11.4){$T_{1/4}$}
\put(30.8,6){$\redden{S_{3/13}}$}
\put(32,11.4){$T_{4/17}$}
\put(39.3,1.3){$\redden{S_{7/30}}$}
\put(42,11.4){$T_{7/30}$}
\end{picture}
\caption{Here $\alpha = 24/103$.  The fundamental group $\pie M_{\alpha}^\circ$ is freely generated by $k_0$ and $k_1$.  Certain elements of $\pie M_{\alpha}^\circ$ visible in the diagram may be identified with favorite elements of $\Gamma$ using the isomorphism $\pie M_{\alpha}^\circ \to \Gamma$. }  \label{fig: SW link II}
\end{figure}

\begin{dfn} \label{def: generators}
Consider the truncated diagram $D_\alpha$ from Definition~\ref{def: diagram}.  Take $k_0,k_1\in\pie M_\alpha^\circ$ corresponding to the two leftmost
arcs shown in Figure~\ref{fig: SW link II}.  The diagram can always be oriented so that the upper arc, labeled $k_0$, is oriented in the clockwise direction and the lower arc is oriented in the counterclockwise direction.  The elements $k_0$ and $k_1$ generate $\pie M_\alpha^\circ$ and are referred to as the {\it canonical
generators} associated to $\alpha$. \index{canonical generators|idxdef}

When $\alpha$ is hyperbolic, the canonical generators may also be viewed as elements of $\pie M_\alpha$, where they continue to generate.
\end{dfn}

\begin{dfn} \label{def: UW}
For $\omega \in \Q_0$, write
\begin{align*}
 d&= d_\calQ(\omega) &Q&= \calQ(\omega) & N&= \calN(\omega) & V&= \calV(\omega)
 \end{align*}
 and define
\[
W_\omega \ =\ \begin{bmatrix} 1&0 \\ \frac{-Q^2}{d}&1 \end{bmatrix} \qquad \text{and} \qquad U_\omega \ =\ \begin{bmatrix} 1&-V \\ 0 & 1 \end{bmatrix}
\]
in $\pslQx$.  
\end{dfn}

Notice that the definition agrees with that of $U_0$ from Section \ref{sec: completeness}, so
\[ U_0^{-1} \ =\ \begin{bmatrix} 1&1\\0&1 \end{bmatrix}.\]
Also,
\[ W_0 \ =\ \begin{bmatrix} 1&0\\-1/x & 1 \end{bmatrix}.\]

\begin{remark} \label{rem: specializing} If $\frz \in \calH(M_\alpha^\circ)$, then $\calZ(e_j)_\frz \in \IH^2$ for each interior edge $e_j$ of $F_\alpha$ and each $\frz \in \calH(M_\alpha^\circ)$.  In particular, both $\dQ(\omega_i)_\frz$ and $\calQ(\omega_i)_\frz$ are nonzero for every hub $\omega_i$ of $F_\alpha$.  Hence, the specialized elements $\left( W_{\omega_i} \right)_\frz$ and $\left( U_{\omega_i} \right)_\frz$ are well-defined elements of $\pslC$.  Furthermore, the same is true for any elements of $\pslQx$ defined as products of elements $U_\omega$ and $W_\omega$. \index{funnel!hub}
\end{remark} 

\begin{dfns} \label{def: generic holonomy}
The {\it generic holonomy group} is defined as 
\[ \Gamma \ =\ \left\langle U_0, W_0 \right\rangle \ < \ \pslQx.\]
It is a standard fact that $\Gamma$ is free of rank two.  The fundamental group of $M_\alpha^\circ$ is also free of rank two, generated by its canonical generators $k_0$ and $k_1$.  The {\it generic holonomy representation} is the isomorphism
\[ \varphi \co \langle k_0, k_1 \rangle \longrightarrow \Gamma < \pslQx\]
given by
\[ k_0 \mapsto U_0^{-1} \qquad \text{and} \qquad k_1 \mapsto W_0.\] 
If $\frz \in \C$ is nonzero, then $\frz$ determines a homomorphism $\Gamma \to \pslC$ by specializing $x$ to $\frz$.  Let $\Gamma_\frz$ denote its image and decorate elements of $\Gamma$ with the subscript $\frz$ to denote their images in $\Gamma_\frz$.  
\end{dfns}
\index{holonomy!generic|idxdef}
\index{holonomy!generic representation|idxdef}

As usual, the group $\Gamma$ acts on $\Q(x) \cup \left\{ \infty \right\}$ by M\"obius transformations and this action specializes to the familiar action of $\Gamma_\frz$ on $\bound \IH^3$.  

The next main goal is to see that, for every $\frz \in \calH(M_\alpha^\circ)$, $\varphi$ specializes to the holonomy representation for $\frz$, pairing the faces of the domain $\Omega(\alpha)_\frz$.  Using the apparatus set up in the first part of this paper, it is also possible to obtain formulas for the elements of $\Gamma_\frz$ which identify the faces of the domain $\Omega(\alpha)_\frz$, edges of the triangles in the lifts of the punctured spheres $\Sigma_j$ to $\Omega(\alpha)_\frz$, and certain useful elements of $\pie M_\alpha^\circ$ which are easily visible in the diagram $D_\alpha$.  Some of these elements are shown in Figure \ref{fig: SW link II}.

In what follows, it is helpful to notice that 
\begin{itemize}
\item every $W_\omega$ fixes zero and every $U_\omega$ fixes $\infty$,
\item $U_1=U_0$ and $W_1=W_0^{-1}$, and
\item for every nonnegative integer $j$,
\begin{align*}  
W_0 \left( \calV \left( \infty \oplus^j 0 \right)\right) &= 
U_0 \left( \calV \left( \infty \oplus^{j+2} 0 \right) \right)
\end{align*}
\end{itemize}

This identity is a special case of a more general relationship between $W_\beta$, $U_\beta$, and the vertices around the Stern--Brocot triangle $\Delta(\omega)$.  This relationship is used to prove the next lemma which concerns the following collection of elements in $\pslQx$.

\begin{dfns} \label{def: S and T}
Given $\omega \in \Q_0$, take $U_\omega, W_\omega \in \pslQx$ as defined in Definition \ref{def: UW}.  Now, define
\[ T_\omega = U_\omega^{-1} W_\omega^{-1} U_\omega \qquad \text{and} \qquad S_\omega = T_\omega^{-1} U_0^{-1}.\] 
The element $S_\omega$ is called the {\it $\omega$-crossing loop} and $T_\omega$ is the {\it $\omega$-meridian}.  
\end{dfns}
\index{crossing loop!holonomy|idxdef} \index{holonomy!meridian|idxdef}

The names given here to $S_\omega$ and $T_\omega$ will be justified in Section \ref{sec: xing circles} by Theorem \ref{thm: diagram}, where their geometric interpretations as crossing loops and meridians become clear.  Refer also to the labels in Figure \ref{fig: SW link II}.  A calculation shows that
\[
T_\omega = \begin{bmatrix} NQ+d & -N^2 \\ Q^2 & -NQ+d \end{bmatrix} \qquad \text{and} \qquad S_\omega= \begin{bmatrix} -NQ+d & N^2-NQ+d \\ -Q^2 & -Q^2+NQ+d \end{bmatrix}
\]
where $d$, $N$, and $Q$ are defined as in Definition \ref{def: UW}. Thus, it is straightforward to find $S_\omega$ and $T_\omega$ without needing to compute elaborate  matrix products.  
 
\begin{figure}
\setlength{\unitlength}{.1in}
\begin{picture}(28,20)
\put(4.5,0) {\includegraphics[width= 2in]{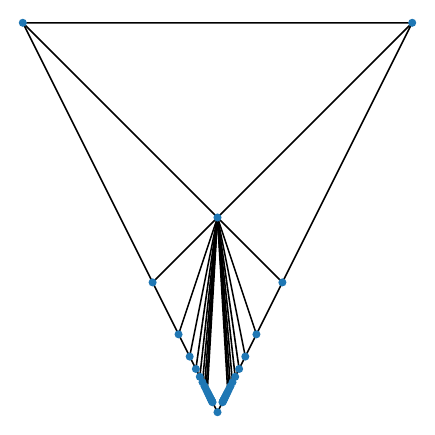}}
\put(24.5,19){$\beta_0=1$}
\put(0,19){$0=\beta_{-1}$}
\put(14.2,11){$\omega$}
\put(4.5,6.5){$1/3=\beta_{-2}$}
\put(5.7,4.0){$2/5=\beta_{-3}$}
\put(18.2,6.7){$\beta_1=2/3$}
\put(17,4.2){$\beta_2=3/5$}
\end{picture}
\caption{The Farey neighbors of $\omega=1/2$ on the boundary of the Stern--Brocot triangle $\Delta(1/2)$.  Their indexing shows the neighbor cycle $\{\beta_j\}_{\Z}$.}  \label{fig: cycle}
\end{figure}

\begin{dfn}\label{def: omega cycle}
Suppose $\omega\in\Q\cap(0,1)$, and let $\beta_0$ and $\beta_{-1}$ be the right and left corners of the Stern--Brocot triangle $\Delta(\omega)$.  For $j\in\Z$, define
\[ \beta_j \ =\  \begin{cases} \beta_0\oplus^j\omega & \text{if } j\geq0,\\ \beta_{-1}\oplus^{-1-j}\omega & \text{if } j<0. \end{cases} \]
The ordered sequence
\[ \ldots,\, \beta_{-3},\, \beta_{-2},\, \beta_{-1},\, \beta_0,\, \beta_1,\, \beta_2,\ldots \]
lists the Farey neighbors of $\omega$ in cyclic order and is called the {\em neighbor cycle} for $\omega$.  By convention, the indices $0$ and $-1$ are reserved for the right and left corners of $\Delta(\omega)$. Elements with negative indices always lie on the left side of $\Delta(\omega)$ and the nonnegative indices on the right side.  See Figure~\ref{fig: cycle}.
\end{dfn}
\index{neighbor cycle|idxdef}

 \begin{lemma} \label{lem: STbeta}
 Suppose $\omega \in \Q \cap (0,1)$, let $\{ \beta_j \}_{\Z}$ be its neighbor cycle, and, for $j \in \Z$, write
 \begin{align*}
 d_j&= d_\calQ(\beta_j) &Q_j&= \calQ(\beta_j) & N_j&= \calN(\beta_j) & V_j&= \calV(\beta_j).
 \end{align*}
 Then
\[ T_\omega\left( \calV(\omega)\right) = \calV(\omega) =  S_\omega\left(-1+\calV(\omega)\right).\]
Also,
\begin{align*}
T_\omega \left( -1+V_j\right) &=  V_{j-2}   &&\text{if } j \leq -1\\
T_\omega \left( V_j\right) &=  V_{j-2}   &&\text{if } j \in \{ 0,1\}\\
T_\omega \left( V_j\right) &= 1+ V_{j-2}   &&\text{if } j \geq 2
\end{align*}
and
\begin{align*}
S_\omega \left( -1+V_j\right) &=  -1+V_{j+2}   &&\text{if } j \leq -3\\
S_\omega \left( -1+V_j\right) &=  V_{j+2}   &&\text{if } j \in \{ -2,-1\}\\
S_\omega \left( V_j\right) &=  V_{j+2}   &&\text{if } j \geq 0.
\end{align*}
 \end{lemma}
  
 \proof
Since $W_\omega$ fixes zero, the first line of equalities follows immediately from the definitions for $S_\omega$ and $T_\omega$.

 Since $S_\omega=T_\omega^{-1}U_0^{-1}$, the last group of equalities are consequences of the second group, so it suffices to prove that
\begin{align*}
W_\omega^{-1}U_\omega U_0\left( V_j\right) &= U_\omega\left( V_{j-2} \right)  &&\text{if } j \leq -1\\
W_\omega^{-1}U_\omega\left( V_j\right) &= U_\omega  \left( V_{j-2} \right)  &&\text{if } j \in \{ 0,1\}\\
W_\omega^{-1}U_\omega\left( V_j\right) &= U_\omega U_0^{-1} \left( V_{j-2} \right)  &&\text{if } j \geq 2.
\end{align*}
These expressions will follow from the Farey recursion relation \ref{eq: defn} and Lemma \ref{lem: det2}.   Specifically, the lemma provides
\[
\begin{vmatrix} N&Q \\ N_j & Q_j \end{vmatrix} \ = \ \begin{cases}
-dQ_{j+1} & \text{if } j\leq -2 \\
-d_jQ_{j+1} & \text{if } j=-1\\
d_jQ_{j-1} & \text{if } j=0 \\
dQ_{j-1} & \text{if } j \geq 1.
\end{cases}
\]
This is useful here because 
\[ U_\omega ( V_j ) = \frac{-1}{QQ_j} \cdot \begin{vmatrix} N&Q \\ N_j & Q_j \end{vmatrix}.\]
It is also useful to notice that, if $A, B \in \Q(x)$, then
\[ W_\omega^{-1} \left( \frac{A}{B} \right) \ =\ \frac{dA}{AQ^2+dB}.\]

Begin by assuming $j\leq -1$.  Then
\[ U_\omega (V_{j-2}) \ = \ \frac{dQ_{j-1}}{QQ_{j-2}}.\]
Now, if $j \leq -2$ then, using the FRF condition,
\begin{align*}
W_\omega^{-1} U_\omega U_0 (V_j) &= W_\omega^{-1} \left( \frac{-QQ_j+dQ_{j+1}}{QQ_j} \right) \\
&= W_\omega^{-1} \left( \frac{ -Q_{j-1}}{QQ_j} \right) \\
&= \frac{dQ_{j-1}}{Q(-dQ_j+QQ_{j-1})} \\
&= \frac{dQ_{j-1}}{QQ_{j-2}}.
\end{align*}
And, if $j=-1$, then
\begin{align*}
W_\omega^{-1} U_\omega U_0 (V_{-1}) &= W_\omega^{-1} \left( \frac{-QQ_{-1}+d_{-1}Q_0}{QQ_{-1}} \right) \\
&= W_\omega^{-1} \left( \frac{-Q_{-2}}{QQ_{-1}}\right) \\
&= \frac{dQ_{-2}}{Q(-dQ_{-1}+QQ_{-2})}\\
&= \frac{dQ_{-2}}{QQ_{-3}}.
\end{align*}
This completes the proof of the $T_\omega$ equality in the case that $j \leq -1$.

Next, consider the case that $j \in \{ 0,1 \}$.  First, notice that
\[ U_\omega(V_{j-2}) \ = \ \begin{cases} \frac{dQ_{-1}}{QQ_{-2}} & \text{if } j=0 \\ \frac{d_{-1}Q_0}{QQ_{-1}} & \text{if } j=1.
\end{cases}\]
Also, if $j=0$ then, using that $d=d_{-1}d_0$,
\begin{align*}
W_\omega^{-1} U_\omega (V_j) &= W_\omega^{-1} \left( \frac{-d_0 Q_{-1}}{QQ_0} \right) \\
&= \frac{dQ_{-1}}{Q(-d_{-1}Q_0+Q_{-1}Q)} \\
&= \frac{dQ_{-1}}{QQ_{-2}}.
\end{align*}
If $j=1$, then
\begin{align*}
W_\omega^{-1} U_\omega (V_j) &= W_\omega^{-1} \left( \frac{-dQ_0}{QQ_1}\right) \\
&= \frac{dQ_0}{Q(QQ_0-Q_1)} \\
&= \frac{d_{-1}Q_0}{QQ_{-1}}.
\end{align*}
This finishes the argument for the $T_\omega$ equality when $j \in \{ 0, 1 \}$.

Finally, assume $j \geq 2$.  Then
\begin{align*}
W_\omega^{-1} U_\omega (V_j) &= W_\omega^{-1} \left( \frac{-dQ_{j-1}}{QQ_j} \right) \\
&= \frac{dQ_{j-1}}{Q(QQ_{j-1}-Q_j)} \\
&= \frac{Q_{j-1}}{QQ_{j-2}}.
\end{align*}
If $j\geq 3$, then
\begin{align*}
U_\omega U_0^{-1} (V_{j-2}) &= \frac{QQ_{j-2}-dQ_{j-3}}{QQ_{j-2}} \\
&= \frac{Q_{j-1}}{QQ_{j-2}}
\end{align*}
and, if $j=2$, then
\begin{align*}
U_\omega U_0^{-1} (V_{j-2}) &= \frac{QQ_0-d_0Q_{-1}}{QQ_0} \\
&= \frac{Q_1}{QQ_0}.
\end{align*}
This completes the proof of the lemma since, as mentioned above, the last group of equalities follow from the expression $S_\omega=T_\omega^{-1}U_0^{-1}$.
 \endproof
 
The next lemma is a consequence of Lemma \ref{lem: STbeta} and provides a relation which will be used later to prove Theorem \ref{thm: generic holonomy}, which says that every $S_\omega$ and every $T_\omega$ belongs to the generic holonomy group $\Gamma$ from Definition \ref{def: generic holonomy}.

\begin{lemma} \label{lem: relation} 
If $e \in \calE$, then 
\[ S_{e^L} = T_{\gamma} T_{e^R}\]
where $\gamma=e^L\oplus e^R$.
\end{lemma}

\proof  
This follows from the fact that a M\"obius function on $\Q(x) \cup \{ \infty\}$ is uniquely determined by its values at three distinct points.  Because $e^L$ and $e^R$ are the left and right corners of $\Delta(\gamma)$, Lemma \ref{lem: STbeta} shows that the relation holds at $-1+\calV\left(e^L\right)$ and $\calV\left(e^R\right)$.  

If $e^R$ is not a corner of $\Delta\left(e^L\right)$, then take $\gamma_0$ so that $\gamma_0 \oplus e^L = e^R$.  Then $e^L$ and $\gamma_0$ are the corners of $\Delta\left(e^R\right)$ so, by Lemma \ref{lem: STbeta}, the relation holds at $\calV\left(e^R\right)$.  

Otherwise, $e^R$ is the right corner of $\Delta\left(e^L\right)$.  So, if $\gamma_0$ is the left corner of $\Delta\left(e^L\right)$, then the relation holds at $-1+\calV(\gamma_0)$.
\endproof

\subsection{Face pairings for $\Omega(\alpha)$} \label{sec: face pairings}

In Definition \ref{def: generic domain}, $\Omega(\alpha)$ is defined as a collection of triples $\delta(e),\delta'(e) \in \Q(x)^3$.  As described in Definition \ref{def: generic tetrahedra}, the pairs of elements of these triples are labeled by rational numbers according to the funnel $F_\alpha$.  This has been arranged so that, if $\frz \in \calH(M_\alpha^\circ)$, then the corresponding developing map preserves labels when it maps tetrahedra from $\wt{M}_\alpha^\circ$ to $\Omega(\alpha)_\frz$.  

Given $\alpha \in \Q \cap (0,1/2)$ with $\alpha \neq 1/3$ and $\frz \in \calH(M_\alpha^\circ)$, let $\omega_1, \ldots, \omega_k$ be the sequence of distinct hubs for $F_\alpha$.  Definition \ref{def: S and T} provides elements $S_{\omega_j}$ and $T_{\omega_j}$ in $\pslQx$.  Recall that
\[ U_0 \ =\ \begin{bmatrix} 1&-1\\0&1 \end{bmatrix}.\]
Following the discussion in Section \ref{sec: completeness} and using Lemma \ref{lem: delta pairs}, $U_0$ takes the pair in $\delta'_{j-1}$ labeled $\omega_j$ to the pair in $\delta_j$ with the same label.  Referring back to the definition of the Sakuma--Weeks triangulation in Definition \ref{def: SW}, this means that $U_0$ is a peripheral face pairing isometry in the holonomy group for $\frz$.  Other face pairing elements of this holonomy group are, so far, unexplored.  In particular, the face pairings of $\Omega(\alpha)_\frz$ which involve the finite faces of the tetrahedra making up $\Omega(\alpha)_\frz$ have not yet been identified.  Also, those that involve the upper face $(\delta_0)_\frz = [0,1,\infty]$ of $\left( \delta_1 \right)_\frz$ and, when $\frz$ is the geometric root for $\calQ(\alpha)$, the lower face $(\delta'_{n+1})_\frz = \left[ \calV\left(e_{n+1}^-\right)_\frz, \calV\left(e_{n+1}^+\right)_\frz, \infty\right]$ of $\left( \delta_n' \right)_\frz$ can, and should, be understood.  \index{fundamental domain!specialized}
 
 \begin{dfn} \label{def: Ai} \index{funnel!hub}
 Given $\alpha$ and $\{ \omega_1, \ldots, \omega_k\}$ as above, define $A_i$ by
\[ A_i = \begin{cases} S_{\omega_i} & \text{if } \omega_i \text{ is on the left} \\ T_{\omega_i} & \text{otherwise.} \end{cases} \] 
\end{dfn}

Notice that $\omega_i$ is on the left if and only if $i$ is odd.  Also, because $\omega_1$ is always $0$, the isometry $A_1$ is always equal to 
\[ S_0 \ =\ \begin{bmatrix} x-1 & x \\ -1 & x \end{bmatrix}.\]

\begin{thm} \label{thm: generic face pairing}
Suppose $\alpha \in \Q \cap (0,1/2)$ and $\alpha \neq 1/3$.  Take $k$, $\omega_i$,  and $A_i$ as above.  If $i \in \{ 1, \ldots, k \}$ and $e_j$ is a spoke for $\omega_i$, then $A_i$ takes the triple $\delta_{j-1} \subset \Q(x)$ to $\delta_j' \subset \Q(x)$, preserving the labeling of pairs by rationals.  Also, $A_1$ takes the triple $\delta_0=\{ 0, 1, \infty\}$ to the finite face of $\delta'_1$ so that the induced map on labels of pairs is 
\[ 1/2 \mapsto 1/2 \quad 0 \mapsto 1 \quad 1 \mapsto 0.\]
\end{thm}

\proof
The first part follows directly from Lemma \ref{lem: STbeta} and the definitions of $\delta_j$ and $\delta_j'$.  The second part is a simple calculation.
\endproof

The next corollary follows immediately.

\begin{cor} \label{cor: holonomy I}
If $\frz \in \calH(M_\alpha^\circ)$, then
\[ \left\{ U_0, (A_1)_\frz, \ldots, (A_k)_\frz \right\} \]
is a generating set for the holonomy group of the developing map for $M_\alpha^\circ$.  Moreover, $(A_1)_\frz$ is the face pairing which completes the upper hairpin fold.
\end{cor}

When $\frz\in\calH(M_\alpha)$, the lower hairpin fold introduces two additional face identifications.  The next corollary shows that these require no additional generators.

\begin{cor} \label{cor: holonomy II}
If $\frz \in \calH(M_\alpha)$, then
\[ \left\{ U_0, (A_1)_\frz, \ldots, (A_k)_\frz \right\} \]
is a generating set for the holonomy group of the hyperbolic structure on $M_\alpha$ determined by $\frz$.
\index{hairpin!fold}
\end{cor}

\proof There are two triangular face identifications needed to accomplish the bottom hairpin fold across the two hairpin edges which are both labeled $\hat{e}_{n+1}$.  The first of these face pairings identifies the faces $F \subset (\delta_n)_\frz$ and $F'\subset (\delta'_n)_\frz$, where
\begin{align*}
F &= \left[ -1+\calV\left(e_n^L\right)_\frz, \ \calV\left(e_n^R\right)_\frz, \ \infty \right] \\
F' &= \left[ \calV\left( e_n^- \oplus e_n^+ \right)_\frz, \ \calV\left( \hat{e}_{n+1} \right)_\frz, \ \infty \right].
\end{align*}
The second face pairing identifies the finite face of $(\delta_n)_\frz$ to the face 
\[ \left[ \calV\left( e_n^- \oplus e_n^+ \right)_\frz, \ \calV\left( e_{n+1}^- \right)_\frz, \ \infty \right].\]

By Theorem \ref{thm: generic face pairing}, $A_k$ takes $\delta_n$ to $\delta'_{n+1}$ and so $(A_k)_\frz$ takes the finite face of $(\delta_n)_\frz$ to
\[ \left[ \calV\left(e^R_{n+1}\right)_\frz, \ \calV\left( e^L_{n+1}\right)_\frz, \  \calV(\alpha)_\frz \right] \ = \ \left[ \calV\left(e^R_{n+1}\right)_\frz, \ \calV\left( e^L_{n+1}\right)_\frz, \  \infty \right].
\]
Since $\left\{ e^L_{n+1}, \, e^R_{n+1} \right\} = \left\{ e^-_n \oplus e^+_n, \, e^-_{n+1} \right\}$, the element $(A_k)_\frz$ of $\Gamma_\frz$ provides the second face identification.

The holonomy element needed for the identification of $F$ to $F'$ depends on the value of $\nu(e_{n+1})$.  First, if $\nu(e_{n+1})=-1$, then 
\[ e_n^L = e_{n+1}^L \quad \text{and} \quad e_n^- \oplus e_n^+ = e_{n+1}^R.\]
So, using Lemma \ref{lem: difference},
\begin{align*}
-1 + \calV\left( e_n^L \right)_\frz &= -1 + \calV\left( e_{n+1}^L \right)_\frz \\
&= \calV\left( e_n^R \right)_\frz \\
&= \calV\left(  e_n^- \oplus e_n^+ \right)_\frz. \\
\end{align*}
Therefore, $F'=F$ and the holonomy element is the identity.  Otherwise, $\nu(e_{n+1})$ is positive and 
\[ e_n^R = e_{n+1}^R \quad \text{and} \quad e_n^- \oplus e_n^+ = e_{n+1}^L.\]
So,
\begin{align*}
\calV\left( e_n^R \right)_\frz &=  \calV\left( e_{n+1}^R \right)_\frz \\
&= -1+\calV\left( e_n^L \right)_\frz \\
&= -1+\calV\left(  e_n^- \oplus e_n^+ \right)_\frz. \\
\end{align*}
Therefore, $F'=F+1$ and $U_0$ accomplishes this part of the hairpin fold.  This, together with Corollary \ref{cor: holonomy I}, establishes this corollary.
\endproof

The generic holonomy group $\Gamma$ (Definition~\ref{def: generic holonomy}) is the subgroup of $\pslQx$ generated by $U_0$ and $W_0$.  The next theorem and corollary justify its name by showing that, for every $\frz\in\calH(M_\alpha^\circ)$, the specialization $\Gamma_\frz$ is the holonomy group of $\dev_\frz$.

\begin{thm} \label{thm: generic holonomy}
For each $\beta \in \Q \cap [0,1/2]$, the face pairing transformations $S_\beta$ and $T_\beta$ belong to $\Gamma$.
\end{thm}

\proof
Recall that
\[ U_1^{-1} \ =\ U_0^{-1} \ = \ \begin{bmatrix} 1&1\\0&1 \end{bmatrix} \qquad \text{and} \qquad
 W_1^{-1} \ = \ W_0 \ =\ \begin{bmatrix} 1&0 \\ -\frac1x&1 \end{bmatrix}.\]
By definition, $T_\beta = U_\beta^{-1} W_\beta^{-1} U_\beta$, so $T_0$ and $T_1$ belong to $\Gamma$.  

It is enough to show that $T_\beta \in \Gamma$ for every $\beta$ because $S_\beta = T_\beta^{-1} U_0^{-1}$.  This is done by induction on the denominator of $\beta$.  The base case, where the denominator is one, is already finished.

Suppose that $\beta$ has denominator larger than one and that the theorem holds for $\gamma$'s with denominator smaller than that of $\beta$.  Then there is an edge $e \in \calE$ with $e^L \oplus e^R=\beta$.  By assumption, $T_{e^L}$ and $T_{e^R}$ belong to $\Gamma$.  The relation of Lemma \ref{lem: relation} now shows that $T_\beta \in \Gamma$.
\endproof

\begin{cor} \label{cor: holonomy specialization} \index{developing map}
If $\frz\in \calH(M_\alpha^\circ)$ then $\Gamma_\frz$ is the holonomy group for the developing map $\dev_\frz$ associated to $\Omega(\alpha)_\frz$.
\end{cor}
\index{holonomy!specialization}

\proof
Theorem \ref{thm: generic holonomy}, together with Corollary \ref{cor: holonomy I}, shows that $\Gamma_\frz$ contains the holonomy group of $\dev_\frz$.

On the other hand, $\Gamma_\frz$ is generated by $(U_0)_\frz$ and $(W_0)_\frz$.  Because
\[ T_0^{-1} = U_0^{-1} \, W_0 \, U_0 \qquad \text{and} \qquad S_0=T_0^{-1}U_0^{-1}\]
$\Gamma_\frz$ is also generated by the face pairing isometries $(U_0)_\frz$ and $(S_0)_\frz$.   Therefore, the holonomy group contains $\Gamma_\frz$.
\endproof

Let $\alpha$ be hyperbolic and consider the geometric root $\frz(\alpha)$ for $\calQ(\alpha)$.  Theorem \ref{thm: complete} implies that the developing map is an isometry $\wt{M}_\alpha \to \IH^3$.  Consequently, the holonomy
map \index{geometric root!2-bridge link}
\[
\pie M_\alpha\longrightarrow\Gamma_{\frz(\alpha)}
\]
is an isomorphism and
\[
M_\alpha\cong\IH^3/\Gamma_{\frz(\alpha)}.
\]

\subsection{Punctured spheres and crossing circles} \label{sec: xing circles}

Assume $\alpha \in \Q \cap (0,1/2)$ with $\alpha \neq 1/3$ and let $\{ \omega_1, \ldots, \omega_k\}$ be the sequence of hubs of the funnel $F_\alpha$.  The goal of this subsection is to relate the topology of $M_\alpha^\circ$ and its truncated diagram $D_\alpha$ to the generic holonomy group $\Gamma$ defined in Definition \ref{def: generic holonomy}.  This is done by looking at the $4$-punctured spheres $\Sigma_j\subset M_\alpha^\circ$ introduced in Definition~\ref{def: quotient complex}.  Recall that each $\Sigma_j$ is a union of faces of the tetrahedral pairs in $M_\alpha^\circ$.

Suppose $\frz \in \calH(M_\alpha^\circ)$.  The first step is to understand which elements of $\Gamma$ specialize to edge pairings for a lift of $\Sigma_j$ to the domain $\Omega(\alpha)_\frz$.  With Figure~\ref{fig: SW link}, Sakuma and Weeks provide a picture of how the triangulation of $\Sigma_j$ sits in the diagram $D_\alpha$.   This will imply that each $\left( S_{\omega_i} \right)_\frz$ is represented by a loop which is freely homotopic to a crossing circle for the twist region associated to the hub $\omega_i$ of $F_\alpha$.  Similarly, it will become easy to see arcs of $D_\alpha$ which correspond to each $T_{e_j^R}$.  

\begin{center}
\begin{figure}[h] 
\setlength{\unitlength}{.1in}
\begin{picture}(26,16)
\put(4.4,0) {\includegraphics[width= 2in]{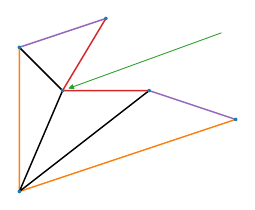}}
\put(13.8,15){$\scriptstyle \varepsilon+\calV\left( \hat{e} \right)$}
\put(5.2,.3){$\scriptstyle \infty$}
\put(23.7,7){$\scriptstyle \calV\left(e^L\right)$}
\put(22.3,14.1){$\scriptstyle \calV\left(e^R\right)$}
\put(16.5,10){$\scriptstyle \calV\left( e^L \oplus e^R \right)$}
\put(0,13){$\scriptstyle -1+\calV\left(e^L\right)$}
\end{picture}
\caption{The labeled simplicial complex $\wt{\Sigma}(e)$.  The vertices are labeled by elements of $\Q(x) \cup \{ \infty \}$ as shown.  The red and orange edges are all labeled by the rational number $e^L$ while the purple edges are labeled $e^R$.  The black edges that share a vertex with a purple edge are labeled $e^L\oplus e^R$ and the last (black) edge is labeled $e^R$.  The purple edges are paired by $S_{e^L}$, the red edges are paired by $T_{e^R}$, and the orange edges are paired by $U_0$.  The quotient of $\wt{\Sigma}(e)$ by these edge pairings is a triangulated 2-sphere.}  \label{fig: Sigma_e}
\end{figure}
\end{center}

The complex $\wt{\Sigma}(e)$ introduced here is designed so that, when $e=e_j$ is an interior edge of $F_\alpha$, the specialization $\wt{\Sigma}(e_j)_\frz$ is a lift of the triangulated $4$-punctured sphere $\Sigma_j$.  For an edge $e \in \calE$, let $\wt{\Sigma}(e)$ be the labeled abstract simplicial complex of four triangles shown in Figure \ref{fig: Sigma_e}.  The number $\varepsilon$ should be chosen so that $\wt{\Sigma}(e)$ contains a triangle labeled by the three elements of $\delta(e)$.   (See Definition \ref{def: generic tetrahedra} to remember how $\varepsilon$ relates to $\delta(e)$.)  Evaluating the functions at the vertices of $\wt{\Sigma}(e)$ at a value $\frz \in \C$ determines a subset $\wt{\Sigma}(e)_\frz$ of $\delta(e)_\frz \cup \delta'(e)_\frz$.  The edges of $\wt{\Sigma}(e)$ are labeled by rational numbers by pulling back the labels from $\delta(e)_\frz \cup \delta'(e)_\frz$.  

\begin{figure}[h] 
   \centering
   \includegraphics[width=5in]{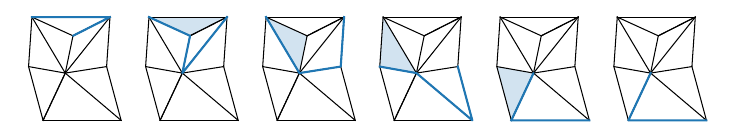} 
   \caption{Here, $\alpha = 5/17$ and $\frz=\frz_1\sim 0.473764 + 0.240160\, i$.  The figure shows a copy of $\Omega(\alpha)_\frz$ for each of the edges $e_0, \ldots , e_{n+1}$ of $F_\alpha$.  In the $j^\text{th}$ copy the subcomplex $\wt{\Sigma}(e_j)_\frz$ is highlighted in blue.  Note that, when $j$ is $0$ or $n+1$,  $\wt{\Sigma}(e_j)_\frz$ collapses from four ideal triangles to two.}  \label{fig: wt Sigma}
\end{figure}

Since $U_0$ fixes $\infty$ and takes $\calV(e^L)$ to $-1+\calV(e^L)$, it can be viewed as a homeomorphism between the orange edges in Figure \ref{fig: Sigma_e} and as an isometry of $\IH^3$ which identifies the corresponding geodesics in $\wt{\Sigma}(e)_\frz$.  Also, from Lemma \ref{lem: STbeta}, 
\begin{align*}
S_{e^L} \left( -1+\calV\left(e^L\right)\right) &= \calV\left(e^L\right) & 
S_{e^L} \left( \epsilon +\calV\left(\hat{e}\right)\right) &= \calV\left(e^L \oplus e^R \right)\\
T_{e^R} \left( \epsilon +\calV\left(\hat{e}\right)\right) &= \calV\left(e^L \oplus e^R \right) &
T_{e^R}\left( \calV\left(e^R\right)\right) &= \calV\left(e^R\right).
\end{align*}
So, $S_{e^L}$ and $T_{e^R}$ can be regarded as homeomorphisms identifying the purple pair of edges and identifying the red pair of edges.  Likewise, the specializations of $S_{e^L}$ and $T_{e^R}$ to $\frz$ are isometries which do the same for the corresponding edges of $\wt{\Sigma}(e)_\frz$.  In the quotient, $\wt{\Sigma}(e)$ folds into a triangulated 2-sphere with four triangles and four vertices.  These four vertices are ideal vertices, so deleting them gives a triangulated $4$-punctured sphere.  When $e=e_j$, the edge labels agree with those on $\Sigma_j$ in Figure~\ref{fig: SW link}.  Hence, the quotient is identified with the triangulated sphere $\Sigma_j$.  In particular, if $\frz \in \calH(M_\alpha^\circ)$ and $e_j$ is the $j^\text{th}$ interior edge of the funnel $F_\alpha$, then the specialization $\wt{\Sigma}(e_j)_\frz$ is a lift of $\Sigma_j$ from $M_\alpha^\circ$ to $\Omega(\alpha)_\frz$.  This is shown in Figure \ref{fig: wt Sigma}.
\index{Sakuma--Weeks complex!4-punctured spheres}

Now, consider the truncated diagram $D_\alpha$ for $M_\alpha^\circ$ as shown in Figure \ref{fig: SW link II}.  By identifying the elements of $\pie M_\alpha^\circ$ with those of $\Gamma$ using the generic holonomy isomorphism, oriented arcs of $D_\alpha$ represent certain elements of $\Gamma$.  Likewise, each twist region determines an element of $\pie M_\alpha^\circ = \Gamma$ as indicated in Figure \ref{fig: SW link II}.  The figure is labeled to show this correspondence and the next theorem describes how this holds in general.  Mentioned in the theorem are the distinguished elements $S_\omega$ and $T_\omega$ in $\Gamma$ from Definition \ref{def: S and T}. \index{funnel!hub}

\begin{thm} \label{thm: diagram}
Take $\alpha \in \Q \cap (0,1/2)$ with $\alpha \neq 1/3$ and identify the elements of $\pie M_\alpha^\circ$ and $\Gamma$ via the generic holonomy isomorphism.  Let $D_\alpha$ be the truncated diagram for $M_\alpha^\circ$ determined by the funnel $F_\alpha$.  The long arc across the top of the diagram (oriented to the left) is $U_0$.  The element $S_{\omega_j}$ is represented by a loop which is freely homotopic to a crossing circle about the $j^\text{th}$ twist region of $D_\alpha$.  If $v_1, \ldots, v_m$ are the vertices on the right side of $F_\alpha$, listed from top to bottom, then the sequence of arcs (oriented to the left) immediately below $U_0$ are $T_{v_1}, \ldots, T_{v_m}$.
\end{thm}
\index{crossing loop}
\index{crossing loop!holonomy}
\index{crossing loop!twist region}

\proof
The preceding construction identifies the three pairs of edges of $\wt{\Sigma}(e)$ by $U_0$, $S_{e^L}$, and $T_{e^R}$.  When $e=e_j$ is an interior edge of $F_\alpha$, $\wt{\Sigma}(e_j)$ gives a lift of the triangulated $4$-punctured sphere $\Sigma_j$.  Thus these edge pairings can be compared directly with the position of $\Sigma_j$ in Figure~\ref{fig: SW link}.

The edges paired by $U_0$ project, as the spheres $\Sigma_j$ are followed through the diagram, to the long arc across the top of
$D_\alpha$.  With the orientation in the statement, this arc therefore represents $U_0$.  Similarly, the edges paired by $T_{e^R}$ project to the arcs immediately below the top arc.  As the edges $e_0,\ldots,e_{n+1}$ of the funnel are followed from top to bottom, their right endpoints run through the vertices $v_1,\ldots,v_m$ on the right side of $F_\alpha$, with repetitions when consecutive edges have the same right endpoint.  Suppressing these repetitions, the corresponding arcs therefore represent $T_{v_1},\ldots,T_{v_m}$.

It remains to identify the crossing loops.  Let $\omega$ be a hub of $F_\alpha$.  The spokes at $\omega$ correspond to the layers of the associated twist region of $D_\alpha$.  If $\omega$ is a left hub, then $e_j^L=\omega$ for each of its spokes.  The corresponding purple edge pairing of $\wt{\Sigma}(e_j)$ is therefore $S_\omega$, and comparison with Figure~\ref{fig: SW link} shows that its associated loop is freely homotopic to a crossing circle about this twist region.

If $\omega$ is a right hub, then $e_j^R=\omega$ for each of its spokes.  The two punctures enclosed by a crossing circle for the corresponding twist region are represented by the arcs $U_0$ and $T_\omega$.  Hence, after choosing an orientation and a basepoint, the crossing circle represents $T_\omega^{-1}U_0^{-1}$, up to conjugacy.  By Definition~\ref{def: S and T}, this element is $S_\omega$.  In either case, $S_\omega$ is represented by a loop freely homotopic to the corresponding crossing circle.
\endproof

\begin{cor} \label{cor: geodesic crossing circles}
Suppose $\alpha$ is hyperbolic and $\frz\in\calH(M_\alpha)$. If $\omega$ is a hub of $F_\alpha$, then $(S_\omega)_\frz$ is loxodromic.  If $g\subset\IH^3$ is its axis, then the image of $g$ under the quotient
\[ \IH^3\longrightarrow \IH^3/\Gamma_\frz\cong M_\alpha \]
is the geodesic representative of the free homotopy class of the crossing circle for the twist region of $D_\alpha$ associated to $\omega$.
\end{cor}

\proof
By Theorem~\ref{thm: diagram}, $(S_\omega)_\frz$ represents the free homotopy class of the corresponding crossing circle.  This loop is nonperipheral, so completeness of the hyperbolic structure on $M_\alpha$ implies that $(S_\omega)_\frz$ is loxodromic and its axis projects to the geodesic representative of this class.
\endproof

\begin{remark} \label{rem: slopes}
The isotopy classes of essential loops on the punctured spheres $\Sigma_j\subset M_\alpha^\circ$ are naturally associated to $\Q\cup\{\infty\}$ as slopes through the Sakuma--Weeks construction (see Definitions~\ref{def: tetrahedral complex} and \ref{def: quotient complex}).  These agree with the rational labels on the edges triangulating the spheres; the boundary of a regular neighborhood of an edge of slope $\omega$ in $\Sigma_j$ is an essential loop of slope $\omega$.

Under the generic holonomy isomorphism $\pie(M_\alpha^\circ)\cong\Gamma$, the element $S_\omega$ corresponds to a loop freely homotopic to an essential loop of slope $\omega$ on $\Sigma_j$.  In particular, when $\omega$ is a hub of $F_\alpha$, Theorem~\ref{thm: diagram} identifies this class with the corresponding crossing circle of $D_\alpha$.

For a general $\frz\in\calH(M_\alpha^\circ)$, even when $M_\alpha^\circ$ is a manifold and $(S_\omega)_\frz$ is loxodromic, its axis in $\IH^3$ need not lie in the image of the developing map.  Thus, in general, it may not be true that the corresponding free homotopy class has a geodesic representative contained in $M_\alpha^\circ$.  When $\frz\in\calH(M_\alpha)$, completeness removes this issue, as in Corollary~\ref{cor: geodesic crossing circles}.

For a complete characterization of peripheral essential simple loops on $2$-bridge spheres in hyperbolic $2$-bridge link complements, see Theorem 2.4 of \cite{LS-loops}. 
\end{remark}

\section{Computer Code} \label{subsec: code III}

Together with the code given in Sections~\ref{subsec: code I} and~\ref{subsec: code II}, the commands listed here can be used to have \lstinline$SageMath$ draw accurate pictures of $\Omega(\alpha)_{\frz(\alpha)}$, where $\alpha$ is hyperbolic and $\frz(\alpha)$ is the geometric root for $\calQ(\alpha)$.  As before, a rational number $\alpha$ is represented by a vector \lstinline$v$.

The function \lstinline$coor(z)$ gives the Cartesian coordinates of a complex number, and \lstinline$Vc(v, z)$ evaluates $\calV$ at a vertex of $\calG$ represented by \lstinline$v$ and specializes the resulting rational function at \lstinline$z$.  The command \lstinline$domain_points(v)$ computes the coordinates of the ideal vertices used to draw $\Omega(\alpha)_{\frz(\alpha)}$, while \lstinline$domain(v)$ constructs the corresponding graphics object.

\begin{lstlisting}[language=Python]
def coor(z):   # z is a complex number.
    return [z.real(), z.imag()]

def Vc(v, z):   # v is a vector and z is in C.
    return coor((N(v)/Q(v)).substitute(x = z))

def domain_points(v):   # v is a vector.
    groot = geometric_root(v)
    fedges = funnel_edges(v)
    return [[Vc(e[0], groot), Vc(e[1], groot)] 
    	for e in fedges]

def domain(v):   # v is a vector.
    dp = domain_points(v)
    L1 = sum([line([vector(e[0]) - vector([1,0]), e[1], e[0]]) 
              for e in dp])
    L2 = line([vector(e[0]) - vector([1,0]) for e in dp])
    L3 = line([e[1] for e in dp])
    L4 = line([e[0] for e in dp])
    return L1 + L2 + L3 + L4
\end{lstlisting}

After defining the previous \lstinline$SageMath$ functions, the following command draws $\Omega(\alpha)_{\frz(\alpha)}$.

\begin{lstlisting}[language=Python]
show(domain(v), aspect_ratio=1, figsize=3, axes=False)
\end{lstlisting}
\chapter{Logarithmic spirals}\label{chap: logs}
\index{logarithmic spiral}

This chapter uses the tools developed in the previous chapters to uncover some geometric properties of the Sakuma--Weeks complexes and their developed domains.  The main goals are 
\index{Sakuma--Weeks complex}
\begin{itemize}
\item to show that the $\omega$-crossing circles $S_\omega \in \Gamma$ are conjugate to the square of the recursion matrix for the Farey recursive functions $\calQ$ and $\calN$ around the Stern--Brocot triangle $\Delta(\omega)$, and 
\index{crossing loop}
\index{crossing loop!holonomy}
\index{Farey recursion}
\index{Farey recursive function}
\index{Stern--Brocot!triangle}
\item to show that the ideal vertices of the domains $\Omega(\alpha)_\frz$ lie on logarithmic spirals for the square roots of $\left(S_\omega\right)_\frz$, where $\omega$ is taken from the hubs of the funnel $F_\alpha$.
\index{fundamental domain}
\index{funnel!hub}
\end{itemize}

It is natural that logarithmic spirals appear; if $B$ is a loxodromic element of $\pslC$ and $p \in \bound \IH^3$ is not fixed by $B$, then the $\langle B \rangle$-orbit of $p$ lies on a unique logarithmic spiral in $\bound \IH^3$.  Certainly, this applies to any ideal point of $\Omega(\alpha)_\frz$ and any loxodromic element of $\Gamma_\frz$.  The surprising aspects are that $S_\omega$ is so closely related to its recursion matrix and that all of the ideal vertices of the domain $\Omega(\alpha)_\frz$ lie on these particular spirals.
 
 To reach these goals, the chapter begins with a discussion of ratios of recursive sequences and their relationships to M\"obius functions.  When interpreted properly, these relationships extend around the corners of the Stern–Brocot triangles $\Delta(\omega)$.  The second part of the chapter discusses this.  As usual, specialized versions hold when $x$ is specialized to values $\frz$ in $\C$.

\section{M\"obius functions and ratios of recursive sequences} \label{sec: ratios}

Let $d$ and $Q$ be elements of the field $\Q(x)$ with $d \neq 0$ and $Q^2 \neq 4d$.  Take $\K$ to be the field extension 
\[ \K \  = \ ^{\displaystyle \Q(x)[y]}\big/_{\displaystyle \left(y^2-(Q^2-4d)\right).}\]
Then, over $\K$, the matrix
\[ X \ = \begin{pmatrix} 0 & 1 \\ -d & Q \end{pmatrix} \]
has eigenvalues $\lambda_{\pm} = \frac12 \, \left( Q \pm y \right)$.  For $n \in \Z$, define
\[ \calA_n \ = \ \lambda_-^n - \lambda_+^n \]
and observe that $\calA_0= 0$,  $\calA_1 = -y$, and $\lambda_+ \lambda_-=d$.  

\begin{lemma} \label{lem: An}
The sequence $\{ \calA_n \}_{n \in \Z}$ is recursive with recursion matrix $X$.
\end{lemma}

\proof
Because
\begin{align*}
Q \calA_{n-1} &= \left( \lambda_+ + \lambda_- \right) \left( \lambda_-^{n-1} - \lambda_+^{n-1} \right) \\
&= \calA_n + \lambda_+ \lambda_- \left( \lambda_-^{n-2}-\lambda_+^{n-2} \right) \\
&= \calA_n + d \calA_{n-2}, 
\end{align*}
it must be true that $\calA_n = -d \calA_{n-2} + Q \calA_{n-1}$.
\endproof

\begin{lemma} \label{lem: X recursion}
If $\{ p_n \}_{n \in \Z}$ is a recursive sequence in $\Q(x)$ with recursion matrix $X$, then
\[ p_n \calA_1  = -d p_0 \calA_{n-1}  + p_1 \calA_n.\]
\end{lemma}
\proof
Since $Q^2\neq 4d$, the eigenvalues $\lambda_\pm$ are distinct and
\[ X^n \ = \ \frac{1}{\calA_1} \, \begin{pmatrix} 1&1 \\ \lambda_+ & \lambda_- \end{pmatrix} \, \begin{pmatrix} \lambda_+^n &0 \\ 0& \lambda_-^n \end{pmatrix} \, \begin{pmatrix} \lambda_- & -1 \\ -\lambda_+ & 1 \end{pmatrix}.\]
The top row of the matrix product on the right-hand side is
\begin{align*}
\left( \lambda_+^n \lambda_- - \lambda_-^n \lambda_+ \quad \lambda_-^n-\lambda_+^n \right) &= \left( -\lambda_+ \lambda_- \calA_{n-1} \quad \calA_n \right)\\
&= \left( -d \calA_{n-1} \quad \calA_n \right).
\end{align*}
Since
\[ X^n \, \begin{pmatrix} p_0 \\ p_1 \end{pmatrix} \ = \ \begin{pmatrix} p_n \\ p_{n+1} \end{pmatrix}\]
the result follows.
\endproof

\begin{lemma} \label{lem: quotients}
Suppose that $\{ p_n\}_{n \in \Z}$ and $\{ q_n \}_{n \in \Z}$ are recursive sequences in $\Q(x)$, both having recursion matrix $X$, and assume that $p_1q_0-p_0q_1\neq0$.  Let
\[Y=\begin{pmatrix} -d p_0 & p_1 \\ -d q_0 & q_1 \end{pmatrix} \] 
and regard $X$ and $Y$ as M\"obius functions on $\Q(x) \cup \infty$.
Then
\[ \frac{p_n}{q_n} \ = \ Y X^{n}\left( \infty \right).\]
The sequence $\left\{ \frac{p_n}{q_n} \right\}_{n \in \Z}$ lies on the logarithmic spiral for $YXY^{-1}$ through $\frac{p_0}{q_0}$. 
\end{lemma}

\proof
The determinant hypothesis implies that $Y$ is invertible.  By Lemma~\ref{lem: An}, $\{\calA_n\}$ has recursion $X$.  Since $X(\infty)=0$,
\[ X^n(\infty)\ =\ X^{n-1}(0)\ =\ [\calA_{n-1}:\calA_n] \]
in the projective line.  Therefore,
\begin{align*}
YX^n(\infty) &= \left[ -dp_0\calA_{n-1}+p_1\calA_n: -dq_0\calA_{n-1}+q_1\calA_n \right] \\
&= [p_n\calA_1:q_n\calA_1] \\
&= [p_n:q_n]
\end{align*}
where the second equality follows from Lemma~\ref{lem: X recursion} and the last uses $\calA_1=-y\neq0$.  Hence,
\[ \frac{p_n}{q_n}\ = \ YX^n(\infty) \ =\  \left(YXY^{-1}\right)^n \left(\frac{p_0}{q_0}\right). \]
\endproof

A related phenomenon appears in Section 10 of \cite{GF}.  In a numerical family of punctured torus bundles, Gu\'eritaud shows that certain ideal vertices are discrete points on a parametrized curve
\[ \big\{ \cot(a+tb) \, \big| \, t \in \R \big\}\]
where $a, b \in \C$.  The addition formula for cotangent shows that these points form an orbit
of the M\"obius transformation 
\[z\ \mapsto \ \frac{z\cos b-\sin b}{z\sin b+\cos b}.\]
After conjugating its fixed points to $0$ and $\infty$, this transformation is multiplication by $e^{2ib}$.  Thus Gu\'eritaud's cotangent parametrized vertices give a concrete instance of the logarithmic spiral phenomenon
described in Lemma~\ref{lem: quotients}.  Compare Figure 11 of \cite{GF} to Figures \ref{fig: spiral 35/86}, \ref{fig: spiral_7_38}, \ref{fig: bottom_spiral I}, \ref{fig: spiral 21/46}, and \ref{fig: spiral 9_41} here.   

The following sections identify, directly from Farey recursion, the M\"obius functions governing the corresponding fans in the Sakuma--Weeks complexes.

\section{Generic spirals} \label{sec: generic spirals} 

Suppose $\omega \in \Q \cap (0,1/2]$ and let $\{ \beta_j \}_{\Z}$ be its neighbor cycle (Definition~\ref{def: omega cycle}). Write $Q=\calQ(\omega)$ and $d=d_\calQ(\omega)$ and
\index{neighbor cycle}
\begin{align*}
 d_j&= d_\calQ(\beta_j) &Q_j&= \calQ(\beta_j) & N_j&= \calN(\beta_j) & V_j&= \calV(\beta_j).
 \end{align*}  

\begin{lemma} \label{lem: 4d}
\[ 4d \neq Q^2.\]
\end{lemma}
\proof
Because $\omega \in \Q_0$, Remark \ref{rem: constant term} states that the constant term of $Q$ is one.  On the other hand, the constant term of $d$ is zero.
\endproof

\begin{dfns} \label{def: XYR}
Define
\begin{align*} 
X_\omega&= \begin{bmatrix} 0&1 \\ -d &Q \end{bmatrix},  & 
Y_\omega&=\begin{bmatrix} -dN_0 & N_1 \\ -dQ_0 & Q_1 \end{bmatrix} 
\end{align*}
and $R_\omega=Y_\omega X_\omega Y_\omega^{-1}$.  Together with the fact that $d_\calQ$ is never zero, Lemma \ref{lem: det2} shows that $X_\omega$, $Y_\omega$, and $R_\omega$ may always be regarded as elements of $\pglQx$.
\end{dfns}

The following theorem is the main result of this section. It shows that the $\omega$-crossing loop $S_\omega$ (from Definition \ref{def: S and T}) is conjugate to the square of the recursion matrix $X_\omega$ and that the sequences $\{ \calV(\beta_j) \}_0^\infty$ and $\{ -1+\calV(\beta_{-j})\}_1^\infty$ are the forward and reverse orbits of $\calV(\beta_0)$ under the M\"obius action of $R_\omega$.  (See Definition \ref{def: calV} for $\calV$.) \index{crossing loop!holonomy}

\begin{thm} \label{thm: spiral} 
If $\omega \in \Q \cap (0,1/2]$, then
\[ R_\omega^2=S_\omega \qquad \text{and} \qquad  R_\omega^j \left( \calV(\beta_0) \right) \ = \ \begin{cases} \calV(\beta_j) & \text{if } j\geq 0 \\
-1 + \calV(\beta_j) & \text{if } j<0.\end{cases}\]
For the matrix representatives in Definition~\ref{def: XYR},
\[S_\omega\ =\ -\frac{1}{d_\calQ(\omega)}R_\omega^2.\]
So, as elements of $\pglQx$, $S_\omega = R_\omega^2$.
\end{thm}

\proof This proof utilizes the notation outlined at the beginning of the current section.  Begin by setting
\[Z\ =\ \begin{bmatrix} 1&-1 \\ 0 & 1 \end{bmatrix} \, \begin{bmatrix} -dN_{-1} & N_{-2} \\ -dQ_{-1} & Q_{-2} \end{bmatrix}.\]
Because the sequences $\{ N_j \}_0^\infty$, $\{ Q_j \}_0^\infty$, $\{ Q_{-j} \}_1^\infty$, and $\{ N_{-j} - Q_{-j}\}_1^\infty$ all have recursion $X_\omega$, Lemma \ref{lem: quotients} gives that
\begin{align} 
V_j &= R_\omega^j \left(V_0 \right) & \text{if } j\geq 0 \label{eq: actionI}\\ 
-1 +V_j &= ZX_\omega^{-j-1}Z^{-1}\left(-1+V_{-1}\right) & \text{if } j<0. \label{eq: actionII} 
\end{align}

The last line of Lemma \ref{lem: STbeta} gives that, for $j\geq0$,
\[ S_\omega^j \left( V_0 \right) \ =\ V_{2j}.\]
Together with \ref{eq: actionI}, this shows that $R_\omega^2=S_\omega$.

To finish the proof, it is enough to show that 
\[ ZX_\omega Z^{-1}=R_\omega^{-1}\qquad \text{and} \qquad R_\omega^{-1} \left( V_0 \right) = -1 + V_{-1}.\]
Lemma \ref{lem: STbeta} shows that, for $j<0$,
\[ S_\omega^{j} \left( -1+V_{-1} \right) \ = \ -1+ V_{2j-1}.\] 
Thus, using \ref{eq: actionII}, we obtain $S_\omega^{-1}=ZX_\omega^{2}Z^{-1}$ and so $R_\omega^{-1}=ZX_\omega Z^{-1}$.

Again from Lemma \ref{lem: STbeta}, 
\[-1+V_{-1} =\ S_\omega^{-1} \left( V_1\right).\]
Since $V_1=R_\omega\left( V_0 \right)$ and $R_\omega^{-1}=S_\omega^{-1}R_\omega$, the second equality holds.
\endproof

Next, consider the special case $\omega=0$.  Here, only the vertices on the right side of the Stern--Brocot triangle $\Delta(\omega)$ arise as vertices of funnels $F_\alpha$ and a slight variation of Theorem \ref{thm: spiral} is desirable.  Recall the involution
\[ \sigma_u \ = \ \begin{bmatrix} -1&1\\0&1 \end{bmatrix} \]
from Definition \ref{def: concrete involutions}.    It acts on $\IH^3$ as the order-2 involution that fixes the points of the geodesic from $1/2$ to $\infty$.

\begin{thm} \label{thm: w=0}
Suppose $\omega=0$ and, for $j \geq 0$, define $\beta_j =  1/j$.  Then 
\[
X_\omega \ = \ \begin{bmatrix} 0&1 \\ -x & 1 \end{bmatrix} \qquad 
\text{and} \qquad Y_\omega \ =\ \begin{bmatrix} -x & 1 \\ 0 & 1 \end{bmatrix}.\] 
Also, $R_\omega^2=S_0$ and, for $j \geq 0$,
\[ \calV(\beta_j) =  R_\omega^j (\beta_0) \qquad \text{and} \qquad \sigma_u \calV(\beta_j) = R_\omega^{-j} (\beta_0).\]
\end{thm}

\proof  Recall that $R_\omega = Y_\omega X_\omega Y_\omega^{-1}$
and let 
\[
\eta = \begin{bmatrix} 0 & x \\ 1 & 0 \end{bmatrix}.
\]
Notice that both $\eta$ and $\sigma_u$ have order two and $R_\omega=\sigma_u \eta$.  Hence, 
\[ R_\omega^{-1} = \eta \sigma_u = \sigma_u \sigma_u \eta \sigma_u = \sigma_u R_\omega \sigma_u.\]
Since $Q(\omega)^2\neq 4\dQ(\omega)$ and $Y_\omega(\infty)=\infty$, Lemma \ref{lem: quotients} applies to give
\[ \calV(\beta_j) = R_\omega^j Y_\omega(\infty) = R_\omega^j (\beta_0)\]
for $j \geq 0$.  Because $\sigma_u$ fixes $\beta_0$,
\[  \sigma_u \calV(\beta_j) = \sigma_u R_\omega^j \sigma_u (\beta_0)  = R_\omega^{-j} (\beta_0).\]
\endproof

This completes the main goal for this section, but it makes sense to pause to take note of a useful formula and a result.

\subsection{Shapes and traces} \label{sec: shapes and traces} Traces in $\pslQx$ are defined only up to sign.  However, Theorems \ref{thm: spiral} and \ref{thm: w=0} show that the representative of the $\omega$-crossing loop $S_\omega$ in $\SL_2\left(\Q(x)\right)$ is conjugate in $\GL_2(\Q(x))$ to
\[ -\frac{1}{d_\calQ(\omega)}\, X_\omega^2. \]
Hence, this choice fixes the sign of $\tr S_\omega$.  Since
\[ \tr\left(X_\omega^2\right)\ = \  \calQ(\omega)^2-2d_\calQ(\omega)\]
a calculation gives
\begin{align}\label{eq: trS one}
\tr S_\omega &= 2-\frac{\calQ(\omega)^2}{\dQ(\omega)}. \end{align}
This trace can also be expressed in terms of the discriminant of the characteristic polynomial of $X_\omega$.  This discriminant will also be relevant later while investigating limits of geometric roots.

\begin{remark} \label{rem: EMS}
The {\it Farey polynomial} for $\omega$ defined in \cite{EMS2} is equal to the polynomial in $\Z[\mu]$ obtained from $\tr S_\omega$ by replacing $x$ with $\mu=-1/x$.  Conjecture 3.8 of \cite{EMS2} concerns the factorization of the polynomial obtained from the Farey polynomial by subtracting $2$.   The conjecture follows immediately from the trace formula \ref{eq: trS one}, because
\[ \tr S_\omega -2 \ = \ -\frac{\calQ(\omega)^2}{\dQ(\omega)}.\]
\end{remark}

\begin{dfn} \label{def: Disc}
The characteristic polynomial of $X_\omega$ is a quadratic polynomial with coefficients in $\Z[x]$.  The {\it discriminant polynomial} $D_\calQ(\omega)$ is the discriminant for this characteristic polynomial.  Specifically,
\[ D_\calQ(\omega) \ =\ \calQ(\omega)^2-4 \dQ(\omega).\]
\end{dfn}
\index{discriminant polynomial|idxdef}

Therefore, it is also true that 
\begin{align} \label{eq: trS two} 
\tr S_\omega \ = \ -2-\frac{D_\calQ(\omega)}{\dQ(\omega)}.
\end{align}

This helps to understand the relationships between the zeros of $\calQ$ and $D_\calQ$ and the specializations of $S_\omega$.

\begin{thm} \label{thm: Sparabolic}
Suppose $\omega \in \Q \cap [0,1/2]$ and $\frz \in \C$.   
\begin{enumerate}
\item $D_\calQ(\omega)_0 \neq 0$.
\item $\calQ(\omega)$ and $D_\calQ(\omega)$ have no common roots.
\item $(S_\omega)_\frz$ is trivial if and only if  $\calQ(\omega)_\frz = 0$. 
\item $(S_\omega)_\frz$ is a nontrivial parabolic if and only if $D_\calQ(\omega)_\frz=0$.
\end{enumerate}
\end{thm}

\proof
Suppose $\frz$ is a root of $D_\calQ(\omega)$. Then $\calQ(\omega)_\frz^2 = 4 \dQ(\omega)_\frz$.  So, if $\frz=0$, then $\calQ(\omega)_\frz =0$.  This would contradict Corollary \ref{cor: degree}.  This proves (1).

Suppose $\frz$ is a root of $\calQ(\omega)$ and of $D_\calQ(\omega)$.  Then 
\[ 0 \ = \ \calQ(\omega)^2_\frz - D_\calQ(\omega)_\frz \ =\ 4 \dQ(\omega)_\frz.\]
But this implies that $\frz=0$, which again contradicts Corollary \ref{cor: degree}.  This proves (2).

If $(S_\omega)_\frz$ is trivial, then the expression for $S_\omega$ given just after Definition \ref{def: S and T} shows that $\calQ(\omega)_\frz=0$.  Conversely, if $\calQ(\omega)_\frz=0$, then $\tr (X_\omega)_\frz =  0$ and $(X_\omega)_\frz$ has order two.  So, using Theorem \ref{thm: spiral}, if $\calQ(\omega)_\frz=0$ then $(S_\omega)_\frz$ is trivial.  This proves (3).

Suppose $(S_\omega)_\frz$ is a nontrivial parabolic.  If $\tr(S_\omega)_\frz=2$, then the trace identity \ref{eq: trS one} gives $\calQ(\omega)_\frz=0$.  So, by part (3), $\tr (S_\omega)_\frz = -2$.  The trace identity \ref{eq: trS two} gives $D_\calQ(\omega)_\frz=0$.  On the other hand, if $D_\calQ(\omega)_\frz=0$, then the same identity gives that $\tr (S_\omega)_\frz = -2$.  Parts (2) and (3) show that $(S_\omega)_\frz$ is nontrivial.
\endproof

\begin{cor} \label{cor: epimorphism}
If $\alpha$ is hyperbolic and $\frz$ is a root of $\calQ(\alpha)$, then there is an epimorphism $\pie M_\alpha \to \Gamma_\frz$.  
\end{cor}

\proof
Recall that the generic holonomy group $\Gamma = \langle U_0, W_0 \rangle$ from Definition \ref{def: generic holonomy} is free of rank two.  By Theorem \ref{thm: generic holonomy}, $S_\alpha \in \Gamma$.  It is well-known that $\pie M_\alpha \cong \Gamma/\langle \! \langle S_\alpha \rangle \! \rangle$.  Hence, the corollary follows from part (3) of Theorem \ref{thm: Sparabolic}.
\index{holonomy!generic}
\endproof

The purpose of much of Section \ref{sec: holonomy} was to make connections between the geometric triangulations determined by a point $\frz$ of $\calH(M_\alpha^\circ)$ and the associated holonomy group.  The next result makes another striking connection between traces and shape parameters.  
\index{geometric triangulation}
\index{holonomy}

Central to this paper are the shape parameter functions 
\[ \calZ \co \calE \longrightarrow \Q(x) \]
from Definition \ref{def: shape parameter function}.  By Theorem \ref{thm: main1}, if $\frz \in \calH(M_\alpha^\circ)$, then the shape parameter for $\Delta_j$, the $j^\text{th}$ pair of tetrahedra in the Sakuma--Weeks complex, is given by $\calZ(e_j)_\frz$, where $e_j$ is the $j^\text{th}$ interior edge of the funnel $F_\alpha$.  The next theorem provides an alternative formula for $\calZ(e)$ in terms of traces of $\omega$-crossing circles.
\index{shape parameter}
\index{shape parameter!function}

\begin{thm} \label{thm: traces and shapes}
\ShapeThm
\end{thm}

\proof
Recall from Definition~\ref{def: shape parameter function} that
\[ \calZ(e) \ =\  -\dQ(\hat e)^{-\nu(e)} \left(\frac{\calQ(e^L)}{\calQ(e^R)}\right)^2. \]
Also, Formula~(\ref{eq: trS one}) gives
\[ 2-\tr S_\omega \ =\  \frac{\calQ(\omega)^2}{\dQ(\omega)}. \]
Therefore,
\[ -\frac{2-\tr S_{e^L}}{2-\tr S_{e^R}} \ =\  -\frac{\dQ(e^R)}{\dQ(e^L)} \left(\frac{\calQ(e^L)}{\calQ(e^R)}\right)^2. \]

By Definition~\ref{def: E and more}, the low vertex satisfies $e^+=e^-\oplus\hat e$.  Since $\dQ$ is multiplicative, $e^L=e^{\nu(e)}$,  and $e^R=e^{-\nu(e)}$, it follows that
\[ \frac{\dQ(e^R)}{\dQ(e^L)} \ =\  \dQ(\hat e)^{-\nu(e)}. \]
Substituting this into the previous expression gives
\[ -\frac{2-\tr S_{e^L}}{2-\tr S_{e^R}} \ =\  -\dQ(\hat e)^{-\nu(e)} \left(\frac{\calQ(e^L)}{\calQ(e^R)}\right)^2 \ =\  \calZ(e)\]
as desired.
\endproof

\section{Specialization} \label{sec: special}
To ensure that the funnel $F_\alpha$ has at least one interior edge and that the Sakuma--Weeks complex for $M_\alpha^\circ$ has at least one pair of ideal tetrahedra, assume throughout this section that $\alpha \in \Q \cap (0,1/2)$ and $\alpha \neq 1/3$.  This section considers the specializations of Theorems \ref{thm: spiral} and \ref{thm: w=0} at values $\frz \in \calH(M_\alpha^\circ)$.  

Recall that the work in Section \ref{sec: ratios} requires $d \neq 0$ and $Q^2 \neq 4d$.  The work in Section \ref{sec: generic spirals} requires the results of Section \ref{sec: ratios}.  By definition of $D_\calQ$, $\calQ(\omega)^2=4\dQ(\omega)$ if and only if $D_\calQ(\omega)=0$.  Therefore, the next lemma is necessary.

\begin{lemma} \label{lem: DQ not zero}
Take $\frz \in \calH(M_\alpha^\circ)$.  Let $\omega$ be a hub of the funnel $F_\alpha$ and assume that $(S_\omega)_\frz$ is not parabolic.
Then $\dQ(\omega)$ and $D_\calQ(\omega)$ are both nonzero at $\frz$.
\end{lemma}

\proof
By definition of $\calH(M_\alpha^\circ)$, $\frz \neq 0$.  Hence $\dQ(\omega)_\frz \neq 0$.  Part (4) of Theorem \ref{thm: Sparabolic} gives that $D_\calQ(\omega)_\frz \neq 0$.
\endproof

Let $\omega$ be a hub of $F_\alpha$.
By Lemma \ref{lem: DQ not zero}, either $(S_\omega)_\frz$ is parabolic or the specialized versions of Theorems \ref{thm: spiral} and  \ref{thm: w=0} are valid.  The results in this section assume that $(S_\omega)_\frz$ is not parabolic and follow immediately from these two theorems.
\index{holonomy!specialization}

\begin{cor} \label{cor: left hub} \index{funnel!hub} \index{neighbor cycle}
Fix $\alpha \in \Q \cap (0,1/2)$ with $\alpha \neq 1/3$, $\frz \in \calH(M_\alpha^\circ)$, and a nonzero hub $\omega$ on the left side of $F_\alpha$.  Assume that $(S_\omega)_\frz$ is not parabolic.  Let $\{ \beta_j\}_{\Z}$ be the neighbor cycle for $\omega$.  Take $k \in \Z_{\geq 1}$ as large as possible under the requirement that $\beta_k$ is a vertex of $F_\alpha$.   
For $j \in \Z$, write 
\[V_j=(R_\omega)^j_\frz \left(\calV(\beta_0)_\frz\right).\]  
Then the points $\{ V_{-1}, \ldots, V_k \}$ are ideal vertices for $\Omega(\alpha)_\frz$.  More specifically, when $0\leq j \leq k-1$,
\[ \delta\left( [ \omega, \beta_j]\right)_\frz \  = \ \left[ -1+\calV(\omega)_\frz, \, V_{j-1}, \, V_j, \, \infty \right] \]
and
\[ \delta'\left( [\omega, \beta_j] \right)_\frz  \  =\  \ \left[ \calV(\omega)_\frz, \, V_j, \, V_{j+1}, \, \infty\right].\]
\end{cor}
\index{fundamental domain!specialized}

\begin{figure}
\setlength{\unitlength}{.1in}
\begin{picture}(46,23)
\put(3,0) {\includegraphics[width= 3.92in]{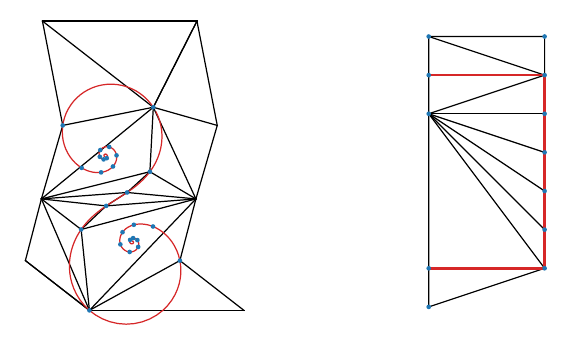}}
\put(4.6,1.2){$\scriptstyle R^5_\omega \calV(\beta_0)$}
\put(16,5.8){$\scriptstyle R^6_\omega \calV(\beta_0)$}
\put(0.6,10.1){$\scriptstyle -1+\calV(\omega)$}
\put(17.2,10.0){$\scriptstyle \calV(\omega)$}
\put(1.6,14.9){$\scriptstyle R^{-1}_\omega \calV(\beta_0)$}
\put(14.2,16.5){$\scriptstyle \calV(\beta_0)$}
\put(27.5,18.4){$\scriptstyle \beta_{-1}=1/3$}
\put(41.2,18.4){$\scriptstyle \beta_0=1/2$}
\put(27.0,5.1){$\scriptstyle \beta_6=13/32$}
\put(41.2,5.1){$\scriptstyle \beta_5=11/27$}
\put(28.6,15.8){$\scriptstyle \omega=2/5$}
\end{picture}
   \caption{\index{fundamental domain!specialized} The fundamental domain $\Omega(\alpha)_{\frz_1}$ for $\alpha=35/86$ with its logarithmic spiral for the hub $\omega=2/5$ of $F_\alpha$.  The red lines shown in the funnel are the lines shared by the funnel and the boundary of the triangle $\Delta(2/5)$ in the Stern--Brocot diagram $\calG$.}  \label{fig: spiral 35/86}
\end{figure}

Compare the corollary to Figure \ref{fig: spiral 35/86}, which illustrates the typical situation.

The next result is a corollary for Theorem \ref{thm: w=0}.  It discusses the special case where the first hub $\omega=0$ is chosen.   Notice that $D_\calQ(0)=-4x+1$, so $(S_0)_\frz$ cannot be parabolic if $\frz \in \IH^2$.  Recall that the upper involution for $M_\alpha^\circ$ acts by isometry on $M_\alpha^\circ$, fixing the points of the upper tunnel.  This isometry lifts to the isometry $\sigma_u$ of $\IH^3$ given in Definition \ref{def: concrete involutions}.

\begin{cor}\label{cor: zero spiral} \index{funnel!hub}
Fix $\alpha \in \Q \cap (0,1/2)$ with $\alpha \neq 1/3$ and take $\frz \in \calH(M_\alpha^\circ)$.  For $j \geq 0$, define $\beta_j=1/j$ and take $k \in \Z_{\geq 1}$ as large as possible under the requirement that $\beta_k$ is a vertex of $F_\alpha$.   
For $j \in \Z$, write 
\[V_j=(R_0)^j_\frz \, (\infty).\]  
The points $\{ V_{-1}, \ldots V_k \}$ are ideal vertices for $\Omega(\alpha)_\frz$ and, when $2\leq j \leq k-1$,
\[ \left( \delta_{j-1} \right)_\frz \ = \ \delta\left( [ 0, 1/j]\right)_\frz \  = \ \left[ 0, \, V_{j-1}, \, V_j, \, \infty \right] \]
and
\[ \left( \delta'_{j-1} \right)_\frz \ = \ \delta'\left( [ 0, 1/j]\right)_\frz \  = \ \left[ 1, \, V_{j}, \, V_{j+1}, \, \infty \right]. \]
If $j \in \Z$, then $\sigma_u (V_j)=V_{-j}$, so the logarithmic spiral for $(R_0)_\frz$ through $\infty$ is invariant under $\sigma_u$.
\end{cor}

\begin{figure}
\setlength{\unitlength}{.1in}
\begin{picture}(48,27)
\put(1,0) {\includegraphics[width= 4.4in]{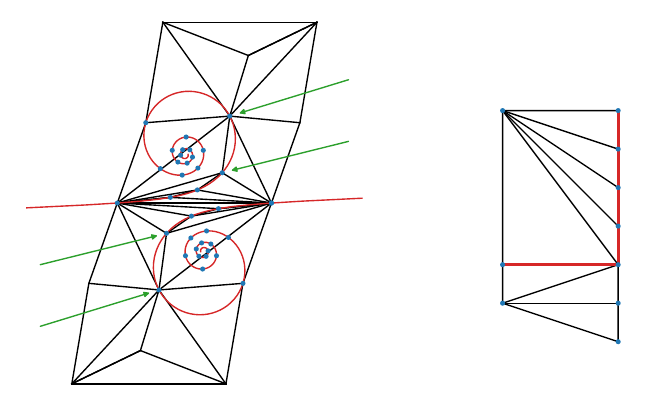}}
\put(19.8,12.7){$\scriptstyle R_\omega \infty$}
\put(5.2,14.5){$\scriptstyle R^{-1}_\omega \infty$}
\put(18,7.8){$\scriptstyle R^6_\omega \infty$}
\put(0.5,5.3){$\scriptstyle R^5_\omega \infty$}
\put(0.5,9.6){$\scriptstyle R^4_\omega \infty$}
\put(25.2,22.3){$\scriptstyle R^{-5}_\omega \infty$}
\put(25.2,18.0){$\scriptstyle R^{-4}_\omega \infty$}
\put(6.9,19.3){$\scriptstyle R^{-6}_\omega \infty$}
\put(44.0,20){$\scriptstyle \beta_1=1$}
\put(31.0,9.4){$\scriptstyle \beta_6=1/6$}
\put(44.0,9.4){$\scriptstyle \beta_5=1/5$}
\put(32.6,20){$\scriptstyle \omega=0$}
\end{picture}
\caption{\index{fundamental domain!specialized} The fundamental domain $\Omega(\alpha)_{\frz_1}$ for $\alpha=7/38$ with its logarithmic spiral for its initial hub $\omega=0$ of $F_\alpha$.}  \label{fig: spiral_7_38}
\end{figure}

Compare the previous corollary to Figure \ref{fig: spiral_7_38} which also illuminates the following discussion.  Assuming the hypotheses and notation of Corollary \ref{cor: zero spiral}, notice that 
\[V_{-1}=0=-1+\calV(0), \qquad V_0=\infty, \qquad \text{and} \qquad V_1=1=\calV(0).\]  Because the lifts of tetrahedra in the Sakuma--Weeks complex are invariant under the action of $(\sigma_u)_\frz$, the spiral lies nicely across both $\Omega(\alpha)_\frz$ and $(\sigma_u)_\frz\big(\Omega(\alpha)_\frz\big)$.   If $2 \leq j \leq k-1$, then
\begin{align*}
(\sigma_u)_\frz \left( \delta_{j-1} \right)_\frz &=   \left[ V_{1}, \, V_{1-j}, \, V_{-j}, \, V_0 \right] \\
(\sigma_u)_\frz \left( \delta'_{j-1} \right)_\frz &=  \left[ V_{-1}, \, V_{-j}, \, V_{1-j}, \, V_0 \right]
\end{align*}
which shows that 
\begin{align*}
(\sigma_u)_\frz(\delta_1)_\frz &= S_0^{-1} (\delta'_1)_\frz  \\ 
(\sigma_u)_\frz(\delta'_1)_\frz &= S_0^{-1} (\delta_1)_\frz \\ 
(\sigma_u)_\frz(\delta_2)_\frz &= S_0^{-1} (\delta_2)_\frz.
\end{align*}
This is not surprising since the upper involution preserves the Sakuma--Weeks complex and inverts the side edges of $\Delta_1$ and the top edges of $\Delta_2$.

Results similar to Corollary~\ref{cor: zero spiral} hold for the last hub of $F_\alpha$.  This first version applies when the last hub lies on the left side of $F_\alpha$. Corollary \ref{cor: bottom spiral right} will provide a right-sided version.

To begin, assume that $\alpha$ is hyperbolic, and let $\frz_1\in\calH(M_\alpha)$ be the geometric root of $\calQ(\alpha)$.  In this setting,  the lower involution on $M_\alpha$ fixes the points of the lower tunnel.   The lower involution lifts to  $(\sigma_\ell)_{\frz_1}$ as given in Definition \ref{def: concrete involutions}.  Recall that $(\sigma_\ell)_{\frz_1}$ fixes $\infty$ and interchanges $\calV\left( e_{n+1}^- \right)_{\frz_1}$ and $\calV\left( e_{n+1}^+ \right)_{\frz_1}$.  Therefore, $(\sigma_\ell)_{\frz_1}$ is conjugate to $(\sigma_u)_{\frz_1}$ by any isometry which fixes $\infty$ and takes this pair of points to $0$ and $1$.  This relationship makes it possible to use Theorem \ref{thm: w=0} to obtain a version of Corollary \ref{cor: zero spiral} for the bottom hub.
\index{geometric root!2-bridge link}

\begin{figure}
\setlength{\unitlength}{.1in}
\begin{picture}(46,25)
\put(0,.5) {\includegraphics[width= 4.23in]{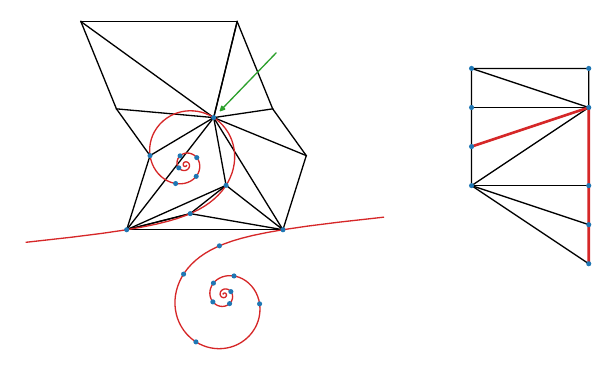}}
\put(18.3,4.7){$\scriptstyle R^{9}_\omega \calV(\beta_0)$}
\put(9,1.2){$\scriptstyle R^{8}_\omega \calV(\beta_0)$}
\put(7.7,6.8){$\scriptstyle R^{7}_\omega \calV(\beta_0)$}
\put(19,8.7){$\scriptstyle R^{5}_\omega \calV(\beta_0)$}
\put(3.8,10.6){$\scriptstyle R_\omega^{3}\calV(\beta_0)$}
\put(4.7,14.7){$\scriptstyle R^{-1}_\omega \calV(\beta_0)$}
\put(19.3,22.2){$\scriptstyle \calV(\beta_0)$}
\put(27.5,15.6){$\scriptstyle \beta_{-1}=2/5$}
\put(41.3,18.4){$\scriptstyle \beta_0=1/2$}
\put(41.3,7){$\scriptstyle \beta_{3}=10/23$}
\put(28.7,13.0){$\scriptstyle \omega=3/7$}
\end{picture}
\caption{\index{fundamental domain!specialized} The fundamental domain $\Omega(\alpha)_{\frz_1}$ for $\alpha=13/30$ with its logarithmic spiral for its last hub $\omega=3/7$ of $F_\alpha$.  Here, $R^4 \,\calV(\beta_0)_{\frz_1} =\infty$.  The spiral is invariant under $(\sigma_\ell)_{\frz_1}$ and lies nicely across $\Omega(\alpha)_{\frz_1} \cup (\sigma_\ell)_{\frz_1} \, \Omega(\alpha)_{\frz_1}$.}  \label{fig: bottom_spiral I}
\end{figure}

For the left-hub case considered below, Corollary~\ref{cor: geodesic crossing circles} shows that
$(S_\omega)_{\frz_1}$ is not parabolic.  See Figure~\ref{fig: bottom_spiral I} for an example.

\begin{cor}\label{cor: bottom spiral left} \index{funnel!hub} \index{neighbor cycle}
Suppose $\alpha$ is hyperbolic and let $\frz_1$ be the geometric root for $\calQ(\alpha)$.  Assume that the last hub $\omega=e_{n+1}^-$ of $F_\alpha$ is a left hub.  Let $\{\beta_j\}_{\Z}$ be the neighbor cycle for $\omega$ and take $k\in\Z_{\geq1}$ as large as possible under the requirement that $\beta_k$ is a vertex of $F_\alpha$.  For $j \in \Z$, write 
\[V_j=(R_\omega)^j_{\frz_1} \, \left(\calV(\beta_0)_{\frz_1}\right).\]  
The points $\{ V_{-1}, \ldots V_k \}$ are ideal vertices for $\Omega(\alpha)_{\frz_1}$ and, when $0\leq j \leq k-1$,
 \[ \delta\left( [ \omega, \beta_j]\right)_{\frz_1} \ =\  \left[ V_{k}, \, V_{j-1}, \, V_j, \, \infty \right]\]
 and
 \[ \delta'\left( [ \omega, \beta_j]\right)_{\frz_1} \ =  \ \left[ V_{k+2}, \, V_{j}, \, V_{j+1}, \, \infty \right]. \]
The involution $(\sigma_\ell)_{\frz_1}$ fixes $V_{k+1}=\infty$ and interchanges $V_j$ and $V_{2k+2-j}$ for each integer $j$, so the logarithmic spiral for $(R_\omega)_{\frz_1}$ through $\infty$ is invariant under $(\sigma_\ell)_{\frz_1}$.
\end{cor}

Here, the equality $V_{k+1}=\infty$ follows from $\beta_{k+1}=\alpha$ and Theorem~\ref{thm: spiral}, since
\[ V_{k+1} \ = \ \calV(\beta_{k+1})_{\frz_1}\ =\  \calV(\alpha)_{\frz_1} \ =\  \infty. \]
Here again, the lifts of tetrahedra in the Sakuma--Weeks complex are invariant under the action of $(\sigma_\ell)_{\frz_1}$ and the spiral lies nicely across both $\Omega(\alpha)_{\frz_1}$ and $(\sigma_\ell)_{\frz_1}\big(\Omega(\alpha)_{\frz_1}\big)$.

The next results are the corresponding corollaries for hubs on the right side of $F_\alpha$.  Only slight adjustments need to be made.  Recall, from Definition \ref{def: UW}, that $U_0$ is the M\"obius function which acts on $\C$ by $z \mapsto z-1$.

\begin{cor} \label{cor: right hub} \index{funnel!hub} \index{neighbor cycle}
Fix $\alpha \in \Q \cap (0,1/2)$ with $\alpha \neq 1/3$, $\frz \in \calH(M_\alpha^\circ)$, and a hub $\omega$ on the right side of $F_\alpha$.  Assume that $(S_\omega)_\frz$ is not parabolic.  Let $\{\beta_j\}_{\Z}$ be the neighbor cycle for $\omega$ and take $k <0$ so that $\beta_k$ is a vertex of $F_\alpha$ and the denominator of $\beta_k$ is as large as possible.   
For $j \in \Z$, write 
\[V_j=(R_\omega)^j_\frz \left(\calV(\beta_0)_\frz\right).\]  
The points $\{ V_0, \ldots V_{k+1} \}$ are ideal vertices for $\Omega(\alpha)_\frz$.  When $k+1 \leq j \leq -1$,
\[ \delta\left( [ \omega, \beta_j]\right)_\frz \  = \ \left[ \calV(\omega)_\frz, \, V_{j+1}, \, V_j, \, \infty \right] \]
and
\[ U_0 \left( \delta'\left( [\omega, \beta_j] \right)_\frz \right) \  =\  \ \left[ -1+\calV(\omega)_\frz, \, V_j, \, V_{j-1}, \, \infty\right].\]
\end{cor}

\begin{figure}
\setlength{\unitlength}{.1in}
\begin{picture}(48,23)
\put(4.2,1) {\includegraphics[width= 4in]{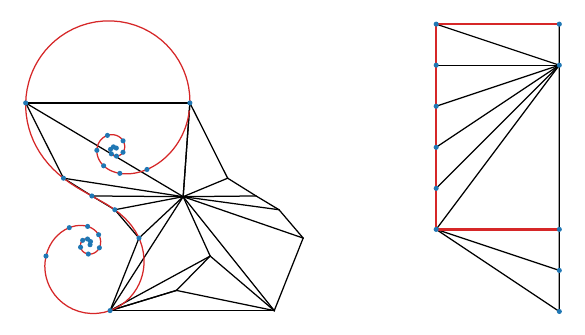}}
\put(9,1){$\scriptstyle R^{-6}_\omega \calV(\beta_0)$}
\put(1.9,6.4){$\scriptstyle R^{-7}_\omega \calV(\beta_0)$}
\put(3.3,10.7){$\scriptstyle R_\omega^{-2}\calV(\beta_0)$}
\put(0.6,16.8){$\scriptstyle R^{-1}_\omega \calV(\beta_0)$}
\put(17.8,16.8){$\scriptstyle \calV(\beta_0)$}
\put(30.1,22){$\scriptstyle \beta_{-1}=0$}
\put(28.5,7.9){$\scriptstyle \beta_{-6}=5/11$}
\put(43.1,22){$\scriptstyle \beta_0=1$}
\put(43.1,7.9){$\scriptstyle \beta_{-7}=6/13$}
\put(43.1,19.2){$\scriptstyle \omega=1/2$}
\end{picture}
   \caption{\index{fundamental domain!specialized} The fundamental domain $\Omega(\alpha)_{\frz_1}$ for $\alpha=21/46$ with its logarithmic spiral for the hub $\omega=1/2$ of $F_\alpha$.  Notice that $1+R_\omega^{-7} \, \calV(\beta_0) = \calV(6/13)$.}  \label{fig: spiral 21/46}
\end{figure}

See Figure \ref{fig: spiral 21/46} for a typical example.  

In the next corollary, as in Corollary \ref{cor: bottom spiral left}, $(S_\omega)_{\frz_1}$ cannot be parabolic.

\begin{cor}\label{cor: bottom spiral right} \index{funnel!hub} \index{neighbor cycle}
Suppose $\alpha$ is hyperbolic and let $\frz_1$ be the geometric root for $\calQ(\alpha)$.  Assume that the last hub $\omega=e_{n+1}^-$ of $F_\alpha$ is a right hub.  Let $\{\beta_j\}_{\Z}$ be the neighbor cycle for $\omega$ and take $k <0$ so that $\beta_k$ is a vertex of $F_\alpha$ and the denominator of $\beta_k$ is as large as possible. For $j \in \Z$, write 
\[V_j=(R_\omega)^j_{\frz_1} \, \left(\calV(\beta_0)_{\frz_1}\right).\]  
The points $\{ V_0, \ldots V_{k+1} \}$ are ideal vertices for $\Omega(\alpha)_{\frz_1}$ and, when $k+1 \leq j \leq -1$,  
\[ \delta\left( [ \omega, \beta_j]\right)_{\frz_1} \ =  \ \left[ V_{k}, \, V_{j+1}, \, V_j, \, \infty \right] \]
and
\[ U_0 \left( \delta'\left( [ \omega, \beta_j]\right)_{\frz_1}\right) \ = \  \left[ V_{k-2}, \, V_{j}, \, V_{j-1}, \, \infty \right].\]
Define the involution
\[ \tau \ =\ U_0 \, (\sigma_\ell)_{\frz_1} \,  U_0^{-1}.\]
This involution  fixes $V_{k-1} = \infty$ and interchanges $V_j$ and $V_{2k+2-j}$ for each integer $j$, so the logarithmic spiral for $(R_\omega)_{\frz_1}$ through $\infty$ is invariant under $\tau$.
\end{cor}

\begin{figure}
\setlength{\unitlength}{.1in}
\begin{picture}(46,19)
\put(0,0) {\includegraphics[width= 4.42in]{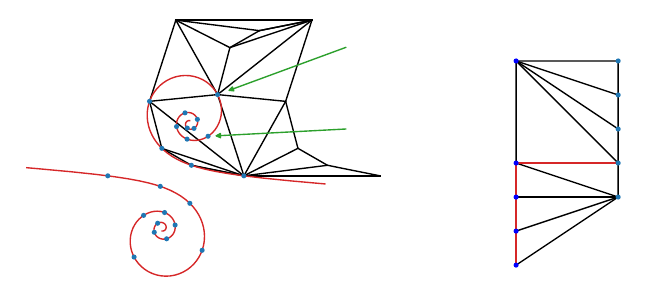}}
\put(14.1,2.3){$\scriptstyle R^{-9}_\omega \calV(\beta_0)$}
\put(4,6.7){$\scriptstyle R^{-6}_\omega \calV(\beta_0)$}
\put(4.6,13){$\scriptstyle R_\omega^{-1}\calV(\beta_0)$}
\put(14.5,6.7){$\scriptstyle R_\omega^{-4}\calV(\beta_0)$}
\put(24,11.2){$\scriptstyle R_\omega \calV(\beta_0)$}
\put(24,16.9){$\scriptstyle \calV(\beta_0)$}
\put(30.4,8.7){$\scriptstyle \beta_{-1}=1/5$}
\put(30,1.3){$\scriptstyle \beta_{-4}=7/32$}
\put(43.1,8.9){$\scriptstyle \beta_0=1/4$}
\put(43.1,6.6){$\scriptstyle \omega=2/9$}
\end{picture}
\caption{\index{fundamental domain!specialized} The fundamental domain $\Omega(\alpha)_{\frz_1}$ for $\alpha=9/41$ with its logarithmic spiral for its last hub $\omega=2/9$ of $F_\alpha$.  Here, $R^{-5} \, \calV(\beta_0)_{\frz_1} = \calV(9/41)_{\frz_1} = \infty$.  The spiral is invariant under the conjugate of $(\sigma_\ell)_{\frz_1}$ by $U_0$, this conjugate is a lift of the lower involution to $\IH^3$.  Both the upper and lower parts of the spiral lie nicely across $\Omega(\alpha)_{\frz_1} \cup U_0 \left( \sigma_\ell \, \Omega(\alpha)_{\frz_1}\right)$.}  \label{fig: spiral 9_41}
\end{figure}

See Figure \ref{fig: spiral 9_41} for an example.

\chapter{Applications and Observations}\label{chap: Applications}

\section{Real loci for shape parameter functions} \label{sec: real loci}

As in Chapters~\ref{chap: SakumaWeeks} and~\ref{chap: FD}, whenever $M_\alpha^\circ$ is considered, assume that $\alpha\in\Q\cap(0,1/2)$ and $\alpha\neq1/3$.  The sets $\calH(M_\alpha^\circ)$ (Definition~\ref{def: hyperbolic structures}) are of central importance.  Using Corollary ~\ref{cor: main1}, $\calH(M_\alpha^\circ)$ is taken as a subset of $\IH^2$  namely
\[ \calH(M_\alpha^\circ)  \ =\ \left\{ \frz \in \IH^2 \, \big| \, \calZ(e_j)_\frz \in \IH^2 \text{ for every interior edge } e_j \text{ of } F_\alpha \right\}.\]
If $\alpha$ is hyperbolic (Definition~\ref{def: hyperbolic alpha}), then $\calQ(\alpha)$ has a geometric root $\frz_1$ (Definition ~\ref{def: geometric root}).  Because each shape parameter for a point in $\calH(M_\alpha^\circ)$ must have positive imaginary part, given an interior edge $e_j$ of the funnel $F_\alpha$, the set $\calH(M_\alpha^\circ)$ is contained in the open region
\[  \left\{ \frz \in \IH^2 \, \big| \,  \calZ(e_j)_\frz \in \IH^2 \right\}.\]

\begin{dfn} \label{dfn: real locus} 
Given a nonhorizontal edge $e$ in the edge set $\calE$, the {\it real shape locus} for $e$ is the real algebraic set
\[  \R \calZ(e) \ = \ \left\{ \frz \in \IH^2 \, \big| \,  \calZ(e)_\frz \in \R \right\}.\]
\end{dfn}
\index{shape parameter!real shape locus|idxdef}
\index{shape parameter!function}

\begin{figure}
\setlength{\unitlength}{.1in}
\begin{picture}(47.7,22)
\put(2.8,0) {\includegraphics[width= 4.5in]{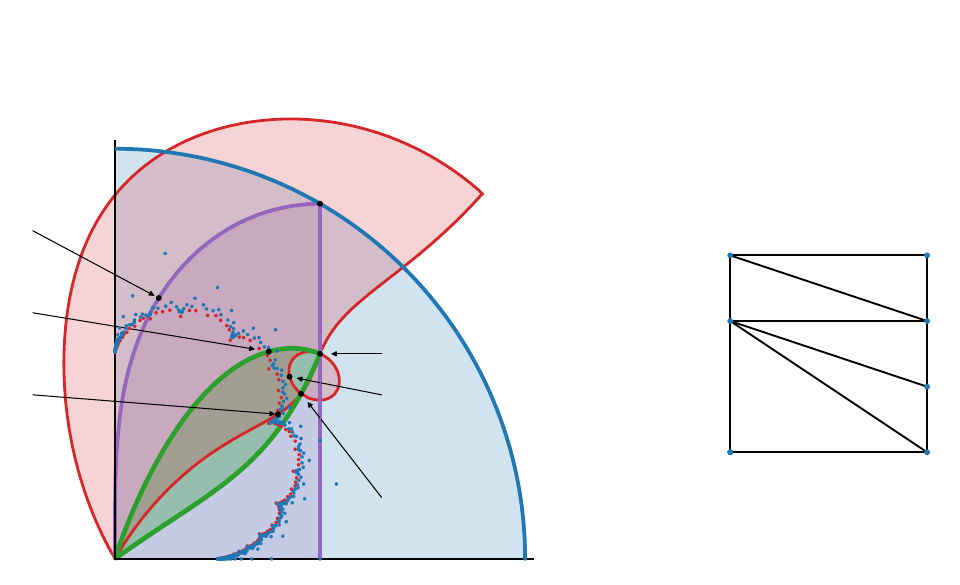}}
\put(18.2,18.4){$\scriptstyle \frz(2/5)$}
\put(21,10.7){$\scriptstyle \frz(3/8)$}
\put(21,8.7){$\scriptstyle \frz(7/19)$}
\put(21,3.5){$\scriptstyle \frz(4/11)$}
\put(0,12.7){$\scriptstyle \frz(21/55)$}
\put(0,8.8){$\scriptstyle \frz(20/57)$}
\put(0.5,16.8){$\scriptstyle \frz(7/16)$}
\put(46.7,5.8){$\scriptstyle 3/8$}
\put(46.7,9.1){$\scriptstyle 2/5$}
\put(46.7,12.2){$\scriptstyle 1/2$}
\put(46.7,15.4){$\scriptstyle 1$}
\put(34.9,12.2){$\scriptstyle 1/3$}
\put(34.4,5.8){$\scriptstyle 4/11$}
\put(36.1,15.4){$\scriptstyle 0$}

\end{picture}
\caption{\index{shape parameter!real shape locus} The diagram on the right shows the combinatorics of a portion of the Stern--Brocot diagram.  The figure on the left shows regions where $\calZ(e)$ takes values with positive imaginary part when $e$ is an interior (or bottom) edge of the diagram.  The region for $[1/2,1/3]$ is blue, for $[1/3,2/5]$ is purple, for $[1/3,3/8]$ is green, and for $[4/11,3/8]$ is red. Labeled points lie on boundary curves as explained by Theorem ~\ref{thm: real loci}.}  \label{fig: RZ}
\end{figure}

\begin{examples} \label{ex: real loci}
When the vertices of $e \in \calE$ have small denominators, it is not difficult to compute $\R \calZ(e)$.  
\begin{enumerate}
\item When $e$ is the edge between $0$ and $1/2$, then $\calZ(e)=x$.  Hence, $\calZ(e)$ is never real on $\IH^2$ and \[\calH\big( M_{1/4}^\circ\big) \ = \ \calH\big( M_{2/5}^\circ\big) \ = \ \IH^2.\]
\item When $e$ is the edge between $0$ and $1/3$, then 
\[\calZ(e) \ =\ \frac{x^2}{(x-1)^2}.\]  
The M\"obius transformation $x \mapsto x/(x-1)$ has real coefficients and preserves $\R\cup\{\infty\}$.  Therefore, as a function on $\IH^2$, the shape parameter function $\calZ(e)$ is real precisely on the image of $i\R \cup \{ \infty \}$ under the map $x \mapsto x/(x-1)$.  Hence, $\R \calZ(e)$ is the intersection of $\IH^2$ with the circle with radius $1/2$ and center $1/2$.  The sets $\calH\big(M_{2/7}^\circ \big)$ and $\calH\big(M_{1/5}^\circ \big)$ are equal to the interior of the corresponding semicircular disk.

\item When $e$ is the edge between $1/2$ and $1/3$, then 
\[\calZ(e) = - \frac{(x-1)^2}{x}.\]  
Here, $\R \calZ(e)$ is the intersection of $\IH^2$ with the unit circle.  This is the blue curve in Figure ~\ref{fig: RZ}.   The sets $\calH\big(M_{3/8}^\circ\big)$ and $\calH(M_{3/7}^\circ\big)$ are the intersection of $\IH^2$ with the open unit disk.
\item Other examples of the real algebraic curves $\R \calZ(e)$ can be computed.  The curves for the edges $[1/3,2/5]$, $[1/3,3/8]$, and $[4/11,3/8]$ are shown in Figure ~\ref{fig: RZ}.
\end{enumerate}
As seen in the figure, certain geometric roots lie on the real shape loci.  This will be explained in Theorem ~\ref{thm: real loci}. 
\end{examples}

\begin{dfn} \label{def: fat}
Suppose $e$ is a nonhorizontal edge in $\calE$ and let $\alpha=e^- \oplus^2 e^+$.  Since the numerator of $e^+$ is at least one, $\alpha$ is hyperbolic and the geometric root $\frz(\alpha)$ exists.  The {\it thick region} for $\calZ(e)$ is the connected component $\calU_e$  of $\IH^2 - \R\calZ(e)$ which contains the point $\frz(\alpha)$.  
\end{dfn}
\index{thick region|idxdef}

Example ~\ref{ex: real loci} and Figure ~\ref{fig: RZ} provide some examples $\calU_e$.  

\section{Reversing continued fractions} \label{sec: reverse}

\begin{dfns} \label{def: reverse}
Suppose that $\alpha\in\Q\cap(0,1/2)$ and write $\alpha=[0;a_1,\ldots,a_k]$ for its standard positive termed continued fraction expansion with $a_k\geq2$. Its {\it reverse} is $\bar{\alpha} = [0; a_k, \ldots, a_1]$.  A number equal to its reverse is called a {\it palindrome}.
\end{dfns}
\index{reverse|idxdef} \index{palindrome|idxdef}

By Definition ~\ref{def: diagram}, the diagrams $D_\alpha$ and $D_{\bar{\alpha}}$ correspond to links that differ, at most, by a reflection.

Assume, for the remainder of the section, that $\alpha$ is given by a positive termed continued fraction $[0; a_1, \ldots, a_k]$, where $k\geq 2$ and $a_1, a_k \geq 2$.  These assumptions imply that $\alpha$ is hyperbolic.  In particular, $F_\alpha$ contains an interior edge, $M_\alpha$ is a hyperbolic $2$-bridge link complement, and its Sakuma--Weeks triangulation contains at least one pair of ideal tetrahedra.  Notice that $\bar \alpha$ also satisfies these standing assumptions.

From Definition ~\ref{def: hubs}, it follows that $e_{n+1}^-$ is the last hub for the funnel $F_\alpha$.  Write $r/s=e_{n+1}^-$.  Work in \cite{Hat} shows that $r/s$ can also be obtained as the truncated continued fraction $[0;a_1, \ldots, a_{k-1}]$.  Sublemma II.3.4 of \cite{SW} gives $rq-ps = (-1)^k$.  Therefore,
\[ A_\alpha \ = \ \begin{bmatrix} -s & r \\ -q &p \end{bmatrix}\]
is an element of $\PGL_2 \Z$ with determinant $(-1)^k$.   The element $A_\alpha$ determines a M\"obius function $\Phi_\alpha \co \wh{\Q} \to \wh{\Q}$.  Since $\det A_\alpha=\pm1$, the map $\Phi_\alpha$ preserves Farey pairs.  In general, $\Phi_\alpha$ does not preserve $\Q\cap[0,1]$, but the next theorem shows that its restriction to the vertices of the funnel $F_\alpha$ extends to a simplicial bijection $F_\alpha \to F_{\bar \alpha}$.

\begin{thm} \label{thm: reverse}
The restriction of $\Phi_\alpha$ to vertices of $F_\alpha$ extends to a simplicial bijection \[ \Phi_\alpha \co F_\alpha \longrightarrow F_{\bar \alpha}.\]
It takes the $j^\text{th}$ interior edge of $F_\alpha$ to the $(n+1-j)^\text{th}$ interior edge of $F_{\bar{\alpha}}$.  It takes hubs to hubs and takes the bottom edge of $F_\alpha$ to the top edge of $F_{\bar{\alpha}}$.  Moreover $\Phi_\alpha(\alpha)=\infty$ and $\Phi_\alpha(\infty) = \bar{\alpha}$.  In particular, $\bar{\alpha}=s/q$.
\end{thm}

\proof
It is immediate that
\[ \Phi_\alpha(\alpha)\ =\ \infty,\qquad \Phi_\alpha(e^-_{n+1})\ =\ 0,\AND \Phi_\alpha(e^+_{n+1})\ =\ 1.\]
The preceding three identities, together with the continued fraction description of the Farey triangle sequence in Section~II.3 of \cite{SW}, show that $\Phi_\alpha(F_\alpha)=F_{\bar\alpha}$ as described.
\endproof

The existence of the reversal symmetry $\Phi_\alpha$ is a consequence of the combinatorial classification in Theorem~II.3.1 of \cite{SW}, which says that a pair of Sakuma--Weeks decompositions are combinatorially equivalent precisely when their continued fraction sequences agree up to reversal.  The direct
description here is kept because the explicit action of $\Phi_\alpha$ on the edges of the funnel is needed below.

If $e_0$ and $e_{n+1}$ are the top and bottom edges of $F_\alpha$ and $e_1, \ldots, e_n$ are the interior edges of $F_\alpha$, listed from top to bottom, then $\Phi_\alpha(e_0)$ and $\Phi_\alpha(e_{n+1})$ are the bottom and top edges of $F_{\bar{\alpha}}$ and $\Phi_\alpha(e_1), \ldots, \Phi_\alpha(e_n)$ are the interior edges of $F_{\bar{\alpha}}$, listed from bottom to top.  If $e_{n+1}^-$ lies on the left side of $F_\alpha$, then $\Phi_\alpha$ preserves the left-right ordering of the endpoints of every $e_j$.  Otherwise, it reverses the orderings.

The simplicial homeomorphism $\Phi_\alpha\co F_\alpha\longrightarrow F_{\bar\alpha}$ induces a simplicial homeomorphism $\phi_\alpha\co M_\alpha\longrightarrow M_{\bar\alpha}$ of the corresponding Sakuma--Weeks triangulations.  The map $\phi_\alpha$ interchanges the top and bottom of the Sakuma--Weeks triangulation and takes the lower tunnel to the upper tunnel (and vice versa).  Its effect on orientation is opposite to the effect of $\Phi_\alpha$ on the oriented Stern--Brocot diagram.  Since $\det A_\alpha=(-1)^k$, the following are equivalent \index{core tunnel}
\begin{enumerate}
\item $\phi_\alpha$ is orientation preserving;
\item $k$ is odd;
\item $\nu(e_{n+1})=-1$.
\end{enumerate} 
This comes from Lemma~II.3.3 and Sublemma~II.3.4 of \cite{SW} and the orientation calculation on page 416 of \cite{SW}.

\begin{remark}\label{rem: SWautomorphisms}
The symmetry $\phi_\alpha$ is part of the complete calculation of the combinatorial automorphism group of the Sakuma--Weeks decomposition
from Theorem~II.3.2 of \cite{SW}.  If $\alpha$ is not a palindrome, then the full combinatorial automorphism group of the Sakuma--Weeks triangulation is the Klein $4$-group $G$.  If $\alpha$ is a palindrome, then $\phi_\alpha$ has order two and provides an additional symmetry.  In this case, the full
combinatorial automorphism group is generated by $G$ and $\phi_\alpha$.
\end{remark}

\begin{figure}
\setlength{\unitlength}{.1in}
\begin{picture}(42,20)
\put(1.7,0) {\includegraphics[width= 3.5in]{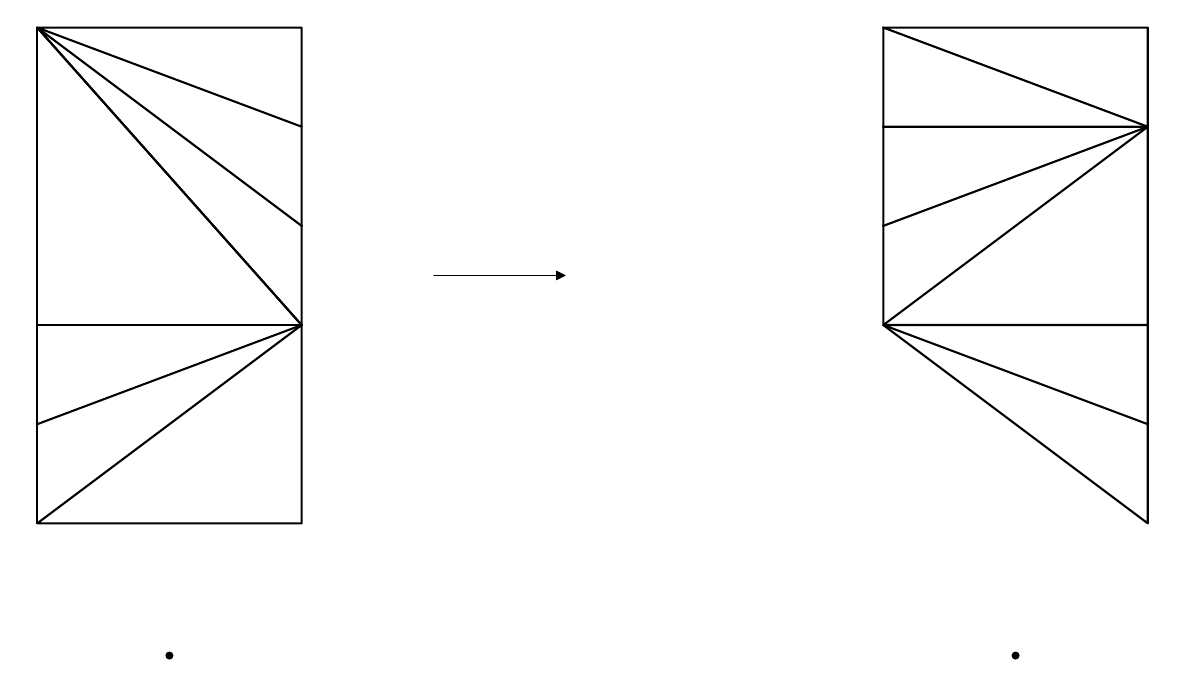}}
\put(1.6,19){$\scriptstyle 0$}
\put(.7,10.2){$\scriptstyle 1/5$}
\put(.7,7.3){$\scriptstyle 2/9$}
\put(36,7.3){$\scriptstyle 7/16=\Phi_\alpha(1/2)$}
\put(0,4.5){$\scriptstyle 3/13$}
\put(36,4.5){$\scriptstyle 10/23=\Phi_\alpha(1)$}
\put(11.2,19){$\scriptstyle 1$}
\put(11.2,16.1){$\scriptstyle 1/2$}
\put(21,16.1){$\scriptstyle \Phi_\alpha(2/9)=1/3$}
\put(36,16.1){$\scriptstyle 1/2=\Phi_\alpha(1/4)$}
\put(11.2,13.2){$\scriptstyle 1/3$}
\put(21,13.2){$\scriptstyle \Phi_\alpha(1/5)=2/5$}
\put(11.2,10.3){$\scriptstyle 1/4$}
\put(22,10.3){$\scriptstyle \Phi_\alpha(0)=3/7$}
\put(36,10.3){$\scriptstyle 4/9=\Phi_\alpha(1/3)$}
\put(11.2,4.5){$\scriptstyle 4/17$}
\put(7.5,.6){$\scriptstyle \alpha=7/30$}
\put(32.3,.6){$\scriptstyle \bar{\alpha} =13/30 =\Phi_\alpha(\infty)$}
\put(16,12.6){$\scriptstyle \Phi_\alpha$}
\put(21.5,19){$\scriptstyle \Phi_\alpha(3/13)=0$}
\put(36,19){$\scriptstyle 1=\Phi_\alpha(4/17)$}
\end{picture}
\caption{If $\alpha = 7/30$ then $\bar{\alpha}=13/30$.  The simplicial bijection $\Phi_\alpha$ takes $F_\alpha$ to $F_{\bar{\alpha}}$. }  \label{fig: flip}
\end{figure}

Let $\frz_1$ be the geometric root of $\calQ(\alpha)$ and write $\frz_j=\calZ(e_j)_{\frz_1}$ for the shape parameter of the $j$th tetrahedral pair in the
geometric Sakuma--Weeks triangulation of $M_\alpha$.  The map $\phi_\alpha$ carries the tetrahedral pair associated to an interior edge $e$ of $F_\alpha$ to the tetrahedral pair associated to $\Phi_\alpha(e)$.  It is natural, then, to compare $\frz_j$ with the parameter for the pair $\phi_\alpha(\Delta_j)$ in the Sakuma--Weeks triangulation of $M_{\bar\alpha}$.

Recall from Section ~\ref{sec: equations} that the Sakuma--Weeks triangulation is canonical in the sense of \cite{EP}.  Together with Mostow--Prasad rigidity, this allows $\phi_\alpha$ to be realized by an isometry between the hyperbolic manifolds $M_\alpha$ and $M_{\bar \alpha}$ which respects the corresponding canonical triangulations.  Henceforth, always take $\phi_\alpha$ to be this isometry.  This map interchanges the top and bottom edges of corresponding tetrahedra.  Recall from Definition~\ref{def: parameters} that the shape parameter is the edge parameter of the top edge and that the opposite bottom edge has the same edge parameter.

Suppose first that $\phi_\alpha$ is orientation preserving. The top edge is taken to the bottom edge of the corresponding
tetrahedron and these opposite edges have the same edge parameter, so the edge parameters of corresponding tetrahedra agree.  Since
$\Phi_\alpha(e_n)$ is the first interior edge of $F_{\bar\alpha}$, the geometric root of $\calQ(\bar\alpha)$ is $\frz_n$.  More
generally,
\[ \frz_j \ =\ \calZ\left( \Phi_\alpha(e_j) \right)_{\frz_n}\]
for every $j$.

\begin{figure}
\setlength{\unitlength}{.1in}
\begin{picture}(45,17)
\put(0,0) {\includegraphics[width= 4.5in]{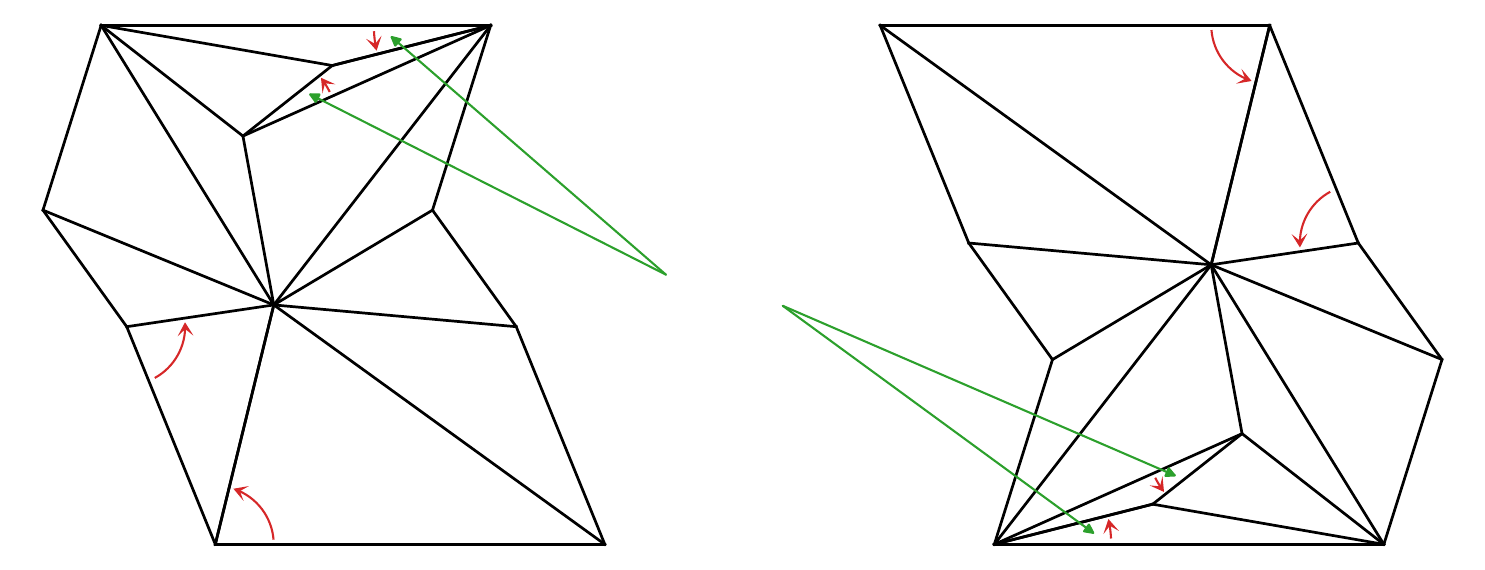}}
\put(21.5,8.4){$\frz_1$}
\put(5.7,6.1){$\scriptstyle \frz_6$}
\put(8.3,2){$\scriptstyle \frz_6$}
\put(38.5,11.0){$\scriptstyle \frz_6$}
\put(35.9,15.1){$\scriptstyle \frz_6$}
\end{picture}
\caption{\index{fundamental domain!specialized} Here $\alpha = 7/30$ and $\bar{\alpha}=13/30$.  The homeomorphism $\phi_\alpha$ is orientation preserving.  Note that $\frz(7/30) \sim 0.406636 + 0.103064 \, i$ and $\frz(13/30) \sim 0.150149 + 0.614630 \, i$. In the figure, $\frz_1$ is the first shape parameter for $7/30$ and $\frz_6$ is the last shape parameter for $7/30$.}  \label{fig: OP}
\end{figure}

\begin{figure}
\setlength{\unitlength}{.1in}
\begin{picture}(45,19)
\put(0,0) {\includegraphics[width= 4.5in]{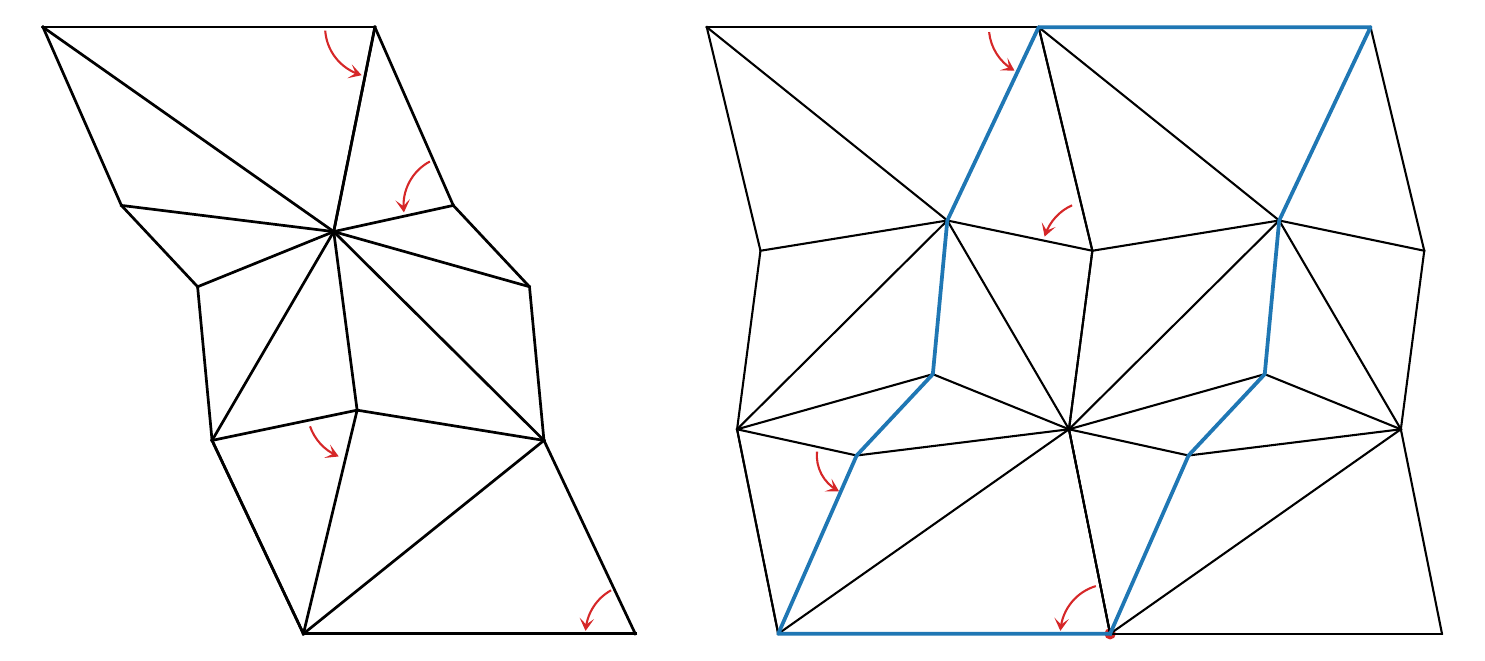}}
\put(9.1,17.9){$\scriptstyle \frz_1$}
\put(11.3,14.7){$\scriptstyle \frz_1$}
\put(8.7,6.0){$\scriptstyle \frz_6$}
\put(16.8,1.8){$\scriptstyle \frz_6$}
\put(28.1,17.8){$\scriptstyle \bar{\frz}_6^{-1}$}
\put(30.6,13.9){$\scriptstyle \bar{\frz}_6^{-1}$}
\put(30.9,1.9){$\scriptstyle \bar{\frz}_1^{-1}$}
\put(23.0,5){$\scriptstyle \bar{\frz}_1^{-1}$}
\end{picture}
\caption{\index{fundamental domain!specialized} Here $\alpha = 17/39$ and $\bar{\alpha}=16/39$.  The homeomorphism $\phi_\alpha$ is orientation reversing.  Calculations provide $\frz(17/39) \sim 0.1240059 + 0.616139 \, i$ and $\frz(16/39) \sim 0.275235 + 0.582608 \, i$.  In the figure, $\frz_1$ is the first shape parameter for $17/39$ and $\frz_6$ is the last shape parameter for $17/39$.}  \label{fig: OR}
\end{figure}

Next, assume instead that $\phi_\alpha$ is orientation reversing.  The edge parameter of the mirror image of an ideal hyperbolic
tetrahedron is obtained from the original by inversion in the unit circle.   Consequently, the geometric root of $\calQ(\bar\alpha)$ is
$\bar{\frz}_n^{-1}$ and, for every $j$,
\[ \bar{\frz}_j^{-1} \ =\ \calZ\left( \Phi_\alpha(e_j) \right)_{\bar{\frz}_n^{-1}}.\]
Figures~\ref{fig: OP} and ~\ref{fig: OR} illustrate the two possibilities.

It is curious that the question of whether the homeomorphism $\phi_\alpha$ is orientation preserving or reversing can be deduced from the location of $\frz_n$.  This curiosity is recorded as the next result.

\begin{thm}\label{thm: phi OP} \index{geometric root!2-bridge link}
Let $\frz_1,\ldots,\frz_n$ be the geometric shape parameters for the Sakuma--Weeks triangulation of $M_\alpha$.  (Recall that $\phi_\alpha$ is orientation preserving when the length of the standard continued fraction expansion for $\alpha$ is odd and orientation reversing when it is even.)

\begin{enumerate}
\item If $\phi_\alpha$ is orientation preserving, then $\frz_n$ is the geometric root of $\calQ(\bar\alpha)$ and, for every $j\in\{1,\ldots,n\}$,
\[ \frz_j \ = \ \calZ\left(\Phi_\alpha(e_j)\right)_{\frz_n}. \]
\item If $\phi_\alpha$ is orientation reversing, then $\bar{\frz}_n^{-1}$ is the geometric root of $\calQ(\bar\alpha)$ and, for every $j\in\{1,\ldots,n\}$,
\[ \bar{\frz}_j^{-1} \ =\ \calZ\left(\Phi_\alpha(e_j)\right)_{\bar{\frz}_n^{-1}}.\]
\item The homeomorphism $\phi_\alpha$ is orientation reversing if and only if
\[ |\frz_n|\ \geq\ 1.\]
\end{enumerate}
\end{thm}

\proof
The first two statements follow from the preceding discussion.

Next, suppose $\phi_\alpha$ is orientation preserving.  Then, using part (1), $\frz_n$ is the geometric root of $\calQ(\bar{\alpha})$.  Remark ~\ref{remark: Riley} (2) implies that $|\frz_n|\leq 1$.  To see that the inequality is strict, suppose that $|\frz_n|=1$.  Remark ~\ref{remark: Riley} (2) then gives $\bar{\alpha}=2/5$.  Since $2/5=[0;2,2]$ is a palindrome, it follows that $\alpha=2/5$. But the continued fraction $[0;2,2]$ has even length, so $\phi_\alpha$ is orientation reversing, a contradiction.  Thus $|\frz_n|<1$.

Now, assume $\phi_\alpha$ is orientation reversing.  Then, by part (2), $\bar{\frz}_n^{-1}$ is the geometric root of $\calQ(\bar{\alpha})$.  Remark ~\ref{remark: Riley} (2) implies that $|\bar{\frz}_n^{-1}|\leq1$, and hence $|\frz_n|\geq1$.
\endproof

Suppose $e \in \calE$ is an edge of the Stern-Brocot diagram whose vertices lie in $[0,1/2]$.  Then $e$ is an edge of the funnel $F_\alpha$ for many choices of $\alpha$.  It can be interesting to find a palindrome $\alpha$ such that
\begin{enumerate}
\item $e$ is an edge of $F_\alpha$ and
\item $\Phi_\alpha(e)=e$.
\end{enumerate}
 The palindromic symmetries of Sakuma--Weeks triangulations are classified in Lemma~II.3.3 of \cite{SW}, and the three possibilities are illustrated in \cite[Figure~II.3.4]{SW}.  The next lemma records them in the notation of this paper.

\begin{lemma} \label{lem: fixed edge}
Suppose $e \in \calE$ is an edge of the Stern--Brocot diagram whose vertices lie in $[0,1/2]$.  Write $e^R=a/b$ and $e^L=c/d$.  Set
\[ \alpha_+ \ =\ \frac{cd-ab}{d^2-b^2} \qquad \alpha_- \ = \ \frac{ab+cd}{b^2+d^2} \qquad \alpha_0 \ = \ \frac{ad+bc}{2bd}.\]
Then
\begin{enumerate}
\item The numbers $\alpha_-$ and $\alpha_0$ lie in $(0,1/2)$. The same is true for $\alpha_+$ unless $e=[0,1/2]$, in which case $\alpha_+=2/3$.  (For the remaining statements, only the numbers $\alpha_\epsilon\in(0,1/2)$ are considered.)
\item A number $\alpha\in\Q\cap(0,1/2)$ is a palindrome for which $e$ is an edge of $F_\alpha$ and $\Phi_\alpha(e)=e$ if and only if $\alpha \in \{ \alpha_\pm, \alpha_0\} \cap (0,1/2)$.
\item The homeomorphism $\Phi_{\alpha_0}\co F_{\alpha_0} \to F_{\alpha_0}$ fixes the endpoints of $e$, while $\Phi_{\alpha_\pm}$ interchange the endpoints of $e$.
\item If $\alpha_-$ is hyperbolic, then the homeomorphism $\phi_{\alpha_-}\co M_{\alpha_-}\to M_{\alpha_-}$ is orientation reversing.  The homeomorphisms $\phi_{\alpha_0}$ and
$\phi_{\alpha_+}$ are orientation preserving whenever the corresponding numbers are hyperbolic.
\end{enumerate}
\end{lemma}

\proof
Since $e^L<e^R$,
\[ \frac cd\ <\ \frac ab \AND  ad-bc\ =\ 1. \]
Notice first that
\[ \alpha_-\ =\ \frac{b^2}{b^2+d^2}\cdot \frac ab+ \frac{d^2}{b^2+d^2}\cdot \frac cd \AND \alpha_0 \ =\  \frac12\left(\frac ab+\frac cd\right).\]
So both $\alpha_-$ and $\alpha_0$ lie strictly between the endpoints of $e$, and belong to $(0,1/2)$.  For $\alpha_+$, direct calculation using $ad-bc=1$ gives
\begin{align}
 \alpha_+ -\frac cd & = \frac{b}{d(b^2-d^2)} \label{eq: alpha1}\\ 
 \alpha_+-\frac ab &=\ \frac{d}{b(b^2-d^2)}. \label{eq: alpha2}
 \end{align}
Suppose first that $d>b$.  Then $b^2-d^2<0$ and, by (\ref{eq: alpha1}), 
\[ \alpha_+ \ < \ \frac{c}{d} \ < \frac{1}{2}.\]
Notice that $c\geq1$ because $d>b$.  Also,
\[ d^2 - b^2 \ =\ (d-b)(d+b) \ \geq \ d+b \ > \ b.\]  
Hence,
\[ c\, (d^2-b^2)-b \ > \ 0.\]
Now, again using (\ref{eq: alpha1}),
\[ cd-ab \  = \ \frac{c(d^2-b^2)-b}{d} \ > \ 0.\]
So, $\alpha_+ \in (0,1/2)$.

Next, assume $d<b$ and $a/b<1/2$.  Then $b-2a \geq 1$ and
\[ \frac{1}{2} -\frac{a}{b} \ = \ \frac{b-2a}{2b} \ \geq \ \frac{1}{2b}.\]
Now, because $b^2-d^2>0$, the equality (\ref{eq: alpha2}) gives
\[ 0 \ < \ \frac{a}{b} \ < \ \alpha_+.\]
Notice,
\[ b^2-d^2 \ =\ (b-d)(b+d) \ \geq \ b+d \ > \ 2d.\]
So, again using (\ref{eq: alpha2}) and the inequality above,
\[ \alpha_+ - \frac{a}{b} \ < \ \frac{1}{2b} \ \leq \ \frac{1}{2} -\frac{a}{b}.\]
Hence, $\alpha_+ \in (0,1/2)$.  

Finally, if $d<b$ and $a/b=1/2$, then the equation $ad-bc=1$ forces $e=[0,1/2]$.  A calculation shows $\alpha_+=2/3$ in this case.  This concludes the argument for (1).

Define
\[ A\ =\ \begin{bmatrix} a&c\\ b&d \end{bmatrix} \quad J_+\ =\ \begin{bmatrix}0&1\\1&0\end{bmatrix} \quad J_-\ =\ \begin{bmatrix}0&-1\\1&0\end{bmatrix} \quad J_0\ =\ \begin{bmatrix}1&0\\0&-1\end{bmatrix}\]
in $\PGL_2\Z$.  Write $J_\epsilon$ to denote an element of the set $\{ J_\pm, J_0 \}$.  Acting as a M\"obius function, $J_\epsilon$ preserves the edge $[0,\infty]$; $J_0$ fixes its endpoints and the $J_\pm$ interchange them.  Since $A(\infty)=e^R$ and $A(0)=e^L$, the M\"obius function $AJ_\epsilon A^{-1}$ preserves $e$; $AJ_0 A^{-1}$ fixes its endpoints and the $AJ_\pm A^{-1}$ interchange them.  Direct calculation shows that the quotient of the first column of $AJ_\epsilon A^{-1}$ is $\alpha_\epsilon$.

Fix $\epsilon$ for which $\alpha_\epsilon\in(0,1/2)$ and write $\alpha_\epsilon = p/q$ in lowest terms.  Since $\tr AJ_\epsilon A^{-1}=0$, 
\[ AJ_\epsilon A^{-1} \ = \ \begin{bmatrix} p&u\\ q&-p \end{bmatrix}\]
for some $u \in \Z$.  Because $\det AJ_- A^{-1}=1$, $p^2 \equiv -1 \pmod q$ in this case.  By Lemma~II.3.3 of \cite{SW}, $\alpha_-$ is a palindrome
and the length of its positive termed continued fraction is even.  Also, the comment preceding Remark ~\ref{rem: SWautomorphisms}, shows that $\phi_{\alpha_-}$ is orientation reversing.  Now, assume $\alpha_\epsilon \in \{ \alpha_+, \alpha_0 \}$.  Then $\det AJ_\epsilon A^{-1}=-1$ and $p^2 \equiv 1 \pmod q$.  If $q$ is odd, this is Case (2) of \cite[Lemma~II.3.3]{SW}.  If $q$ is even, then
\[ p^2\ \equiv \ 1\pmod{2q} \qquad\text{or}\qquad p^2\ \equiv \ q+1\pmod{2q}\]
which gives Case (3) or Case (2), respectively.  Thus $\alpha_\epsilon$ is again a palindrome, and in either case its
positive termed continued fraction has odd length.  Again, using the comment just before Remark ~\ref{rem: SWautomorphisms}, $\phi_{\alpha_\epsilon}$ must be orientation preserving.  (This proves (4).)

In fact, 
\[ A_{\alpha_\epsilon} \ =\ -AJ_\epsilon A^{-1} \]
in all cases.   To see this, take relatively prime positive integers $r$ and $s$ so that $r/s$ is the high vertex $e_{n+1}^-$ of the bottom edge of the funnel $F_{\alpha_\epsilon}$.  Since $\alpha_\epsilon$ is a palindrome, Theorem~\ref{thm: reverse} gives $s/q=\bar\alpha_\epsilon=p/q$ and so $s=p$. (The reciprocal integer cases excluded by the standing assumptions are addressed in Remark~\ref{rem: reciprocals} below.) The above parity and determinant conditions force $r=-u$ to give $A_{\alpha_\epsilon} = -AJ_\epsilon A^{-1}$ as claimed.

It remains to show that $e$ is an edge of $F_{\alpha_\epsilon}$. By Lemma~II.3.3 and Figure~II.3.4 of \cite{SW}, the Farey triangle sequence associated to a palindromic continued fraction has a two-fold symmetry which preserves an edge of that sequence. Consequently, $\Phi_{\alpha_\epsilon}$ preserves some edge of $F_{\alpha_\epsilon}$.  Because $[0,\infty]$ is the unique edge preserved by $J_\epsilon$, $e$ is the unique edge preserved by $\Phi_{\alpha_\epsilon}$.  Therefore, the invariant edge of $F_{\alpha_\epsilon}$ supplied by \cite{SW} must be $e$.

Conversely, suppose that $\alpha\in\Q\cap(0,1/2)$ is a palindrome, $e$ is an edge of $F_\alpha$, and $\Phi_\alpha(e)=e$.  Then $A^{-1}A_\alpha A$ is a nontrivial element of $\PGL_2\Z$ which preserves $[0,\infty]$.  The set of nontrivial elements of $\PGL_2 \Z$ which preserves $[0,\infty]$ setwise is $\{ J_\pm, J_0\}$ so there is a $J_\epsilon$ so that $A_\alpha = AJ_\epsilon A^{-1}$.  The quotient of the first column of $A_\alpha$ is $\alpha$, so $\alpha = \alpha_\epsilon$.  This proves (2), and the actions of the matrices $J_\epsilon$ on $[0,\infty]$ prove (3).
\endproof

\begin{remark} \label{rem: reciprocals}
Some of the numbers in Lemma~\ref{lem: fixed edge} are reciprocal integers and therefore do not satisfy the standing assumptions for this section.  These cases can be described explicitly.

Suppose first that $\alpha_\epsilon=1/q$.  In the notation of the proof,
\[ AJ_\epsilon A^{-1} \ = \ \begin{bmatrix} 1&u\\ q&-1 \end{bmatrix}. \]
The case $J_\epsilon=J_-$ cannot occur, since $\det(AJ_-A^{-1})=1$ would imply that $q\mid2$, contrary to $\alpha_-\in(0,1/2)$.  Otherwise, the determinant is $-1$, so $u=0$. In the $J_+$ case, $a=c=1$ and $d=b+1$.  Hence
\[e \ =\ \left[\frac1{b+1},\frac1b\right] \AND \alpha_+\ =\ \frac1{2b+1}. \]
In the $J_0$ case, $c=0$ and $a=d=1$. Therefore
\[ e\ =\ \left[0,\frac1b\right] \AND \alpha_0\ =\ \frac1{2b}. \]
The conclusions of Lemma ~\ref{lem: fixed edge} are easily verified in these special cases.
\end{remark}

The $\alpha_+$ case of Lemma~\ref{lem: fixed edge} identifies several distinguished points at which the shape function $\calZ(e)$ takes real values, including a geometric point on the real shape locus $\R\calZ(e)$. (Recall Definition~\ref{dfn: real locus}.)  In what follows, for an edge $e=[c/d, a/b]$, write
\[ \alpha_e\ =\ \alpha_+\ =\ \frac{cd-ab}{d^2-b^2}.\] \index{shape parameter!real shape locus}

\begin{thm} \label{thm: real loci}
\RealThm
\end{thm}

\proof
Recall from Definition~\ref{def: shape parameter function} that
\[ \calZ(e) \ =\  -d_{\calQ}(\hat e)^{-\nu(e)} \left( \frac{\calQ(e^L)}{\calQ(e^R)} \right)^2. \]
Since $e^L$ and $e^R$ are Farey partners, Lemma~\ref{lem: adjacent zeros} implies that $\calQ(e^L)$ and $\calQ(e^R)$ have no common zero.
Consequently,
\[ \calZ(e)_{z_L}\ =\ 0 \AND \calZ(e)_{z_R} \ =\ \infty. \]

Recalling that $\zeta(w) = \frac{1}{1-w}$, the formulas in Lemma~\ref{lem: zetas} give
\[ 1-\calZ(e) \ = \  \begin{cases} \displaystyle \frac{ \calQ(e^-\oplus e^+)\calQ(\hat e) }{ d_{\calQ}(\hat e)\calQ(e^-)^2 } &
\text{if $\nu(e)=1$},
\\[4mm]
\displaystyle \frac{ \calQ(e^-\oplus e^+)\calQ(\hat e) }{ \calQ(e^+)^2 } & \text{if $\nu(e)=-1$}.
\end{cases} \]
 Lemma~\ref{lem: adjacent zeros}, together with Remark~\ref{rem: constant term} and Example~\ref{ex: dQ}, ensures that the relevant denominators do not vanish at roots of $\calQ(\hat e)$ or $\calQ(e^-\oplus e^+)$.  Hence,
\[ \calZ(e)_{z_{\hat e}}\ =\ \calZ(e)_{z_\oplus}\ = \ 1.\]

It remains to consider the geometric root for $\calQ(\alpha_e)$.  By Lemma~\ref{lem: fixed edge}, $\alpha_e$ is a
palindrome between $0$ and $1/2$, $e$ is an edge of $F_{\alpha_e}$, $\Phi_{\alpha_e}$ preserves $e$ while interchanging its endpoints, and $\phi_{\alpha_e}\co M_{\alpha_e}\longrightarrow M_{\alpha_e}$ is an orientation preserving isometry.  Equivalently, this is the endpoint interchanging, orientation preserving palindromic symmetry appearing in the second diagram of \cite[Figure~II.3.4]{SW}.

Since the isometry preserves the tetrahedral pair corresponding to $e$ and interchanges the labels
$e^L$ and $e^R$, the corresponding geometric holonomies have equal trace
\[ \tr\left(S_{e^L}\right)_{\frz(\alpha_e)} \ =\  \tr\left(S_{e^R}\right)_{\frz(\alpha_e)}.\]
Theorem~\ref{thm: traces and shapes} gives
\[\calZ(e) \ = \ -\frac{2-\tr S_{e^L}}{2-\tr S_{e^R}}.\]
Evaluating at $\frz(\alpha_e)$ therefore gives
\[ \calZ(e)_{\frz(\alpha_e)}\ =\ -1 \]
as claimed.
\endproof

Refer back to Figure ~\ref{fig: RZ} to see some examples.

\begin{figure}[h] 
   \centering
   \includegraphics[width=4.2in]{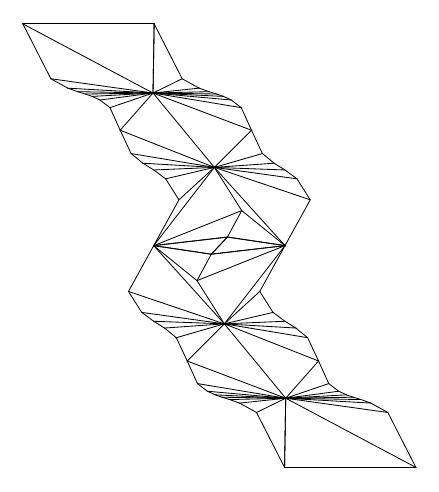} 
\caption{\index{fundamental domain!specialized} If $e=[133/281, 195/412]$, then $\alpha_e= 42967/90783$.  The knot $L_{\alpha_e}$ has $D(2,8,1,6,5,6,1,8,2)$ as an alternating 4-plat diagram, so $L_{\alpha_e}$ is a 39 crossing knot.  The polynomial $\calQ(\alpha_e)$ has degree $45391$ and the degree of $p_e$ is 411.  The geometric root $\frz$ is approximately  $0.008782828575501701 + 0.5267009878219225 \, i$.  The fundamental domain $\Omega(\alpha_e)_\frz$ consists of $72$ tetrahedra with a total volume of approximately $33.2723826156$.}     \label{fig: 42967 90783}
\end{figure}

\begin{remarks} \label{rem: real loci}
\[\]
\vspace{-.4in}
\begin{enumerate}
\item The $\alpha_-$ case of Lemma~\ref{lem: fixed edge} gives an orientation reversing analogue of the last statement of Theorem~\ref{thm: real loci}.  Suppose that $\alpha_-$ is hyperbolic, so that the geometric root $\frz(\alpha_-)$ exists.  The isometry $\phi_{\alpha_-}\co M_{\alpha_-}\longrightarrow M_{\alpha_-}$ preserves the tetrahedral pair corresponding to $e$ and interchanges the labels $e^L$ and $e^R$.  Since this isometry is orientation reversing, it follows that
\[ \tr\left(S_{e^L}\right)_{\frz(\alpha_-)} \AND \tr\left(S_{e^R}\right)_{\frz(\alpha_-)} \]
are complex conjugates.  This is the first of the three palindromic configurations in \cite[Figure~II.3.4]{SW}.

\item Suppose that $e$ and $\alpha_e$ are as in the last statement of Theorem~\ref{thm: real loci}.  Consider the polynomial
\[ p_e \ =\  d_{\calQ}(\hat e)\calQ(e^-)^2-\calQ(e^+)^2. \]
By Equation~(\ref{eq: trS one}), the polynomial $p_e$ vanishes precisely when
\[ \tr S_{e^L}\ =\ \tr S_{e^R}. \]
The proof of Theorem~\ref{thm: real loci} shows that
\[ p_e\bigl(\frz(\alpha_e)\bigr)\ =\ 0. \]
Consequently, the irreducible factor of $\calQ(\alpha_e)$ having $\frz(\alpha_e)$ as a root divides $p_e$.

The degree reduction can be dramatic.  For example, if $e=\left[133/281,\,195/412\right]$, then $\alpha_e=42967/90783$.  In this case, $p_e$ has degree $411$, while $\calQ(\alpha_e)$ has degree $45391$.  Thus the palindromic symmetry produces an exact polynomial equation of much lower degree satisfied by the geometric root of $\calQ(\alpha_e)$.  The corresponding geometric fundamental domain is shown in Figure~\ref{fig: 42967 90783}.

Similar algebraic properties have been previously observed for palindromic $2$-bridge links.  For example, compare Riley's discussion of reducibility of the corresponding character variety in \cite{RRalg}.
\end{enumerate}
\end{remarks}

\begin{example}[Golden links] \label{ex: golden}

\begin{figure}[h]
   \centering
   \includegraphics[width=3.5in]{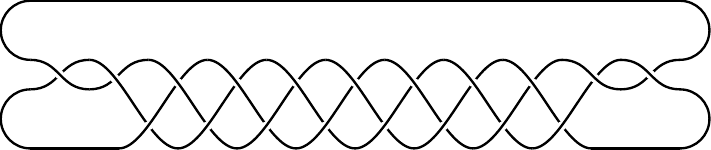} 
   \caption{The diagram $D_{\gamma(15)}$ shows the $19$-crossing golden link for $\gamma(15) = [0;2,1,\ldots, 1,2] = 2584/6765$.}  \label{fig: goldenlink}
\end{figure}

Let $F_n$ be the $n^\text{th}$ Fibonacci number where $F_0=0$ and $F_1=1$.  The successive quotients of the Fibonacci numbers have nice continued fraction expansions
\[ \frac{F_{n+1}}{F_n} \ =\ \left[ 1; 1, \ldots, 1\right] \]
where the expression on the right contains $n$ ones.  Also, if $\varphi = (1/2) \left( 1+\sqrt{5}\right)$ is the golden ratio, then
\[ \lim_{n \to \infty} \left(\frac{F_{n+1}}{F_{n}} \right)\ = \ \varphi.\]
Now define
\[ \gamma(n) \ =\ \frac{F_{n+3}}{F_{n+5}}.\]
Then 
\[ \{ \gamma(n) \}_0^8 \ =\ \left\{ \frac{2}{5}, \frac{3}{8}, \frac{5}{13}, \frac{8}{21}, \frac{13}{34}, \frac{21}{55}, \frac{34}{89}, \frac{55}{144}, \frac{89}{233} \right\} \]
and 
\[\gamma(n) \ = \ 2 - \frac{F_{n+6}}{F_{n+5}}.\]  
It follows that the 2-bridge link $L_{\gamma(n)}$ is equivalent to the link for $F_{n+6}/F_{n+5}$.  For this reason, the link $L_{\gamma(n)}$ is called the $n^\text{th}$ {\it golden link}.  (The $15^\text{th}$ golden link is shown in Figure ~\ref{fig: goldenlink}.)  The number $\gamma(n)$ has a nice continued fraction expansion
\[ \gamma(n) \ =\ [0;2, 1, \ldots, 1, 2]\] 
where there are $n$ ones sandwiched between the two twos.  

\begin{figure}
   \centering
   \includegraphics[width=4.85in]{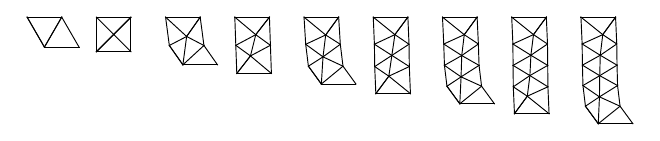} 
   \caption{\index{fundamental domain!specialized} Fundamental domains $\Omega(\gamma(n))_{\frz_1}$ where $n \leq 8$. The first two golden links are the figure-eight knot and the Whitehead link.}  \label{fig: smallgoldens}
\end{figure}

The expected degrees of the polynomials $\calQ(\gamma(n))$ grow rapidly with $n$.  Indeed, since the denominator of $\gamma(n)$ is $F_{n+5}$, the expected degree of $\calQ(\gamma(n))$ is
\[ \left\lfloor \frac{F_{n+5}-1}{2}\right\rfloor \]
(see Remark~\ref{rem: conj QandL}).  Thus, exact symbolic calculations with $\calQ(\gamma(n))$ can quickly involve polynomials of very large degree.  Figure~\ref{fig: smallgoldens} shows the geometric fundamental domains for the first several members of this family.

\begin{figure}
   \centering
   \includegraphics[width=4.85in]{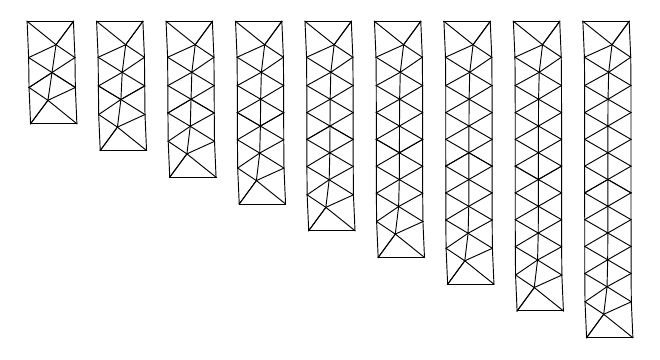} 
   \caption{\index{fundamental domain!specialized} Fundamental domains $\Omega(\alpha)_{\frz_1}$ where $\alpha$ is in $\{ 21/55,
 55/144,
 144/377,
 377/987,
 987/2584,
 2584/6765,
 6765/17711,
 17711/46368,
 46368/121393 \}$. }  \label{fig: oddgoldens}
\end{figure}

For the odd-indexed members of this family, Theorem~\ref{thm: real loci} and part~(2) of Remark~\ref{rem: real loci} give a striking algebraic degree reduction.  Let $e(n)$ be the edge $[\gamma(n), \gamma(n+2)]$ in the Stern--Brocot diagram $\calG$.  A direct calculation from the definition gives
\[\alpha_{e(n)}\ =\ \gamma(2n+7).\]
Let $p_n=p_{e(n)}$ be the polynomial from Remark~\ref{rem: real loci} (2).

Since Theorem~\ref{thm: real loci} gives
\[ \calZ(e(n))_{\frz(\alpha_{e(n)})}\ =\ -1\]
the proof of part (2) of Remark~\ref{rem: real loci} shows that $\frz(\alpha_{e(n)})$ is a root of $p_n$.  Thus the geometric root of $\calQ(\gamma(2n+7))$ satisfies the substantially lower-degree polynomial $p_n$.  For instance, the expected degree of $\calQ(\gamma(21))$ is $60696$, while the degree of $p_{7}$ is only $376$. \index{expected degree}

Figure~\ref{fig: oddgoldens} shows geometric fundamental domains for several of these odd-indexed golden links.  No analogous algebraic degree reduction is apparent for the even-indexed golden links, since the polynomials  $\calQ(\gamma(2n))$ often appear to be irreducible over $\Q$.

Nevertheless, the numerical Newton method described in Section~\ref{subsec: code II} makes it possible to compute the geometries of much larger golden links.  Figure~\ref{fig: longgoldens} shows examples of both parities.

\begin{figure}[h]
   \centering
   \includegraphics[width=5in]{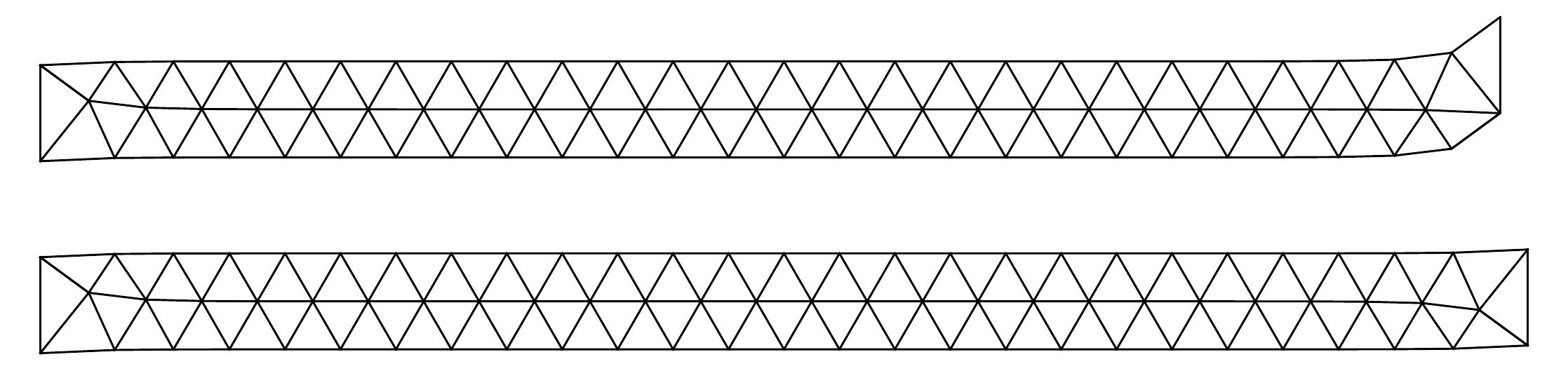}
   \caption{\index{fundamental domain!specialized} Geometric fundamental domains for the golden links corresponding to $\gamma(50)$ and $\gamma(51)$.}
   \label{fig: longgoldens}
\end{figure}

\begin{conj} \label{conj: golden} \index{geometric root!2-bridge link}
For $n \in \Z_{\geq 0}$, let $\frg_n$ be the geometric root for $\calQ(\gamma(n))$.  The roots $\frg_n$ converge to a point $\hat{\frg}$ for which $\Gamma_{\hat{\frg}}$ is a geometrically infinite Kleinian group.   Moreover, if $e_j$ is the $j^\text{th}$ edge of the infinite funnel for $2-\varphi = (1/2)(3-\sqrt{5})$, then
\[ \lim_{j \to \infty} \calZ(e_j)_{\hat{\frg}} \ =\ e^{i\pi/3}\]
and $\hat{\frg}$ is the unique point in the intersection
\[ \bigcap_{j \geq 1} \left( \calU_{e_j} \cap (\C^\ast - \calR) \right)\]
where $\calR$ is the Riley slice defined in Definition ~\ref{def: Riley slice} and $\calU_{e_j}$ is the thick region for $\calZ(e_j)$ from Definition ~\ref{def: fat}. \index{Riley slice} \index{thick region}
\end{conj}

More generally, suppose that $\{ \alpha_j \}_1^\infty$ is a sequence of distinct hyperbolic rational numbers that converge to an irrational real number $\lambda \in (0,1/2)$.  It is expected that the geometric roots $\frz(\alpha_j)$ will converge to a complex number $\hat{\frz}$ which is the unique element in the intersection 
\[ \bigcap_{j \geq 1} \left( \calU_{e_j} \cap (\C^\ast - \calR) \right)\]
where $e_j$ is the $j^\text{th}$ interior edge of the infinite funnel for $\lambda$.

\end{example}

\section{Cusp shapes} \label{sec: cusps}

\begin{dfns} \label{def: cusp field}
Suppose that $P \subset \C$ is a parallelogram with $0$ as a vertex and let $\lambda, \mu \in \C$ be the vertices of $P$ adjacent to $0$.  Let $T$ be the torus obtained by identifying the opposite edges of $P$.  The {\it shape of $T$ relative to $(\lambda, \mu)$ } is the quotient $\lambda/\mu$.  The {\it shape } of $T$ is the orbit of $\lambda/\mu$ under the M\"obius action of $\text{PGL}_2(\Q)$ together with complex conjugation.  It is well known that the modulus of $T$ is a commensurability invariant for $T$ (see, for example, Lemma 4.9 of \cite{CD}).  The {\it cusp field } of a cusp in a hyperbolic 3-manifold is the number field generated by any representative of its shape.    
\end{dfns}
\index{cusp shape}
\index{cusp field}

\begin{thm} \label{thm: cuspfield} \index{geometric root!2-bridge link}
Suppose $\alpha$ is hyperbolic.  Let $\frz$ be the geometric root of $\calQ(\alpha)$.  The cusp field of any cusp of $M_\alpha$ is $\Q(\xi)$, where $\xi$ is either element of 
\[ \left\{ \calV\left( e_{n+1}^- \right)_{\frz}, \,  \calV\left( e_{n+1}^+ \right)_{\frz} \right\}.\]
\end{thm}

\proof  If $M_\alpha$ has more than one cusp, an involution in the Klein 4-group $G$ interchanges them.  Hence, the cusp fields for all cusps are identical.

From the proof of Theorem ~\ref{thm: complete}, $M_\alpha$ has a cusp whose Euclidean holonomy is generated by the parabolics
\[ z \mapsto z+1 \qquad \text{and} \qquad z \mapsto z+2 \, \calV\left( e_{n+1}^- \right)_\frz. \]
The shape for a cross-section relative to this pair is $2 \, \calV\left( e_{n+1}^- \right)_\frz$ and so the cusp field is $\Q(\xi)$, where $\xi=\calV\left( e_{n+1}^- \right)_\frz$.

Lemma ~\ref{lem: difference} shows that alternatively $\xi$ can be taken to be $\calV\left( e_{n+1}^+ \right)_\frz$.
\endproof

Valuable invariants of Kleinian groups include their trace fields.  It is well known that the trace field of the 2-parabolic generator group $\Gamma_\frz$ (Definition ~\ref{def: generic holonomy}) is the extension of $\Q$ obtained by  adjoining the traces of $(W_0)_\frz$ and $(U_0 W_0)_\frz$.  It follows that the trace field for $\Gamma_\frz$ is $\Q(\frz)$.  \index{trace field}

A related invariant is the shape field, which is obtained from a geometric ideal triangulation of a hyperbolic manifold $M$ by adjoining to $\Q$ its shape parameters.  If $\alpha$ is hyperbolic and $\frz=\frz(\alpha)$ is the geometric root, then Corollary~\ref{cor: main1} shows that the shape field of the Sakuma--Weeks triangulation of $M_\alpha$ is $\Q(\frz)$.  Thus, in this case, the shape field coincides with the trace field of $\Gamma_\frz$.  More generally, this holds for all hyperbolic links in $S^3$, as shown in Corollary 2.3 and Theorem 2.4 of \cite{NR}. \index{shape field}

There are many known examples of hyperbolic 3-manifolds whose cusp fields are properly contained in their trace fields. However, such examples are rare among hyperbolic knot complements in $S^3$. Only a small number of such examples are currently known; these include the two dodecahedral knots of Aitchison and Rubinstein \cite{AR}, the 12-crossing Boyd knot $12n706$ \cite{GHH}, the knot 
$15n132539$ identified by Nathan Dunfield, and four prism knots identified by Deblois, Gharagozlou, and Hoffman in \cite{hidden}, which exhibit hidden symmetries and contain closed geodesic surfaces.  

In the context of knot complements, Question 2 of \cite{NR} asks when the cusp field equals the trace field.   In spite of the fact that this question has held interest in the field for over 30 years and that there are very few known knot complements with properly contained cusp fields, there are also only a few known infinite families of knots where equality is known.  In \cite{NR}, Neumann and Reid prove equality for the twist knots, $\alpha = \frac{m}{2m+1}$.  Inspired by their argument, the same is proven in \cite{CDHMMW} for the $(1,m)$-Dehn fillings of one cusp of the link $L_{3/10}$. (These fillings are knots but are not, in general, 2-bridge knots.)  Finally, Hoste and Shanahan show in \cite{HS} that the cusp fields and trace fields for the knot complements $M_\alpha$, with $\alpha = \frac{2m}{1+6m}$, are also equal.  

\begin{figure}
   \centering
   \includegraphics[width=4.85in]{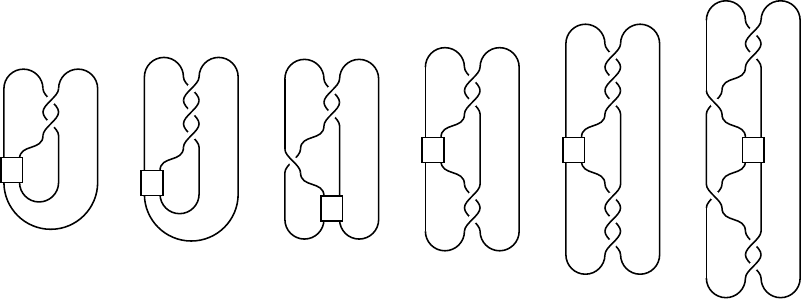} 
   \caption{The six families of knots and links covered by Corollary ~\ref{cor: field equality}.}  \label{fig: CuspLinks}
\end{figure}

\begin{cor}\label{cor: field equality}
\CuspCor
\end{cor}

The statement of the corollary includes the hyperbolic twist knots ($\alpha = \frac{2}{2m+1}$ and $\alpha=\frac{m}{2m+1}$) which, as mentioned above, were already known from \cite{NR}.  In the proof below, a different argument, using Theorem ~\ref{thm: cuspfield}, is given.  The statement also includes the links with $\alpha =\frac{3}{3m+1}$ and $\alpha=\frac{m}{3m+1}$.  These links may have either one or two components, those with one component are covered by the Hoste/Shanahan result in \cite{HS}.  Again, the argument given here is different and uses Theorem ~\ref{thm: cuspfield}. 

\proof[Proof of Corollary ~\ref{cor: field equality}]
First, observe that
\[ \calV(1/2) \ =\ -1+x \qquad \text{and} \qquad \calV(1/3) \ = \ \frac{1-2x}{1-x}. \]
Hence, $x$ can be expressed as a rational function in either $\calV(1/2)$ or $\calV(1/3)$.  This means that, if the cusp field of a hyperbolic link complement $M_\alpha$ is $\Q(\xi)$, where $\xi \in \left\{ \calV(1/2)_{\frz}, \, \calV(1/3)_{\frz} \right\}$, then the cusp field is $\Q(\frz)$.  This observation will be used in all cases for this theorem.

If $\alpha=(m+1)/(2m+3)$ and $e_{n+1}$ is the bottom edge of $F_\alpha$, then $e_{n+1}^-=1/2$.  So, by Theorem ~\ref{thm: cuspfield} and the observation above, the cusp and trace fields for $M_\alpha$ are the same.  Because the link complement $M_{\bar{\alpha}}$ for the reverse (Definition ~\ref{def: reverse}) of $\alpha$ is isometric to $M_\alpha$ by $\phi_\alpha$, the cusp and trace fields for $M_{\bar{\alpha}}$ are equal.  Note that, if $\alpha = (m+1)/(2m+3)$, then its reverse is $\bar{\alpha} = 2/(2m+3)$.

Suppose $\alpha$ is of the form $m/(3m+1)$ or $(m+1)/(3m+2)$.  If $\alpha = 2/5$, then the result is already known by the previous argument.  Otherwise, $e_{n+1}^- = 1/3$ and, as before, Theorem ~\ref{thm: cuspfield} and the observation provide the result.  Since the reverse of $m/(3m+1)$ is $3/(3m+1)$ and the reverse of $(m+1)/(3m+2)$ is $3/(3m+2)$, these cases are also covered.
 
 In the remaining cases, $\alpha$ is a palindrome and symmetries of the cusp triangulations are used.  As in the proof of Theorem ~\ref{thm: complete}, take $\Upsilon(\alpha)_\frz$ to be the triangulated intersection of the domain $\Omega(\alpha)_\frz$ with a horosphere centered at $\infty$.
 
\begin{figure}
   \centering
   \includegraphics[width=3.5in]{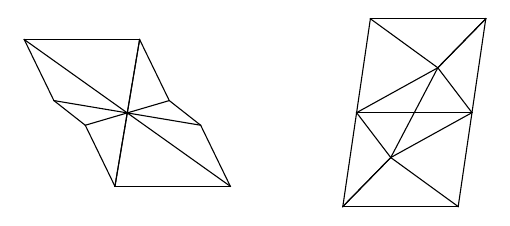} 
   \caption{\index{fundamental domain!specialized} On the left is the cusp triangulation for $\alpha = 7/16$ and on the right is the triangulation for $\alpha =11/30$.  In both cases, the geometric root is used to draw the figure.}  \label{fig: cusps}
\end{figure}

The continued fraction expansion for $\alpha=\frac{2m+1}{4m+4}$ is the palindrome $[0; 2,m,2]$.   The last hub $e_{n+1}^-$ for $F_\alpha$ is $[0;2,m]=m/(2m+1)$.  The formula for $A_\alpha$, provided  just after  Definition ~\ref{def: reverse}, becomes 
\[ A_\alpha \ =\ \begin{bmatrix} -2m-1 & m \\ -4m-4 & 2m+1 \end{bmatrix}.\]
By Theorem~\ref{thm: reverse}, this matrix determines a Mo\"bius function $\Phi_\alpha$ which acts as a simplicial isomorphism on the funnel $F_\alpha$.  This isomorphism fixes the vertex $1/2$ of $F_\alpha$.  As discussed in Section ~\ref{sec: reverse}, $\Phi_\alpha$ induces an orientation preserving isometry $\phi_\alpha: M_\alpha \to M_\alpha$.  This, in turn, lifts to an order-2 hyperbolic isometry which restricts to a self-isometry of $\Upsilon(\alpha)_\frz$ fixing $\calV(1/2)_\frz$.  This implies that $\calV(e_{n+1}^-)_{\frz}=2 \cdot \calV(1/2)_{\frz}$.  Theorem ~\ref{thm: cuspfield} and the observation finish the case.  For reference, the cusp triangulation in the case $m=3$ is shown on the left in Figure ~\ref{fig: cusps}.

The case $\alpha=\frac{3m+1}{9m+6}$ is similar.  Here, $\alpha$ is the palindrome $[0;3,m,3]$ and $\Phi_\alpha$ acts on $F_\alpha$ fixing $1/3$.  It follows that $\calV(e_{n+1}^-)_{\frz}=2 \cdot \calV(1/3)_{\frz}$, which finishes the case as before.  If $m=3$ then $\alpha=10/33$, which is illustrated in Figure ~\ref{fig: domain}.

Finally, the number $\alpha = \frac{3m+2}{9m+3}$ expands as the palindrome $[0; 2, 1, m-1, 1, 2]$.  In this case, $\Phi_\alpha$ acts on $F_\alpha$ fixing $1/3$ and the result follows as before.  For a more visual argument, let $A=\calV(1/3)_{\frz}$ and $\xi=\calV(e_{n+1}^-)_{\frz}$ and keep an eye on the righthand image in Figure ~\ref{fig: cusps}.  The right edges of $\Upsilon(\alpha)_{\frz_1}$ have vertices $1$, $A$, and $\xi$ while the vertices on the left side are translates of these by $U_0$.  Hence, the order-2 rotational isometry of $\Upsilon(\alpha)_{\frz_1}$ must fix the point $-1/2+A=\xi/2$. Therefore, the cusp field is $\Q(A)$ and the result follows as usual.  
\endproof

\section{Riley polynomials} \label{sec: riley}

In \cite{RR}, Riley studies irreducible representations of 2-bridge knot groups into $\pslC$ under which peripheral elements map to parabolics.  Such representations are called {\it p-reps}.  Given $\alpha$, let $k_0$ and $k_1$ be the associated canonical generators for $\pie M_\alpha$.  Every p-rep for $L_\alpha$ may be conjugated so that \index{canonical generators}
\[
k_0 \mapsto
\begin{bmatrix}
1&1\\
0&1
\end{bmatrix}
\qquad \text{and} \qquad
k_1 \mapsto
\begin{bmatrix}
1&0\\
\mathfrak{Z}&1
\end{bmatrix}
\]
for some $\mathfrak{Z} \in \C$.  In Theorem 2 of \cite{RR}, Riley proved that, for knots, this can be done precisely when $\mathfrak{Z}$ satisfies the Riley polynomial $\Lambda(\alpha) \in \Z[Z]$.
\index{Riley polynomial}

In \cite{C}, the definition of Riley polynomials is extended from knots to all 2-bridge links using Farey recursion.  Following the work there, let
\index{Farey recursion}
\[\calP \co \wh{\Q}_0 \longrightarrow \Z[i,z]\]
be the FRF with constant determinant $1$, trace $\calP$, and initial conditions
\[
\calP(0)=iz,
\qquad
\calP(\infty)=0,
\qquad
\calP(1)=z,
\]
and define $f \co \wh{\Q}_0 \to \Z[i,z]$ by
\[ f(p/q) \ =\ \begin{cases}
iz & \text{if $p$ is even and $q$ is odd},\\
z & \text{if $p$ and $q$ are both odd},\\
iz^2 & \text{otherwise}.
\end{cases} \]
It is shown in \cite{C} that, if $Z=z^2$, then the image of $\frac{\calP}{f}$ lies in $\Z[Z]$ and agrees with $\Lambda$ on rational numbers with odd denominator.  So, set
\[ \Lambda\ =\ \frac{\calP}{f}. \]
A nonzero root $Z_0$ of $\Lambda(\alpha)$ corresponds to an irreducible p-rep for the 2-bridge link $L_\alpha$ in the same way as above.  A root $Z_0=0$ does not correspond to an irreducible p-rep; its meaning will become clear in Theorem~\ref{thm: QandL}.

\begin{dfns} \label{def: deg}
For $\alpha=p/q\in\Q_0$, written in lowest terms with $q>0$, define 
\[ \deg(\alpha) \ =\ \deg \calQ(\alpha) \AND \expect(\alpha) \ =\  \left\lfloor\frac{q-1}{2}\right\rfloor. \]
An induction argument shows that $\deg(\alpha) \leq\expect(\alpha)$.  Since equality holds in many cases, we call $\expect(\alpha)$ the
{\it expected degree} of $\calQ(\alpha)$.  It is shown in \cite{C} that the degree of $\Lambda(\alpha)$ is equal to $\expect(\alpha)$. Define the {\it defect} of $\alpha$ by
\[ \mu(\alpha) \ =\  \expect(\alpha)-\deg(\alpha).\]
Thus $\mu(\alpha)\geq0$, with equality precisely when $\calQ(\alpha)$
has its expected degree.  Finally, define the reversed polynomial
\[ \overline{\calQ}(\alpha) \ =\  x^{\deg(\alpha)}\calQ(\alpha)_{-x^{-1}}.\]
\end{dfns}
\index{expected degree|idxdef} \index{defect|idxdef}

The next theorem gives an explicit relationship between the Riley polynomial $\Lambda$ and the Farey recursive polynomial $\calQ$.  Remark 5.3.13 of \cite{ASWY} observes that, after removing a parity dependent monomial factor, the Markoff polynomials considered there become hyperbolicity polynomials for the $2$-bridge links.  The theorem below gives a precise version of this relationship in terms of $\calQ$ and the Riley polynomial $\Lambda$.
\index{Farey recursive function}

\begin{thm}\label{thm: QandL}
Suppose $\alpha \in\Q_0$.  Then
\[ \Lambda(\alpha)_x \ =\  \pm x^{\mu(\alpha)} \, \overline{\calQ}(\alpha).\]
\end{thm}

\proof
Write $\alpha=p/q$ in lowest terms, with $q>0$, and introduce an auxiliary variable $z$.  Let
\[ c(\alpha)\ =\ (-1)^p i^{p+q}z^q \AND \calP^\ast(\alpha)\ =\  c(\alpha)\, \calQ(\alpha)_{-z^{-2}}.\]
Set $\calP^\ast(\infty)=0$. The first step is to show that $\calP^\ast=\calP$.

Begin by taking $\alpha=p/q$ and $\gamma=r/s$ so that they form a Farey pair. Set
\[\beta\ =\ \alpha\oplus\gamma \AND \delta=\gamma\oplus^2\alpha.\]
The Farey recursion for $\calQ$ gives
\[ \calQ(\delta) \ =\  -\dQ(\alpha)\, \calQ(\gamma) + \calQ(\alpha)\, \calQ(\beta).\]
After specializing at $x=-z^{-2}$,
\[ \dQ(\alpha)_{-z^{-2}} \ =\  (-1)^{p+q}z^{-2q}\]
and a direct calculation gives
\[ c(\delta)\, \dQ(\alpha)_{-z^{-2}} \ =\  c(\gamma) \AND  c(\delta)=c(\alpha)\, c(\beta).\]
Hence,
\[ \calP^\ast(\delta) \ =\  -\calP^\ast(\gamma) + \calP^\ast(\alpha)\calP^\ast(\beta)\]
and $\calP^\ast$ satisfies the same Farey recursion as $\calP$.  Because
\[ \calP^\ast(0)=iz, \qquad \calP^\ast(\infty)=0, \qquad \calP^\ast(1)=z \]
$\calP^\ast=\calP$ as desired.  Also
\[ \calP(\alpha) \ =\  c(\alpha) \,  \calQ(\alpha)_{-z^{-2}}.\]

Recall that $\Lambda=\calP/f$.  A check of the three parity cases in the definition of $f$ gives
\[
\frac{c(\alpha)}{f(\alpha)}
=
\pm z^{2\expect(\alpha)}.
\]
Consequently,
\[ \Lambda(\alpha)_{z^2} \ =\  \pm z^{2\expect(\alpha)} \calQ(\alpha)_{-z^{-2}}. \]
Replacing $z^2$ by $x$ gives
\[ \Lambda(\alpha)_x \ =\  \pm x^{\expect(\alpha)} \, \calQ(\alpha)_{-x^{-1}}. \]
Finally, Definition~\ref{def: deg} gives
\[ \expect(\alpha) \ =\  \mu(\alpha)+\deg(\alpha) \AND \overline{\calQ}(\alpha) \ =\  x^{\deg(\alpha)} \, \calQ(\alpha)_{-x^{-1}}.\]
Therefore,
\[ \Lambda(\alpha)_x \ =\  \pm x^{\mu(\alpha)}\overline{\calQ}(\alpha).\]
\endproof

Theorem~\ref{thm: QandL} makes the relationship between the roots of $\calQ$ and $\Lambda$ explicit.  Since $\overline{\calQ}(\alpha)$ has nonzero constant term, $\mu(\alpha)$ is the order of the zero at zero for $\Lambda(\alpha)$.  The remaining roots of $\Lambda(\alpha)$ are precisely $-\frz^{-1}$, where $\frz$ ranges over the roots of $\calQ(\alpha)$, with multiplicity preserved.  In particular, if $\frz(\alpha)$ is the geometric root of $\calQ(\alpha)$, then Corollary~\ref{cor: holonomy specialization} shows that $-\frz(\alpha)^{-1}$ is a root of $\Lambda(\alpha)$ corresponding to a discrete faithful p-rep.

\begin{figure}
   \centering
   \includegraphics[width=5in]{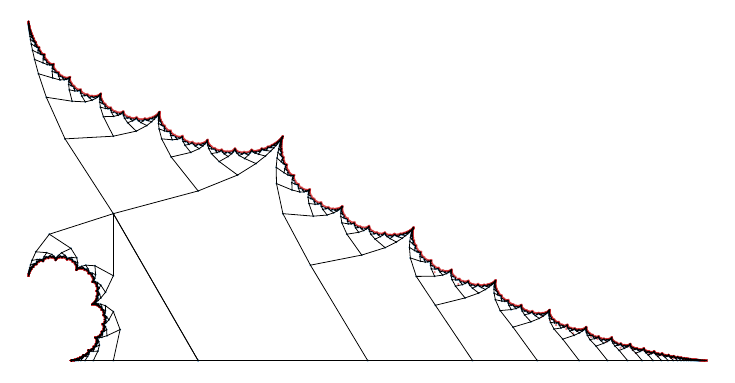}
   \caption{Under the transformation $\frz\mapsto-1/\frz$, the tree of geometric roots for $\calQ$ becomes the familiar picture of Riley polynomial roots.  See, for instance, Figure 3 of \cite{Ober}.}
   \label{fig: inversion}
\end{figure}

Thus, after inversion through the unit sphere, the tree of geometric roots for $\calQ$ shown in Figure~\ref{fig: roots} becomes the familiar picture of geometric Riley roots; compare Figure~1 of \cite{KS} and Figures~0.2a and~0.2b of \cite{ASWY}.  This relationship is shown in Figure~\ref{fig: inversion}.

It is also worth noting that $\Lambda$ is not Farey recursive.

\begin{remarks}\label{remark: Riley}
Suppose $\alpha$ is hyperbolic (Definition~\ref{def: hyperbolic alpha}).
\begin{enumerate}

\item
There are exactly two roots of $\Lambda(\alpha)$ for which the corresponding p-rep is discrete and faithful, and these roots are complex conjugates.  If $\frz_1$ is the geometric root of $\calQ(\alpha)$, Theorem~\ref{thm: QandL} identifies them as
\[ -\frac{1}{\frz_1} \AND -\frac{1}{\overline{\frz_1}}. \]
Thus, the discrete faithful Riley roots may be found by locating the geometric root of $\calQ(\alpha)$.  As described in Section~\ref{subsec: code II}, this is the unique root $\frz_1$ for which
\[ \calZ(e_j)_{\frz_1}\ \in \ \IH^2 \]
for every interior edge $e_j$ of $F_\alpha$.

\item
By Corollary~\ref{cor: holonomy specialization}, the group $\Gamma_{\frz(\alpha)}$ from Definition~\ref{def: generic holonomy} is discrete and not free. As outlined in Section~2.2 of \cite{EMS}, a result of Brenner \cite{Bre} gives
\[ |\frz(\alpha)|\ >\ \frac14,\]
while results of Shimizu \cite{Shi} and Leutbecher \cite{Leu} give
\[ |\frz(\alpha)|\ \leq\ 1. \]
Adams \cite{Ad3} showed that equality in the latter bound occurs
only when $\alpha=2/5$. These bounds are visible in Figures~\ref{fig: roots} and~\ref{fig: inversion}.

\end{enumerate}
\end{remarks}

\begin{remark}\label{rem: conj QandL}
Theorem~\ref{thm: QandL} also raises a different, arithmetic question about the defect.  Namely, when is
\[ \mu(p/q)\ >\ 0? \]

Computations using the Farey-recursive algorithms from Section~\ref{subsec: code II} suggest a strong $2$-adic pattern. For reduced $p/q\in\Q_0$ with $q\leq300$, one finds
\[ \mu(p/q)=0 \qquad\text{whenever}\qquad 8\nmid q, \]
while for $q\leq1000$, the condition $\mu(p/q)\geq2$ occurs only when $64\mid q$.  Larger defects do occur; for example,
\[ \mu(1503/4096)\ =\ \mu(1567/4096)=4.\]

It would be interesting to characterize $\mu(p/q)$ arithmetically. In particular, is it true that
\[ \mu(p/q)\ =\ 0 \qquad\text{whenever}\qquad 8\nmid q? \]
More generally, what divisibility conditions on $q$ are forced by $\mu(p/q)\geq r$?
\end{remark}

\section{Limits for 2-bridge links} \label{sec: limits}

This section studies limits of geometric roots and the corresponding Kleinian groups.  For instance, it shows that accumulation points of geometric roots are distinguished roots of the discriminant polynomials $D_\calQ(\omega)$ (Definition~\ref{def: Disc}), and identifies these points with cusps of the Riley slice.  It also examines the cusps associated to $\omega=1/3$ and $1/4$ in detail, describing the corresponding algebraic limits, convex cores, and geometric limits of $2$-bridge link groups (See, for instance, Figures~\ref{fig: limit} and \ref{fig: geometric}). 
\index{discriminant polynomial}

\subsection{Geometric roots and Riley cusps} \index{Riley slice!cusp}

\begin{example} \label{ex: torus limit}
Consider the numbers $\alpha_m=1/m$ with $m> 3$.  The funnel $F_{\alpha_m}$ has only one hub, namely $\omega=0$.  As $D_\calQ(0) = -4x+1$, the only root of $D_\calQ(\omega)$ is $\hat{\frz}=1/4$.

Remark ~\ref{rk: chebyshev} explains that the roots of $\calQ(\alpha_m)$ are the values $(1/4) \sec^2(k\pi/m)$.  Although $\alpha_m$ is not hyperbolic, and Definition~\ref{def: geometric root} does not assign a geometric root to $\calQ(\alpha_m)$, it is convenient in this example to designate
\[ \frz(\alpha_m)\ =\ \frac14\sec^2\left(\pi/m \right). \]
(This definition is particularly natural in light of Figure~\ref{fig: roots}.)  It is quick to check that $\frz(\alpha_m) \to \hat{\frz}$ as $m \to \infty$.

The links $L_{\alpha_m}$ are $(2,m)$-torus links and are not hyperbolic.  The diagrams $D_{\alpha_m}$ have only one twist region and, evidently, the link complement $M_{\alpha_m}$ can be obtained by Dehn filling on the vertical component of one of the links shown in Figure ~\ref{fig: Aug Torus}. 

\begin{figure}
   \centering
   \includegraphics[width=2in]{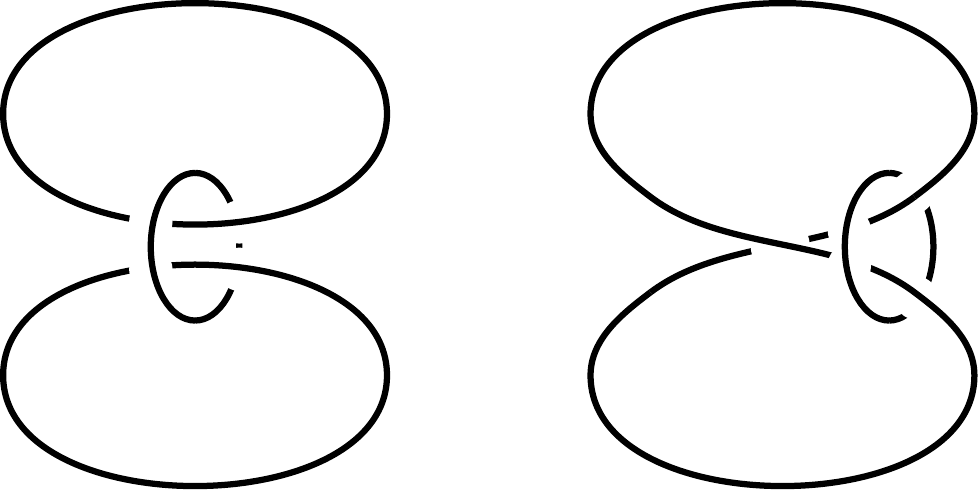} 
   \caption{Doing $(1,j)$-Dehn filling on vertical components yields 2-bridge torus links.}  \label{fig: Aug Torus}
\end{figure}
\end{example}

This same limiting behavior occurs when the number of crossings in the first twist region of an arbitrary hyperbolic $2$-bridge link is allowed to grow.

\begin{thm} \label{thm: convergence I}
Suppose that $\alpha$ is hyperbolic.  The funnel $F_\alpha$ provides a positive termed continued fraction expansion $\alpha = [0;a_1, \ldots, a_k]$ with each of $a_1$, $a_k$, and $k$ at least two.  For each integer $m \geq 2$, define $\alpha_m = [0; m, a_{2}, \ldots]$.  Then $\calH(M_{\alpha_m})$ consists of a single point $\frz(\alpha_m)$ and
\[ \lim_{m \to \infty} \frz(\alpha_m) \ =\ \frac{1}{4}. \]
\end{thm}

\proof
Because $\alpha$ is hyperbolic, Theorem~\ref{thm: complete} implies that $\calH(M_\alpha)$ consists only of the geometric root $\frz$ of $\calQ(\alpha)$.  Regard $M_\alpha$ as the hyperbolic manifold determined by $\frz$.

The paragraph just after Definition~\ref{def: diagram} explains how $F_\alpha$ provides a positive termed continued fraction expansion $\alpha=[0;a_1,\ldots,a_k]$ with $a_1,a_k\geq2$.  Since $\alpha$ is hyperbolic, we also have $k\geq2$. It follows that, for every $m\geq2$, $\alpha_m=[0;m,a_2,\ldots,a_k]$ is hyperbolic.  So, again by Theorem~\ref{thm: complete}, $\calH(M_{\alpha_m})$ consists of the single geometric root $\frz(\alpha_m)$ of $\calQ(\alpha_m)$.

As always, the first hub of $F_\alpha$ is $\omega_1=0$.  Let $C$ be the crossing circle for the first twist region of the diagram $D_\alpha$, and let
\[ N\ =\ S^3-(L_\alpha\cup C) \]
be the corresponding augmented link complement.  
\index{augmented link}
\index{crossing loop}
\index{crossing loop!twist region}

The main result of \cite{Ad1} (see also \cite{Ad2}) shows that $N$ admits a complete hyperbolic structure with finite volume.  Suppose first that $a_1-m$ is even.  Then $M_{\alpha_m}$ is obtained from $N$ by Dehn filling the crossing circle cusp.  Let $\lambda_C$ be the peripheral slope on this cusp which is parallel to the crossing circle $C$.  Under the filling epimorphism
\[ \pi_1N\longrightarrow\pi_1M_{\alpha_m} \]
the element represented by $\lambda_C$ maps to the free homotopy class of the crossing circle in $M_{\alpha_m}$.  By Theorem~\ref{thm: diagram}, this is the conjugacy class represented by $(S_0)_{\frz(\alpha_m)}$.

Thurston's hyperbolic Dehn filling theorem implies that, as $j\to\infty$, the loxodromic isometries $(S_0)_{\frz(\alpha_{a_1+2j})}$ converge to a parabolic isometry.  The geometric roots lie in a bounded subset of $\C$, so there is an accumulation point $\hat{\frz}$ for the sequence $\{\frz(\alpha_{a_1+2j})\}$.  Since $\tr(S_0)_\frz=2-1/\frz$ and the corresponding conjugacy classes converge to a parabolic class, it must be true that $\hat{\frz}\neq0$.  Passing to a convergent subsequence,
\[ (S_0)_{\hat{\frz}} \ =\  \frac{1}{\hat{\frz}} \, \begin{bmatrix} \hat{\frz}-1 & \hat{\frz} \\ -1 & \hat{\frz} \end{bmatrix} \]
is parabolic.  This occurs precisely when $\hat{\frz}$ satisfies $-4\hat{\frz}+1=0$. Hence, every accumulation point is $1/4$ and
\[ \lim_{j\to\infty}\frz\bigl(\alpha_{a_1+2j}\bigr)\ =\ \frac14. \]

The parity of $a_1-m$ changes only the augmented link used in the Dehn filling argument.  Replacing $a_1$ by $a_1+1$ gives the same argument for the other parity.   Consequently,
\[ \lim_{m\to\infty}\frz(\alpha_m)\ =\ \frac14. \]
\endproof

More generally, if the number of crossings in the $j^\text{th}$ twist region of a diagram $D_\alpha$ increases, then the geometric roots converge to a root of the discriminant $D_\calQ(\omega)$, where $\omega$ is given by truncating the continued fraction for $\alpha$ at its $j^\text{th}$ term.  Before stating this more general theorem, some useful connections to the geometry of the Stern--Brocot diagram $\calG$ are reviewed.

\begin{dfn} \label{def: twist}
Suppose $\alpha \in \Q_0$ and let $[0;a_1, \dots, a_k]$ be the positive termed continued fraction expansion for $\alpha$ with $a_1,a_k \geq 2$.  Take $j \in \{ 1, \ldots, k\}$.  Given $m \in \Z$, the {\it level-$j$ $m$-twisted adjustment} to $\alpha$ is the rational number given by the continued fraction obtained by replacing the $j^\text{th}$ term in $[0;a_1, \ldots, a_k]$ with the number $m$.
\end{dfn}
\index{twisted adjustment|idxdef}

\begin{figure}
\setlength{\unitlength}{.1in}
\begin{picture}(25.3,25)
\put(.3,.3) {\includegraphics[width= 2.5in]{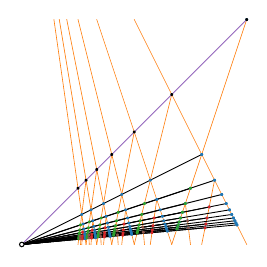}}
\put(0,1.0){$\scriptstyle (0,0)$}
\put(24.2,24.0){$\scriptstyle [0;2]=1/2$}
\put(11.5,16.5){$\scriptstyle [0;3]=1/3$}
\put(20.1,11){$\scriptstyle [0;2,2]=2/5$}
\put(21.2,8.5){$\scriptstyle [0;2,3]=3/7$}
\end{picture}
\caption{This figure shows how the level-1 $m$-twisted adjustments to a number $\alpha$ converge to the origin.  For reference, the lines of the Stern--Brocot diagram are drawn faintly in orange.  The level-1 adjustments to $1/2$ are black and converge down the purple line from $(1/2,1/2)$ to the origin.   The blue points are points of the form $[0;m,n]$.   Green points are of the form $[0;m,1,n]$.  Red points are of the form $[0;m,2,n]$.  Brown points are of the form $[0;m,1,2,n]$. } \label{fig: to origin}
\end{figure}

The numbers $\alpha_m$ described in Theorem ~\ref{thm: convergence I} are the level-1 $m$-twisted adjustments to $\alpha$.  In that case, $m$ was taken to be at least two in order to guarantee that $\alpha_m \in (0,1/2)$.  As shown in Figure ~\ref{fig: to origin}, if $\alpha$ is fixed and $m \geq 2$ is allowed to vary, the Stern--Brocot vertices for the level-1 $m$-twisted adjustments to $\alpha$ converge monotonically along a Euclidean line towards $(0,0) \in \R^2$.

\begin{figure}
   \centering
   \includegraphics[width=5in]{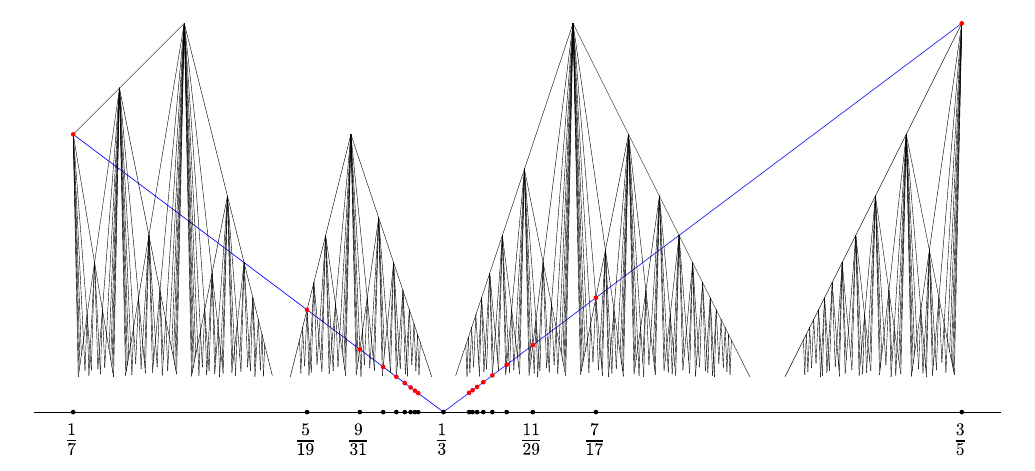}
   \caption{Take $\alpha = 9/31$.  The level-2 $m$-adjustments of $\alpha$ are marked with red points.  They lie on the blue lines $\ell^\pm$ across $\calG$. (This is the same as Figure 5 of \cite{Nav}.)} \label{fig: vee}
\end{figure}

\begin{remark} \label{rk: extension}
This remark provides relevant content from the main theorem of \cite{Nav}.  

Suppose that $\alpha$ is hyperbolic.  Let $[0;a_1, \dots, a_k]$ be the positive termed continued fraction expansion for $\alpha$ with $a_1,a_k \geq 2$ and take $j \in \{ 1, \ldots, k\}$.  For $m \in \Z$, let $\alpha_m$ be the level-$j$ $m$-twisted adjustment to $\alpha$ and define $\omega$ to be the number $[0;a_1, \ldots, a_{j-1}]$.  Let $\ell^+ \subset \R^2$ be the Euclidean line through $(\omega, 0)$ and the Stern--Brocot vertex for $\alpha$.  Let $\ell^-$ be the reflection of $\ell^+$ across the $x$-axis.   Write $v_m$ to denote the Stern--Brocot vertex for $\alpha_m$.
\begin{enumerate}
\item The 2-bridge link $L_{\alpha_m}$ has a 4-plat diagram $D(a_1, \ldots, m, \ldots, a_k)$ as described following Definition ~\ref{def: diagram}.  The diagram is alternating if and only if $m\geq 0$.
\item For every $m\in \Z$, $v_m \in \ell^+ \cup \ell^- \cup \{ \infty \}$.
\item If $m \leq -1$, then $v_m \in \ell^-$.
\item If $m \geq 0$, then $v_m \in \ell^+$.
\item The quantities
\[
\big| v_m - (\omega, 0) \big| \AND \big| v_{-m} - (\omega, 0) \big|
\]
converge monotonically to zero as $m \to \infty$.
\item Assume $j\geq 2$ and, if $a_1=j=2$, assume also that $m \geq 0$.  If $|m|$ is large enough, then $\alpha$ is hyperbolic and $\omega$ is a hub for the funnel $F_{\alpha_m}$.
\end{enumerate}
See Figures ~\ref{fig: to origin} and ~\ref{fig: vee}.
\end{remark}

When $j>1$, an analogous result (and argument) to Theorem~\ref{thm: convergence I} holds for the level-$j$ $m$-twisted adjustments to a number $\alpha$.  In this more general case, the discriminant polynomial may have several roots.  The next theorem identifies the possible accumulation points.  Afterwards, Theorem~\ref{thm: DQ uniqueness} will show that only one of them can occur.

\begin{thm} \label{thm: convergence II}
Suppose that $\alpha$ is hyperbolic and let $\alpha=[0;a_1,\ldots,a_k]$ be the positive termed continued fraction expansion provided by the funnel $F_\alpha$.  Take $j\in\{2,\ldots,k\}$ and define
\[ \omega\ =\ [0;a_1,\ldots,a_{j-1}]. \]
For $m\in\Z$, let $\alpha_m$ be the level-$j$ $m$-twisted adjustment to $\alpha$.  Whenever $|m|$ is sufficiently large and the conclusion of Remark~\ref{rk: extension} (6) applies, let $\frz(\alpha_m)$ be the geometric root of $\calQ(\alpha_m)$.  Then every accumulation point of $\left\{ \frz(\alpha_m) \right\}$ is a root of $D_\calQ(\omega)$.
\end{thm} \index{discriminant polynomial}

\begin{figure}
   \centering
   \includegraphics[width=2.25in]{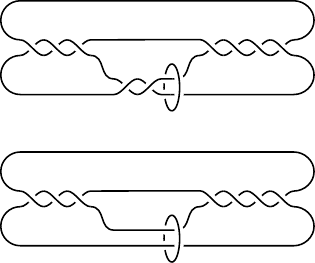}
   \caption{\index{crossing loop!twist region} On the top is the augmented link $N$, described in the proof of Theorem ~\ref{thm: convergence II}, in the case where $\alpha = 9/31$ and $j=2$.  The two link complements shown are homeomorphic to each other via a homeomorphism which twists along the disk bounded by the crossing circle.  A diagram for the augmented link $N'$, also discussed in the proof of Theorem ~\ref{thm: convergence II}, looks similar except that there is an odd number of half-twists on the pair of strands passing through the crossing circle.  Note that, using this paper's techniques, the geometries of $N$ and $N'$ are computed explicitly later in this subsection.}  \label{fig: Augment}
\end{figure}

Notice that, because every geometric root $\frz(\alpha)$ has positive real and imaginary parts, every accumulation point arising in Theorem~\ref{thm: convergence II} must have nonnegative real and imaginary parts.  

\proof
Throughout the proof, $m$ ranges over the admissible values from the statement, with $|m|$ sufficiently large that the properties in Remark~\ref{rk: extension} (6) hold.  In particular, each $\alpha_m$ is hyperbolic and $\omega$ is a hub of $F_{\alpha_m}$.

Let $C$ be the crossing circle for the $j^\text{th}$ twist region of the diagram $D_\alpha$, and let
\[
N\ =\ S^3-(L_\alpha\cup C)
\]
be the corresponding augmented link complement.  Let $\alpha'=\alpha_{a_j+1}$ and define $N'$ analogously by augmenting the $j^\text{th}$ twist region of $D_{\alpha'}$.  By Theorem~\ref{thm: diagram}, the element $S_\omega$ represents the free homotopy class of the crossing circle in each case.  By \cite{Ad1}, both $N$ and $N'$ admit complete finite volume hyperbolic structures.

Using (1) of Remark~\ref{rk: extension}, if $m-a_j$ is even, then $M_{\alpha_m}$ is obtained by Dehn filling the crossing circle cusp of $N$, while if $m-a_j$ is odd, it is obtained by Dehn filling the corresponding cusp of $N'$.  Arguing as in the proof of Theorem~\ref{thm: convergence I}, Thurston's hyperbolic Dehn filling theorem implies that, along each such filling sequence, $(S_\omega)_{\frz(\alpha_m)}$ converges to a parabolic isometry.

The geometric roots lie in a bounded subset of $\C$, so they have accumulation points.  Let $\hat{\frz}$ be any such accumulation point.  After passing to a subsequence of fixed parity, the corresponding conjugacy classes converge to a parabolic class.  The trace formula~(\ref{eq: trS two}) shows that $\hat{\frz}\neq0$, and hence $(S_\omega)_{\hat{\frz}}$ is a nontrivial parabolic.  Therefore, Theorem~\ref{thm: Sparabolic} implies that
\[ D_\calQ(\omega)_{\hat{\frz}}\ =\ 0. \]
\endproof

To prove that the accumulation point is unique, the next lemma will be useful.  From Definition~\ref{def: generic holonomy}, recall that the generic holonomy group $\Gamma$ is freely generated by $U_0$ and $W_0$ in $\pslQx$.

\begin{lemma} \label{lem: boundary groups} \index{holonomy!generic}
Suppose that $\{ \alpha_m \}$ is a sequence of pairwise distinct hyperbolic rational numbers.  Assume also that the geometric roots $\frz(\alpha_m)$ converge to a number $\hat{\frz} \in \C$. Then $\hat{\frz} \neq 0$, $\Gamma_{\hat{\frz}}$ is a nonelementary Kleinian group, and the evaluation homomorphism
\[\rho \co \Gamma \longrightarrow \Gamma_{\hat{\frz}}\]
is an isomorphism.  In particular, $\Gamma_{\hat{\frz}} = \big\langle (U_0)_{\hat{\frz}}, (W_0)_{\hat{\frz}} \big\rangle$ is free of rank two.
\end{lemma}

\proof
Remark~\ref{remark: Riley} (2) implies that $|\hat{\frz}|\geq 1/4>0$.  Hence, $\Gamma_{\hat{\frz}}$ is well-defined as the image of the evaluation $\rho$.  This proof first shows that $\Gamma_{\hat{\frz}}$ is a nonelementary Kleinian group and then proves that $\rho$ is injective.

To start, pass to a subsequence for which the numbers $\frz(\alpha_m)$ are pairwise distinct.  To see that this is possible, suppose $\frz(\alpha_m)=\frz(\alpha_n)=\frz$.  By Theorem~\ref{thm: complete} and Corollary~\ref{cor: holonomy specialization}, both $M_{\alpha_m}$ and $M_{\alpha_n}$ are identified with $\IH^3/\Gamma_\frz$.  Moreover, under these identifications, their canonical generators are both carried to
\[(U_0^{-1})_\frz \AND (W_0)_\frz \]
(see Definitions~\ref{def: generators} and~\ref{def: generic holonomy}).  The resulting homeomorphism $M_{\alpha_m}\to M_{\alpha_n}$ preserves the meridional filling slopes and extends over the meridional fillings to a homeomorphism of pairs $(S^3,L_{\alpha_m})\to(S^3,L_{\alpha_n})$.  Schubert's classification of $2$-bridge links \cite{Sch} implies that, for fixed $\alpha_m$, there are only finitely many possibilities for $\alpha_n$.  Thus each value of $\frz(\alpha_m)$ occurs only finitely many times.  Since the $\alpha_m$ are pairwise distinct, the sequence $\{\frz(\alpha_m)\}$ takes infinitely many distinct values, and the desired subsequence exists.

The defining generators of $\Gamma_{\frz(\alpha_m)}$ converge to those of $\Gamma_{\hat{\frz}}$, so these groups converge algebraically to $\Gamma_{\hat{\frz}}$ in the sense of \cite{Mar}. Each $\Gamma_{\frz(\alpha_m)}$ is a nonelementary Kleinian group, with quotient the hyperbolic manifold $M_{\alpha_m}$.  Theorem~4.1.1 of \cite{Mar} shows that $\Gamma_{\hat{\frz}}$ is also a nonelementary Kleinian group.  The same theorem shows that, for all sufficiently large $m$, the assignments
\[ (U_0)_{\hat{\frz}}\mapsto(U_0)_{\frz(\alpha_m)} \AND (W_0)_{\hat{\frz}}\mapsto(W_0)_{\frz(\alpha_m)} \]
determine a homomorphism $\varphi_m\co\Gamma_{\hat{\frz}} \to \Gamma_{\frz(\alpha_m)}$.  So, for large $m$, 
\[\varphi_m\, \rho(A)\ =\ (A)_{\frz(\alpha_m)}\] 
for every $A\in\Gamma$.
\index{algebraic limit}

It remains to prove that $\rho$ is injective.  Take $A\in\ker\rho$ and choose a representative
\[ \begin{pmatrix}a&b\\c&d\end{pmatrix}\in\SLtwo{\Q(x)} \]
for $A$.  Since $\rho(A)$ is trivial in $\PSL_2(\C)$, $(A)_{\frz(\alpha_m)}=\varphi_m\rho(A)$ is represented by $\pm I$.
Hence
\[ 0 \ =\  b(\frz(\alpha_m)) \ =\  c(\frz(\alpha_m)) \ =\  (a-d)(\frz(\alpha_m)). \]
The numbers $\frz(\alpha_m)$ are pairwise distinct, while a nonzero rational function has only finitely many zeros.  So,  $b=c=a-d=0$ in $\Q(x)$.  Since the chosen representative has determinant one, $a^2=1$, and $A$ is trivial in $\PSL_2(\Q(x))$.  This shows that $\rho$ is injective.
\endproof

\begin{dfns} \label{def: Riley slice}
If $\Gamma_\frz$ is discrete, its {\it regular set} $\calD(\Gamma_\frz)$ is the subset of the Riemann sphere on which $\Gamma_\frz$ acts properly discontinuously.  The {\it Riley slice} is the set
\[ \calR \ =\ \left\{ \frz \in \C^\ast \, \big| \, \Gamma_\frz \text{ is discrete and } \calD(\Gamma_\frz)/\Gamma_\frz \text{ is a 4-punctured sphere}\right\}. \]
The closure $\overline{\calR}$ of $\calR$ in $\C^\ast$ is called the {\it extended} Riley slice.  Its {\it boundary} is the set $\bound \calR = \overline{\calR}-\calR$.  A point $\frz \in \bound \calR$ for which $\calD(\Gamma_\frz)/\Gamma_\frz$ is a pair of 3-punctured spheres is called a {\it cusp}.
\end{dfns}
\index{regular set|idxdef}
\index{Riley slice|idxdef}
\index{Riley slice!boundary|idxdef}
\index{Riley slice!cusp|idxdef}

\begin{facts} \label{fact: Riley} The following items are well-established.
\begin{enumerate}
\item The extended Riley slice $\overline{\calR}$ is contained in the unit disk  (second point in Remark ~\ref{remark: Riley}).
\item The boundary of the Riley slice $\bound \calR$  is a single Jordan curve, \cite{OM}.
\item From \cite{OM},
\[ \overline{\calR} \ =\ \left\{ \frz \in \C^\ast \, \big| \, \Gamma_\frz \text{ is discrete and free with rank 2} \right\}. \]
\item The cusps are dense in $\bound \calR$, \cite{CCHS} and  \cite{McM}.
\end{enumerate}
\end{facts}

\begin{thm}\label{thm: DQ uniqueness}
Suppose $\omega\in\Q\cap(0,1/2]$.  There is a unique root $\hat{\frz}(\omega)$ of $D_\calQ(\omega)$ which can occur as an accumulation point in Theorem~\ref{thm: convergence II}.  In particular, the geometric roots in that theorem converge to this well-defined $\hat{\frz}(\omega)$.  Together with the point $\hat{\frz}(0)=1/4$ from Theorem~\ref{thm: convergence I}, these points constitute the full set of cusps for $\bound\calR$ in the first quadrant of $\C$.
\end{thm}

\proof
Take $\omega$ and $\{\alpha_m\}$ as in Theorem~\ref{thm: convergence II}, and let $\hat{\frz}$ be an accumulation point of the geometric roots  $\frz(\alpha_m)$.  Pass to a subsequence to assume that $\frz(\alpha_m)\to \hat{\frz}$.  By Theorem~\ref{thm: convergence II}, $\hat{\frz}$ is a root of
$D_\calQ(\omega)$.  Theorem~\ref{thm: Sparabolic} implies that $\hat{\frz}\neq0$ and that $(S_\omega)_{\hat{\frz}}$ is a nontrivial parabolic.

By Remark~\ref{rk: extension}, after discarding finitely many terms, the numbers $\alpha_m$ are pairwise distinct hyperbolic rational numbers.  Lemma~\ref{lem: boundary groups} applies and $\Gamma_{\hat{\frz}}$ is a nonelementary Kleinian group which is free of rank two.  By Fact~\ref{fact: Riley} (3),
$\hat{\frz}\in\overline{\calR}$. On the other hand, each $\Gamma_{\frz(\alpha_m)}$ is a hyperbolic $2$-bridge link group and is not free, so the same fact gives
$\frz(\alpha_m)\notin\overline{\calR}$.  Therefore, $\hat{\frz}\in\bound\calR$.

The elements 
\[ \big(U_0\big)_{\hat{\frz}},\qquad \big(W_0\big)_{\hat{\frz}},\AND \big(S_\omega\big)_{\hat{\frz}}\]
are nontrivial parabolics.  Since $\Gamma_{\hat{\frz}}$ is free, it contains no subgroup isomorphic to $\Z^2$, so each of the maximal parabolic subgroups containing these parabolics is cyclic.  Lemma~\ref{lem: boundary groups} says that the specialization
\[ \rho_{\hat{\frz}}\co \Gamma\longrightarrow\Gamma_{\hat{\frz}} \]
is an isomorphism.  Since $U_0$ and $W_0$ freely generate $\Gamma$,
neither is a proper power.  Hence, the maximal parabolic subgroups containing $\big( U_0\big)_{\hat{\frz}}$ and $\big( W_0\big)_{\hat{\frz}}$ are
\[ \Big\langle \big( U_0\big)_{\hat{\frz}} \Big\rangle \AND \Big\langle \big( W_0\big)_{\hat{\frz}} \Big\rangle \]
These two subgroups are not conjugate, because $\Gamma_{\hat{\frz}}$ is freely generated by $\big( U_0\big)_{\hat{\frz}}$ and $\big( W_0\big)_{\hat{\frz}}$.

Now, the trace identity~(\ref{eq: trS two}) and the definition of $D_\calQ$ imply that
\[ \big(\tr S_\omega\big)^2-4 \ =\ \frac{\calQ(\omega)^2D_\calQ(\omega)}{\dQ(\omega)^2}.\]
By Theorem~\ref{thm: Sparabolic}, this expression is nonzero in $\Q(x)$.  On the other hand, 
\[ \big(\tr U_0\big)^2 \ =\ \big(\tr W_0 \big)^2 \ = \ 4 \]
in $\Q(x)$.  This shows that $S_\omega$ is conjugate into neither $\big\langle U_0 \big\rangle$ nor $\big\langle W_0 \big\rangle$.  Because $\rho_{\hat{\frz}}$ is an isomorphism, $\Gamma_{\hat{\frz}}$ has at least three conjugacy classes of maximal parabolic subgroups.  In fact, it has exactly three, by Theorem~4.9 of \cite{MT}, which states that three is the maximum possible number for a free rank-two Kleinian group.  Together with the remark following the proof of Theorem~4.9 in \cite{MT}, this implies that $\Gamma_{\hat{\frz}}$ is geometrically finite and that the boundary of its convex core consists of two totally geodesic $3$-punctured spheres.  Therefore, $\hat{\frz}$ is a cusp point.

It is shown in \cite{KMS} that, for each $\omega$, there is exactly one cusp point for which the specialization of $S_\omega$ is parabolic.  Therefore, every accumulation point arising in Theorem~\ref{thm: convergence II} is the same point $\hat{\frz}(\omega)$.  Since the geometric roots form a bounded sequence, they consequently converge to $\hat{\frz}(\omega)$.

For $\omega=0$, Theorem~\ref{thm: convergence I} gives $\hat{\frz}(0)=1/4$.  Finally, at any cusp of the Riley slice the third maximal parabolic class corresponds to an essential simple closed curve on the $4$-punctured sphere.  By the slope
parametrization from Remark~\ref{rem: slopes}, and the usual symmetries of the slope parametrization, the cusps in the first quadrant are represented by $S_\omega$ with $\omega\in\Q\cap[0,1/2]$. The uniqueness result of \cite{KMS} shows that these are precisely the points $\hat{\frz}(\omega)$.
\endproof

\begin{figure}
\setlength{\unitlength}{.1in}
\begin{picture}(50,24)
\put(0,0) {\includegraphics[width= 4.8in]{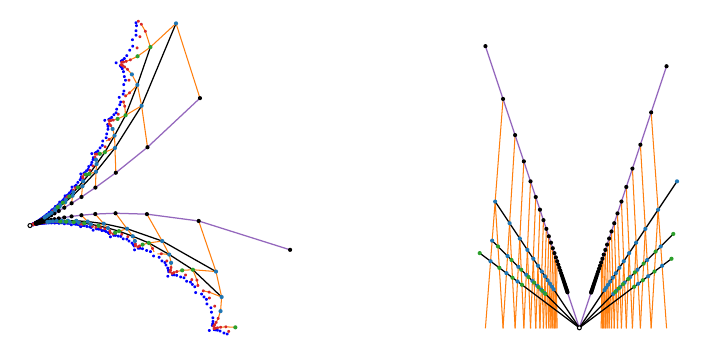}}
\put(0,7.5){$\scriptstyle \hat{\frz}(1/3)$}
\put(12.2,23.1){$\scriptstyle 9/25$}
\put(14.1,17.8){$\scriptstyle 5/14$}
\put(20.3,7){$\scriptstyle 4/13$}
\put(13.6,9.5){$\scriptstyle 5/16$}
\put(9.9,13){$\scriptstyle 6/17$}
\put(15.3,5.5){$\scriptstyle 9/29$}
\put(15.5,3.6){$\scriptstyle 13/42$}

\put(30.2,7.8){$\scriptstyle 13/42$}
\put(37.8,0.4){$\scriptstyle (1/3,0)$}
\put(46.7,12.1){$\scriptstyle 9/25$}
\put(42.7,20.4){$\scriptstyle 5/14=[0;2,1,4]$}
\put(30.5,21.8){$\scriptstyle 4/13=[0;3,4]$}
\put(31.5,17.3){$\scriptstyle 5/16$}
\put(45.0,16.5){$\scriptstyle 6/17$}
\put(31.1,10.6){$\scriptstyle 9/29$}
\end{picture}
\caption{This figure compares convergence in the Stern--Brocot diagram with convergence of the corresponding geometric roots near the cusp for $\omega=1/3$.  In the right-hand image, several families of twisted adjustments converge to $(1/3,0)$.  The black families include $[0;3,m]$ and $[0;2,1,m]$; the blue and green families are obtained by appending fixed tails, for example $[0;3,m,r]$, $[0;2,1,m,r]$, $[0;3,m,1,r]$, and $[0;2,1,m,1,r]$, for the values shown in the figure.  The left-hand image shows the corresponding geometric roots, using the same colors.  Every family converges to $\hat{\frz}(1/3)$, marked by the small open circle.} \label{fig: cusp13 w diag}
\end{figure}

\begin{dfn} \label{def: D geometric root}
Given $\omega\in\Q\cap[0,1/2]$, Theorem~\ref{thm: convergence I} and Theorem~\ref{thm: DQ uniqueness} provide a well-defined root $\hat{\frz}(\omega)$ of the discriminant polynomial $D_\calQ(\omega)$.  This distinguished root is referred to as the {\it geometric root} of $D_\calQ(\omega)$.
\end{dfn}
\index{geometric root!discriminant polynomial|idxdef}
\index{discriminant polynomial}

Theorems~\ref{thm: convergence I}, \ref{thm: convergence II}, and~\ref{thm: DQ uniqueness} show that each geometric root $\hat{\frz}(\omega)$ arises as the limit of many different sequences of geometric roots associated to complete hyperbolic $2$-bridge link complements.  In the case of $\omega=1/3$, this can be seen in Figure ~\ref{fig: cusp13 w diag}.  Observe that
\index{geometric root!2-bridge link}
\[ D_\calQ(1/3)\ =\ (4x^2-3x+1)(x+1) \AND
\hat{\frz}(1/3) \ =\  \frac18\left(3+i\sqrt7\right). \]

The next theorem shows that the Riley slice description separates the roots of the $\calQ$-polynomials and their discriminants from the interior of the Riley slice.

\begin{thm}\label{thm: R is empty}
\ButterThm
\end{thm}

\proof
Suppose first that $\frz$ is a root of $\calQ(\alpha)$. Remark~\ref{rem: constant term} shows that $\frz\neq0$. Since $\calQ(\alpha)$ is not the zero polynomial, Theorem~\ref{thm: Sparabolic} (3) implies that $S_\alpha$ is nontrivial in the generic free group $\Gamma$.  On the other hand, the same theorem gives that $S_\alpha$ is in the kernel of the evaluation $\Gamma \to \Gamma_\frz$.  Hence, $\Gamma_\frz$ is not free of rank two and Fact~\ref{fact: Riley} (3) implies that $\frz\notin\overline{\calR}$.

Now, suppose that $\hat{\frz}$ is a root of $D_\calQ(\alpha)$.  Theorem~\ref{thm: Sparabolic} shows that $\hat{\frz}\neq0$ and that $(S_\alpha)_{\hat{\frz}}$ is a nontrivial parabolic.  If $\hat{\frz}\in\calR$, then $\calD(\Gamma_{\hat\frz})/\Gamma_{\hat\frz}$ is a $4$-punctured sphere.  By Remark~\ref{rem: slopes}, $S_\alpha$ represents the essential, nonperipheral simple loop of slope $\alpha$.  Such an element is loxodromic in the interior of
the Riley slice \cite{KS}.  This contradicts the fact that $(S_\alpha)_{\hat{\frz}}$ is parabolic.  Hence, $\hat{\frz}\notin\calR$.
\endproof

This shows that the empty region in the center of Figure \ref{fig: allroots} truly contains no roots. 

It is time to revisit Figure~\ref{fig: roots}.  The figure shows the quadrant of $\C$ with nonnegative real and imaginary parts.  The trivalent blue vertices are the geometric roots of the polynomials $\calQ(\alpha)$ and the red vertices are the geometric roots of the discriminant polynomials $D_\calQ(\alpha)$.  Here, the numbers $\alpha$ are taken in $\Q \cap (0,1/2)$. The black edges are drawn between blue vertices according to the Stern--Brocot tree shown in Figure ~\ref{fig: tree}.  The full Riley slice is not confined to this quadrant; it is symmetric under negation and complex conjugation.  Hence, from Figure ~\ref{fig: roots}, it is easy to imagine the cauliflower picture for the full Riley slice.  Because the red points are cusp points and they are dense in $\bound \calR$, they mark the Jordan curve of $\bound \calR$.  After deleting the origin, the open region bounded by this curve is the Riley slice.  The blue points correspond to the 2-bridge links and they lie outside the Jordan curve and converge back to $\bound \calR$.

\subsection{Algebraic and geometric limits at the cusps} \index{algebraic limit} \index{geometric limit}

To help understand the limiting geometry, it will be helpful to locate the fixed point of the parabolic crossing loop. \index{crossing loop!holonomy}

\begin{lemma} \label{lem: fixed point parabolic S}
Suppose $\omega \in \Q \cap [0,1/2]$ and $\hat{\frz}$ is a root of $D_\calQ(\omega)$.  The isometry $(S_\omega)_{\hat{\frz}}$ is parabolic  and fixes
\[ \hat{p} \ =\ -\frac12 + \calV(\omega)_{\hat{\frz}}. \]
Its action on $\IH^3$ restricts to an action on the geodesic plane over $\R+\calV(\omega)_{\hat{\frz}}$.
\end{lemma}

\proof
By Theorem ~\ref{thm: Sparabolic}, $(S_\omega)_{\hat{\frz}}$ is a nontrivial parabolic.  Let
\[ d = \dQ(\omega), \quad Q=\calQ(\omega), \AND N=\calN(\omega). \]
Then, using the formula given just after Definition ~\ref{def: S and T},
\[ S_\omega \ =\ \begin{bmatrix} -NQ+d & N^2-NQ+d \\ -Q^2 & -Q^2+NQ+d \end{bmatrix}. \]
Matrix multiplication provides
\[ S_\omega \, \begin{bmatrix} -Q+2N \\ 2Q \end{bmatrix} \ =\ \begin{bmatrix} -NQ^2 + 2Nd + Qd \\ -Q(Q^2 - 2d)\end{bmatrix}. \]
Because $\hat{\frz}$ is a root of $D_\calQ(\omega)$,
\[ \calQ(\omega)_{\hat{\frz}}^2 \ =\  4\dQ(\omega)_{\hat{\frz}}.\]
Hence,
\[ (S_\omega)_{\hat{\frz}} \, \begin{bmatrix} -Q+2N \\ 2Q \end{bmatrix}_{\hat{\frz}}
\ =\ \begin{bmatrix} -d (-Q+2N) \\ -2dQ\end{bmatrix}_{\hat{\frz}}. \]
Since $\hat{p} = \left( \frac{-Q+2N}{2Q} \right)_{\hat{\frz}}$, it follows that $\hat{p}$ is the fixed point for $(S_\omega)_{\hat{\frz}}$.

Lemma ~\ref{lem: STbeta} gives that $S_\omega(-1+\calV(\omega)) = \calV(\omega)$ and a quick computation shows that $(S_\omega)_{\hat{\frz}} (\infty) = -1/4+\calV(\omega)_{\hat{\frz}}$.  This finishes the argument that the action of $(S_\omega)_{\hat{\frz}}$ restricts to an isometry of the geodesic plane over $\R+\calV(\omega)_{\hat{\frz}}$.
\endproof

The remainder of this subsection focuses on a concrete example which illustrates both the algebraic and geometric convergence described in Theorems~\ref{thm: convergence I} and~\ref{thm: convergence II}. The Sakuma--Weeks triangulations and the machinery developed in this paper, are used to explicitly describe the algebraic limits, convex cores, and geometric limits associated with their associated augmented link complements.  The terminology concerning algebraic and geometric limits of Kleinian groups follows \cite{Mar}.  The example anticipates the more general results of \cite{CEP}.

Let
\[ \alpha\ =\ \frac9{31}=[0;3,2,4], \]
and, for $m\geq1$, consider the level-2 $m$-adjustments $\alpha_m=[0;3,m,4]$ of $\alpha$.  For example, $\alpha_8=[0;3,8,4]=33/103$.  The numbers $\alpha_m$ lie on the left-hand line $\ell^+$ shown in Figure~\ref{fig: vee}.  As discussed in Section~\ref{sec: reverse}, the reversed numbers
\[ \bar\alpha_m\ =\ [0;4,m,3] \]
determine hyperbolic manifolds $M_{\bar\alpha_m}$ isometric to $M_{\alpha_m}$.

\begin{figure}
   \centering
   \includegraphics[width=5in]{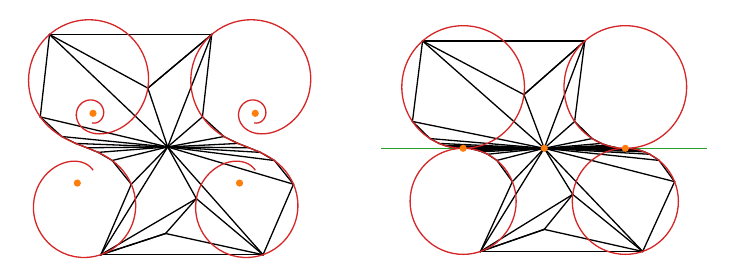}
\caption{\index{fundamental domain!specialized}\index{geometric root!discriminant polynomial} Take $\alpha_8=[0;3,8,4]=33/103$ and $\omega=1/3$.  On the  left is the fundamental domain $\Omega(\alpha_8)_{\frz(\alpha_8)}$ for the complete hyperbolic structure on $M_{\alpha_8}$.  The orange points are the fixed points for $(S_\omega)_{\frz(\alpha_8)}$ and its horizontal conjugate.  The picture on the right is an amalgamation of the domains shown in Figure~\ref{fig: twolimits}.  Later, it will be identified with the geometric limit of the domains $\Omega(\alpha_{2m})_{\frz(\alpha_{2m})}$ as $m\to\infty$. It describes the geometry of the augmented link shown in Figure~\ref{fig: Augment}.} \label{fig: limit}
\end{figure}

\subsubsection{Algebraic limits} \label{sec: algebraic}
\index{algebraic limit}
 
By Theorems~\ref{thm: convergence II} and~\ref{thm: DQ uniqueness}, the geometric roots $\frz(\alpha_m)$ converge to the geometric root
\[ \hat{\frz}(1/3) \ =\ \frac{1}{8}\left( 3+i\sqrt{7} \right) \]
of $D_\calQ(1/3)$.  Likewise, the geometric roots $\frz(\bar{\alpha}_m)$ converge to the geometric root $\hat{\frz}(1/4)$ of
\[ D_\calQ(1/4) \ = \ 4x^4 + 4x^2 - 4x + 1. \]
Here, the number $\hat{\frz}(1/4)$ is the unique root of $D_\calQ(1/4)$ whose real and imaginary parts are both positive and
\[ \hat{\frz}(1/4) \ \sim \ 0.39307568887871164 + 0.1360098247570345i. \]

By the observation following Corollary~\ref{cor: holonomy specialization}, the holonomy maps associated to $\frz(\alpha_m)$ and $\frz(\bar\alpha_m)$ are isomorphisms, and
\[ M_{\alpha_m} \ \cong\  \IH^3/\Gamma_{\frz(\alpha_m)} \AND M_{\bar\alpha_m} \ \cong \ \IH^3/\Gamma_{\frz(\bar\alpha_m)}. \]
The generators for $\Gamma_{\frz(\alpha_m)}$ and $\Gamma_{\frz(\bar\alpha_m)}$ converge to those of $\Gamma_{\hat{\frz}(1/3)}$ and $\Gamma_{\hat{\frz}(1/4)}$.  Therefore,
\[ \Gamma_{\frz(\alpha_m)} \longrightarrow \Gamma_{\hat{\frz}(1/3)} \AND \Gamma_{\frz(\bar\alpha_m)} \longrightarrow \Gamma_{\hat{\frz}(1/4)}\]
algebraically.

\subsubsection{Convex cores in the limit}
\index{convex core}

In what follows, it will be convenient to consider the numbers $[0;3,m]$ and $[0;4,m]$ obtained by truncating the continued fractions for $\alpha_m$ and $\bar{\alpha}_m$.  Assume $q\in\{3,4\}$ and write $\omega=1/q$ and $\hat{\frz}=\hat{\frz}(\omega)$.  Theorems~\ref{thm: convergence II} and~\ref{thm: DQ uniqueness} imply that the geometric roots $\frz([0;q,m])$, for $m\geq2$, converge to $\hat{\frz}$ as $m\to\infty$.  Likewise, the nonelementary Kleinian groups $\Gamma_{\frz([0;q,m])}$, for $m\geq2$, converge algebraically to $\Gamma_{\hat{\frz}}$.

\begin{figure}
\setlength{\unitlength}{.1in}
\begin{picture}(50,10)
\put(0,0) {\includegraphics[width= 5in]{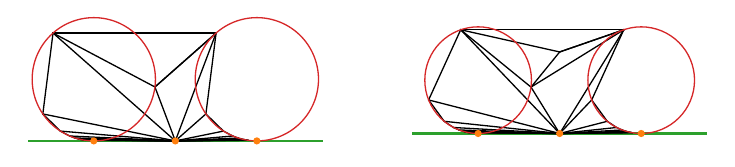}}
\put(10.5,.5){$\scriptstyle \calV(\omega)$}
\put(3.5,.5){$\scriptstyle -\frac{1}{2}+\calV(\omega)$}
\put(11,5.0){$\scriptstyle  V_{0}$}
\put(1.2,9.2){$\scriptstyle  V_{-1}$}
\put(.5,3.2){$\scriptstyle  V_{-2}$}
\put(36.6,1){$\scriptstyle \calV(\omega)$}
\put(29.6,1){$\scriptstyle -\frac{1}{2}+\calV(\omega)$}
\put(36.7,4.8){$\scriptstyle  V_{0}$}
\put(28.8,9.4){$\scriptstyle  V_{-1}$}
\put(27.0,4.1){$\scriptstyle  V_{-2}$}
\end{picture}
\caption{\index{fundamental domain!specialized}\index{geometric root!discriminant polynomial} The picture on the left shows $\wh{\Omega}_3$ and
the picture on the right shows $\wh{\Omega}_4$.}  
\label{fig: twolimits}
\end{figure}

\begin{figure}
\setlength{\unitlength}{.1in}
\begin{picture}(31.5,23.3)
\put(1.2,1.8) {\includegraphics[width= 3in]{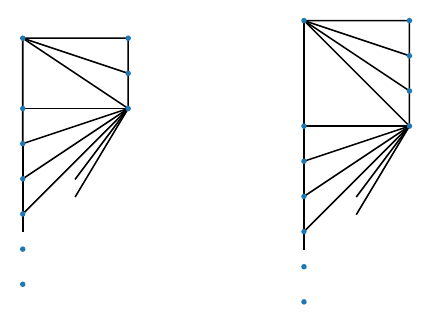}}
\put(3.2,1.6){$\scriptstyle \vdots$}
\put(3.4,4){$\scriptstyle 6/19=[0;3,6]$}
\put(3.4,6.5){$\scriptstyle 5/16=[0;3,5]$}
\put(0,9){$\scriptstyle 4/13$}
\put(0,11.4){$\scriptstyle 3/10$}
\put(0.5,13.8){$\scriptstyle 2/7$}
\put(0.5,16.3){$\scriptstyle 1/4$}
\put(1.6,21.3){$\scriptstyle 0$}
\put(30.1,22.5){$\scriptstyle 1$}
\put(30.1,20.0){$\scriptstyle 1/2$}
\put(30.1,17.6){$\scriptstyle 1/3$}
\put(30.1,15.2){$\scriptstyle 1/4$}
\put(10.6,21.3){$\scriptstyle 1$}
\put(10.6,19.0){$\scriptstyle 1/2$}
\put(10.6,16.3){$\scriptstyle 1/3$}
\put(21.3,22.5){$\scriptstyle 0$}
\put(20.1,15.2){$\scriptstyle 1/5$}
\put(20.1,12.7){$\scriptstyle 2/9$}
\put(19.6,10.2){$\scriptstyle 3/13$}
\put(19.6,7.9){$\scriptstyle 4/17$}
\put(22.9,5.3){$\scriptstyle 5/21=[0;4,5]$}
\put(22.9,2.8){$\scriptstyle 6/25=[0;4,6]$}
\put(22.9,.2){$\scriptstyle \vdots$}
\end{picture}
\caption{On the left is the infinite funnel ending with the numbers $[0;3,m]$.  On the right is the infinite funnel ending with the numbers $[0;4,m]$. } \label{fig: InfiniteFunnels}
\end{figure}

Define $\Omega\big( [0;q,m]\big)_{\hat{\frz}}$ as in Definition ~\ref{def: generic domain}.  Then there is a chain of inclusions
\[
\Omega\big( [0;q,1]\big)_{\hat{\frz}}
\ \subset \
\Omega\big( [0;q,2]\big)_{\hat{\frz}}
\, \subset \ \cdots \ \subset \
\Omega\big( [0;q,m]\big)_{\hat{\frz}}
\ \subset \ \cdots .
\]
Define
\[
\wh{\Omega}_q \ =\ \bigcup_1^\infty \ \Omega\big( [0;q,m]\big)_{\hat{\frz}} \ \subset \ \IH^3.
\]
See Figure ~\ref{fig: twolimits}.  Let $\wh{F}_q$ be the infinite funnel
\[
\wh{F}_q \ = \ \bigcup_1^\infty F_{[0;q,m]}.
\]
Both possibilities are pictured in Figure ~\ref{fig: InfiniteFunnels}.  As in Definition ~\ref{def: funnel}, take $\{ e_j \}_1^\infty$ to be the interior edges of $\wh{F}_q$ indexed from top to bottom.  For $j \geq 1$, write $\delta_j=\delta(e_j)_{\hat{\frz}}$ and $\delta'_j=\delta'(e_j)_{\hat{\frz}}$, where $\delta$ and $\delta'$ are defined as in Definition ~\ref{def: generic tetrahedra}.  Then
\[
\wh{\Omega}_q \ =\ \bigcup_{j \geq 1} \Big( \delta_j \cup  \delta'_j \Big).
\]

Following the notation from Corollary ~\ref{cor: right hub}, let $\beta_0=1/(q-1)$.  Recall from Theorem ~\ref{thm: Sparabolic}, that $(S_\omega)_{\hat{\frz}}$ is parabolic.  Recall $R_\omega$ from Definition~\ref{def: XYR}.  By Theorem ~\ref{thm: spiral}, $R^2_\omega = S_\omega$.  Therefore, $(R_\omega)_{\hat{\frz}}$ is also parabolic.  Define
\[ V_j \ =\  (R_\omega)^j_{\hat{\frz}} \, \left(\calV(\beta_0)_{\hat{\frz}}\right) \]
so then
\[ V_j \ =\ \lim_{m \to \infty} (R_\omega)^j_{\frz([0;q,m])} \, \left(\calV(\beta_0)_{\frz([0;q,m ])}\right). \]

Because $(R_\omega)_{\hat{\frz}}$ is parabolic, the points $\{V_j \}_\Z$ lie on a circle through the fixed point for $(R_\omega)_{\hat{\frz}}$.  By Lemma ~\ref{lem: fixed point parabolic S}, this fixed point is
\[ \hat{p}_q \ = \ -1/2+V_\omega \]
where $V_\omega = \calV(\omega)_{\hat{\frz}}$.  Also by Lemma ~\ref{lem: fixed point parabolic S}, the line $\R+V_\omega$ is invariant under the action of $(R_\omega)_{\hat{\frz}}$.  Therefore, the line and the circle are tangent at $\hat{p}_q$.  This geometry is evident in Figure ~\ref{fig: twolimits}.

By Definition ~\ref{def: generic tetrahedra} and Corollary ~\ref{cor: right hub},
\[ \delta_j \ =\ \Big[ V_\omega, \, V_{-j+1}, \, V_{-j+2}, \, \infty\Big] \AND \delta'_j \ =\ \Big[ V_\omega, \, 1+V_{-j}, \, 1+V_{-j+1}, \, \infty\Big]. \]
  
\begin{lemma} \label{lem: fat}
For each $j\geq 1$, the imaginary part of $\calZ(e_j)_{\hat{\frz}}$ is strictly positive.  Hence,
\[ \hat{\frz}\ \in\ \calH\left(M_{[0;q,m]}^\circ\right) \]
for every $m\geq1$.
\end{lemma}

\proof
It is easy to verify by hand that $\calZ(e_j)_{\hat{\frz}}$ has positive imaginary part for $j \in \{ 1,2,3 \}$.  Also, the value $\calZ(e_j)_{\frz([0;q,m])}$ has positive imaginary part for all large integers $m$, so the imaginary part of the limiting value $\calZ(e_j)_{\hat{\frz}}$ is at least zero. This reduces the problem to showing that, if $i \in \Z_{\geq 2}$, then $V_\omega$, $V_{-i}$, and $V_{-i-1}$ are not collinear in $\C$. 

As described above, the points $\{V_j\}_{j\in\Z}$ lie on a circle $C$ tangent to the line $\R+V_\omega$ at $\hat{p}_q$.  The two lines from $V_\omega$, tangent to $C$, divide $C$ into two arcs.  Let $A$ denote the arc closest to $V_\omega$.   A calculation shows that
\[ V_{-1} \ =\ 0 \ \in \ C-A \AND V_0 \ \in \ A.\]
This is evident in Figure~\ref{fig: twolimits}.  Since $\big( R_\omega^{-1} \big)_{\hat{\frz}}$ is a parabolic which fixes $\hat{p}_q$ and takes $V_{j}$ to $V_{j-1}$, it must move points counterclockwise around $C$, limiting towards $\hat{p}_q$.  So,  
\[ V_{-j} \ \in \ C-A\]
for every $j\geq 1$. 

Now, fix a particular $i \geq 2$.  Because $V_{-i} \in C-A$, the line between $V_{-i}$ and $V_\omega$ meets $C$ in exactly two points; namely $V_{-i}$ and one point in $A$.  But $V_{-i-1} \in C-A$, so $V_\omega$, $V_{-i}$, and $V_{-i-1}$ are not collinear.
\endproof

For each $m\geq 2$, endow $M_{[0;q,m]}^\circ$ with the hyperbolic structure given by $\hat{\frz}$.  By Lemma \ref{lem: fat}, there is a developing map
\[ \dev_m \co \wt{M}_{[0;q,m]}^\circ \longrightarrow \IH^3 \] \index{developing map}
with associated holonomy
\[ H_m \co \pie M_{[0;q,m]}^\circ \longrightarrow \pslC.\]
By Theorem \ref{thm: domain}, $\Omega([0;q,m])_{\hat{\frz}}$ is the image of a fundamental domain under $\dev_m$.  Corollary \ref{cor: holonomy specialization} shows that $\Gamma_{\hat{\frz}}$ is the image of $H_m$.  Recall from Definition \ref{def: generators} that $\pie M_{[0;q,m]}^\circ$ is freely generated by its canonical generators $k_0$ and $k_1$.  As mentioned in Definition \ref{def: generic holonomy}, this group is identified with the generic holonomy group $\Gamma$ by the isomorphism determined by $k_0 \mapsto U_0^{-1}$ and $k_1 \mapsto W_0$.   When this is done, the holonomy map is given by the evaluation homomorphism $\rho \co \Gamma \to \Gamma_{\hat{\frz}}$ discussed in Lemma~\ref{lem: boundary groups}.  By that lemma, the holonomy map is an isomorphism onto $\Gamma_{\hat{\frz}}$.  Hence the
\index{fundamental domain}
developing map descends to a local isometry
\[ \overline{\dev}_m\co M_{[0;q,m]}^\circ\longrightarrow \IH^3/\Gamma_{\hat{\frz}}. \]
In fact, this map is an embedding.  Any failure of injectivity would give an open overlap between two developed tetrahedra.  Any overlap would persist after specializing to $\frz([0;q,n])$ for all sufficiently large $n$.  This would contradict the complete Sakuma--Weeks triangulation of $M_{[0;q,n]}$.  

The image of $M_{[0;q,m]}^\circ$ is triangulated by the image of $\Omega([0;q,m])_{\hat{\frz}}$ and 
\[ \overline{\dev}_2 M_{[0;q,2]}^\circ \ \subset \ \overline{\dev}_3 M_{[0;q,3]}^\circ \, \subset \ \cdots \ \subset \ \overline{\dev}_m M_{[0;q,m]}^\circ \ \subset \ \cdots \ \subset \ \wh{\Omega}_q/\Gamma_{\hat{\frz}}.\]
By construction, the inclusions respect the labeled triangulations of the tangle complements $M_{[0;q,m]}^\circ$.  This provides a hyperbolic ideal triangulation of the subset $\wh{\Omega}_q/\Gamma_{\hat{\frz}}$ of $\IH^3/\Gamma_{\hat{\frz}}$.  The edges of the triangulations are labeled with the rational numbers from the vertices of $\hat{F}_q$.  Recalling Definition \ref{def: quotient complex}, this also determines embedded 4-punctured spheres $\Sigma_j$ in $\IH^3/\Gamma_{\hat{\frz}}$.  As in the manifolds $M_{[0;q,m]}^\circ$, these spheres are triangulated by faces of tetrahedra.  See Figure \ref{fig: SW link} again. 

Each $\Sigma_j$ has a triangular face which is the quotient of a face $[V_\omega, \, V_{-j+1}, \, \infty]$ between $\delta_j$ and $\delta_{j+1}$.  (Here, indices are chosen to match the case $q=3$;  for $q=4$, they need to be adjusted by one.)  Similarly, each $\Sigma_j$ has a face which is the quotient of a triangle $[V_\omega, \, 1+V_{-j}, \, \infty]$ between $\delta'_j$ and $\delta'_{j+1}$.  Since $V_{-j} \to \hat{p}_q$ as $j \to \infty$, these pairs of triangles converge to the square $[ \hat{p}_q, \, V_\omega, \, 1+\hat{p}_q, \, \infty ]$ which lies in the hyperplane over $\R + V_\omega$.  This square intersects $\wh{\Omega}_q$ only in the edge $[V_\omega, \, \infty]$ labeled $\omega$.   Because
\[ U_0 \cdot \Big[ 1+\hat{p}_q, \, \infty \Big] \ = \ \Big[ \hat{p}_q, \, \infty\Big] \qquad \text{and} \qquad (T_\omega)_{\hat{\frz}} \cdot \Big[ \hat{p}_q, \, V_\omega \Big] \ = \ \Big[ 1+\hat{p}_q, \, V_\omega\Big],\] 
the image of the square in $\IH^3/\Gamma_{\hat{\frz}}$ is a totally geodesic 3-punctured sphere in the closure of $\wh{\Omega}_q/\Gamma_{\hat{\frz}}$ which intersects $\wh{\Omega}_q/\Gamma_{\hat{\frz}}$ precisely in the image of the geodesic $[V_\omega, \infty]$.

\begin{figure}
\setlength{\unitlength}{.1in}
\begin{picture}(25,14)
\put(0,0) {\includegraphics[width= 2.5in]{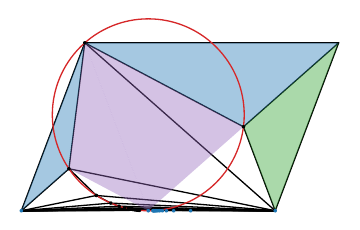}}
\put(10,0){$\scriptstyle \hat{p}$}
\put(3.6,13.2){$\scriptstyle V_{-1}$}
\put(18,6.8){$\scriptstyle V_0$}
\put(5.6,4.5){$\scriptstyle V_{-2}$}
\end{picture}
\caption{\index{fundamental domain!specialized}\index{geometric root!discriminant polynomial} The black lines indicate a fundamental domain for $\wh{\Omega}_3$.  The tetrahedra $\delta_j$ are shown as usual, but the tetrahedra $\delta'_j$ have been replaced with their horizontal translates $U_0 \cdot \delta'_j$.  The blue faces are paired by $U_0 S_0$. As always, the triangle $[0,1,\infty]$ is paired with the finite face (green) of $\delta'_1$ by $S_0$.  The quotient of the purple square by $\Gamma_{\hat{\frz}}$ is a totally geodesic 3-punctured sphere on the boundary of the closure of $\wh{\Omega}_q/\Gamma_{\hat{\frz}}$.}   \label{fig: othersquare} 
\end{figure}

The other pairs of triangles in the triangulated spheres $\Sigma_j$'s also limit to a 3-punctured sphere in the closure of $\wh{\Omega}_q$, but this is slightly more difficult to see.  Figure \ref{fig: othersquare}, which illustrates the $q=3$ case, may be helpful in the following discussion.  

First, recall that $R^2_\omega=S_\omega$ and $V_j  =  (R_\omega)^j_{\hat{\frz}} \, V_0$.    This means that any edge $[V_{-j}, \, V_{-j+1}]$, with $j\geq 1$, can be moved by a nonnegative power of $(S_\omega)_{\hat{\frz}}$ to either $[V_{-1}, \, V_0]$ or $[V_{-2}, \, V_{-1}]$.  If $j$ is large, the power will also be large.  This means that each tetrahedron in $\{ \delta_j \}_2^\infty \cup \{ U_0 \cdot \delta'_j \}_2^\infty$ can be moved by a power of $(S_\omega)_{\hat{\frz}}$ so that one of its edges (labeled $\omega$) coincides with $[V_{-1}, \, V_0]$ or $[V_{-2}, \, V_{-1}]$.  Again, if $j$ is large, then so must be the power.

The edge of $\delta_j$ (or $U_0 \cdot \delta'_j$) opposite to this edge is also labeled $\omega$.  These edges are $[\hat{p}_q+1/2, \infty ]$ in the $\delta_j$ case and $[\hat{p}_q-1/2, \infty ]$ in the case of $U_0 \cdot \delta'_j$.  In either case, both ends of the edge lie on the extended line $(\R+V_\omega) \cup \{ \infty \}$ which is invariant under $(S_\omega)_{\hat{\frz}}$.  In particular, if $z$ is an endpoint of one of these edges, then $\{ (S_\omega)^j_{\hat{\frz}} (z)\}_0^\infty$ accumulates to $\hat{p}_q$, the fixed point for $(S_\omega)_{\hat{\frz}}$.  (The blue points in Figure \ref{fig: othersquare} are the points of this accumulating set.)  The conclusion is that the image of the totally geodesic square $[\hat{p}_q, V_0, \, V_{-1}, \, V_{-2}]$ is in the closure of $\wh{\Omega}_q/\Gamma_{\hat{\frz}}$.  The isometry $(S_\omega)_{\hat{\frz}}$ pairs the edges $[\hat{p}_q, \, V_{-2}]$ and $[\hat{p}_q, \, V_0]$ fixing $\hat{p}_q$.  Likewise, $(W_0)_{\hat{\frz}}$ pairs $[V_{-2}, \, V_{-1}]$ and $[V_0, \, V_{-1}]$ fixing $V_{-1}$.   

Some selected results from the discussion above are summarized in the theorem below.

\begin{thm} \label{thm: core}
The closure of $\wh{\Omega}_q/\Gamma_{\hat{\frz}}$ in $\IH^3/\Gamma_{\hat{\frz}}$ is the convex core for $\IH^3/\Gamma_{\hat{\frz}}$.  Its boundary is the union of a pair of totally geodesic 3-punctured spheres.  One of these is the quotient of the totally geodesic ideal square $ [ \hat{p}_q, \, V_\omega, \, 1+\hat{p}_q, \, \infty ]$ and the other is the quotient of the square $[\hat{p}_q, \, V_0, \, V_{-1}, \, V_{-2}]$.  The manifold $\IH^3/\Gamma_{\hat{\frz}}$ has three annular cusps corresponding to the parabolic elements, $U_0$, $(W_0)_{\hat{\frz}}$, and $(S_\omega)_{\hat{\frz}}$.
\end{thm} 
\index{convex core}

\begin{remark} \label{rem: core}
Theorem~\ref{thm: core} is similar to the rational pinching phenomena seen in Theorem 10.1 of \cite{Gue09}, which describes geometric triangulations of convex cores of punctured torus groups.  As seen there, rational pleating curves may be pinched so that components of the convex core boundary become $3$-punctured spheres.  The examples in Section~\ref{sec: Heck} below suggest a related picture for the cone manifold deformations along extended pleating rays.
\end{remark}

\subsubsection{Geometric limits.}
\index{geometric limit}

\begin{figure}
\setlength{\unitlength}{.1in}
\begin{picture}(30,8)
\put(0,0) {\includegraphics[width=3in]{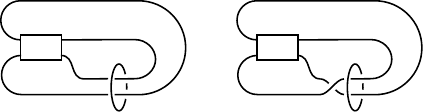} }
\put(2.6,4.4){$\scriptstyle q$}
\put(19.4,4.4){$\scriptstyle q$}
\end{picture}
   \caption{The geometric limits $N_q$ and $N_q'$ for the link complements $M_{[0;q,m]}$.  See Theorem \ref{thm: augment example}.}  \label{fig: augment short}
\end{figure}

First, using the notation from above, consider the sequence $\{ M_{[0;q,2m]}\}$ of hyperbolic knot complements.  As in the proof of Theorem \ref{thm: convergence II}, take $N_q$ to be the complement of the link obtained by augmenting the last twist region in the diagram $D(q,2)$ with a crossing circle, see Figure \ref{fig: augment short}.  Then $N_q$ is hyperbolic and $M_{[0;q,2m]}$ is obtained from $N_q$ by a hyperbolic Dehn filling on the crossing circle cusp.  Using the version of Thurston's hyperbolic Dehn filling theorem stated on page 214 of \cite{Mar}, the groups $\Gamma_{\frz([0;q,2m])}$ converge geometrically to their envelope $\Gamma_\infty$, which is a nonelementary Kleinian group with $\IH^3/\Gamma_\infty = N_q$.

The manifold $N_q$ can also be obtained by correctly identifying the 3-punctured spheres on the boundary of the hyperbolic manifold given in Theorem \ref{thm: core} as the closure of the quotient of  $\wh{\Omega}_q$ by $\Gamma_{\hat{\frz}}$.  In particular, there is an isometry $\mu_q \in \pslC$ which identifies the $\Gamma_{\hat{\frz}}$-orbits of the ideal squares listed in Theorem \ref{thm: core} and such that $\Gamma_\infty = \left\langle \mu_q, \, \Gamma_{\hat{\frz}} \right\rangle$.   The isometry $\mu_q$ can be chosen to be a parabolic fixing $\hat{p}_q$ and taking $V_{-1}$ to $\infty$.  This determines 
\[ \mu_q \ =\ \begin{bmatrix} 2 & -\hat{p}_q \\ 1/\hat{p}_q & 0 \end{bmatrix}\]
and 
\[ \mu_q(V_0) \ = \ V_\omega \qquad \text{and} \qquad \mu_q(V_{-2}) \ = \ -1+V_\omega.\]
Thus, $\mu_q$ acts on $\IH^3$ identifying the lifts of the two 3-punctured spheres.  In particular, since $V_{-1}$ and $\infty$ are ideal points for two
distinct cusps of $\IH^3/\Gamma_{\hat{\frz}}$, this gluing identifies those cusps, as expected from the augmented-link description of $N_q$.

Now, take $N'_q$ to be the complement of the link obtained by augmenting the last twist region in the diagram $D(q,3)$ with a crossing circle.  Then, as above, the groups $\Gamma_{\frz([0;q,2m+1])}$ converge geometrically to a nonelementary Kleinian group $\Gamma'_{\infty}$ with $\IH^3/\Gamma'_{\infty}=N'_q$.  The group $\Gamma'_\infty$ is generated by $\Gamma_{\hat{\frz}}$ together with an isometry $\eta_q$ which differs from $\mu_q$ by a lift of an isometry of one of the 3-punctured spheres which fixes the cusp for $\hat{p}_q$ and interchanges the other two.  Define
\[ \tau_q \ = \ \begin{bmatrix} 1-2\hat{p}_q & 2\hat{p}_q^2 \\ -2 & 1+2\hat{p}_q \end{bmatrix}.\]
Then $\tau_q$ is a parabolic isometry that fixes $\hat{p}_q$ and satisfies
\[ \tau_q \begin{pmatrix} \infty \\ V_\omega \\ -1+V_\omega  \end{pmatrix} \ = \ \begin{pmatrix} -1+V_\omega \\ \infty \\ -3/4+V_\omega  \end{pmatrix}.\]
This is a lift of the desired isometry and 
\[\eta_q \ = \ \tau_q  \mu_q \ = \ \begin{bmatrix} 2\big(1-\hat{p}_q\big) & \hat{p}_q\big( -1+2\hat{p}_q \big) \\ -2+1/\hat{p}_q & 2\hat{p}_q \end{bmatrix}.\]

The results of this discussion are summarized in the next theorem.

\begin{thm} \label{thm: augment example}
Let $N_q$ and $N'_q$ be the hyperbolic link complements obtained by augmenting the diagrams $D(q,2)$ and $D(q,3)$ at the last twist region.  With $\mu_q$ and $\eta_q$ in $\pslC$ as above,
\[ N_q \ =\  \IH^3/\left\langle \mu_q,\Gamma_{\hat{\frz}}\right\rangle \AND N'_q \ =\  \IH^3/\left\langle \eta_q,\Gamma_{\hat{\frz}}\right\rangle. \]
Moreover, these groups are the geometric limits of $\Gamma_{\frz([0;q,2m])}$ and $\Gamma_{\frz([0;q,2m+1])}$, respectively.
\index{augmented link}
\end{thm}

Recall that $\alpha_m=[0;3,m,4]$.  The last project for this section is to understand the geometric limits for $\Gamma_{\frz(\alpha_{2m})}$ and $\Gamma_{\frz(\alpha{2m+1}}$.

Take $N$ to be the complement of the hyperbolic link obtained by augmenting the second twist region in the diagram $D(\alpha_2)$ with a crossing circle.
Define $N'$ analogously using $D(\alpha_1)$.  Figure \ref{fig: Augment} shows a pair of diagrams, both for $N$.  Thurston's hyperbolic Dehn filling theorem shows that these are the geometric limits of the even and odd subsequences, respectively.

\begin{figure}
\setlength{\unitlength}{.1in}
\begin{picture}(47,17)
\put(0,0) {\includegraphics[width= 2.2in]{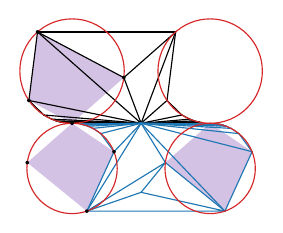}}
\put(27,0) {\includegraphics[width= 2.2in]{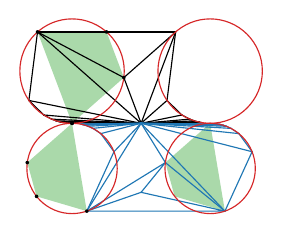}}
\end{picture}
\caption{\index{fundamental domain!specialized} Fundamental domains for the augmented links that arise as the geometric limits of $\Gamma_{\frz(\alpha_{2m})}$ and $\Gamma_{\frz(\alpha_{2m+1})}$. The domains are the same, although the limit groups are different.  For even subscripts, the limit group contains the isometry $\xi_1$, which fixes $\hat{p}_3$ and identifies the two leftmost purple ideal squares.  For odd subscripts, the limit contains $\xi_2$, which also fixes $\hat{p}_3$, but identifies the two leftmost green ideal squares.}    \label{fig: geometric} 
\end{figure}

By Theorem \ref{thm: diagram} (and Figure \ref{fig: SW link II}), there is a topological embedding of the convex core for $\Gamma_{\hat{\frz}(1/3)}$ into both $N$ and $N'$.  Its boundary 3-punctured spheres can be seen in the link diagrams as a pair of 2-punctured disks bounded by the crossing circle and that separate the link complement into two pieces.  The second piece is the image of a topological embedding of the convex core for $\Gamma_{\hat{\frz}(1/4)}$.  By realizing $N$ (or $N'$) by identifying these convex cores correctly by an isometry of their boundaries, a complete hyperbolic structure is induced on the augmented link complement.  Mostow rigidity implies that this agrees with the original geometric structure.  Therefore, the topological embeddings of the convex cores for $\Gamma_{\hat{\frz}(1/3)}$ and $\Gamma_{\hat{\frz}(1/4)}$ are, in fact, isometric embeddings that meet exactly along their boundaries.  So, using Lemma \ref{lem: fat} and Theorem \ref{thm: core}, the complements of the two 3-punctured spheres in $N$ and $N'$ are filled with geometric ideal tetrahedra according to $\wh{\Omega}_q$, $q \in \{ 3,4\}$.  These tetrahedra can be lifted to $\IH^3$ to agree exactly with $\wh{\Omega}_{3}$ and to agree with a copy of $\wh{\Omega}_{4}$ developed across the ideal square 
$[\hat{p}_3, \, V_{1/3}, \, 1+\hat{p}_3, \, \infty]$.  This is shown on the right in Figure \ref{fig: limit} and in both pictures of Figure \ref{fig: geometric}.

Let $\sigma$ be the order-2 elliptic which fixes $\infty$ and the point halfway between $\hat{p}_{3}$ and $\hat{p}_{4}$.  Then
\[ \sigma \ = \ \begin{bmatrix} -i & i\left(\hat{p}_3+\hat{p}_4\right) \\ 0 & i \end{bmatrix}\]
and $\sigma$ moves $\wh{\Omega}_4$ into the position below $\wh{\Omega}_3$ as in Figure \ref{fig: geometric}.  In the figure, $\wh{\Omega}_3$ is drawn in black and the image of $\wh{\Omega}_4$ is blue.  Define
\[ \xi_1 \ = \ \sigma \, \mu_4^{-1} \, \sigma \, \mu_3 \qquad \text{and} \qquad \xi_2 \ = \ \sigma \, \mu_4^{-1} \, \sigma \, \eta_3.\]
Then $\xi_1$ and $\xi_2$ are both parabolic isometries that fix $\hat{p}_3$.  The leftmost purple ideal squares on the left side Figure \ref{fig: geometric} are identified by $\xi_1$ and the leftmost green ideal squares on the right side of the figure are identified by $\xi_2$.  In particular,
\[ \xi_1 \, \begin{pmatrix} \calV(1/2)_{\hat{\frz}(1/3)} \\ 0 \\  -1+\calV(1/4)_{\hat{\frz}(1/3)}\end{pmatrix} \ =\  \sigma \, \begin{pmatrix} -1+\calV(1/5)_{\hat{\frz}(1/4)} \\ 0 \\ \calV(1/3)_{\hat{\frz}(1/4)} \end{pmatrix}\]
and 
\[ \xi_2 \, \begin{pmatrix} \calV(1/2)_{\hat{\frz}(1/3)} \\ 0 \end{pmatrix} \ =\  \sigma \, 
\begin{pmatrix}  0 \\ \calV(1/3)_{\hat{\frz}(1/4)} \end{pmatrix}.\]

The last theorem for this example follows.

\begin{thm} \label{thm: geometric limits}
Let $N$ and $N'$ be the hyperbolic link complements obtained by augmenting the diagrams $D(\alpha_2)$ and $D(\alpha_1)$ with crossing circles at their second twist regions.  Take $\sigma$, $\xi_1$, and $\xi_2$ in $\pslC$ as above.  Then
\[ \Gamma_1 \ = \ \Big\langle  \xi_1, \, \Gamma_{\hat{\frz}(1/3)}, \, \sigma \Gamma_{\hat{\frz}(1/4)} \sigma \Big\rangle \qquad \text{and} \qquad \Gamma_2 \ = \ \Big\langle  \xi_2, \, \Gamma_{\hat{\frz}(1/3)}, \, \sigma \Gamma_{\hat{\frz}(1/4)} \sigma \Big\rangle\]
are Kleinian groups with
\[ N \ = \ \IH^3/ \Gamma_1\qquad \text{and} \qquad N' \  = \ \IH^3/\Gamma_2.\]
Moreover,
\[ \Gamma_{\frz(\alpha_{2m})}\longrightarrow\Gamma_1 \AND \Gamma_{\frz(\alpha_{2m+1})}\longrightarrow\Gamma_2 \]
geometrically.
\end{thm}

\section{Heckoid orbifolds} \label{sec: Heck}

This paper closes by turning to Heckoid orbifolds, which arise naturally in the classification of Kleinian groups generated by two parabolics and play a central role in the geometry underlying the Riley slice.

A natural project is to classify all Kleinian groups generated by a pair of parabolic isometries. Since every such group is, up to conjugation in $\pslC$, inversion of generators, and complex conjugation, equal to $\Gamma_\frz$ for some $\frz$ in the first quadrant of $\C$, this classification problem is highly relevant to this paper.  Recent work in \cite{AHOPSY}, motivated by ideas in \cite{Agol}, completes this classification. Previously, Adams showed in \cite{Adams_2gen} that if $\Gamma_\frz$ is torsion-free and $\IH^3/\Gamma_\frz$ has finite volume, then $\Gamma_\frz$ is a uniformizing group for a hyperbolic $2$-bridge link. Later, as stated in Fact~\ref{fact: Riley}~(3), \cite{OM} established that $\Gamma_\frz$ is a free Kleinian group if and only if $\frz \in \overline{\calR}$. Finally, the recent work mentioned above proved that $\Gamma_\frz$ is a nonfree Kleinian group if and only if it is a hyperbolic 2-bridge link group or a Kleinian Heckoid group (Theorem ~\ref{thm: AHOPSY}). This important result calls for a deeper understanding of these groups and their associated orbifolds.

The Heckoid orbifolds can be viewed as generalizations of 2-bridge link complements and the first part of this section reviews their construction from \cite{AHOPSY}, \cite{CMS}, and \cite{Lee-Sakuma}.  As mentioned in the introduction to this paper, all Kleinian groups are considered to be nonelementary.

First, for $\alpha \in \Q$ and $n \in \frac12 \Z$ with $n \geq 3/2$, there is a cone manifold $M(\alpha)_n$ obtained from $M_\alpha$ by assigning a cone angle of $2\pi/n$ to the lower core tunnel $\tau_-$ of $M_\alpha$.  When $\alpha$ is hyperbolic, this arc is the lower core tunnel of Definition~\ref{def: tunnel}.  When $\alpha$ is not hyperbolic, the upper and lower core tunnels can instead be described diagrammatically as the standard unknotting arcs in the link diagram $D_\alpha$.
\index{core tunnel}

The standard Klein $4$-group $G$ of symmetries of the $2$-bridge link $L_\alpha$ acts by homeomorphisms on $M_\alpha$.  When $\alpha$ is hyperbolic, this is the Klein $4$-group described in Definition~\ref{def: Klein} and Section~\ref{sec: equations}.  The involutions in $G$ can be taken to preserve the upper and lower core tunnels.  When this is arranged, they either fix a given tunnel pointwise or invert it.  Consequently, the involutions in $G$ also act by orbifold homeomorphisms on the cone manifold $M(\alpha)_n$.  Regardless of $\alpha$ and $n$, there is a unique nontrivial element of $G$ which fixes $\tau_-$ pointwise.

\begin{dfn} \label{def: Heckoid orbs}
Suppose $\alpha \in \Q$ and $n \in \frac12 \Z$ with $n \geq 3/2$.  If $n \in \Z$, then the cone manifold $M(\alpha)_n$ is called an {\it even Heckoid orbifold} and is denoted as $S(\alpha; n)$.  If $n \notin \Z$, take $\sigma$ to be the nontrivial element of the Klein 4-group $G$ which fixes the points of the lower tunnel $\tau_-$ under its action on $M(\alpha)_n$.  In this case, the quotient cone manifold $M(\alpha)_n/\sigma$ is called an {\it odd Heckoid orbifold} and is denoted as $S(\alpha; n)$.
\end{dfn}
\index{Heckoid orbifold|idxdef}
\index{Heckoid orbifold!even|idxdef}
\index{Heckoid orbifold!odd|idxdef}

Assume $\alpha \in \Q$ and $n \in \frac12 \Z$ with $n \geq 3/2$.  Theorem 2.2 of \cite{Lee-Sakuma} shows that $S(\alpha; n)$ admits the structure of a geometrically finite hyperbolic orbifold and its uniformizing group is generated by a pair of parabolic isometries.  Standard arguments show that this structure is uniquely determined.    

\begin{figure}
   \centering
   \includegraphics[width=2.5in]{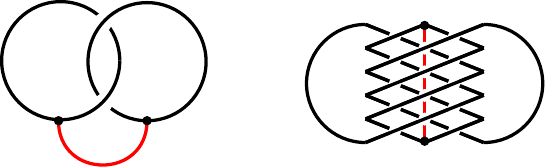} 
   \caption{These are the weighted graphs which determine the even Heckoid orbifolds $S(1/2; n)$ and $S(2/5; n)$.  The red arcs are weighted by $n$ and the black arcs are weighted by $\infty$.}  \label{fig: Even}
\end{figure}

As described in \cite{AHOPSY}, the Heckoid orbifolds $S(\alpha; n)$ are determined by a weighted trivalent graph $\Sigma$ in $S^3$.  If $n \in \Z$, then the graph for the even orbifold $S(\alpha; n)$ is the union of the link $L_\alpha \subset S^3$ together with its lower tunnel.    The lower tunnel is weighted with the integer $n$ and the remaining arcs of $\Sigma$ are weighted by $\infty$.  The arcs weighted by $\infty$ correspond to the cusps of $S(\alpha; n)$ and the arc labeled $n$ indicates its singular set.  Figure ~\ref{fig: Even} shows the weighted graphs for the even orbifolds $S(1/2; n)$ and $S(2/5; n)$.  The black arcs should be labeled $\infty$ and the red arcs labeled $n$.  

\begin{figure}
   \centering
   \includegraphics[width=3in]{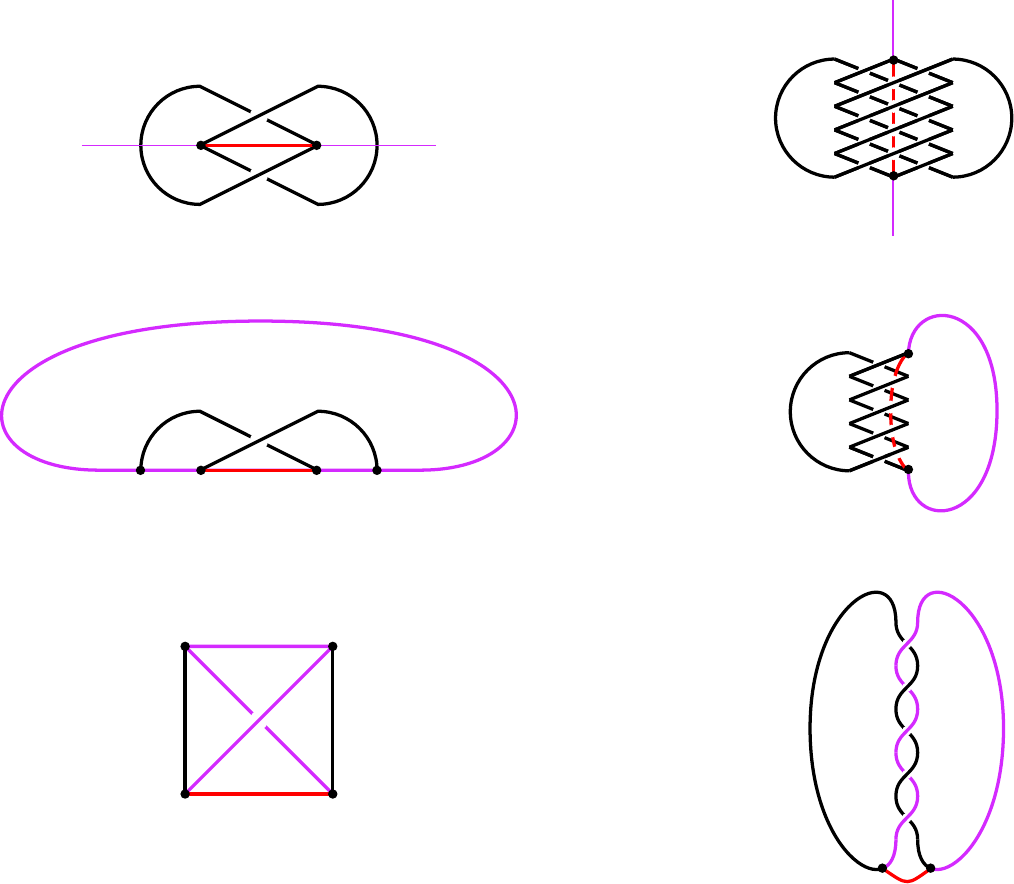} 
   \caption{The bottom images are the weighted graphs which determine the odd Heckoid orbifolds $S(1/2; m/2)$ and $S(2/5; m/2)$.  The red arcs are weighted by $m$, the black arcs are weighted by $\infty$, and the purple arcs are weighted by $2$.  The upper images are meant to show how the Heckoid orbifolds arise as quotients of the cone manifolds $M(1/2)_{m/2}$ and $M(2/5)_{m/2}$.  The bottom two rows differ only by isotopy.}  \label{fig: Odd}
\end{figure}

It is less obvious that the odd Heckoid orbifolds also have weighted graphs made up of 2-bridge links and their lower tunnels.  It turns out that this is the case, although the link is not typically $L_\alpha$.  This is explained in detail in Section 5 of \cite{Lee-Sakuma} and nicely illustrated in Figures 5 and 6 there.  Here, Figure ~\ref{fig: Odd} shows weighted graphs for the odd orbifolds $S(1/2;n)$ and $S(2/5;n)$.  The top images show the links $L_\alpha$ (black), their lower tunnels (red), and the axis of symmetry for $\sigma$ (purple).  After assigning cone angle $2\pi/n$ to the red arc, the top image describes the cone manifold $M(\alpha)_n$.  The lower images show the quotients by $\sigma$.  For the weighted graphs, the black arcs should be labeled $\infty$, the red arc should be labeled with the odd integer $m=2n$, and the purple arcs should be labeled $2$.  Notice that the graph for the odd orbifolds $S(1/2;n)$ is made up of the 2-bridge link $L_{1/1}$ and its lower tunnel.  The graph for the odd orbifolds $S(2/5;n)$ is made up of the 2-bridge link $L_{1/5}$ and its lower tunnel.

\begin{dfn} \label{def: Heckoid groups}
Suppose $\alpha \in \Q$ and $n \in \frac12  \Z$ with $n \geq 3/2$.  As mentioned at the start of this section, a  uniformizing group for the hyperbolic orbifold $S(\alpha; n)$ can be normalized to coincide with a group $\Gamma_\frz$, where $\frz$ lies in the first quadrant of $\C$.  The group $\Gamma_\frz$ is called a {\it Kleinian Heckoid group} and is denoted as $G(\alpha; n)$. 
\end{dfn}
\index{Heckoid group|idxdef}

It seems reasonable to restate the main theorem of \cite{AHOPSY} here.  This is the key classification result: every nonfree 2-parabolic Kleinian group is either a hyperbolic $2$-bridge link group or a Kleinian Heckoid group.

\begin{thm}[{\cite[Theorem~1.1]{AHOPSY}, in the present normalization}] \label{thm: AHOPSY}
A nonfree Kleinian group $\Gamma$ is generated by a pair of parabolic isometries if and only if (up to conjugation and complex conjugation)
\begin{enumerate}
\item There exists a hyperbolic $\alpha$ such that  $\Gamma = \Gamma_\frz$, where $\frz$ is the geometric root for $\calQ(\alpha)$, or
\item there exists $\alpha \in \Q \cap [0,1/2]$ and $n \in \frac12  \Z$ with $n \geq 3/2$, where 
$\Gamma = G(\alpha; n)$.
\end{enumerate}
\end{thm}

Recall Definition ~\ref{def: Riley slice} for the Riley slice $\calR$ and its cusps.  Also, by Theorem ~\ref{thm: DQ uniqueness}, the first quadrant cusps of $\calR$ are exactly the geometric roots for the polynomials $\{ D_\calQ(\alpha) \, | \, \alpha \in \Q \cap [0,1/2]\}$.

\begin{figure}
   \centering
   \includegraphics[width=3in]{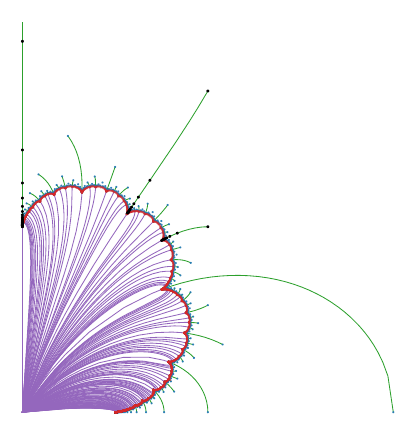} 
   \caption{ Here, Figure ~\ref{fig: pleats0} is shown again.  The red points here are the geometric roots of the discriminant polynomials $D_\calQ(\alpha)$.  They are dense in the Jordan curve $\bound \calR$.  Those blue points in $\C-\overline{\calR}$ which are not real numbers are geometric roots for hyperbolic $\alpha \in \Q$.  The purple curves are the pleating rays of Definition ~\ref{def: KSray} and the green curves are extensions of the pleating rays into the $\C-\calR$ which end at the 2-bridge geometric roots.  Finally, the three sequences of black points on the pleating extensions to the roots for $1/2$, $2/5$, and $3/8$ are the points corresponding to the (even and odd) Heckoid groups $G(1/2;n)$, $G(2/5;n)$, and $G(3/8;n)$ of Definition ~\ref{def: Heckoid groups}.}  \label{fig: pleats}
\end{figure}

\begin{dfn} \label{def: KSray}
Keen and Series show in \cite{KS} that, for $\alpha \in \Q \cap [0,1/2]$,  the set
\[ \left\{ \frz \in \overline{\calR} \, \big| \, \tr (S_\alpha)_\frz \in \R \right\}\]
contains a nonsingular arc ${\mathfrak  p}_\alpha$  from $0$ to the geometric root $\hat{\frz}(\alpha)$ for $D_\calQ(\alpha)$.  This arc is called the {\it pleating ray} for $\alpha$.  These are the purple arcs in Figure ~\ref{fig: pleats}.
\end{dfn}
\index{pleating ray|idxdef}

The preface of \cite{ASWY} describes how the pleating rays should extend as smooth arcs from the cusps of the Riley slice into its exterior $\C-\calR$.  In particular, for every $\alpha\in\Q\cap(0,1/2]$ there should be a pleating ray extension from the geometric root $\hat{\frz}$ of $D_\calQ(\alpha)$ to a
distinguished endpoint $\frz_1$.  If $\alpha$ is hyperbolic, then $\frz_1$ is the geometric root of $\calQ(\alpha)$.  If $\alpha=1/t$, take
\[ \frz_1\ =\ \frac14 \, \sec^2\left(\frac{\pi}{t}\right)\]
where $\frz_1=\infty$ when $t=2$.  The points on the extension for $\alpha$ are parametrized by the trace of $S_\alpha$ (recall Section ~\ref{sec: shapes and traces}) on the interval $[-2,2]$ with $\tr S_\alpha$ taking the value $-2$ at $\hat{\frz}$ and $2$ at $\frz_1$.  Alternatively, the arc is parametrized by a number $\theta_\frz \in [0, 2\pi]$, where 
\begin{align} \label{eq: param} \tr \left( S_\alpha \right)_\frz & =  - 2 \cos \frac{\theta_\frz}{2}. 
\end{align}
As in \cite{ASWY}, the point $\frz$ should correspond to a hyperbolic cone manifold structure on the manifold $M_\alpha$, where the singular set consists only of the lower tunnel and the cone angle at points of the tunnel is $\theta_\frz$.  In particular, if $n \in \Z_{\geq 2}$ and $\theta_\frz=2\pi/n$, then $\Gamma_\frz$ is the even Heckoid group $G(\alpha; n)$ and $\IH^3/\Gamma_\frz$ is the hyperbolic Heckoid orbifold $S(\alpha; n)$.  Alternatively, if $m \geq 3$ is odd, $n=m/2$, and $\theta_\frz=2 \pi/n$, then the associated hyperbolic cone manifold 
is not a hyperbolic orbifold, but its 2-fold quotient is the odd Heckoid orbifold $S(\alpha;n)$, which is hyperbolic.  Here, $\Gamma_\frz$ is the odd Heckoid group $G(\alpha; n)$ and $\IH^3/\Gamma_\frz = S(\alpha;n)$.  Extensions of the pleating rays can be seen in Figure ~\ref{fig: pleats} as green curves.  On each of the three curves for $\alpha \in \{ 1/2, \,  2/5, \,  3/8\}$, the black dots mark the points for the Heckoid groups $G(\alpha;n)$.

\begin{remark} \label{rem: triangle}
The weighted graph for an even Heckoid orbifold $S(0;n)$ is a theta-graph whose arcs are labeled $\infty$, $\infty$, and $n$.  The weighted graph for an odd Heckoid orbifold $S(0;n)$ is a theta-graph whose arcs are labeled $\infty$, $2n$, and $2$.  Hence, the even and odd Heckoid groups $G(0;n)$ are Fuchsian triangle groups.
\end{remark}

The remainder of this section is devoted to computing the geometries of the Heckoid orbifolds $S(\alpha;n)$ for
$\alpha\in\{1/2,2/5\}$ and some small values of $n$.  These examples may be viewed as cone manifold analogues of Gu\'eritaud's convex core constructions mentioned in Remark~\ref{rem: core}.  They also provide specific examples which illustrate some of the more general results of \cite{CEP}.

\begin{example} \label{ex: Heckoid25}
Consider the even and odd Heckoid orbifolds $S(2/5;n)$.  In the rank-2 free group $\langle f, g \rangle$, let 
\[ w \ =\ fgfGFgfgFG \]
where $F=f^{-1}$ and $G=g^{-1}$.   As defined, the isomorphism $\langle f, g \rangle \to \Gamma$ given by 
\[f \mapsto U_0^{-1} \AND g \mapsto W_0\]
takes $w$ to $S_{2/5}$.  (The element $w \in \langle f, g \rangle$ is sometimes referred to as a {\it Farey word} for $2/5$.)  It follows that, if $n \in \Z_{\geq 2}$, then
\[ G(2/5;n) \ = \ \big\langle f, g \, \big| \, w^n \big\rangle.\]
It is shown in \cite{CMS}, that if $k \in \Z_{\geq 1}$, $m=2k+1$, and $n=m/2$, then
\[  G(2/5;n) \ = \ \Big\langle f, g \, \Big| \, w^m, \, \left( fw^k \right)^2, \, Gzfz \Big\rangle\]
where $z=FGw^kfgf$. In these presentations, the elements $f$ and $g$ are both peripheral for the orbifolds $S(2/5;n)$.  As mentioned above, \cite{Lee-Sakuma} shows that these groups can be represented as points in the exterior of the Riley slice $\C-\calR$.  

Given $n \in \frac12 \Z$ with $n \geq 3/2$, let $\frz(n)$ be the unique point in the first quadrant of $\C-\calR$ for which $\Gamma_{\frz(n)}$ is a Kleinian group and the function 
\[ f \ \mapsto \ \left(U_0\right)_{\frz(n)}^{-1} \qquad g \ \mapsto \ \left( W_0 \right)_{\frz(n)} \]
extends to a group isomorphism $G(2/5;n) \to \Gamma_{\frz(n)}$. (Recall Definition ~\ref{def: UW}.)   Then 
\[w \ \mapsto \ (S_{2/5})_{\frz(n)}.\] 
Also, set $\frz(1)$ to be the geometric root $e^{i\pi/3}$ of $\calQ(2/5)$ and extend the notation for the Heckoid orbifolds by setting
\[ S(2/5;1) \ =\ M_{2/5}.\]

Let ${\mathfrak h}_{2/5}$ be the extension of the pleating ray ${\mathfrak p}_{2/5}$ into $\C-\calR$ from the geometric root $\hat{\frz}$ for $D_\calQ(2/5)$ to the geometric root $e^{ i\pi/3}$ for $\calQ(2/5)$.  As mentioned before, some early points $\frz(n)$ can be seen in Figure ~\ref{fig: pleats};  these are the black points beginning with the geometric root $\frz(1)=e^{\pi i/3}$ accumulating along ${\mathfrak  h}_{2/5}$ to $\hat{\frz}$.

Define $p(x) = 2  x^{5} - x^{4} + 2  x^{3} - 3  x^{2} + 2  x - 1$.  Then the trace formula ~\ref{eq: trS one} of Section ~\ref{sec: shapes and traces} gives
\[ \tr S_{2/5} \ = \ x^{-5} \, p(x)\]
and $(S_{2/5})_\frz$ has order-2 in $\pslC$ if and only if $p(\frz)=0$.

Define 
\[ X \ = \ \begin{pmatrix} 0 & 1 \\ -x^{10} & p(x) \end{pmatrix} \]
and let $\{ f(j)\}_0^\infty$, $\{ f_e(j)\}_0^\infty$, and $\{ f_o(j)\}_0^\infty$ be the $X$-recursive sequences with initial conditions
\begin{align*}
f(0)&=0 & f(1)&=1 \\
f_e(0)&=-x^{-5} & f_e(1)&=1 \\
f_o(0)&=x^{-5} & f_o(1)&=1 \\
\end{align*}

The next lemma follows from standard arguments.

\begin{lemma} \label{lem: order polys 2/5}
Suppose $\frz$ lies in the first quadrant of $\C$ and write $S_\frz = (S_{2/5})_\frz$.  Take $k \in \Z_{\geq 1}$.
\begin{enumerate}
\item $S_\frz^k$ is trivial if and only if $f(k)_\frz = 0$.
\item $U_0^{-1}S_\frz^k$ has order two if and only if $f_o(k+1)_\frz=0$.
\item \[ f(2k+1) \ =\ f_o(k+1) \, f_e(k+1)\]
\end{enumerate}
\end{lemma}

This lemma has the following consequences.
\begin{enumerate}
\item If $n \in \Z_{\geq 2}$ then the Heckoid point $\frz(n)$ satisfies $f(n)$.  If $n=2k+1$ for $k \in \Z_{\geq 1}$, then $\frz(n)$ satisfies $f_e(k+1)$.
\item If $k \in \Z_{\geq 1}$, $m=2k+1$, and $n=m/2$, then $\frz(n)$ satisfies $f_o(k+1)$.
\item Given $k \in \Z_{\geq 1}$, there is a polynomial $\hat{f}(k) \in \Z[x]$ such that $S_\frz$ has order $k$ if and only if $\hat{f}(k)_\frz=0$.  This polynomial can be obtained by dividing $f(k)$ by polynomials $f(\ell)$ where $\ell$ divides $k$.  
\end{enumerate}

These facts make it possible to easily find polynomials which must be satisfied by the numbers $\frz(n)$.  At least for the first twelve values of $n$, these polynomials are irreducible over $\Z$.

\begin{itemize}
\item $n=1$.  The knot complement $S(2/5;1)=M_{2/5}$ has cone angle $2\pi$ about the lower tunnel.  The number $\frz(1)$ has been defined to be the geometric root $e^{i \pi/3}$ for $\calQ(2/5)$.

\item $n=3/2$.  Here, if $k=1$, then $m=2k+1$ and $n=m/2$.  The relations for the odd Heckoid group $G(2/5;3/2)$ are $w^3$ and $(fw)^2$.  It must be true, then, that $f_o(2)_{\frz(n)}=0$.  A computation shows
\[ f_o(2) \ =\ (x^{4} + 2  x^{2} - x + 1)(x - 1).\]
Since $G(2/5;3/2)$ is not Fuchsian, $\frz(3/2)\neq 1$ and so, it must satisfy the quartic factor.  This factor has only one root in quadrant one.  An approximation is
\[ \frz(3/2) \ \sim \ 0.343814597201477 + 0.6253578117826761\, i.\]

\item $n=2$.  Recall that
\[ f(2) \ = \ 2  x^{5} - x^{4} + 2  x^{3} - 3  x^{2} + 2  x - 1\]
which is satisfied by $\frz$ if and only if $S_\frz$ has order two.  This polynomial has exactly one root which lies in the open first quadrant.  An approximation is
\[ \frz(2) \ \sim \ 0.3127418007960793 + 0.5800061630078577\, i.\]

\item $n=5/2$. Take $k=2$ and $m=2k+1$ so that $n=m/2$.  The number $\frz(5/2)$ must satisfy $f_o(3)$ which has degree 10 and is irreducible over $\Z$.  This polynomial has exactly three roots in the open first quadrant.  Only one of these lies on the ray ${\mathfrak  h}_{2/5}$ between the geometric root $\hat{\frz}$ for $D_\calQ(2/5)$ and $\frz(2)$.  An approximation is
\[ \frz(5/2) \ \sim \ 0.30041815765580926 + 0.56194811320874\, i.\]

\item $n=3$.  Since $f(3)=f_o(2) \cdot f_e(2)$ and $f_o(2)$ is not satisfied by $\frz(3)$, it must be true that $f_e(2)_{\frz(3)}=0$.
In this case, $f_e(2)$ has degree five and exactly one root in the open first quadrant
\[ \frz(3) \ \sim \ 0.29416839128790606 + 0.5527615885489652\, i.\]

\item $n=7/2$.  Take $k=3$ and $m=2k+1$ so that $n=m/2$.  
The polynomial $f_o(4)$ is irreducible over $\Z$ and has degree $15$.  As in the $n=5/2$ case, it has exactly one root in the right position on ${\mathfrak h}_{2/5}$ and
\[ \frz(7/2) \ \sim \ 0.29053749887933267 + 0.5474138904196743\, i.\]

\item $n=4$.  Here $f(4)=f(2) \cdot \hat{f}(4)$.  Since $S_{\frz(4)}$ has order four, $\frz(4)$ must satisfy $\hat{f}(4)$.  This polynomial is irreducible over $\Z$ and has degree 10.  It has only one root in the correct position
\[ \frz(4) \ \sim \ 0.28823341902406935 + 0.5440159268950905\, i.\]

\item $n=9/2$.  Take $k=4$ and $m=2k+1$ so that $n=m/2$.  
The number $\frz(9/2)$ must satisfy $f_o(5)$ but not its factor $f_o(2)$.   The other factor is irreducible over $\Z$ and of degree 15.  It has one root in the correct position
\[ \frz(9/2) \ \sim \ 0.2866768619956519 + 0.5417183268012936\, i.\]
\end{itemize}

Similar arguments provide
\begin{align*}
\frz(5) &\sim 0.2855747510036132 + 0.5400904905296778 \, i \\
\frz(11/2) &\sim 0.2847652796426727 + 0.5388943295171682 \, i\\
\frz(6) &\sim  0.2841529699387652 + 0.5379891963994432 \, i\\
\frz(13/2) &\sim 0.28367844077452586 + 0.537287542979274 \, i\\
\frz(7) &\sim 0.2833031488608303 + 0.5367325054434339 \, i
\end{align*}
Also,
\begin{itemize}
\item $\frz(5)$ satisfies $f_e(3)$,
\item $\frz(11/2)$ satisfies $f_o(6)$,
\item $\frz(6)$ satisfies $\hat{f}(6) = \frac{f(6)}{f(2) f(3)}$, and
\item $\frz(13/2)$ satisfies $f_o(7)$.
\end{itemize}

\begin{remark} \label{rem: unique heckoid roots} As shown in \cite{CMS}, results of \cite{Knapp} imply that 
\[\tr S_{\frz(n)}\ =\ -2 \cos \frac{\pi}{n}.\] 
(This also follows from the parametrization ~\ref{eq: param}.)  By counting the degree of the polynomial used to find $\frz(n)$, it is possible to see that this polynomial has exactly one root in $\IH^2$ which meets this condition.
\end{remark}

In Definition ~\ref{def: fat}, the thick region associated to a nonhorizontal edge $e$ of the Stern–-Brocot diagram is the open connected subset $\calU_e \subset \IH^2$ where $\calZ(e)$ has positive imaginary part and which contains the geometric root of $e^- \oplus^2 e^+$.  Let $e_1=[0,1/2]$ and $e_2=[1/3,1/2]$.  Then, as discussed in Example ~\ref{ex: real loci}, $\calU_{e_1}  =  \IH^2$ and $\calU_{e_2}$ is the intersection of $\IH^2$ with the open unit disk.  Except for the point $e^{ i\pi/3}$, the ray ${\mathfrak h}_{2/5}$ is a subset of $\calU_{e_2}$.  So, recalling Definition ~\ref{def: hyperbolic structures},
\[ \Big\{ \frz(m/2) \Big\}_{m=3}^\infty \ \subset \ {\mathfrak h}_{2/5} - \{ \frz(1) \} \ \subset \  \calH\left(M^\circ_{3/7}\right).\] 
This means that, if $n\in \Z$, then, with the geometry determined by $\frz(n)$, there is an isometric embedding $M_{3/7}^\circ \to S(2/5;n)$.  If $n \notin \Z$, there is an isometric embedding $M_{3/7}^\circ/\sigma \to S(2/5;n)$.   Note that $\Omega(3/7)_{\frz(n)}$ is made up of two pairs of ideal hyperbolic tetrahedra.  Under the specializations for $\frz(n)$ with $n$ up to $7$, these four tetrahedra are shown in blue in Figures ~\ref{fig: heck1}, ~\ref{fig: row1}, ~\ref{fig: row2}, and ~\ref{fig: row3}.

Next, let
\[ \wt{\Sigma}_n\ =\ \wt{\Sigma}(e_2)_{\frz(n)} \qquad \text{where} \qquad e_2=[1/3,1/2].\]
In other words, $\wt{\Sigma}_n$ is the specialization of $\wt{\Sigma}(e_2)$ at $\frz(n)$. Explicitly, $\wt{\Sigma}_n$ is the union of the four ideal triangles
\begin{align*}
t_1&=\big[ 0, \, -1+\calV(1/3), \, \calV(1/2) \big]_{\frz(n)} & 
t_2&=\big[  -1+\calV(1/3), \, \calV(1/2), \, \infty \big]_{\frz(n)} \\
t_3&=\big[  \calV(1/2), \, \calV(2/5), \, \infty \big]_{\frz(n)} &
t_4&=\big[  \calV(1/3), \, \calV(2/5), \, \infty \big]_{\frz(n)}.
\end{align*}
In particular, under the geometry determined by $\frz(n)$, $\wt{\Sigma}_n$ is a lift of the triangulated $4$-punctured sphere $\Sigma_2=\bound M_{3/7}^\circ$ to $\IH^3$.

\begin{figure}
\setlength{\unitlength}{.1in}
\begin{picture}(40,15)
\put(0,0) {\includegraphics[width= 4in]{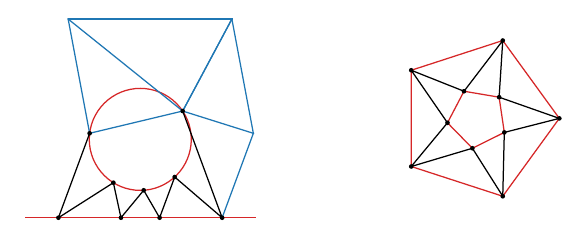}}
\put(1,0){$\scriptstyle -1+\calV(2/5)$}
\put(14.1,0){$\scriptstyle \calV(2/5)$}
\put(0,7.4){$\scriptstyle -1+\calV(1/3)$}
\put(13.2,8.8){$\scriptstyle \calV(1/2)$}
\put(10,0){$\scriptstyle S \cdot \infty$}
\end{picture}
\caption{\index{convex core} \index{fundamental domain!specialized} The domain on the left is $\Omega(3/7)_{\frz(5/2)}$ together with the $5$-drum $\calD_{5/2}$.  Their union projects to the convex core of $S(2/5; 5/2)$.  The faces of the drum lie on the geodesic planes whose boundaries are the red circle and horizontal line.  On the right, an isometry has been applied to the drum to make it easier to see.  The isometry has taken the axis of symmetry for the drum to the geodesic $[0,\infty]$.}   \label{fig: heck1} 
\end{figure}

The main goal for this example is to show that there is an ideal hyperbolic drum which fits nicely against $\Omega(3/7)_{\frz(n)}$ sharing faces of $\wt{\Sigma}_n$ and so that their union projects onto the convex core of $S(2/5;n)$ under the covering map.  The action of $\Gamma_{\frz(n)}$ on $\Omega(3/7)_{\frz(n)}$ is already well understood and it will follow that the action on the drum is also easy to understand.  This will provide a good understanding of the geometry of $S(2/5;n)$ in terms of four ideal tetrahedra and an ideal drum.
\index{fundamental domain!specialized}

\begin{figure}
\setlength{\unitlength}{.1in}
\begin{picture}(50,23)
\put(0,1) {\includegraphics[width= 5in]{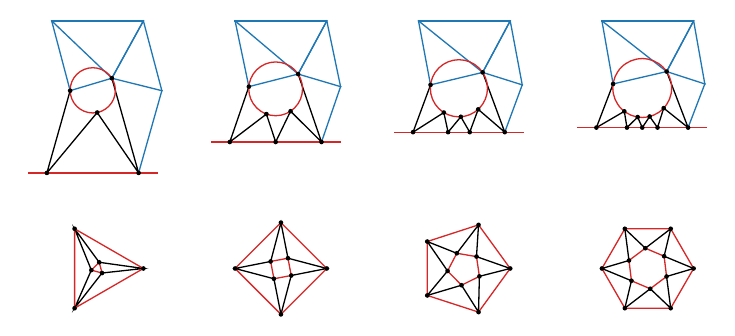}}
\put(5,0){$\scriptstyle n=3/2$}
\put(18,0){$\scriptstyle n=2$}
\put(30,0){$\scriptstyle n=5/2$}
\put(43,0){$\scriptstyle n=3$}
\end{picture}
\caption{\index{fundamental domain!specialized} Similar to Figure ~\ref{fig: heck1}, this shows pictures for the Heckoid orbifolds $S(2/5;n)$ for the given values of $n$.  The upper images show $\Omega(3/7)_{\frz(n)} \cup \calD_n$ and the lower images show images of the drums under isometries that take their axes to $[0,\infty]$. }   \label{fig: row1} 
\end{figure}

\begin{dfns} \label{def: drum}
Given $k \in \Z_{\geq 3}$, let $P$ be a regular, ideal, geodesic  $k$-gon in $\IH^3$ and let $\gamma$ be the unique geodesic perpendicular to $P$ and passing through its center.  Suppose $\phi$ is a loxodromic isometry which acts as a nontrivial translation along $\gamma$.   The convex hull $\calD$ of $P$ and $\phi P$ is called a {\it hyperbolic $k$-drum}.  The $k$-gons $P$ and $\phi P$ are called the {\it heads} of $\calD$ and the complement of the heads in $\bound \calD$ is called the {\it shell}.  The edges of the shell formed where geodesic faces meet are the {\it laces} of $\calD$.

If $\phi$ has no rotation component, the resulting drum is a symmetric square-faced drum.  It is also possible to choose $\phi$ so that the resulting drum is a {\it regular} drum whose dihedral angles are all $\pi/2$.  Both of these especially symmetric types of drums are discussed in Chapter 6 of  \cite{Th_notes}.
\end{dfns}
\index{drum|idxdef} \index{drum!head|idxdef} \index{drum!shell|idxdef} \index{drum!laces|idxdef} 

\begin{figure}
\setlength{\unitlength}{.1in}
\begin{picture}(50,23)
\put(0,1) {\includegraphics[width= 5in]{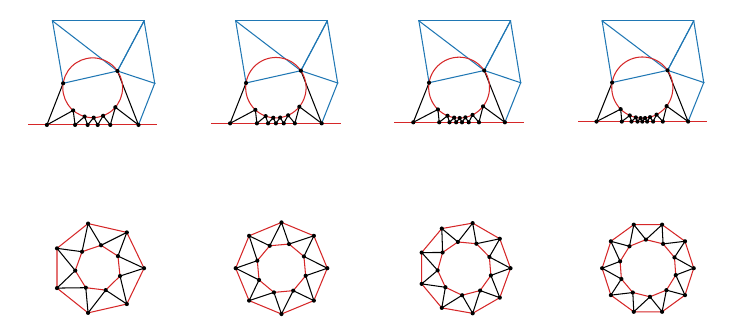}}
\put(5,0){$\scriptstyle n=7/2$}
\put(18,0){$\scriptstyle n=4$}
\put(30,0){$\scriptstyle n=9/2$}
\put(43,0){$\scriptstyle n=5$}
\end{picture}
\caption{\index{fundamental domain!specialized} This shows pictures for the Heckoid orbifolds $S(2/5;n)$ for the given values of $n$.}   \label{fig: row2} 
\end{figure}

To describe the drum needed for $S(2/5;n)$, it is helpful to recall the element $R_{2/5}$ defined in Definition~\ref{def: XYR}.  For $\frz \in \C$, write $R_\frz = (R_{2/5})_\frz$.  By Theorem ~\ref{thm: spiral}, $R_\frz^2 = S_\frz$.  Let $\{\beta_j\}_{\Z}$ be the neighbor cycle for $\omega=2/5$ from Definition~\ref{def: omega cycle}.  Thus, $\beta_0=1/2$ and $\beta_{-1}=1/3$.  Then, again using Theorem ~\ref{thm: spiral}, 
\[R_{2/5}^j \left( \calV(\beta_0) \right) \ = \ \begin{cases} \calV(\beta_j) & \text{if } j\geq 0 \\
-1 + \calV(\beta_j) & \text{if } j<0.\end{cases}\]
By definition of $\frz(n)$, the order of $S_{\frz(n)}$ is $n$ if $n \in \Z$ and is $2n$ if not.  In either case, the order of $R_{\frz(n)}$ is $2n$.   Let $\calD_n$ be the convex hull of the points
\[ \big\langle R_{\frz(n)} \big\rangle \cdot \big\{ \calV(1/2)_{\frz(n)}, \, \infty \big\} \ \subset \ \bound \IH^3.\]
The  $\langle R_{\frz(n)}\rangle$-orbit of $\calV(1/2)_{\frz(n)}$ is a regular ideal geodesic $2n$-gon, and the same holds for the orbit of $\infty$.  The axis $\gamma$ of $S_{\frz(n)}$ passes orthogonally through the center of each polygon, and $\calD_n$ forms a hyperbolic $2n$-drum with $\gamma$ as its axis of symmetry.    \index{neighbor cycle}

\begin{figure}
\setlength{\unitlength}{.1in}
\begin{picture}(50,23)
\put(0,1) {\includegraphics[width= 5in]{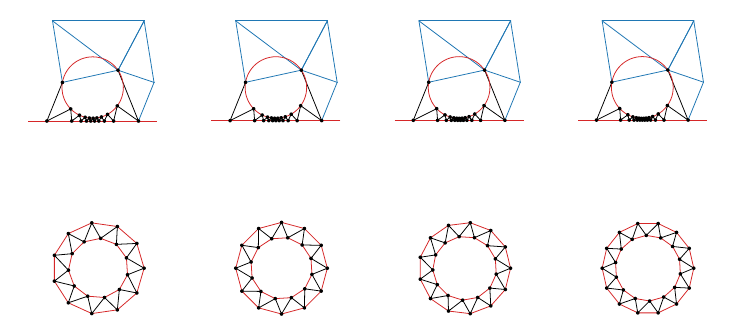}}
\put(5,0){$\scriptstyle n=11/2$}
\put(18,0){$\scriptstyle n=6$}
\put(30,0){$\scriptstyle n=13/2$}
\put(43,0){$\scriptstyle n=7$}
\end{picture}
\caption{\index{fundamental domain!specialized} This shows pictures for the Heckoid orbifolds $S(2/5;n)$ for the given values of $n$.}   \label{fig: row3} 
\end{figure}

The unions $\Omega(3/7)_{\frz(n)} \cup \calD_n$ for $n \in \frac12 \Z$ with $3/2 \leq n \leq 7$ are shown in Figures  ~\ref{fig: heck1}, ~\ref{fig: row1}, ~\ref{fig: row2}, and ~\ref{fig: row3}.  In the figures, the geodesic edges at the intersections of the heads and the shell are shown in red and the laces are black.  In each case, 
\[ \calD_n \, \cap \, \Omega(3/7)_{\frz(n)}\ = \ t_2 \cup t_3 \ \subset \ \wt{\Sigma}_n.\] 
The triangle $t_4$ is taken by $U_0$ to the face $\big[  -1+\calV(1/3), \, -1+\calV(2/5), \, \infty \big]_{\frz(n)}$ of $\calD_n$ and, by Lemma ~\ref{lem: STbeta},  the triangle $t_1$ is taken by $(T_{1/2})_{\frz(n)}$ to the face $\big[\calV(1/2), \,\calV(2/5), \, \calV(3/7) \big]_{\frz(n)}$ of $\calD_n$.  So, because $U_0$ and $(T_{1/2})_{\frz(n)}$ are elements of $\Gamma_{\frz(n)}$, every triangle of $\wt{\Sigma}_n$ shares its $\Gamma_{\frz(n)}$-orbit with a triangular face on the shell of $\calD_n$.  Because $S_{\frz(n)} \in \Gamma_{\frz(n)}$, every triangle in the shell of $\calD_n$ shares its $\Gamma_{\frz(n)}$-orbit with a triangle in $\wt{\Sigma}_n$.

\begin{figure}
   \centering
   \includegraphics[width=3in]{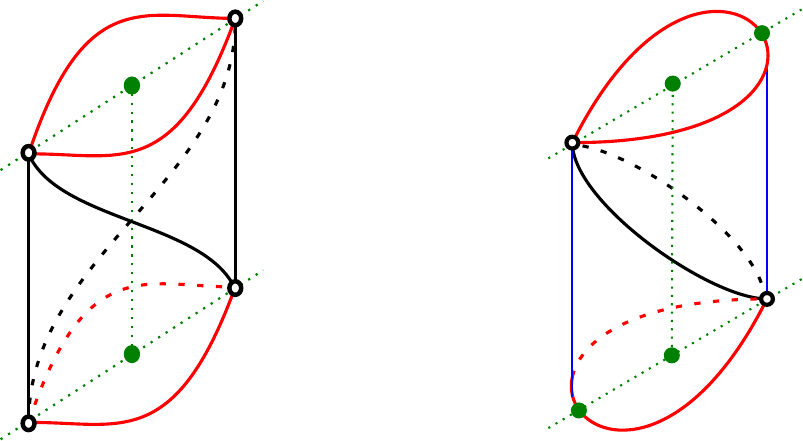} 
   \caption{This shows the orbifold quotients $\calD_n/S_{\frz(n)}$.  The edges are color coded red and black to match Figures ~\ref{fig: heck1}, ~\ref{fig: row1}, ~\ref{fig: row2}, and ~\ref{fig: row3}.  On the left is the case that $n$ is an integer and on the right is the case that $n=m/2$ for an odd integer $m$.  On the right, the annular quotient of the shell of $\calD_n$ is subdivided into a pair of triangles by the two black edges shown.  (The blue lines are not part of the triangulation of this annulus.)}  \label{fig: singulartubes}
\end{figure}

Suppose, for a moment, that $n$ is an integer.  Then $S_{\frz(n)}$ has order $n$ and the orbifold quotient $\calD_n/S_{\frz(n)}$ is shown on the left in Figure ~\ref{fig: singulartubes}.  The four triangles $(T_{1/2})_{\frz(n)}(t_1)$, $t_2$, $t_3$, and $U_0(t_4)$ make up a fundamental domain for the action of $\langle S_{\frz(n)} \rangle$ on the shell of $\calD_n$, so the quotient of the shell is an annulus, with two ideal points on each boundary component, which is triangulated by four ideal triangles.  The heads of $\calD_n$ descend to totally geodesic ideal bigons with order-$n$ singular points at their centers.  Together with the annular quotient of the shell, these bigons bound a ball with an order-$n$ cone singularity along the image of $\gamma$.  Take $P_n$ to be a wedge shaped fundamental domain for $\langle S_{\frz(n)} \rangle$ which contains the four triangles $(T_{1/2})_{\frz(n)}(t_1)$, $t_2$, $t_3$, and $U_0(t_4)$ as well as the axis $\gamma \cap \calD_n$.  Then $\Gamma_{\frz(n)}$ contains the face pairings for the polyhedron $\Omega(3/7)_{\frz(n)} \cup P_n$ that fold it into its image in $S(2/5;n)$ under the covering map. Because the image of this polyhedron is convex, it must be the convex core for $S(2/5;n)$.  Notice that the images of the heads of $\calD_n$ constitute the boundary of the convex core.  These boundary components are totally geodesic $(n, \infty, \infty)$-turnovers.  They can be visualized in Figure ~\ref{fig: singulartubes} by identifying the edges of each bigon across the diagonal dotted green lines.  These turnovers can also be seen as boundaries of small neighborhoods of the vertices in the weighted graph on the right in Figure ~\ref{fig: Even}.  
\index{convex core}
\index{fundamental domain}

At this point, it seems worth mentioning the case $n=1$.  Here, the cone angle around the lower tunnel is $2\pi$ and there is no cone singularity.  The Heckoid orbifold $S(2/5;1)$ coincides with the hyperbolic link complement $M_{2/5}$ and so $\frz(1)$ is the geometric root $e^{\pi i/3}$ for $\calQ(2/5)$.  A calculation shows that $R_{\frz(1)}$ is the order-2 isometry
\[ R_{\frz(1)} \ = \ \begin{bmatrix} -i & i+\sqrt{3} \\ 0 & i \end{bmatrix}.\]
So, $S_{\frz(1)}$ is trivial, as expected.  Also, $R_{\frz(1)}$ fixes both $\calV(1/2)_{\frz(1)}$ and $\infty$.  So, as defined above, the heads of the drum $\calD_1$ are just the points $\calV(1/2)_{\frz(1)}$ and $\infty$.  The drum has collapsed to a geodesic line.  Moreover, 
\begin{align*}
\calV(1/2)_{\frz(1)} & = -1+\calV(1/3)_{\frz(1)} & \calV(2/5)_{\frz(1)} &= \infty.
\end{align*}
Hence, $\Omega(3/7)_{\frz(1)}$ is the same as $\Omega(2/5)_{\frz(1)}$ which, by Theorem~\ref{thm: domain} and Corollary~\ref{cor: holonomy specialization}, is a fundamental domain for $G(2/5;1)=\Gamma_{\frz(1)}$.

Moving on to the odd case, assume that $n=m/2$ where $m \in \Z_{\geq 3}$ is odd.  In this case, both $S_{\frz(n)}$ and $R_{\frz(n)}$ have order $m$.   Here, the triangles $(T_{1/2})_{\frz(n)}(t_1)$ and $t_3$ make up a fundamental domain for the action of $\langle S_{\frz(n)} \rangle$ on the shell of $\calD_n$.  As shown on the right of Figure ~\ref{fig: singulartubes}, the quotient of the shell is an annulus, with one ideal point on each boundary component.  The image of the shell is triangulated by two ideal triangles, the images of $(T_{1/2})_{\frz(n)}(t_1)$ and $t_3$.  The heads of $\calD_n$ descend to totally geodesic ideal monogons with order-$m$ singular points at their centers.  As before, with the annular quotient of the shell, these monogons cobound a ball with an order-$m$ cone singularity along the image of $\gamma$.  Now, let $P_n$ be a wedge shaped fundamental domain for $\langle S_{\frz(n)} \rangle$ which contains the triangles $(T_{1/2})_{\frz(n)}(t_1)$ and $t_3$ as well as the axis $\gamma \cap \calD_n$.  As in the integer case, the face pairings in $\Gamma_{\frz(n)}$ take $\Omega(3/7)_{\frz(n)} \cup P_n$ to the convex core of $S(2/5;n)$.  The totally geodesic boundary of this convex core is the image of the heads of $\calD_n$ and consists of a pair of $(2,m,\infty)$-turnovers.  They can be visualized on the right in Figure ~\ref{fig: singulartubes} by identifying the two half edges on the boundary of the monogon across the diagonal dotted green lines.  They can also be seen on the bottom right in Figure ~\ref{fig: Odd} as boundaries of vertex neighborhoods in the weighted graph pictured there.

\end{example}

\begin{example} \label{ex: Heckoid12}
This example investigates the even and odd Heckoid orbifolds $S(1/2;n)$.  In the rank-2 free group $\langle f, g \rangle$, let $w$ be the commutator $[f,g]$.  As in Example ~\ref{ex: Heckoid25}, if $n \in \Z_{\geq 1}$, then
\[ G(1/2;n) \ = \ \big\langle f, g \, \big| \, w^n \big\rangle.\]
In \cite{CMS} (see also eg. \cite{MM1}), it is shown that if $k \in \Z_{\geq 1}$, $m=2k+1$, and $n=m/2$, then
\[  G(1/2;n) \ = \ \Big\langle f, g \, \Big| \, w^m, \, \left( fw^k \right)^2, \, (Gw^k)^2, \, (Gw^kf)^2 \Big\rangle.\]
Again, $f$ and $g$ are both peripheral in this presentation.  

For $n \in \frac12 \Z$ with $n \geq 3/2$, let $\frz(n)$ be the unique point in the first quadrant of $\C-\calR$ for which $\Gamma_{\frz(n)}$ is a Kleinian group and the function 
\[ f \ \mapsto \ U_0^{-1} \qquad g \ \mapsto \ \left( W_0 \right)_{\frz(n)} \]
extends to a group isomorphism $G(1/2;n) \to \Gamma_{\frz(n)}$.  Here, 
\[w \ \mapsto \ (S_{1/2})_{\frz(n)}.\]
Some points $\frz(n)$ can be seen in Figure ~\ref{fig: pleats} as the black points beginning with $\frz(3/2)=i$ and accumulating along the extension of the pleating ray ${\mathfrak  p}_{1/2}$ to the geometric root $\hat{\frz}=i/2$ for $D_\calQ(1/2)$.  \index{pleating ray}

\begin{figure}
\setlength{\unitlength}{.1in}
\begin{picture}(15,17)
\put(0,0) {\includegraphics[width= 1.5in]{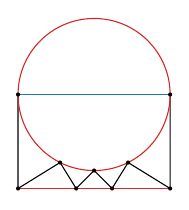}}
\put(0,0){$-\frz$}
\put(12.2,0){$1-\frz$}
\put(0,9){$0$}
\put(14.2,9){$1$}
\put(7.9,0){$\scriptstyle S\cdot \infty$}
\end{picture}
\caption{\index{fundamental domain!specialized} The right angled drum for the orbifold $S(1/2;5/2)$.  The faces of the drum include $(\delta_0)_{\frz(5/2)} = [0,1,\infty]$ and $(\delta'_0)_{\frz(5/2)}=[1,1-\frz(5/2),\infty]$.}   \label{fig: sqheck} 
\end{figure}

Considerations similar to those described in Example ~\ref{ex: Heckoid25} provide analogous recursive polynomials $f(j)$, $f_e(j)$, and $f_o(j)$ that give rise to a set of polynomials in $\Z[x]$ which are satisfied by the numbers $\frz(n)$.  A comparison of these polynomials to those discussed in Remark ~\ref{rk: chebyshev} provides
\[ \frz(n) \ =\ \frac{i}{2} \cdot \sec\left(\frac{\pi}{2n}\right).\]

\begin{figure}
\setlength{\unitlength}{.1in}
\begin{picture}(50,24)
\put(0,1) {\includegraphics[width= 5in]{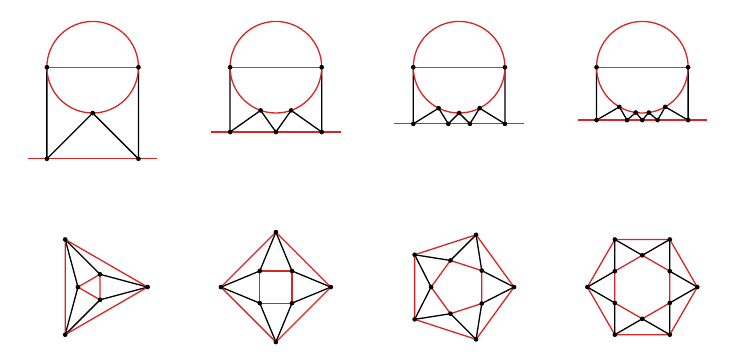}}
\put(5,0){$\scriptstyle n=3/2$}
\put(18,0){$\scriptstyle n=2$}
\put(30,0){$\scriptstyle n=5/2$}
\put(43,0){$\scriptstyle n=3$}
\end{picture}
\caption{\index{fundamental domain!specialized} This shows the pictures for the Heckoid orbifolds $S(1/2;n)$ for the given values of $n$.  Here, the drums are fully symmetric right angled drums and the remaining portion of the cell decomposition is a pair of triangular faces on the drum.}   \label{fig: row_for_12} 
\end{figure}

Recall from Definition~\ref{def: generic tetrahedra} that, for the horizontal edge between $0$ and $1$, the degenerate tetrahedra
\[ (\delta_0)_\frz=[0,1,\infty] \AND (\delta'_0)_\frz=[1,\calV(1/2)_\frz,\infty] \]
are ideal triangles.  These are the two distinguished triangular faces used in the construction.  By Definitions~\ref{def: S and T} and~\ref{def: XYR},
\[ R_{1/2} \ = \ \frac{i}{x} \cdot \begin{bmatrix} 1-x & x \\ 1 & x\end{bmatrix} \qquad \text{and} \qquad S_{1/2} \ =\ -\frac{1}{x^2} \cdot \begin{bmatrix} x^2-x+1 & x \\ 1 & x(x+1)\end{bmatrix}.\]
As before, write $R_\frz = (R_{1/2})_{\frz}$ and $S_\frz = (S_{1/2})_\frz$.  By Theorem ~\ref{thm: spiral} (or by direct calculation), $R_\frz^2 = S_\frz$.  Let $\{\beta_j\}_{\Z}$ be the neighbor cycle for $\omega=1/2$ from Definition~\ref{def: omega cycle}.  Thus, $\beta_0=1$ and $\beta_{-1}=0$.  Then, again by Theorem ~\ref{thm: spiral}, 
\[R_{1/2}^j \left( \calV(\beta_0) \right) \ = \ \begin{cases} \calV(\beta_j) & \text{if } j\geq 0 \\
-1 + \calV(\beta_j) & \text{if } j<0.\end{cases}\]
As in Example ~\ref{ex: Heckoid25}, the order of $S_{\frz(n)}$ is $n$ if $n \in \Z$ and is $2n$ if not.  The order of $R_{\frz(n)}$ is always $2n$.  Let $\calD_n$ be the hyperbolic $2n$-drum whose ideal vertices are the $\langle R_{\frz(n)} \rangle$-orbits of $0$ and $\infty$ and let $\gamma$ be the fixed axis for $R_{\frz(n)}$.  Since 
\[R_{\frz(n)}\cdot 0 \ = \ 1 \qquad \text{and} \qquad R_{\frz(n)}\cdot \infty \ = \ 1-\frz(n)\ =\ \calV(1/2)_{\frz(n)}\]
 the two ideal triangles $(\delta_0)_{\frz(n)}$ and $(\delta'_0)_{\frz(n)}$ are faces of $\calD_n$, as are their images
 \[ \left(W_0 \cdot \delta_0 \right)_{\frz(n)} \ =\ \big[0,\, -x, \, -1+\calV(1/3)\big]_{\frz(n)} \quad \text{and} \quad \left( U_0 \cdot  \delta'_0\right)_{\frz(n)} \ =\ [0, \, -x, \, \infty]_{\frz(n)}\] 
 under the generators $(U_0)_{\frz(n)}$ and $(W_0)_{\frz(n)}$ of $\Gamma_{\frz(n)}$.  So, as in Example ~\ref{ex: Heckoid25}, the convex core of $S(1/2;n)$ is the image, under the covering map, of $\calD_n$.  Moreover, this convex core can be constructed by taking the cyclic quotient of $\calD_n$ by the group $\langle S_{\frz(n)} \rangle$ (see Figure ~\ref{fig: singulartubes}) and then identifying the triangles in the image of the shell using the isometries induced by $(U_0)_{\frz(n)}$ and $(W_0)_{\frz(n)}$. \index{neighbor cycle}
 
 In this case, the drums $\calD_n$ are regular drums.  Recall that this means that all of the dihedral angles of $\calD_n$ are right angles.  Due to the rotational symmetry $R_{\frz(n)}$ of $\calD_n$, this is quick to check.
 \begin{enumerate}
 \item The head spanned by the vertices $R^j_{\frz(n)} \cdot 0$ meets the shell at right angles because the geodesic plane bounded by the circle with radius $1/2$ and center $1/2$ meets the triangle $[0,1,\infty]$ in a right angle.
 \item The other head of $\calD_n$ meets the shell at right angles because the geodesic plane bounded by the line $-\frz(n)+\R$ meets the triangle $[0,  -\frz(n), \infty]$ in a right angle.
 \item The dihedral angles at the laces are $\pi/2$ because the triangle $[0,1,\infty]$ meets the triangles $[0,  -\frz(n), \infty]$ and $[1, 1-\frz(n), \infty]$ in right angles.
 \end{enumerate}

Attributed to Thurston and Milnor from Chapters 6 and 7 of \cite{Th_notes}, the paper \cite{ABVE} provides a formula for the volume of the regular ideal drums.  This, in turn, makes it easy to provide analytic formulas for the convex cores of the Heckoid orbifolds $S(1/2;n)$.  

These volumes are expressed in terms of the Lobachevsky function
\[ L(x) \ = \ - \int_0^x \log |2 \sin t |\, dt.\]
Their formulas give
\[ \text{Vol} \left( \calD_n \right) \ = \ 4n \cdot 
\left( L \left( \frac{\pi}{4} + \frac{\pi}{4n} \right) + L\left( \frac{\pi}{4} - \frac{\pi}{4n} \right) \right)
 \]
where $\calD_n$ is the $2n$-drum defined for $S(1/2;n)$.  Combined with the work above, this shows that, if $n \in \Z$,  the volume of the convex core of $S(1/2;n)$ is
\[ 4 \cdot 
\left( L \left( \frac{\pi}{4} + \frac{\pi}{4n} \right) + L\left( \frac{\pi}{4} - \frac{\pi}{4n} \right) \right)\]
and, if $n=m/2$ where $m=2k+1$ and $k \in \Z_{\geq 1}$, then the volume of the convex core of $S(1/2;n)$ is
\[2 \cdot 
\left( L \left( \frac{\pi}{4} + \frac{\pi}{4n} \right) + L\left( \frac{\pi}{4} - \frac{\pi}{4n} \right) \right) .\]
 
\end{example}

\appendix
\chapter{Polynomials} \label{chap: list}

\begin{center}
\begin{tabular}{ | c | c | c | }
\hline
 $\scriptstyle \alpha$ & $\scriptstyle \calQ(\alpha)$ & $\scriptstyle \calN(\alpha)$ \\
\hline
$\scriptstyle \infty$ &$\scriptstyle 0$ &$\scriptstyle 1$ \\
$\scriptstyle  1$ &$\scriptstyle 1$ &$\scriptstyle 1$ \\
$\scriptstyle  1/2$ &$\scriptstyle 1$ &$\scriptstyle -x+1$ \\
$\scriptstyle  1/3$ &$\scriptstyle -x + 1 $ &$\scriptstyle -2x + 1$ \\
$\scriptstyle 1/4$ &$\scriptstyle  -2x + 1$ &$\scriptstyle x^2 - 3x + 1$ \\
$\scriptstyle  1/5$ &$\scriptstyle x^2 - 3x + 1 $ &$\scriptstyle (3x - 1)(x - 1)$ \\
$\scriptstyle  2/5$ &$\scriptstyle x^2 - x + 1$ &$\scriptstyle (x - 1)^2$ \\
$\scriptstyle  1/6$ &$\scriptstyle (3x - 1)(x - 1)$ &$\scriptstyle -x^3 + 6x^2 - 5x + 1$ \\
$\scriptstyle  1/7$ &$\scriptstyle -x^3 + 6x^2 - 5x + 1$ &$\scriptstyle -(2x^2 - 4x + 1)(2x - 1)$ \\
$\scriptstyle  2/7$ &$\scriptstyle x^3 + 2x^2 - 3x + 1$ &$\scriptstyle (2x - 1)^2$ \\
$\scriptstyle  3/7$ &$\scriptstyle -x^3 + 2x^2 - x + 1$ &$\scriptstyle -2x^3 + 2x^2 - 2x + 1$ \\
$\scriptstyle  1/8$ &$\scriptstyle -(2x^2 - 4x + 1)(2x - 1)$ &$\scriptstyle (x^3 - 9x^2 + 6x - 1)(x - 1)$ \\
$\scriptstyle  3/8$ &$\scriptstyle 2x^2 - 2x + 1$ &$\scriptstyle -(x^3 + x^2 - 2x + 1)(x - 1)$ \\
$\scriptstyle  1/9$ &$\scriptstyle (x^3 - 9x^2 + 6x - 1)(x - 1)$ &$\scriptstyle (5x^2 - 5x + 1)(x^2 - 3x + 1)$ \\
$\scriptstyle  2/9$ &$\scriptstyle x^4 - 2x^3 + 7x^2 - 5x + 1$ &$\scriptstyle (x^2 - 3x + 1)^2$ \\
$\scriptstyle  4/9$ &$\scriptstyle x^4 - 2x^3 + 3x^2 - x + 1$ &$\scriptstyle x^4 - 4x^3 + 3x^2 - 2x + 1$ \\
$\scriptstyle  1/10$ &$\scriptstyle (5x^2 - 5x + 1)(x^2 - 3x + 1)$ &$\scriptstyle -x^5 + 15x^4 - 35x^3 + 28x^2 - 9x + 1$ \\
$\scriptstyle  3/10$ &$\scriptstyle -(3x^2 - 3x + 1)(x^2 + x - 1)$ &$\scriptstyle x^5 - 3x^4 - 3x^3 + 8x^2 - 5x + 1$ \\
$\scriptstyle  1/11$ &$\scriptstyle -x^5 + 15x^4 - 35x^3 + 28x^2 - 9x + 1$ &$\scriptstyle -(x^2 - 4x + 1)(3x - 1)(2x - 1)(x - 1)$ \\
$\scriptstyle  2/11$ &$\scriptstyle x^5 + 3x^4 - 13x^3 + 16x^2 - 7x + 1$ &$\scriptstyle (3x - 1)^2(x - 1)^2$ \\
$\scriptstyle  3/11$ &$\scriptstyle -x^5 - x^4 - 3x^3 + 8x^2 - 5x + 1$ &$\scriptstyle -(x^4 + 4x^2 - 4x + 1)(2x - 1)$ \\
$\scriptstyle  4/11$ &$\scriptstyle x^5 - x^4 - x^3 + 4x^2 - 3x + 1$ &$\scriptstyle (2x^3 + x^2 - 2x + 1)(x - 1)^2$ \\
$\scriptstyle  5/11$ &$\scriptstyle -x^5 + 3x^4 - 3x^3 + 4x^2 - x + 1$ &$\scriptstyle -2x^5 + 3x^4 - 6x^3 + 4x^2 - 2x + 1$ \\
$\scriptstyle  1/12$ &$\scriptstyle -(x^2 - 4x + 1)(3x - 1)(2x - 1)(x - 1)$ &$\scriptstyle x^6 - 21x^5 + 70x^4 - 84x^3 + 45x^2 - 11x + 1$ \\
$\scriptstyle  5/12$ &$\scriptstyle -(2x^2 - x + 1)(x^2 + 1)(x - 1)$ &$\scriptstyle x^6 - 3x^5 + 4x^4 - 6x^3 + 5x^2 - 3x + 1$ \\
$\scriptstyle  1/13$ &$\scriptstyle x^6 - 21x^5 + 70x^4 - 84x^3 + 45x^2 - 11x + 1$ &$\scriptstyle (7x^3 - 14x^2 + 7x - 1)(x^3 - 6x^2 + 5x - 1)$ \\
$\scriptstyle  2/13$ &$\scriptstyle x^6 - 3x^5 + 22x^4 - 40x^3 + 29x^2 - 9x + 1$ &$\scriptstyle (x^3 - 6x^2 + 5x - 1)^2$ \\
$\scriptstyle  3/13$ &$\scriptstyle x^6 - 5x^5 + 6x^4 - 16x^3 + 17x^2 - 7x + 1$ &$\scriptstyle 3x^6 - 6x^5 + 14x^4 - 28x^3 + 23x^2 - 8x + 1$ \\
$\scriptstyle  4/13$ &$\scriptstyle x^6 + 5x^5 - 6x^4 - 4x^3 + 9x^2 - 5x + 1$ &$\scriptstyle -x^6 + 8x^5 - 4x^4 - 10x^3 + 13x^2 - 6x + 1$ \\
$\scriptstyle  5/13$ &$\scriptstyle x^6 - x^5 + 2x^4 - 4x^3 + 5x^2 - 3x + 1$ &$\scriptstyle x^6 + 2x^4 - 6x^3 + 7x^2 - 4x + 1$ \\
$\scriptstyle  6/13$ &$\scriptstyle x^6 - 3x^5 + 6x^4 - 4x^3 + 5x^2 - x + 1$ &$\scriptstyle x^6 - 6x^5 + 6x^4 - 8x^3 + 5x^2 - 2x + 1$ \\
$\scriptstyle  1/14$ &$\scriptstyle (7x^3 - 14x^2 + 7x - 1)(x^3 - 6x^2 + 5x - 1)$ &$\scriptstyle -(x^4 - 24x^3 + 26x^2 - 9x + 1)(x^2 - 3x + 1)(x - 1)$ \\
$\scriptstyle  3/14$ &$\scriptstyle -(x^3 + 6x^2 - 5x + 1)(x^3 - 2x^2 + 3x - 1)$ &$\scriptstyle (x^5 + x^4 - 6x^3 + 11x^2 - 6x + 1)(x^2 - 3x + 1)$ \\
$\scriptstyle  5/14$ &$\scriptstyle -(x^3 - 4x^2 + 3x - 1)(x^3 - x + 1)$ &$\scriptstyle -(x^6 + 3x^5 - 5x^4 - x^3 + 6x^2 - 4x + 1)(x - 1)$ \\
$\scriptstyle  1/15$ &$\scriptstyle -(x^4 - 24x^3 + 26x^2 - 9x + 1)(x^2 - 3x + 1)(x - 1)$ &$\scriptstyle -(2x^4 - 16x^3 + 20x^2 - 8x + 1)(2x^2 - 4x + 1)(2x - 1)$ \\
$\scriptstyle  2/15$ &$\scriptstyle x^7 + 4x^6 - 34x^5 + 86x^4 - 91x^3 + 46x^2 - 11x + 1$ &$\scriptstyle (2x^2 - 4x + 1)^2(2x - 1)^2$ \\
$\scriptstyle  4/15$ &$\scriptstyle (x^4 + 2x^2 - 3x + 1)(x^3 + 4x^2 - 4x + 1)$ &$\scriptstyle (2x^4 + 4x^2 - 4x + 1)(2x - 1)^2$ \\ 
\hline
\end{tabular}
\end{center}

\backmatter
\chapter*{List of symbols}

\noindent
The following list records some recurring notation used throughout the paper.

\bigskip
\begingroup
\renewcommand{\arraystretch}{1.08}
\setlength{\LTpre}{0pt}
\setlength{\LTpost}{0pt}
\begin{longtable}{@{}p{0.19\textwidth}p{0.50\textwidth}p{0.25\textwidth}@{}}
\textbf{Symbol} & \textbf{Meaning} & \textbf{First use}\\
\hline
\endfirsthead
\textbf{Symbol} & \textbf{Meaning} & \textbf{First use}\\
\hline
\endhead

\\

\multicolumn{3}{@{}l}{\textit{Basic sets, groups, and conventions}}\\[2pt]
$\C^\ast$ & Nonzero complex numbers. & Ch.~\ref{chap: intro}\\
$\Q_0$ & $\Q\cap[0,1]$. & Sec.~\ref{sec: SBD}\\
$\widehat{\Q}_0$ & $\Q_0\cup\{\infty\}$. & Sec.~\ref{sec: SBD}\\
$[v_1,\ldots,v_k]$ & Geodesic convex hull. & Def.~\ref{def: parameters}\\

\\[-2pt]
\multicolumn{3}{@{}l}{\textit{Farey and Stern--Brocot notation}}\\[2pt]
$\oplus$ & Farey sum. & Def.~\ref{def: Farey pair}\\
$\gamma\oplus^j\alpha$ & Iterated Farey sum with $\alpha$. & Def.~\ref{def: Farey pair}\\
$\calG$ & Stern--Brocot diagram. & Def.~\ref{def: SternBrocot}\\
$\Delta(\alpha)$ & Stern--Brocot triangle centered at $\alpha$. & Def.~\ref{def: triangle}\\
$\calE$ & Edges of $\calG$ with endpoints in $\Q_0$. & Def.~\ref{def: E and more}\\
$\nu(e)$ & Sign of the slope of $e \in \calE$. & Def.~\ref{def: E and more}\\
$\widehat e$ & Opposite vertex for $e \in \calE$. & Def.~\ref{def: E and more}\\
$e^-,\, e^+$ & High and low endpoints of $e\in \calE$. & Def.~\ref{def: E and more}\\
$e^L,\, e^R$ & Left and right endpoints of $e \in \calE$. & Def.~\ref{def: E and more}\\
$F_\alpha$ & Funnel determined by $\alpha$. & Def.~\ref{def: funnel}\\
$s_0,\ldots,s_n$ & Triangles of $F_\alpha$. & Def.~\ref{def: funnel}\\
$e_0,\ldots,e_{n+1}$ & Distinguished edges of $F_\alpha$. & Def.~\ref{def: funnel}\\
$\{\beta_j\}_{\Z}$ & Neighbor cycle for a hub $\omega$. & Def.~\ref{def: omega cycle}\\

\\[-2pt]
\multicolumn{3}{@{}l}{\textit{Farey-recursive functions and polynomials}}\\[2pt]
$d_\calQ$ & Determinant $d_\calQ(p/q)=(-1)^p x^q$. & Ex.~\ref{ex: dQ}\\
$\calQ$ & Major Farey recursive function. & Ex.~\ref{ex: Q}\\
$\calN$ & Secondary Farey recursive function. & Ex.~\ref{ex: N}\\
$\calZ$ & Shape parameter function. & Def.~\ref{def: shape parameter function}\\
$\calV$ & Vertex function $\calN/\calQ$. & Def.~\ref{def: calV}\\
$D_\calQ$ & Discriminant polynomial function. & Def.~\ref{def: Disc}\\
$\Lambda$ & Riley polynomial function. & Sec.~\ref{sec: riley}\\

\newpage
\\
\multicolumn{3}{@{}l}{\textit{2-bridge links and Sakuma--Weeks complexes}}\\[2pt]
$L_\alpha$ & 2-bridge link of slope $\alpha$. & Def.~\ref{def: 2 bridge link}\\
$M_\alpha$ & Complement $S^3-L_\alpha$. & Def.~\ref{def: 2 bridge link}\\
$\X$ & Punctured plane $\C-\Z[i]$. & Def.~\ref{def: tetrahedral complex}\\
$\X_j$ & Copy of $\X$ associated to $s_j$. & Def.~\ref{def: tetrahedral complex}\\
$\Theta$ & Half turn group acting on $\X$. & Def.~\ref{def: tetrahedral complex}\\
$\widehat K$ & Infinite layered tetrahedral complex. & Def.~\ref{def: tetrahedral complex}\\
$K$ & Quotient $\widehat K/\Theta$. & Def.~\ref{def: quotient complex}\\
$\Sigma_j$ & 4-punctured sphere in $K$. & Def.~\ref{def: quotient complex}\\
$\Delta_j$ & Tetrahedral pair associated to $e_j$. & Def.~\ref{def: quotient complex}\\
$M_\alpha^\circ$ & Complex after the upper hairpin fold. & Def.~\ref{def: SW}\\
$D_\alpha$ & Alternating 4-plat diagram for $L_\alpha$. & Def.~\ref{def: diagram}\\
$\widehat\Theta$ & Half turn group at half lattice points. & Def.~\ref{def: Klein}\\
$G$ & Klein 4-group $\widehat\Theta/\Theta$. & Def.~\ref{def: Klein}\\
$O_\alpha^\circ,\, O_\alpha$ & Klein group quotient complexes. & Def.~\ref{def: orbifolds}\\

\\[-2pt]
\multicolumn{3}{@{}l}{\textit{Hyperbolic structures and developed domains}}\\[2pt]
$\frz$ & Complex specialization parameter. & Ch.~\ref{chap: intro}\\
$\zeta(x)$ & Shape transform $1/(1-x)$. & Def.~\ref{def: parameters}\\
$\calH(M_\alpha^\circ)$ & Positive $G$-invariant gluing solutions. & Def.~\ref{def: hyperbolic structures}\\
$\calH(M_\alpha)$ & Solutions satisfying the lower fold. & Def.~\ref{def: hyperbolic structures}\\
$\dev_{\mathfrak Z},\, \dev_\frz$ & Developing map. & Def.~\ref{def: hyperbolic structures}\\
$\rho_{\mathfrak Z},\, \rho_\frz$ & Holonomy representation. & Def.~\ref{def: hyperbolic structures}\\
$\delta(e),\, \delta'(e)$ & Generic tetrahedra for $e \in \calE$. & Def.~\ref{def: generic tetrahedra}\\
$\delta_j,\, \delta'_j$ & $\delta(e_j),\, \delta'(e_j)$. & Def.~\ref{def: generic tetrahedra}\\
$\Omega(\alpha)$ & Generic developed domain. & Def.~\ref{def: generic domain}\\
$\frz(\alpha),\, \frz_1$ & Geometric root of $\calQ(\alpha)$. & Def.~\ref{def: geometric root}\\
$\sigma_u,\, \sigma_\ell$ & Upper and lower involutions. & Def.~\ref{def: concrete involutions}\\

\\[-2pt]
\multicolumn{3}{@{}l}{\textit{Generic holonomy and face pairings}}\\[2pt]
$k_0,\, k_1$ & Canonical generators of $\pi_1(M_\alpha^\circ)$. & Def.~\ref{def: generators}\\
$U_\omega,\, W_\omega$ & Basic generic holonomy matrices. & Def.~\ref{def: UW}\\
$\Gamma$ & Generic holonomy group $\langle U_0,\, W_0\rangle$. & Def.~\ref{def: generic holonomy}\\
$\Gamma_\frz$ & Specialized holonomy group. & Def.~\ref{def: generic holonomy}\\
$T_\omega$ & $\omega$-meridian. & Def.~\ref{def: S and T}\\
$S_\omega$ & $\omega$-crossing loop. & Def.~\ref{def: S and T}\\
$A_i$ & Generic face pairing element. & Def.~\ref{def: Ai}\\
$\widetilde\Sigma(e)$ & Generic lifted 4-punctured sphere. & Sec.~\ref{sec: xing circles}\\

\newpage
\\
\multicolumn{3}{@{}l}{\textit{Logarithmic spirals}}\\[2pt]
$R_\omega$ & Spiral Möbius transformation. & Def.~\ref{def: XYR}\\

\\[-2pt]
\multicolumn{3}{@{}l}{\textit{Real loci, reversals, and special families}}\\[2pt]
$\R\calZ(e)$ & Real shape locus for $e \in \calE$. & Def.~\ref{dfn: real locus}\\
$\calU_e$ & Thick region for $\calZ(e)$. & Def.~\ref{def: fat}\\
$\bar\alpha$ & Reverse of $\alpha$. & Def.~\ref{def: reverse}\\
$A_\alpha$ & Projective matrix for funnel reversal. & Sec.~\ref{sec: reverse}\\
$\Phi_\alpha$ & Möbius map reversing the funnel. & Sec.~\ref{sec: reverse}\\
$\phi_\alpha$ & Reversal map of Sakuma--Weeks complexes. & Sec.~\ref{sec: reverse}\\
$\frz_j$ & Geometric shape of the $j$th tetrahedral pair. & Sec.~\ref{sec: reverse}\\

\\[-2pt]
\multicolumn{3}{@{}l}{\textit{Riley slice and limiting constructions}}\\[2pt]
$\calD(\Gamma_\frz)$ & Regular set of $\Gamma_\frz$. & Def.~\ref{def: Riley slice}\\
$\calR$ & Riley slice. & Def.~\ref{def: Riley slice}\\
$\widehat\frz(\omega)$ & Geometric root of $D_\calQ(\omega)$. & Def.~\ref{def: D geometric root}\\

\\[-2pt]
\multicolumn{3}{@{}l}{\textit{Heckoid orbifolds and drums}}\\[2pt]
$M(\alpha)_n$ & Cone manifold with cone angle $2\pi/n$. & Sec.~\ref{sec: Heck}\\
$S(\alpha;n)$ & Heckoid orbifold. & Def.~\ref{def: Heckoid orbs}\\
$G(\alpha;n)$ & Heckoid group. & Def.~\ref{def: Heckoid groups}\\
${\mathfrak p}_\alpha$ & Pleating ray of slope $\alpha$. & Def.~\ref{def: KSray}\\
$\calD_n$ & Drum used for $S(2/5;n)$ or $S(1/2;n)$. & Ex.~\ref{ex: Heckoid25}\\
$L(x)$ & Lobachevsky function. & Ex.~\ref{ex: Heckoid12}\\

\end{longtable}
\endgroup

\printindex

\bibliographystyle{amsalpha}
\bibliography{MF}

\end{document}